\documentclass[10pt]{amsart}

\newif\ifdebug
\debugfalse

\usepackage{fullpage}
\usepackage{mathtools}
\usepackage{amsmath,amsthm,amscd,graphicx}
\usepackage[cal=cm, scr=rsfs]{mathalfa}
\usepackage{amssymb}
\usepackage[dvipsnames]{xcolor}
\usepackage{bm}
\usepackage{enumerate}
\usepackage{epsfig}
\usepackage{float}
\usepackage{amsxtra}
\usepackage{pdfpages}
\usepackage{multicol}
\usepackage[normalem]{ulem}
\usepackage{psfrag}
\usepackage{tikz}
\usepackage[T1]{fontenc}

\usepackage{indentfirst}
\usepackage{tikz-cd}
\usepackage{url}
\usepackage[all,pdf]{xy}
\usepackage{color}
\usepackage{verbatim}

\RequirePackage{xspace}
\usepackage{mathabx}
\makeatletter
\def\thm@space@setup{\thm@preskip=2ex \thm@postskip=2ex
}
\makeatother

\usepackage{bbm}

\numberwithin{equation}{subsubsection}

\usepackage{hyperref}
\usepackage{aliascnt}
\usepackage{cleveref}

\theoremstyle{plain}

\newtheorem{thm}[subsubsection]{Theorem~}
\newtheorem{lem}[subsubsection]{Lemma~}
\newtheorem{prop}[subsubsection]{Proposition~}

\newtheorem{cor}[subsubsection]{Corollary~}

\theoremstyle{remark}
\newtheorem{rmk}[subsubsection]{Remark~}
\newtheorem{examp}[subsubsection]{Example~}
\newtheorem{ass}[subsubsection]{Assumption~}

\theoremstyle{definition}
\newtheorem{defn}[subsubsection]{Definition~}

\newcommand\m[1]{\mathrm{#1}}

\newcommand{\Q}{\mathbb{Q}}

\newcommand{\A}{\mathbb{A}}
\newcommand{\Z}{\mathbb{Z}}
\newcommand{\R}{\mathbb{R}}
\newcommand{\C}{\mathbb{C}}

\newcommand{\ov}{\overline}
\newcommand{\wt}{\widetilde}

\newcommand{\gr}{\mathrm{gr}\,}

\newcommand{\adm}{\mathrm{adm}}

\newcommand{\hol}{\mathrm{hol}}

\newcommand{\GL}{\mathrm{GL}}
\newcommand{\ad}{\mathrm{ad}}

\newcommand{\SL}{\mathrm{SL}}
\newcommand{\SO}{\mathrm{SO}}

\newcommand{\Gal}{\mathrm{Gal}}

\newcommand{\Mat}{\mathrm{Mat}}
\newcommand{\Sym}{\mathrm{Sym}}
\newcommand{\Rep}{\mathrm{Rep}}

\newcommand{\End}{\mathrm{End}}
\newcommand{\sEnd}{{\underline{\mathrm{sEnd}}^\circ}}
\newcommand{\sEndint}{\underline{\mathrm{sEnd}}}
\newcommand{\Endint}{\underline{\mathrm{End}}}

\newcommand{\Spa}{\mathrm{Spa}\,}
\newcommand{\Spd}{\mathrm{Spd}\,}
\newcommand{\Spf}{\mathrm{Spf}\,}
\newcommand{\supp}{\mathrm{supp}\,}

\newcommand{\vol}{\mathrm{vol}\,}
\newcommand{\simlr}{\xrightarrow{\sim}}
\newcommand{\Spec}{\mathrm{Spec}\,}
\newcommand{\Frob}{\mathrm{Frob}}
\newcommand{\QuadHyper}{\mathcal{Q}}

\DeclareMathOperator{\length}{length}
\DeclareMathOperator{\Lie}{Lie}
\DeclareMathOperator{\Hom}{Hom}
\DeclareMathOperator{\Homint}{\underline{\mathrm{Hom}}}

\DeclareMathOperator{\Res}{Res}
\DeclareMathOperator{\Fil}{Fil}

\DeclareMathOperator{\tr}{tr}
\DeclareMathOperator{\Tr}{tr}

\newcommand{\charfun}{\mathbbm{1}}
\newcommand{\Sh}{\mathscr{S}}
\newcommand{\ShG}{\mathrm{Sh}}
\newcommand{\Ch}{\mathrm{Ch}}

\newcommand{\dR}{\mathrm{dR}}
\newcommand{\Betti}{\mathrm{B}}

\newcommand{\aux}{\mathrm{aux}}

\newcommand{\Gr}{\mathrm{Gr}}
\newcommand{\bpi}{\boldsymbol{\pi}}
\newcommand{\DR}{\mathrm{dR}}

\newcommand{\perf}{\mathrm{perf}}
\newcommand{\pdivgroup}{\mathrm{BT}}
\newcommand{\Diff}{\mathrm{Diff}}

\newcommand{\pDen}{\partial\mathrm{Den}}
\newcommand{\Int}{\mathrm{Int}}
\newcommand{\pEis}{\partial\mathrm{Eis}}

\newcommand{\disc}{\mathrm{disc}}
\newcommand{\basic}{\mathrm{basic}}
\newcommand{\Has}{\epsilon}
\newcommand{\val}{\mathrm{val}}
\newcommand{\rank}{\mathrm{rank}}
\newcommand{\sAut}{\underline{\mathrm{Aut}}}
\newcommand{\der}{\mathrm{der}}
\newcommand{\betti}{\mathrm{B}}
\newcommand{\cris}{\mathrm{cris}}

\newcommand{\Bcris}{B_{\cris}}
\newcommand{\loc}{\mathrm{loc}}
\newcommand{\Fb}{\breve{\FLocal}}
\newcommand{\kb}{\breve{k}}
\newcommand{\et}{{\acute{\mathrm{e}}\mathrm{t}}}
\newcommand{\OFb}{O_{\breve{\FLocal}}}

\DeclareMathOperator{\Aut}{Aut}
\DeclareMathOperator{\Sp}{Sp}

\DeclareMathOperator{\Pic}{Pic}

\DeclareMathOperator{\im}{im}

\DeclareMathOperator{\Or}{O}

\newcommand{\id}{\ensuremath{\mathrm{id}}\xspace}
\DeclareMathOperator{\GSp}{GSp}
\DeclareMathOperator{\Isom}{\underline{Isom}}

\DeclareMathOperator{\Vect}{Vect}
\DeclareMathOperator{\GSpin}{GSpin}
\DeclareMathOperator{\Spin}{Spin}

\DeclareMathOperator{\Mp}{Mp}
\DeclareMathOperator{\FJ}{FJ}
\DeclareMathOperator{\an}{an}
\DeclareMathOperator{\realization}{real}

\newcommand{\red}{{\m{red}}}

\DeclareMathOperator{\Hilb}{Hilb}
\DeclareMathOperator{\Mor}{Mor}
\DeclareMathOperator{\Loc}{Loc}

\newcommand{\Gm}{\BG_\mathrm{m}}
\newcommand{\bark}{\bar{k}}
\newcommand{\barFp}{\bar{\BF}_p}

\usepackage{relsize}
\usepackage[bbgreekl]{mathbbol}
\usepackage{amsfonts}
\DeclareSymbolFontAlphabet{\mathbbl}{bbold}
\newcommand{\Prism}{{\mathlarger{\mathbbl{\Delta}}}}

\newcommand{\fka}{\ensuremath{\mathfrak{a}}\xspace}

\newcommand{\fkd}{\ensuremath{\mathfrak{d}}\xspace}

\newcommand{\fkf}{\ensuremath{\mathfrak{f}}\xspace}

\newcommand{\sA}{\ensuremath{\mathscr{A}}\xspace}

\newcommand{\sP}{\ensuremath{\mathscr{P}}\xspace}

\newcommand{\sS}{\ensuremath{\mathscr{S}}\xspace}

\newcommand{\sV}{\ensuremath{\mathscr{V}}\xspace}

\newcommand{\sX}{\ensuremath{\mathscr{X}}\xspace}

\newcommand{\fkB}{\ensuremath{\mathfrak{B}}\xspace}

\newcommand{\fkS}{\ensuremath{\mathfrak{S}}\xspace}

\newcommand{\BA}{\ensuremath{\mathbb{A}}\xspace}

\newcommand{\BC}{\ensuremath{\mathbb{C}}\xspace}
\newcommand{\BD}{\ensuremath{\mathbb{D}}\xspace}

\newcommand{\BF}{\ensuremath{\mathbb{F}}\xspace}
\newcommand{\BG}{\ensuremath{\mathbb{G}}\xspace}
\newcommand{\BH}{\ensuremath{\mathbb{H}}\xspace}

\newcommand{\BL}{\ensuremath{\mathbb{L}}\xspace}
\newcommand{\BM}{\ensuremath{\mathbb{M}}\xspace}

\newcommand{\BQ}{\ensuremath{\mathbb{Q}}\xspace}
\newcommand{\BR}{\ensuremath{\mathbb{R}}\xspace}
\newcommand{\BS}{\ensuremath{\mathbb{S}}\xspace}

\newcommand{\BZ}{\ensuremath{\mathbb{Z}}\xspace}

\newcommand{\bH}{\ensuremath{\mathbf{H}}\xspace}

\newcommand{\bu}{\ensuremath{\mathbf{u}}\xspace}

\newcommand{\bx}{\ensuremath{\mathbf{x}}\xspace}

\newcommand{\CA}{\ensuremath{\mathcal{A}}\xspace}
\newcommand{\CB}{\ensuremath{\mathcal{B}}\xspace}

\newcommand{\CE}{\ensuremath{\mathcal{E}}\xspace}
\newcommand{\CF}{\ensuremath{\mathcal{F}}\xspace}
\newcommand{\CG}{\ensuremath{\mathcal{G}}\xspace}
\newcommand{\CH}{\ensuremath{\mathcal{H}}\xspace}
\newcommand{\CI}{\ensuremath{\mathcal{I}}\xspace}

\newcommand{\CM}{\ensuremath{\mathcal{M}}\xspace}

\newcommand{\CO}{\ensuremath{\mathcal{O}}\xspace}
\newcommand{\CP}{\ensuremath{\mathcal{P}}\xspace}
\newcommand{\CQ}{\ensuremath{\mathcal{Q}}\xspace}

\newcommand{\CT}{\ensuremath{\mathcal{T}}\xspace}
\newcommand{\CU}{\ensuremath{\mathcal{U}}\xspace}
\newcommand{\CV}{\ensuremath{\mathcal{V}}\xspace}
\newcommand{\CW}{\ensuremath{\mathcal{W}}\xspace}
\newcommand{\CX}{\ensuremath{\mathcal{X}}\xspace}
\newcommand{\CY}{\ensuremath{\mathcal{Y}}\xspace}
\newcommand{\CZ}{\ensuremath{\mathcal{Z}}\xspace}

\newcommand{\RT}{\ensuremath{\mathrm{T}}\xspace}

\newcommand{\sz}{\mathsf{z}}
\newcommand{\sx}{\mathsf{x}}
\newcommand{\sy}{\mathsf{y}}

\newcommand{\RZname}{Rapoport--Zink\xspace}
\newcommand{\AbSpecialEnd}{\sEndint_2^\circ}

\newcommand{\RealPlace}{\nu}
\newcommand{\commonG}{G^{\mathrm{cmp}}}
\newcommand{\commonCG}{\CG^{\mathrm{cmp}}}
\newcommand{\commonX}{X^{\mathrm{cmp}}}
\newcommand{\commonh}{h^{\mathrm{cmp}}}
\newcommand{\commoniota}{\iota^{\mathrm{cmp}}}
\newcommand{\commonE}{E^{\mathrm{cmp}}}
\newcommand{\commonv}{v^{\mathrm{cmp}}}
\newcommand{\commonx}{\mathsf{x}^{\mathrm{cmp}}}
\newcommand{\commonKernel}{\kernelalphaG^{\mathrm{cmp}}}
\newcommand{\commonK}{K^{\mathrm{cmp}}}
\newcommand{\commonKprimeToP}{K^{\mathrm{cmp},p}}
\newcommand{\commonmu}{\mu^{\mathrm{cmp}}}
\newcommand{\commonW}{W^{\mathrm{cmp}}}
\newcommand{\commonPsi}{\Psi^{\mathrm{cmp}}}
\newcommand{\commonCV}{\CV^{\mathrm{cmp}}}
\newcommand{\commonA}{\CA^{\mathrm{cmp}}}
\newcommand{\commonSpecialEnd}{{\underline{\mathrm{sEnd}}^{\circ,\mathrm{cmp}}}}

\ifdebug
\newcommand{\colorNumberOne}{\color{blue}}
\newcommand{\colorNumberTwo}{\color{cyan}}
\usepackage[notref,notcite]{showkeys}
\usepackage{seqsplit}

\else
\newcommand{\colorNumberOne}{}
\newcommand{\colorNumberTwo}{}
\fi

\newcommand{\dimofV}{{\colorNumberOne m}}
\newcommand{\nInLZorth}{{\colorNumberOne n}}

\newcommand{\Vlocal}{{\colorNumberTwo \mathrm{V}}}
\newcommand{\Vspecial}{{\colorNumberOne \mathrm{V}'}}
\newcommand{\Vspecialp}{{\colorNumberOne \mathrm{V}_p'}}
\newcommand{\VspecialGlobal}{{\colorNumberOne V'}}
\newcommand{\Vincoherent}{{\colorNumberOne \mathbb{V}}}
\newcommand{\Vcoherent}{{\colorNumberOne V}}
\newcommand{\Vadl}{{\colorNumberOne V}}
\newcommand{\VlattiveGlobal}{{\colorNumberOne \Lambda}}
\newcommand{\VlattiveGlobalp}{{\colorNumberOne \Lambda_{(p)}}}
\newcommand{\VlattiveGlobalpcomplete}{{\colorNumberOne \Lambda_{p}}}
\newcommand{\VlattiveGlobalv}{{\colorNumberOne \Lambda_{v}}}
\newcommand{\VlattiveGlobalvv}{{\colorNumberOne \Lambda_{v'}}}
\newcommand{\VGlobalp}{{\colorNumberOne V_{p}}}
\newcommand{\VGlobalv}{{\colorNumberOne V_{v}}}
\newcommand{\VlatticeLocal}{{\colorNumberTwo \Lambda}}
\newcommand{\VlattiveLocalo}{{\colorNumberTwo \Lambda_{v_0}}}
\newcommand{\VlattiveLocalv}{{\colorNumberTwo \Lambda_{v}}}

\newcommand{\VlattiveLocaltot}{{\colorNumberTwo \Lambda}}
\newcommand{\FLocal}{{\colorNumberTwo \mathrm{F}}}
\newcommand{\FGlobal}{{\colorNumberOne F}}
\newcommand{\FGloball}{\FGlobal}

\newcommand{\reflexField}{{\colorNumberOne E_\Phi}}

\newcommand{\RZG}{\CM_{\CG, b, \mu}}

\newcommand{\RZGS}{\CM_{\GSpin^\flat}}
\newcommand{\RZSOflat}{\CM_{\SO^\flat}}

\newcommand{\Lder}{\mathbb{L}}
\newcommand{\Lattice}{L}

\newcommand{\GSpinRZ}{\mathcal M_{\mathrm{GSpin}}}
\newcommand{\SORZ}{\mathcal M_{\mathrm{SO}}}

\newcommand{\SORZcomp}{\mathcal N_{\mathrm{fix}}}

\newcommand{\BTstratum}[1]{\mathcal V_{#1}}
\newcommand{\SOShim}{\ShG_{K_2}}
\newcommand{\SOIntModel}{\sS_{K_2}}

\newcommand{\valFLocal}{\val_\FLocal}
\newcommand{\reflexHodge}{E}
\newcommand{\reflexAbel}{E_2}
\newcommand{\reflexHodgeLocal}{{\colorNumberOne \mathrm{E}}}
\newcommand{\reflexAbelLocal}{{\colorNumberOne \mathrm{E}_2}}

\newcommand{\centralHodge}{{\colorNumberOne \CT}}
\newcommand{\centralAbel}{{\colorNumberOne \CT_2}}
\newcommand{\kernelalpha}{{\colorNumberOne \CT_\alpha}}
\newcommand{\centralHodgeG}{{\colorNumberOne T}}
\newcommand{\centralAbelG}{{\colorNumberOne T_2}}
\newcommand{\kernelalphaG}{{\colorNumberOne T_\alpha}}
\newcommand{\kernelalphaGp}{{\colorNumberOne T_{\alpha,p}}}

\newcommand{\Sha}{\operatorname{Sha}}

\crefname{thm}{theorem}{theorems}
\Crefname{thm}{Theorem}{Theorems}
\crefname{lem}{lemma}{lemmas}
\Crefname{lem}{Lemma}{Lemmas}
\crefname{prop}{proposition}{propositions}
\Crefname{prop}{Proposition}{Propositions}
\crefname{cor}{corollary}{corollaries}
\Crefname{cor}{Corollary}{Corollaries}
\crefname{defn}{definition}{definitions}
\Crefname{defn}{Definition}{Definitions}
\crefname{rmk}{remark}{remarks}
\Crefname{rmk}{Remark}{Remarks}
\crefname{ass}{assumption}{assumptions}
\Crefname{ass}{Assumption}{Assumptions}
\title{On the arithmetic Siegel--Weil formula for orthogonal Shimura varieties over a general totally real field}

\author[X. Li]{Xuanyou Li}
\address{Qiuzhen College, Tsinghua University, Beijing 100084, P. R. China}
\email{lixuanyo21@mails.tsinghua.edu.cn}

\date{\today}

\usepackage[alphabetic]{amsrefs}

\begin{document}
	\bibliographystyle{amsalpha}

	\begin{abstract}
		Let $\FLocal/\Q_p$ be a finite extension with $p\ne 2$, and
		let $\Vlocal$ be a quadratic space of dimension $\dimofV\ge2$ containing
		a self-dual lattice.  We prove the local arithmetic Siegel--Weil formula for GSpin \RZname spaces formulated in \cite{LZ-orth}. 
		We use the global-to-local strategy of
		\cite{MihatschZhang2024,Luo2025KudlaRapoport}. Using a generalization of the prismatic realization functor in \cite{imai2025prismatic}, we construct and study 
		derived special cycles on integral models of SO Shimura varieties and their
		Rapoport--Zink uniformizations. We derive the local intersection formula from the modularity of special cycles on the generic fiber.  As a first application, we prove the
		corresponding semi-global arithmetic Siegel--Weil identity for orthogonal 
		Shimura varieties over a totally real field $\FGlobal\neq \Q$.
	\end{abstract}
	
	\maketitle
	
	\tableofcontents
	
	\section{Introduction}
	The classical Siegel--Weil formula \cite{Weil1965SiegelFormula} relates
	theta integrals to special values of Siegel Eisenstein series, connecting
	the arithmetic of quadratic forms with automorphic forms.  The work of
	Kudla and Millson \cite{KudlaMillson1990} and Kudla
	\cite{Kudla1997AlgebraicCycles} gives this relation a geometric
	interpretation: the cohomology classes of special cycles on orthogonal
	Shimura varieties occur as Fourier coefficients of Siegel modular forms,
	and their degrees are related to Fourier coefficients of Eisenstein
	series.
	
	Kudla \cite{Kudla1997CentralDerivatives,Kudla2004} proposed an arithmetic
	refinement in which arithmetic degrees and height pairings of special
	cycles are governed by derivatives of Eisenstein series.  In the
	orthogonal setting, the arithmetic Siegel--Weil formula predicts that
	arithmetic intersections of special cycles are related to central
	derivatives of incoherent Siegel Eisenstein series.  The work of Kudla,
	Rapoport and Yang \cite{Kudla-Rapoport-Yang} establishes this picture for
	Shimura curves.  In higher dimensions, Bruinier and Yang
	\cite{BruinierYang2021Arithmetic} prove cases of the formula and establish
	local identities at both archimedean and nonarchimedean places, while
	Garcia and Sankaran \cite{garcia_sankaran_2019_green_forms} construct Green
	forms and obtain the archimedean height identities for compact orthogonal
	and unitary Shimura varieties.
	
	A central local problem in this program is to identify derived
	intersection numbers of special cycles on Rapoport--Zink spaces with
	derivatives of local representation densities of quadratic forms.  These
	densities describe the local Whittaker coefficients on the analytic side,
	while Rapoport--Zink uniformization relates the local intersections to
	arithmetic degrees on integral models of Shimura varieties.  In the
	unitary setting, this relation was formulated by Kudla and Rapoport
	\cite{Kudla2011,KR14} and proved by Li and Zhang \cite{LZ21}.
	In the orthogonal setting, Li and Zhang
	\cite{LZ-orth} proved a local arithmetic Siegel--Weil formula for GSpin
	Rapoport--Zink spaces over $\Q_p$, for odd $p$, and deduced a semi-global
	formula.  He, Zhang and Zhu \cite{orth-local-Almost-Self-Dual} treat
	almost self-dual parahoric level by introducing weighted special cycles
	on GSpin Rapoport--Zink spaces.  Their local arithmetic Siegel--Weil
	formulas relate the arithmetic intersections of these cycles to
	derivatives of local Whittaker functions, extending the picture beyond
	self-dual level.
	The relation with automorphic $L$-derivatives is developed through
	arithmetic theta lifting in Liu \cite{Liu11a,Liu11b} and Li--Liu
	\cite{LL21,LL22} for unitary groups.
	
	The purpose of this paper is to establish the corresponding
	local identity over an arbitrary finite extension of $\Q_p$ of odd residue
	characteristic, and to obtain a semi-global formula over a general totally
	real field of degree greater than one.  
	
	\subsection{Main result}
	Let $p$ be an odd prime and let $\FLocal/\Q_p$ be a finite extension.
	Let $\breve{\FLocal}$ be the completion of the maximal unramified
	extension of $\FLocal$.
	Let $(\Vlocal,Q)$ be a nondegenerate quadratic space over $\FLocal$
	of dimension $\dimofV=\nInLZorth+1\ge2$ admitting a self-dual
	$O_{\FLocal}$-lattice $\VlattiveLocaltot$, and put
	$\epsilon:=\chi(\Vlocal)\in\{\pm1\}$. 
	We assume\footnote{Otherwise there is no nonzero special endomorphism (cf. \Cref{prop:V_b-SO-2}) and so problem is trivial, exactly as in the unitary case at a split prime.} $(\dimofV,\epsilon)\neq (2,1)$. 
	Associated to these data is a (basic) SO \RZname space $\SORZ$,
	constructed using \cite{PappasRapoport2026IntegralLocalShimura}.
	It is a formal scheme over $\Spf O_{\breve{\FLocal}}$, formally smooth
	of relative dimension $\dimofV-2$; see
	\Cref{defn:SO-RZ-space,subsec:dim-two}.
	
	Let $\Vspecial$ be the $\FLocal$-space of special quasi-endomorphisms
	of the basic framing object on a chosen GSpin cover. For $\dimofV\ge3$, it is the quadratic
	space of dimension $\dimofV$ with $\chi(\Vspecial)=\epsilon$ and Hasse
	invariant $-1$ described in \Cref{prop:V_b-SO}.
	In dimension two, the same description holds for $\epsilon=-1$, whereas
	$\Vspecial=0$ for $\epsilon=1$; see \Cref{prop:V_b-SO-2}.
	
	For $u\in\Vspecial\setminus\{0\}$, the special divisor
	$\CZ_{\mathrm{SO}}(u)\subseteq\SORZ$ is defined by descending the locus
	on the GSpin \RZname space where $u$ deforms to an endomorphism of the
	universal $p$-divisible group; see
	\Cref{defn:SO-local-special-cycle,subsec:dim-two}.
	The Cartier divisor property is established in
	\Cref{prop:local-special-cycle-properties,prop:n-two-special}.

	For a nonsingular tuple
	$\bu=(u_1,\ldots,u_{\dimofV-1})\in\Vspecial^{\dimofV-1}$,
	we consider a geometric intersection number and an analytic density
	derivative. The arithmetic intersection number $\Int^\epsilon(\bu)$
	is the Euler--Poincar\'e characteristic of the derived intersection
	\[
	\mathcal O_{\CZ_{\mathrm{SO}}(u_1)}
	\otimes^{\BL}_{\mathcal O_{\SORZ}}\cdots
	\otimes^{\BL}_{\mathcal O_{\SORZ}}
	\mathcal O_{\CZ_{\mathrm{SO}}(u_{\dimofV-1})},
	\]
	computed on one connected component of $\SORZ$, as in
	\Cref{defn:derived-special-cycles-RZ}. 
	On the analytic side, the normalized local density polynomial
	$\operatorname{Den}^\epsilon(X,\sum_i O_\FLocal u_i)$ is defined in
	\Cref{defn:local-density}, following \cite{LZ-orth}. Its derivative
	defines the derived local density
	\[
	\pDen^\epsilon(\bu)
	=-\left.\frac{d}{dX}\operatorname{Den}^\epsilon
	\!\left(X,\sum_i O_\FLocal u_i\right)\right|_{X=1}.
	\]
	Our local arithmetic Siegel--Weil formula identifies these two invariants.
	
	\begin{thm}[\Cref{thm:main2}]\label{thm:main}
		Let
		\(\bu=(u_1,\ldots,u_{\dimofV-1})\in\Vspecial^{\dimofV-1}\)
		have a nonsingular Gram matrix. 
		Then
		\[
		\Int^\epsilon(\bu)=\pDen^\epsilon(\bu).
		\]
	\end{thm}
	
	We next describe a semi-global application of the local formula.
	Let $\FGlobal$ be a totally real field with
	$[\FGlobal:\Q]>1$ and $\dimofV\ge3$.
	Choose a quadratic space $\Vcoherent/\FGlobal$ of dimension $\dimofV$,
	with signature $(\dimofV-2,2)$ at one real place $\RealPlace_0$ and
	positive definite at every other real place.
	Fix an even positive integer $\fkd$ and a sufficiently small neat
	compact open product level $K_2=\prod_{w<\infty}K_{2,w}$ whose local factors away from $\fkd$ are the
	parahoric subgroups defined by self-dual lattices.
	The theory of \cite{KisinPappasZhou2024ParahoricII} gives a smooth
	integral model $\SOIntModel/\mathcal O_\FGlobal[1/\fkd]$ of the
	orthogonal Shimura variety, which is proper because $\FGlobal\ne\Q$;
	see \Cref{subsec:orth-global-models,prop:properness}.

	The special-endomorphism scheme of
	\Cref{prop:global-special-endomorphism-model} gives global special
	cycles with prescribed Gram matrix and finite adelic level conditions.
	For a nonsingular $T\in\Sym_{\dimofV-1}(\FGlobal)$ and an admissible
	Schwartz function
	$\phi\in\sS((\Vcoherent\otimes_\FGlobal\BA_{\FGlobal,f})^{\dimofV-1})^{K_2}$,
	the class $\CZ_T^{\BL}(\phi)$ is their derived intersection class,
	weighted by $\phi$, as defined in
	\Cref{defn:derived-special-cycles-Sh}.
	At a finite place $v\nmid\fkd$, its degree
	$\deg_v\CZ_T^{\BL}(\phi)$ is the alternating sum of cohomology lengths
	with the convention of \Cref{eq:deg-v-convention}.

	On the analytic side, let $\Vincoherent$ be the incoherent adelic
	quadratic space with the same finite localizations as $\Vcoherent$
	and positive definite at every real place. 
	The set $\Diff(T,\Vincoherent)$ consists of the
	places where $\Vincoherent$ fails to represent $T$, as in
	\Cref{defn:Eisenstein-difference-set}.
	For a finite place $v\nmid\fkd$, the nearby global space $\VspecialGlobal^{(v)}$
	has the same dimension and discriminant as $\Vincoherent$, the opposite
	Hasse invariant at $v$, and the same localizations at every other place.
	Write $\epsilon_v=\chi(\Vincoherent_v)$ and let $K_2^v$ denote the
	level away from $v$. For a tuple $\mathbf x$ in $\VspecialGlobal^{(v)}$,
	$\operatorname{Orb}_{\mathbf x}(\phi^v)$ denotes the orbital integral
	of the factor of $\phi$ away from $v$; see \Cref{thm:global-int}.
	The local identity and the local--global intersection formula of
	\Cref{thm:global-int} give the following semi-global result.
	\begin{thm}[\Cref{thm:semi-global-ASW}]
		Let $v\nmid\fkd$
		be a finite place and let
		$T\in\Sym_{\dimofV-1}(\FGloball)$ be nonsingular with
		$\Diff(T,\Vincoherent)=\{v\}$.  Choose
		\(\mathbf x\in(\VspecialGlobal^{(v)})^{\dimofV-1}\) with
		\(T(\mathbf x)=T\).  For an admissible factorizable Schwartz
		function $\phi$, one has
		\begin{equation*}
			\deg_v\CZ_T^{\BL}(\phi)
			=\frac{2}{\vol(K_2^v)}\,
			\pDen_v^{\epsilon_v}(\mathbf x)\,
			\operatorname{Orb}_{\mathbf x}(\phi^v).
		\end{equation*}
	\end{thm}
	\subsection{Proof strategy}
	The proof follows the global-to-local approach of
	\cite[\S10]{MihatschZhang2024} and its use in the local arithmetic Siegel--Weil formula 
	by Luo \cite{Luo2025KudlaRapoport}.  
	Basic uniformization relates intersection numbers on Shimura varieties
	to those on \RZname spaces. We then show that the differences between
	the two sides of the desired identity occur as coefficients of a
	modular form, using the modularity of special cycles and Green
	functions established in \cite{YuanZhangZhang2009,Sankaran2025ArithmeticSpecialCycles}.
	The archimedean contribution is computed using
	\cite{garcia_sankaran_2019_green_forms}. To prove that this modular
	form vanishes, we apply the boundary-vanishing criterion of
	\Cref{lem:boundary-vanishing}, which plays the role of a global
	analogue of the ``uncertainty principle'' used in
	\cite{LZ-orth,LZ21}. The required boundary vanishing is
	deduced from the lower-dimensional case by the cancellation law
	in \Cref{cor:cancellation}.

	The main new geometric ingredient is the construction of special
	endomorphisms on abelian-type integral models of orthogonal Shimura
	varieties and the proof that the deformation loci of nonzero special
	endomorphisms are \'etale locally effective Cartier divisors. The usual  constructions in
	\cite{MadapusiPera2016,Madapusi2022Derived}	require $\FGlobal/\Q$ to be unramified at $p$. We extend these
	constructions to general $\FGlobal/\Q$ using the theory of
	prismatic $\CG$-crystals developed in \Cref{sec:prism}. More precisely, we prove the following generalization of 
	\cite[Theorem~2.29 and Corollary~2.31]{IKY26}, from reductive group
	schemes over $\Z_p$ to connected parahoric group schemes. 
	On the $p$-adic completion of a Hodge-type cover, this theorem supplies
	the crystalline realizations used to define special endomorphisms;
	see \Cref{prop:cris-tensor}.
	\begin{thm}[\Cref{thm:para-equiv}]\label{thm:para-equiv-intro}
		Let $X$ be a smooth $p$-adic formal scheme over a complete discrete valuation ring $O_K$ of mixed characteristic and perfect residue field. Let $\CG/\Z_p$ be a parahoric group scheme. 
		Let $\omega\in\CG\text{-}\Loc_{\Z_p}(X_\eta)$, and fix a faithful
		representation $\rho:\CG\hookrightarrow\GL(\Lambda)$. If
		$\omega(\Lambda)$ has prismatic good reduction, there exist an exact
		$\Z_p$-linear symmetric monoidal functor
		\[
		\nu:\Rep_{\Z_p}(\CG)\longrightarrow\Vect^\varphi(X_\Prism,\CO_\Prism)
		\]
		and a monoidal natural isomorphism $\alpha:T_{\et}\circ\nu\simeq\omega$.
	\end{thm}

	For $x\in\Vspecial\setminus\{0\}$, we define a special locus $\CZ_{\mathrm{SO}}(x)$ on the
	orthogonal \RZname space by requiring $x$ to extend to an endomorphism
	of the universal $p$-divisible group on a GSpin cover; see
	\Cref{defn:SO-local-special-cycle}.
	A priori, this gives only a closed subscheme. Basic uniformization
	and the global results show that $\CZ_{\mathrm{SO}}(x)$ is a Cartier divisor.
	When $x$ has unit norm, the same method shows that $\CZ_{\mathrm{SO}}(x)$ is
	smooth over $O_{\breve{\FLocal}}$. Its associated diamond then
	determines the formal scheme. For $\dimofV\ge4$, we identify $\CZ_{\mathrm{SO}}(x)$
	with a lower-dimensional orthogonal \RZname space. For $\dimofV=3$,
	it is the disjoint union of the two rank-two spaces corresponding to
	the two Galois-conjugate cocharacters. These identifications give
	the cancellation law in \Cref{cor:cancellation}.

	\subsection{Further directions}
	The results of this paper may be applicable to questions relating Chow
	groups and central derivatives of automorphic $L$-functions in the
	orthogonal setting, in the spirit of the work of Chao Li and Yifeng Liu
	\cite{LL21,LL22} for unitary groups. Two difficulties are the case $p=2$
	and the absence of an orthogonal analogue of the split places of
	the CM extension used in the unitary case. The cancellation argument
	in \Cref{cor:cancellation} requires odd residue characteristic:
	a unimodular even lattice over $\Z_2$ has no orthogonal direct summand
	of rank~$1$.
	
	It would also be interesting to extend the present method to almost
	self-dual level over general finite extensions of $\Q_p$, generalizing \cite{LZ-orth} and \cite{orth-local-Almost-Self-Dual}.
	
	\subsection{Organization of the paper}
	\Cref{part1} develops the geometric constructions in a general abelian-type
	setting.  In \Cref{sec:special-cycles-abelian-type} we descend special
	cycles from Hodge-type covers and prove independence of the auxiliary
	Hodge-type datum.  In \Cref{sec:RZ} we construct the
	corresponding cycles on \RZname spaces and compare them with the global
	cycles via basic uniformization.
	
	\Cref{part2} specializes these constructions to orthogonal Shimura varieties.
	Section~\ref{sec:GSpin} develops the relevant $\SO$ Shimura varieties
	and their GSpin covers. Section~\ref{sec:GSpin-RZ} considers the associated \RZname spaces.
	Section~\ref{sec:Int-number} defines the local and global derived
	intersection numbers and establishes the local--global intersection
	formula.
	
	\Cref{part3} supplies the analytic input and completes the global-to-local
	argument. Section~\ref{sec:Weil} recalls the Weil representations,
	Whittaker functions and local densities used on the analytic side of
	\Cref{thm:main}. Section~\ref{sec:analytic} introduces the Eisenstein
	series. Section~\ref{sec:height} treats the archimedean
	height comparison, and Section~\ref{sec:modularity} constructs the modular
	error series.  In Section~\ref{sec:proof} we globalize the
	prescribed local datum and use the above results to deduce the main theorem. 
	
	\Cref{sec:prism} establishes a generalization of
	\cite[Theorem~2.29]{IKY26} to connected parahoric group schemes,
	which is used to construct the required crystalline realizations. 
	
	\subsection{Acknowledgments}The author would like to thank Yichao Tian for introducing him to this area and for his constant support and encouragement. He would like to thank Keerthi Madapusi for some helpful discussions.
	
	\subsection{Notation and conventions}\label{nutation}
	\begingroup
\setlength{\leftmargini}{1em}

\subsubsection{Notation on local quadratic spaces and lattices}
\label{sec:notat-quadr-latt-1}
\begin{enumerate}[(L1)]
\setlength{\itemsep}{3pt}
\item Let $(V, (\ ,\ ))$ be a quadratic space of dimension $\dimofV$.
Set $Q(x):=\frac12(x,x)$.
For a tuple
$\mathbf x=(x_1,\ldots,x_s)$ in a quadratic space, the
\emph{Gram matrix} is defined by
\begin{equation}\label{eq:Gram-matrix-convention}
	T(\mathbf x):=\bigl((x_i,x_j)\bigr)_{1\le i,j\le s}.
\end{equation}

\item Let $\FLocal$ be a non-archimedean local field
of odd residue characteristic $p$. Write $O_\FLocal$ for its ring of
integers, $k_\FLocal=\mathbb F_q$ for its residue field, and $\varpi$ for a
uniformizer.  We normalize the valuation and absolute value as
\[
\val=\valFLocal:\FLocal\to\BZ\cup\{\infty\},
\qquad
|\cdot|:\FLocal\to\mathbb{R}_{\ge0},\qquad
\val(\varpi)=1,\quad |\varpi|=q^{-1},\quad |a|=q^{-\val(a)}.
\]
\item Let $\Fb$ be the completion of the
maximal unramified extension of $\FLocal$, with ring of integers $\OFb$
and residue field $\kb$. Let
$\sigma_\FLocal\in\Aut(\OFb)$ denote the lift of the absolute $q$-Frobenius on
$\kb$.

\item A quadratic $O_\FLocal$-lattice
$L$ of rank $r$, with symmetric bilinear form $(\ ,\ )$, is
\emph{nondegenerate} if the induced form on
$L_\FLocal:=L\otimes_{O_\FLocal}\FLocal$ is nondegenerate.
All quadratic lattices are assumed nondegenerate unless stated otherwise.
The dual lattice is
\[
L^\vee:=\{x\in L_\FLocal:(x,L)\subseteq O_\FLocal\},
\]
and $L$ is \emph{integral} if $L\subseteq L^\vee$.

\item For an integral lattice $L$, its \emph{fundamental invariants}
are the unique integers
\[
0\le a_1\le\cdots\le a_r
\quad\text{such that}\quad
L^\vee/L\simeq\bigoplus_{i=1}^r O_\FLocal/\varpi^{a_i}.
\]
Its valuation and type are
$
\val(L):=\sum_{i=1}^r a_i$, 
$
t(L):=\#\{i:a_i\ne0\}.
$

\item For a vector $x$ in a quadratic space
over $\FLocal$, set $\val(x):=\val((x,x))$. With $\val(0)=\infty$ and
$|0|=0$, every isotropic vector has valuation $\infty$, and
$\val(ax)=2\val(a)+\val(x)$ for $a\in\FLocal^\times$.

\item For an integral lattice $L$, a \emph{standard orthogonal basis}
is an orthogonal $O_\FLocal$-basis $\{e_1,\ldots,e_r\}$ satisfying
$\val(e_i)=a_i$; such a basis exists for every integral lattice.

\item A lattice $L$ is
\emph{self-dual} if
$L=L^\vee$, and it is \emph{minuscule}, or a \emph{vertex lattice}, if
\[
L\subseteq L^\vee\subseteq\varpi^{-1}L.
\]
For a vertex lattice $\Lambda$ of rank $r$, its type is
\(t(\Lambda):=\dim_{k_\FLocal}(\Lambda^\vee/\Lambda).\)
Equivalently, a vertex lattice of type $t$ has fundamental invariants
$(0^{(r-t)},1^{(t)})$.

\item For an
$O_\FLocal$-basis $\{x_1,\ldots,x_r\}$ of $L$, define
\[
\det(L):=\det\bigl((x_i,x_j)_{i,j=1}^r\bigr)
\in\FLocal^\times/(O_\FLocal^\times)^2
\]
and
\[
\disc(L):=(-1)^{\binom r2}\det(L)
\in\FLocal^\times/(O_\FLocal^\times)^2.
\]
Then $\val(L)=\val(\disc(L))$.

\item Let
\[
\chi=\left(\frac{\cdot}{\varpi}\right)_\FLocal:
\FLocal^\times/(\FLocal^\times)^2\longrightarrow\{\pm1,0\}
\]
be the quadratic residue symbol, and put
\(\chi(L):=\chi(\disc(L))\).
Then $\chi(L)=0$ exactly when $\val(L)$ is odd.  If $L$ is self-dual,
its discriminant is a unit and
\[
\chi(L)\in\{\pm1\}
\simeq O_\FLocal^\times/(O_\FLocal^\times)^2.
\]

\item Let $V$ be a nondegenerate
quadratic space over $\FLocal$ of dimension
$\dimofV$.  We use the analogous invariants
\[
\det(V)\in\FLocal^\times/(\FLocal^\times)^2,\qquad
\disc(V)\in\FLocal^\times/(\FLocal^\times)^2,\qquad
\chi(V)\in\{\pm1,0\}.
\]
Thus $\chi(V)=+1$ if and only if $\disc(V)$ is a square.  For every
full-rank $O_\FLocal$-lattice $L\subset V$, one has
\(\chi(L)=\chi(V)\).
If $\{e_1,\ldots,e_m\}$ is an orthogonal basis and
\(u_i=(e_i,e_i)\), define the Hasse invariant by
\[
\Has(V):=\prod_{1\le i<j\le \dimofV}(u_i,u_j)_\FLocal\in\{\pm1\},
\]
where $(\ ,\ )_\FLocal$ denotes the Hilbert symbol.

\item Quadratic spaces over $\FLocal$ are
classified by their
dimension, discriminant, and Hasse invariant.  For $\dimofV\ge3$, the
discriminant and Hasse invariant are arbitrary.  For $\dimofV=2$, the pair
\(\disc(V)=1\), \(\Has(V)=-1\) does not occur; for $\dimofV=1$, the Hasse
invariant is necessarily $+1$.  A quadratic space $V$ admits a self-dual
lattice if and only if
\[
\Has(V)=+1
\qquad\text{and}\qquad
\chi(V)\ne0.
\]

\item\label{defn:integral-cone} Define
\[
V^\circ:=\{x\in V:(x,x)\in O_\FLocal\}.
\]

\item For
$\varepsilon\in\{\pm1\}$ and $r\ge1$, let $H_r^\varepsilon$ denote
the self-dual $O_\FLocal$-lattice of rank $r$ with
$\chi(H_r^\varepsilon)=\varepsilon$.

\item For a nondegenerate
quadratic space $U_0$ over $k_\FLocal$, its discriminant
character $\chi$ is the quadratic character of $k_\FLocal^\times$ evaluated
on $(-1)^{\binom{\dim U_0}{2}}\det(U_0)$; we use the empty determinant
$1$ and set $\chi(0)=1$ for the zero-dimensional space.
\item For an integral lattice $L$, we use
$\operatorname{sgn}(L)$ and $\operatorname{sgn}_m(L)$ from
\cite[\S3.2--\S3.4]{LZ-orth}. Set
$L_{k_\FLocal}=L\otimes_{O_\FLocal}{k_\FLocal}$ and
$U_0=L_{k_\FLocal}/\operatorname{rad}(L_{k_\FLocal})$, where
$\operatorname{rad}(L_{k_\FLocal})$ is the radical of its bilinear form. Then
\begin{equation*}
	\operatorname{sgn}(L)=
	\begin{cases}
		\chi(U_0),& \dim U_0\equiv0\pmod 2,\\
		0,& \dim U_0\equiv1\pmod 2,
	\end{cases}
\end{equation*}
and
\begin{equation}\label{eq:LZ-sign-convention}
	\operatorname{sgn}_m(L)=
	\begin{cases}
		\chi(U_0),& \dim U_0\equiv m\pmod 2,\\
		0,& \dim U_0\not\equiv m\pmod 2.
	\end{cases}
\end{equation}
\end{enumerate}

\subsubsection{Notation on global quadratic spaces}
\label{sec:notat-quadr-spac}

\begin{enumerate}[(G1)]
\setlength{\itemsep}{3pt}
\item Let $\FGlobal$ be a totally real
number field and
\(\BA_\FGlobal=\BA\otimes_\Q \FGlobal\) its ring of adeles.  Let $\Vincoherent/\BA_\FGlobal$ be an adelic
quadratic space of rank \(\dimofV=\nInLZorth+1\ge3\), with local components
\[
\Vincoherent_v:=\Vincoherent\otimes_{\BA_\FGlobal}\FGlobal_v.
\]
We call $\Vincoherent$ \emph{totally positive definite} if
$\Vincoherent_v$ has signature $(m,0)$ at every real place.

\item The adelic determinant, discriminant, and Hasse invariant are defined by
\begin{align*}
	\det(\Vincoherent)&:=(\det(\Vincoherent_v))_v
	\in\BA_\FGlobal^\times/(\BA_\FGlobal^\times)^2,\\
	\disc(\Vincoherent)&:=(\disc(\Vincoherent_v))_v
	\in\BA_\FGlobal^\times/(\BA_\FGlobal^\times)^2,\\
	\Has(\Vincoherent)&:=\prod_v\Has(\Vincoherent_v)\in\{\pm1\}.
\end{align*}
\item The quadratic character associated with $\Vincoherent$ is
\[
\chi_\Vincoherent((a_v)_v)
:=
\prod_v(a_v,\disc(\Vincoherent_v))_{\FGlobal_v},
\qquad
(a_v)_v\in\BA_\FGlobal^\times.
\]
If the class $\disc(\Vincoherent)\in\BA_\FGlobal^\times/(\BA_\FGlobal^\times)^2$ lies in the image of
\[
\FGlobal^\times/(\FGlobal^\times)^2
\longrightarrow
\BA_\FGlobal^\times/(\BA_\FGlobal^\times)^2,
\]
then $\chi_\Vincoherent$ factors through
$\BA_\FGlobal^\times/\FGlobal^\times$.
\item Assume that $\disc(\Vincoherent)$ lies in the image of
$\FGlobal^\times/(\FGlobal^\times)^2$. We call $\Vincoherent$ \emph{coherent} when
$\Has(\Vincoherent)=+1$ and \emph{incoherent} when
$\Has(\Vincoherent)=-1$.  In the coherent case, $\Vincoherent$ is isometric to the adelization of a
quadratic space over $\FGlobal$.
In the incoherent case, the product formula prevents such a global
realization.

\item Let
$\Vincoherent/\BA_\FGlobal$ be an incoherent quadratic space of dimension
$\dimofV$, positive definite at every real place.
At the distinguished real place $\RealPlace_0$, define the quadratic space
$\Vcoherent=\Vincoherent^{\RealPlace_0}$ over $\FGlobal$ by requiring signature
$(\dimofV-2,2)$ there, signature $(\dimofV,0)$
at every other real place, and the same finite localizations as
$\Vincoherent$, i.e., $
\Vcoherent_v\simeq\Vincoherent_v$ for $v<\infty$.  These local conditions satisfy the product
formula and determine a unique global isometry class. 
We thus have a finite adelic isometry
\begin{equation}\label{eq:Vincoherent-Vcoherent-bridge}
\Vcoherent_f:=\Vcoherent\otimes_\FGlobal\BA_{\FGlobal,f}
\simeq
\Vincoherent_f:=\Vincoherent\otimes_{\BA_\FGlobal}\BA_{\FGlobal,f}.
\end{equation}
We use this isometry to identify their finite localizations and their
finite adelic Schwartz spaces throughout.
\item For a finite place
$v$, we write $
	\VspecialGlobal^{(v)}$ 
for the global quadratic space with the same dimension and discriminant,
the opposite Hasse invariant at $v$, and the same local isometry class
as $\Vincoherent$ at every other place. In particular it is positive
definite at all real places. These local conditions satisfy the product
formula and determine a unique global isometry class.

\item For a nonsingular symmetric
matrix
$T\in\Sym_\nInLZorth(\FGlobal)$, define
\begin{equation*}
	\Diff(T,\Vincoherent)
	:=
	\{v:\ \Vincoherent_v\text{ does not represent }T_v\}.
\end{equation*}
When $v$ is finite, this is equivalent to the Hasse
invariant criterion \eqref{eq:diff}.

\item At a finite place $v$,
a local Schwartz function $\phi_v\in\sS(\Vcoherent_v^r)$ is called
\emph{non-degenerate} if every tuple in
its support has a nonsingular Gram matrix:
\begin{equation}\label{eq:Schwartz-nondegenerate}
	\operatorname{supp}(\phi_v)
	\subseteq
	\{\mathbf x:\det T(\mathbf x)\ne0\}.
\end{equation}
A factorizable global Schwartz function is non-degenerate at two
places if \eqref{eq:Schwartz-nondegenerate} holds at two distinct
finite places. All Schwartz functions on the finite adelic space
$(\Vcoherent\otimes_\FGlobal\BA_{\FGlobal,f})^r$ are taken to be
$\ov{\Q}$-valued. We denote their space by
$\sS\bigl((\Vcoherent\otimes_\FGlobal\BA_{\FGlobal,f})^r\bigr)$.

\item Fix a positive integer $\fkd$ and a lattice
$\VlattiveGlobal\subseteq \Vcoherent$.
	We say that a Schwartz function $\phi\in \sS((\Vcoherent\otimes_\FGlobal \BA_{\FGlobal,f})^{r})^K$ is \emph{admissible} if it is of the form
	$\phi_{\fkd}\charfun_{(\VlattiveGlobal\otimes \widehat{\Z}^\fkd)^{r}}$ for some $\phi_{\fkd}\in \sS(\Vcoherent_{\fkd}^{r})^{K_\fkd}$. Here $\widehat{\Z}^\fkd=\prod_{\ell\nmid \fkd}\Z_\ell\subset \BA^{f,\fkd}$, $\Vcoherent_{\fkd}=\prod_{v\mid \fkd}\Vcoherent_{v}$.

\end{enumerate}

\subsubsection{Notation on \texorpdfstring{$K$}{K}-groups and arithmetic degrees}
\label{subsec:nota-K-gp}
\begin{itemize}
\setlength{\itemsep}{3pt}
\item We follow
\cite[\S2.3]{LZ-orth} and \cite[Appendix~B]{Zhang2019WeilRA}.
Let $X$ be a noetherian scheme or formal scheme, and let
$i:Z\hookrightarrow X$ be a closed immersion with ideal sheaf $\mathcal J$.
Write $\operatorname{Coh}_Z(X)$ for the category of coherent
$\CO_X$-modules locally annihilated by a power of $\mathcal J$.
For a finite morphism $i:Z\to X$, define $\operatorname{Coh}_Z(X)=\operatorname{Coh}_{\im(Z\to X)}(X)$.
\item Let $\operatorname{Perf}_Z(X)$ be the
full triangulated subcategory of
perfect $\CO_X$-complexes whose cohomology sheaves lie in
$\operatorname{Coh}_Z(X)$, and put
\[
K_0^Z(X):=K_0\bigl(\operatorname{Perf}_Z(X)\bigr),\qquad
K'_0(Z):=K_0\bigl(\operatorname{Coh}(Z)\bigr).
\]
Here $K_0$ denotes the Grothendieck group.

\item The subgroup $\Fil_n K'_0(Z)$ is
generated by classes of bounded coherent complexes $E$ on $Z$ satisfying
$\dim\supp\CH^i(E)\le n$ for every $i$. The subgroup $\Fil^n K_0^Z(X)$ is
generated by classes of $E\in\operatorname{Perf}_Z(X)$ satisfying
$\operatorname{codim}_X\supp\CH^i(E)\ge n$ for every $i$.
Here $\CH^i(E)$ denotes the $i$-th cohomology sheaf of $E$.
For a bounded coherent complex, $[E]=\sum_i(-1)^i[\CH^i(E)]$ in $K'_0(Z)$.

If $X$ is a regular scheme of pure dimension that is either of finite type over a field
or proper over a Dedekind domain\footnote{So we have a good dimension theory for cycles.}, and $Z\subseteq X$ is a closed subscheme, 
coherent resolution and d\'evissage give isomorphisms
\begin{equation*}
\begin{aligned}
K_0^Z(X)&\xrightarrow{\ \sim\ }K'_0(Z),\\
\Fil^n K_0^Z(X)&\xrightarrow{\ \sim\ }\Fil_{\dim X-n} K'_0(Z).
\end{aligned}
\end{equation*}
The cycle map on the usual Chow group gives
\[
\Ch_n(Z)_\Q\xrightarrow{\ \sim\ }
\Fil_nK'_0(Z)_\Q/\Fil_{n-1}K'_0(Z)_\Q,
\qquad [C]\longmapsto[\CO_C].
\]
Here $C\subseteq Z$ is an integral closed subscheme of dimension $n$;
see \cite[\href{https://stacks.math.columbia.edu/tag/0FEW}{Tag~0FEW}]{stacks-project}.
\item For schemes $Z_1$ and $Z_2$ finite over a regular scheme $X$,
the derived tensor product induces
\[
\Fil^n K_0^{Z_1}(X)_\Q\times\Fil^m K_0^{Z_2}(X)_\Q
\longrightarrow\Fil^{n+m} K_0^{Z_1\times_X Z_2}(X)_\Q,
\]
where the subscript $\Q$ denotes tensoring with $\Q$.

\item Let $R$ be a complete discrete
valuation ring with uniformizer $\pi$.
For a bounded $R$-complex $C$ with finite-length cohomology, set
\[
\chi_R(C):=\sum_j(-1)^j\operatorname{length}_R H^j(C).
\]
Suppose that $X$ is over $\Spec R$ or $\Spf R$, and that
$Z\hookrightarrow X$ is an ordinary scheme proper over $\Spec R/\pi^a$ for some $a\geq1$.
For $E\in\operatorname{Perf}_Z(X)$, $R\Gamma(X,E)$ has bounded finite-length cohomology and we define
\[
\chi_X:K_0^Z(X)\longrightarrow\Z,\qquad
[E]\longmapsto\chi_R\bigl(R\Gamma(X,E)\bigr).
\]

\item For a finite place $v$ of $\FGlobal$, put
$R_v:=O_{\breve{\FGlobal_v}}$ and
$q_v:=\#(O_{\FGlobal_v}/\mathfrak m_v)$, where $\mathfrak m_v$ is the
maximal ideal of $O_{\FGlobal_v}$.
For a bounded $R_v$-complex $C_v$ with finite-length cohomology, set
\begin{equation}\label{eq:deg-v-convention}
	\deg_v(C_v)
	:=
	\sum_i(-1)^i
	\operatorname{length}_{R_v}H^i(C_v).
\end{equation}
For a perfect complex on a scheme or formal scheme, $\deg_v$ is computed
by base change to $R_v$ along the chosen place above $v$, followed by
derived global sections and $\chi_{R_v}$, whenever the cohomology is bounded
and of finite length. We extend $\deg_v$ linearly to $K$-classes with coefficients.
\end{itemize}

\subsubsection{Notation on Hilbert modular forms}
\label{subsec:nota-hilbert}
\begin{enumerate}[(H1)]
\setlength{\itemsep}{3pt}
\item Let $\FGlobal/\Q$ be a totally real
field with $\FGlobal\ne\Q$. Fix a nontrivial global additive character
\[
\psi:\FGlobal\backslash\BA_{\FGlobal}\longrightarrow\BC^\times.
\]
For the globalizations used below, its components at the distinguished
$p$-adic places will be chosen unramified.  For every real place \(\RealPlace\), write
\[
\psi_\RealPlace(x)=\exp(2\pi i\kappa_\RealPlace x),
\qquad \kappa_\RealPlace>0,
\]
and use the same exponential formula to evaluate \(\psi_\RealPlace\) on complex
arguments.

\item Set
\[
\mathbb H_{\FGlobal}
:=\prod_{\RealPlace\mid\infty}\mathbb H
=\left\{
\tau=(\tau_\RealPlace)_{\RealPlace\mid\infty}:
\tau_\RealPlace=x_\RealPlace+i y_\RealPlace,\ y_\RealPlace>0
\right\}.
\]
\item For an analytic genus-one index
\(a\in\FGlobal\) and
\(\tau\in\mathbb H_{\FGlobal}\), our Fourier monomial is
\begin{equation*}
	q^a
	:=\psi_\infty(a\tau)
	=\prod_{\RealPlace\mid\infty}\psi_\RealPlace(\RealPlace(a)\tau_\RealPlace)
	=\prod_{\RealPlace\mid\infty}
	\exp\!\left(2\pi i\kappa_\RealPlace\,\RealPlace(a)\tau_\RealPlace\right).
\end{equation*}
In particular, a geometric norm \(\xi=(x,x)=2Q(x)\) contributes the
monomial \(q^{\xi/2}\).  
\item At the chosen cusp $\infty$, write
\[
f(\tau)=\sum_{a\in\FGlobal}c(a,y)q^a
\]
for a genus-one Hilbert modular form $f$, where
$y=(y_\RealPlace)_{\RealPlace\mid\infty}$. Its $a$-th Fourier coefficient is
$[f]_{q^a}:=c(a,y)$; it is independent of $y$ when $f$ is holomorphic.
We write
\(\xi>0\), or equivalently \(\xi\in\FGlobal_{>0}\), when \(\xi\) is
totally positive. 

\item For each place $v$ of
$\FGlobal$, the local metaplectic group
$\Mp_2(\FGlobal_v)$ is a double cover of $\SL_2(\FGlobal_v)$
with kernel $\{(1,\varepsilon_v):\varepsilon_v\in\{\pm1\}\}$,
where $1$ denotes the identity in $\SL_2(\FGlobal_v)$.
\item The adelic metaplectic double cover is realized as
\begin{equation}\label{eq:adelic-metaplectic-cover}
\Mp_2(\BA_{\FGlobal})
\simeq
{\displaystyle\prod_{v\in \Sigma_\FGlobal}'\Mp_2(\FGlobal_v)}\Big/
\left\{
((1,\varepsilon_v))_v:
{
\varepsilon_v\in\{\pm1\}\text{ for every }v, 
\varepsilon_v=1\text{ for almost all }v, 
\prod_v\varepsilon_v=1
}
\right\}.
\end{equation}
The restricted product requires components to lie
in the canonical split images of $\SL_2(O_{\FGlobal_v})$ at almost all
finite places of odd residue characteristic. 

\item For $k\in\tfrac12\Z$,
weight $k$ refers to the character of the
inverse image of $\SO(2,\R)$ in $\Mp_2(\R)$ that sends a lifted
rotation of angle $\theta\in\R/4\pi\Z$ to $e^{ik\theta}$.

\item Write $(1,\varepsilon)$, with
$\varepsilon\in\{\pm1\}$, for the elements of
the kernel of $\Mp_2(\BA_{\FGlobal})\to\SL_2(\BA_{\FGlobal})$.
Under \eqref{eq:adelic-metaplectic-cover}, a kernel tuple
$((1,\varepsilon_v))_v$ has image $(1,\prod_v\varepsilon_v)$.
Every nontrivial local kernel element therefore maps to $(1,-1)$.
Parallel weight $k$ prescribes the covering-kernel character
\begin{equation}\label{eq:hilbert-covering-kernel-character}
f\bigl(g(1,\varepsilon)\bigr)=\varepsilon^{2k}f(g),
\qquad \varepsilon\in\{\pm1\}.
\end{equation}
This character is genuine for $k\in\tfrac12+\Z$; for integral $k$, it
is trivial and $f$ descends to $\SL_2(\BA_{\FGlobal})$.

\item The space
\[
\CA_\infty(\SL_2(\FGlobal)\backslash\Mp_2(\BA_\FGlobal),k)
\]
consists of smooth functions on $\Mp_2(\BA_\FGlobal)$ invariant under
left $\SL_2(\FGlobal)$ and right translation by a compact open subgroup of $\Mp_2(\BA_{\FGlobal,f})$, of parallel weight
$k$ under the inverse images of the real maximal compact groups, and
satisfying \eqref{eq:hilbert-covering-kernel-character}.
No finiteness condition under the center of the universal enveloping
algebra is imposed.

\item\label{H10} The subspace
\[
\CA_{\hol}(\SL_2(\FGlobal)\backslash\Mp_2(\BA_\FGlobal),k)\subset \CA_\infty(\SL_2(\FGlobal)\backslash\Mp_2(\BA_\FGlobal),k)
\]
consists of automorphic forms that are annihilated at every real place by 
$
\frac12\begin{psmallmatrix}i&1\\1&-i\end{psmallmatrix}
$ 
in the complexified Lie algebra of $\Mp_2(\BR)$.
This differential operator commutes with right $\Mp_2(\BA_{\FGlobal,f})$  translation.
By strong approximation, it is equivalent to require the holomorphy on $\Mp_2(\BA_{\FGlobal,\infty})\times\{(1)_{v<\infty}\}$. 
The Koecher principle for $\FGlobal\ne\Q$ gives
holomorphy at every cusp. 
For any fixed compact open subgroup of $\Mp_2(\BA_{\FGlobal,f})$, the corresponding holomorphic space is finite-dimensional over $\BC$.
The $\ov{\BQ}$-structure on the full holomorphic space is
$
\CA_{\hol}(\SL_2(\FGlobal)\backslash\Mp_2(\BA_\FGlobal),k)_{\ov{\BQ}}
$, 
the subspace consisting of forms whose Fourier coefficients all lie in
$\ov{\BQ}$.

\item\label{defn:Fourier-transform} For a place \(v\), a Haar measure
\(dm\) on \(\Vcoherent_v\), and \(\phi\in\sS(\Vcoherent_v)\), define
the \emph{Fourier transform} by
\[
\widehat\phi^{\psi_v}(b_1)
=\int_{\Vcoherent_v}\phi(b_2)
\psi_v\!\left((b_1,b_2)\right)
\,dm(b_2),\qquad b_1\in\Vcoherent_v.
\]
Then \(\widehat\phi^{\psi_v}\in\sS(\Vcoherent_v)\).
\item\label{defn:Haar-self-dual} A Haar measure \(dm\) on
\(\Vcoherent_v\) is \emph{self-dual} with respect to the pairing and
\(\psi_v\) if
$
\widehat{\widehat\phi}^{\psi_v}(x)=\phi(-x)$ for every $\phi\in\sS(\Vcoherent_v)$ and  $x\in\Vcoherent_v$. 
Such a Haar measure exists and is unique. At a finite place where
\(\psi_v\) is unramified, this is equivalent to
\(\vol(L)\vol(L^\vee)=1\) for every lattice
\(L\subseteq\Vcoherent_v\), with both volumes taken using \(dm\).

\item Let $X$ be a smooth scheme over an archimedean or non-archimedean local field $\FLocal$ with a top form $\alpha$. Then $\alpha$ defines a measure on $X(\FLocal)$; see \cite[Chapter~I, no.~3]{Weil1965SiegelFormula}. We call $\alpha$ the \emph{density} of this measure. 
\end{enumerate}

\subsubsection{Notation list}
\label{subsec:frequent-symbols}
\begin{itemize}
\setlength{\itemsep}{3pt}
	\item $\dimofV,\nInLZorth$: $\dimofV=\dim V$ and the integer $\nInLZorth=\dimofV-1$.
	\item $\FLocal,O_\FLocal,k_\FLocal,q,\varpi$: a local
	field, its ring of integers, its residue field, the cardinality of that
	residue field, and a uniformizer.
	
	\item $\valFLocal$: the valuation $\valFLocal:\FLocal\to \Z\cup \{+\infty\}$, normalized so that $\valFLocal(\varpi)=1$.  
	
	\item $\FGlobal,\BA_\FGlobal$: a totally real global field and its ring of adeles.
	
	\item $O_{\FGlobal}, O_{\FGlobal,(p)}, O_{\FGlobal,p}, O_{\FGlobal,v}$: the ring of integers of $\FGlobal$, $O_{\FGlobal,(p)}=O_{\FGlobal}\otimes_{\Z}\Z_{(p)}$, $O_{\FGlobal,p}=O_{\FGlobal}\otimes_{\Z}\Z_p$ and $O_{\FGlobal,v}$ is the completed local ring at a finite place $v$. 
	\item $\Sym_r(R)$: the $R$-module of symmetric $r\times r$ matrices,
	for a commutative ring $R$ and an integer $r\ge1$.
	\item $\Sym_r(\BR)_{>0}$: the cone of positive definite real
	symmetric matrices.
	For $T\in\Sym_r(\BR)$, we write $T>0$ when ${}^txTx>0$ for every
	nonzero column vector $x\in\BR^r$, and $T\ge0$ when ${}^txTx\ge0$
	for every column vector $x\in\BR^r$. Here ${}^tx$ denotes the transpose
	of $x$.
	\item $\Sym_r(F)_{>0}$: for a totally real number
	field $F$, the set of matrices $T\in\Sym_r(F)$ for which $\sigma(T)$
	is positive definite for every
	real embedding $\sigma:F\hookrightarrow\BR$. The matrix $\sigma(T)$
	is obtained by applying $\sigma$ to each entry of $T$. We write
	$T>0$ for this condition, and $T\ge0$ when $\sigma(T)$ is positive
	semidefinite for every real embedding $\sigma$.
	\item $\Vlocal,\VlatticeLocal$: the quadratic space defining the local
	GSpin datum and its fixed self-dual lattice.
	\item $\epsilon$: the local sign $\epsilon=\chi(\Vlocal)$. At a varying
	finite place, the corresponding sign is denoted by $\epsilon_v$.
	
	\item $\Vincoherent,\Vcoherent,\VspecialGlobal^{(v)}$: the incoherent
	adelic quadratic space, its nearby coherent space at the distinguished real
	place, and its nearby coherent space at a finite place $v$.
	\item $\Vspecial$: the quadratic space of special quasi-endomorphisms of
	the basic framing object on a GSpin cover.  It is the nearby space of $\Vlocal$ for
	$(\dimofV,\epsilon)\ne(2,+1)$; in the split binary case
	$(\dimofV,\epsilon)=(2,+1)$, we use $\Vspecial=0$.
	\item $T(\mathbf x)$: the Gram matrix of a tuple, defined in
	\eqref{eq:Gram-matrix-convention}.
	\item $\Int^\epsilon(\mathbf x),\pDen^\epsilon(\mathbf x)$: the derived
	degree on one connected component of the SO Rapoport--Zink space and the
	derivative of the normalized local density. For $\dimofV\ge3$, the
	first quantity equals one half of the total derived degree by
	\Cref{lem:SO-Frob-components}; see
	\Cref{defn:derived-special-cycles-RZ,defn:local-density}.
	These quantities depend only on $T(\mathbf x)$; we also
	write $\Int^\epsilon(T)$ and $\pDen^\epsilon(T)$.  

	\item $\operatorname{Orb}_T(\phi^S)$: the finite adelic orbital integral
	with the measure $dg_f^S$ of \Cref{conv:haar-global-constant}, where
	$S$ is the set of omitted finite places.  If $T=T(\mathbf x)$, the synonymous notation
	$\operatorname{Orb}_{\mathbf x}(\phi^S)$ is used when a representative tuple
	is part of the data.
	\item $K,\vol(K)$: a sufficiently small compact open level and its finite
	adelic Haar volume, normalized by
	\eqref{eq:GL-disintegration} and \Cref{defn:Haar-measure}.  This
	volume measures the finite adelic level; the height formula also uses an
	archimedean orbifold volume.
	\item $a_{\infty,\psi}^{\operatorname{tr}}$: the explicit positive
	archimedean factor in \Cref{lem:global-Haar-phase-a}.
	\item $S_\psi,\fkd$: the finite set of places where the global additive
	character is ramified and an even positive integer divisible by all rational
	primes designated as bad, respectively.  In the globalization used for the
	main theorem, $(\fkd,p)=1$ and $S_\psi$ lies over $\fkd$.
	\item $T_F$, $T_{O_F}$: the Weil restrictions $T_F=\Res_{F/\Q}\Gm$, $T_{O_F}=\Res_{O_F/\Z}\Gm$  for a global field
	$F/\Q$, and $T_F=\Res_{F/\Q_p}\Gm$, $T_{O_F}=\Res_{O_F/\Z_p}\Gm$  for a local field $F/\Q_p$.
	
	\item$\ShG_K, \Sh_K$: given a Shimura datum $(G,X)$ with level structure $K\subseteq G(\BA_f)$, we will abbreviate the Shimura variety $\ShG_K(G,X)$ and its integral model $\Sh_K(G,X)$ to $\ShG_K$ and $\Sh_K$.  
	\item $Z(G)$: For a group scheme $G$ we denote by $Z(G)$ the center of $G$. 
\end{itemize}
\endgroup

	\part{General theory}\label{part1}
	\section{Special cycles on Shimura varieties of abelian type}
\label{sec:special-cycles-abelian-type}
Throughout this section $p>2$. 
In this section we construct special endomorphisms and special cycles on
integral models of Shimura varieties of abelian type with parahoric level, assuming the local model is smooth.  
We use the KPZ integral models in \cite{KisinPappasZhou2024ParahoricII} which are canonical by
\cite[Theorem~B]{DanielsYoucis2025} and, in the Hodge-type case with a chosen Hodge embedding, by the
work of \cite{PappasRapoport2024PadicShtukas}.  They have local model diagrams, so their singularities are the same as those of the local models. 

The
construction of special endomorphisms on a Hodge-type integral model is modeled on \cite{MadapusiPera2016}, and we use the twisting of abelian schemes in \cite[\S4.4]{KisinPappas2018} to descend the special endomorphism sheaves to the abelian-type integral model. 

The verification of the properties of special endomorphisms relies on the crystalline local system associated to a representation. 
Since the reflex field may be ramified at $p$, we use prismatic--\'etale comparison, in the form
developed in \cites{du2022completed,guo2022prismatic,AnschutzLeBras2023,IKY26}, to construct crystalline realizations,  
rather than relying only on crystalline--de Rham comparison as in \cite{MadapusiPera2016}. 
Under \Cref{ass:d=1}, we show that the local deformation loci of non-zero special endomorphisms are relative (i.e. flat over the ring of integers of the
reflex field) Cartier divisors in \Cref{prop:tangent-space-divisor}. 

\subsection{Conditions on Shimura data}
We recall the Hodge embedding that will
be used. 

\begin{defn}\label{defn:fundamental-emb}
	Let $\FGlobal$ be a totally real field. For a Hodge-type Shimura datum $(G,X)$, let $G^c$ be the subgroup generated by
	$G^{\der}$ and the maximal compact subtorus of the center, as in
	\cite[\S7.2.4]{KisinPappasZhou2024ParahoricII}. We assume that
	\begin{equation}\label{eq:KPZ-dagger}
		(\dagger)\qquad
		G^c\simeq \Res_{\FGlobal/\Q}H^c
	\end{equation}
	for a reductive group $H^c/\FGlobal$. 
	Under \eqref{eq:KPZ-dagger}, a Hodge
	embedding
	\(\iota:(G,X)\longrightarrow (\GSp(W,\psi),S^{\pm})\)
	is \emph{fundamental} if $W$ is an $\FGlobal$-vector space and the restriction of
	$\iota$ to $G^c$ is obtained by restriction of scalars from an $\FGlobal$-morphism
	$H^c\to\GL_\FGlobal(W)$, following
	\cite[Definition~7.2.7]{KisinPappasZhou2024ParahoricII}.  Moreover, let $H'$ be the subgroup of $\GL_\FGlobal(W)$ generated by $H^c$ and the homotheties  $\FGlobal^\times$; we then have an inclusion $G \subseteq  G':=\Res_{\FGlobal/\Q}H'$, and the embedding $G \to \GSp(W)$ extends
	to an embedding $G'\to \GL(W)$.
	
	After changing the polarization within
	its natural cone, one may arrange
	\begin{equation}\label{eq:fundamental-trace-polarization}
		\psi=\Tr_{\FGlobal/\Q}\circ\Psi
	\end{equation}
	for an $\FGlobal$-bilinear alternating form $\Psi$ on $W$, and
	$H'$ acts through $\GSp_\FGlobal(W,\Psi)$; see
	the argument of \cite[Lemma~7.2.8]{KisinPappasZhou2024ParahoricII}, which
	uses \cite[\S\S1.1.10 and~1.1.18(a)]{Deligne1979Shimura}.  
\end{defn}

Now we introduce  the assumptions used in the construction of integral models. 
\begin{defn}
	We recall the definition of $R$-smoothness. For a torus $T/\Q_p$,
	choose a finite Galois splitting field $L/\Q_p$, put
	$T_1=\Res_{L/\Q_p}T_L$, and require the Zariski closure of $T$ in the
	locally finite type N\'eron model $T_1^{\mathrm{lft}}$ to be smooth over
	$\Z_p$; a reductive group is $R$-smooth if one (equivalently every)
	maximally $\breve{\Q}_p$-split maximal torus is $R$-smooth. By \cite[Proposition~2.1.5(2)]{KisinPappasZhou2024ParahoricII}, the $R$-smoothness condition is stable under extensions of tori.
\end{defn}

\begin{ass}\label{ass:Shimura-data}
	Let $(G_2,X_2)$ be an abelian-type Shimura datum, and let
	$\CG_2$ be the Bruhat--Tits stabilizer attached to a point of
	$\CB(G_2,\Q_p)$ which is generic in its facet. Let
	\(\alpha:(G,X)\longrightarrow(G_2,X_2)\)
	be a central epimorphism of Shimura data,
	inducing an isomorphism
	$(G^{\ad},X^{\ad})\simeq(G_2^{\ad},X_2^{\ad})$.
	Write $\reflexHodge=E(G,X)$ and $\reflexAbel=E(G_2,X_2)$ for the reflex fields.
	Fix a place $v\mid p$ of $\reflexHodge$ and write $v_2$ for its restriction to $\reflexAbel$.
	Let $k=\kappa(v)$, with a fixed embedding $k\hookrightarrow\barFp$.
	In the global integral-model constructions, $O_{\reflexHodge}$ and $O_{\reflexAbel}$
	mean $O_{\reflexHodge,(v)}$ and $O_{\reflexAbel,(v_2)}$, respectively.
	Their completed local rings are denoted by $O_{\reflexHodge,v}$ and $O_{\reflexAbel,v_2}$;
	we also write $\widehat O_{\reflexHodge}=O_{\reflexHodge,v}$ and
	$\widehat O_{\reflexAbel}=O_{\reflexAbel,v_2}$.
	When both places are denoted by $v$, the reflex field specifies the place.
	For a local reflex field, $O_E$ denotes its complete ring of integers.
	We assume that $(G,X)$ is of Hodge type and has a 
    fundamental Hodge embedding
    \begin{equation}\label{eq:F-linear-Hodge-embedding}
        \iota:(G,X)\hookrightarrow(\GSp_\Q(W,\psi),S^{\pm}).
    \end{equation}
	Choose an $O_\FGlobal$-lattice $W_{O_\FGlobal}\subset W$.  We write $W_\Z$ for its underlying
	$\Z$-lattice and set
	\(W_{\Z_{(p)}}:=W_\Z\otimes_\Z\Z_{(p)}\) and
	\(W_{\Z_p}:=W_\Z\otimes_\Z\Z_p\). 
	Throughout, we require every integral Hodge realization used to define
	or compare special endomorphisms below to
	be self-dual at \(p\) for its \(\FGlobal\)-bilinear form \(\Psi\):
	\[
	W_{O_{\FGlobal,(p)}}:=
	W_{O_\FGlobal}\otimes_\Z\Z_{(p)}
	=W_{O_{\FGlobal,(p)}}^\vee:=
	\{w\in W:\Psi(w,W_{O_{\FGlobal,(p)}})
	\subseteq O_{\FGlobal,(p)}\}.
	\]
	Put
	$\kernelalphaG:=\ker(\alpha)$.  We require:
	\begin{enumerate}
		\item\label{ass:Shimura-data:kottwitz}
		The Shimura datum $(G_2,X_2)$ satisfies the Kottwitz condition; that is,
		$\rank_\Q Z(G_2)=\rank_\R Z(G_2)$.
		
		\item\label{ass:Shimura-data:reflex}
		The canonical inclusion $\reflexAbel\subseteq\reflexHodge$ is unramified at all
		places above $p$.
		
		\item\label{ass:Shimura-data:pi_1} $p\nmid |\pi_1(G^{\der})|$. 
		
		\item\label{ass:Shimura-data:level}
		$K_p=\CG(\Z_p)$ and $K_{2,p}=\CG_2(\Z_p)$, where $\CG$ is the
		Bruhat--Tits stabilizer attached to a point $\mathbf x$ in the extended building $\fkB(G, \Q_p)$, and $\CG_2$ is the
		Bruhat--Tits stabilizer attached to the image $\mathbf x_2\in \fkB(G_2, \Q_p)$.

		\item\label{ass:Shimura-data:group-closed} The Hodge embedding
		induces a closed immersion
		$\CG\hookrightarrow\GL_{\Z_p}(W_{\Z_p})$.

		\item\label{ass:Shimura-data:reductive}
		There is a reductive group scheme
		$\mathcal H/\mathcal O_{\FGlobal,(p)}$ such that
		$G_2\simeq\Res_{\FGlobal/\Q}\mathcal H_{\FGlobal}$ and, compatibly with this identification,
		$\CG_2\simeq
		(\Res_{\mathcal O_{\FGlobal,(p)}/\Z_{(p)}}\mathcal H)_{\Z_p}$.
		
		\item\label{ass:Shimura-data:R-smooth}
		Let $\chi:\GSp_\Q(W,\psi)\to\BG_m$ be the similitude character and
		$\kernelalphaG^c:=\ker((\chi\circ\iota)|_{\kernelalphaG})$.
		These tori fit into the exact sequence over $\Q$
		\[1\to \kernelalphaG^c\to \kernelalphaG\xrightarrow{\chi\circ \iota} \BG_{m}\to 1.\]
		We require that $\kernelalphaGp^c=\kernelalphaG^c\times_{\mathbb Q}\mathbb Q_p$ is an $R$-smooth torus over $\mathbb Q_p$, and
		$X_*(\kernelalphaGp^c)_I$ is torsion-free, where
		$I\subset\Gal(\ov{\Q}_p/\Q_p)$ is the inertia group. This implies $\kernelalphaGp$ is an $R$-smooth torus and $X_*(\kernelalphaGp)_I$ is torsion-free\footnote{The inertia action on the cocharacter lattices factors through a finite group. Its first homology with coefficients in $\Z$ is torsion, so the connecting homomorphism to $X_*(\kernelalphaGp^c)_I$ vanishes.}.

		\item\label{ass:local-model-smooth}
		If $\{\mu\}$ denotes the $p$-adic conjugacy class of cocharacters
		attached to $X$, then the canonical local model
		$\BM^{\loc}_{\CG,\{\mu\}}$ is smooth over 
		$O_{\reflexHodge,v}$.
		
		\item\label{ass:Shimura-data:good}
		The $p$-adic completion of \eqref{eq:F-linear-Hodge-embedding} is a
		\emph{good} integral Hodge embedding in the sense of
		\cite[Definition 3.4.4]{KisinPappasZhou2024ParahoricII}; that is, it
		induces a closed immersion on local models.  By
		\cite[Corollary~5.3.16]{KisinPappasZhou2024ParahoricII} and (\ref{ass:local-model-smooth}), it is then a
		\emph{very good} integral Hodge embedding in the sense of
		\cite[Definition~5.2.5]{KisinPappasZhou2024ParahoricII};
	\end{enumerate}
\end{ass}
The base change $O_{\reflexAbel}\to O_{\reflexHodge}$ is faithfully flat
and, by condition (\ref{ass:Shimura-data:reflex}), essentially \'etale.
The special-endomorphism sheaf constructed over $O_{\reflexHodge}$ is
descended along this map in \Cref{rmk:descent-base}.

The Kottwitz condition (\ref{ass:Shimura-data:kottwitz}) allows us to associate an automorphic local system to a representation of $G_2$ (cf. \cite[II, Proposition~3.3; III, Remark~6.1]{Milne1990CanonicalModels}). 

In what follows, we write $\centralHodgeG, \centralAbelG, \kernelalphaG$ for the center of $G$, the center of $G_2$ and the kernel of $\alpha$.  We define $\centralHodge, \centralAbel, \kernelalpha$ to be the kernels of the morphisms of parahoric group schemes $\CG\to \CG^\ad$, $\CG_2\to \CG^\ad$ and $\CG\to \CG_2$, respectively. 

\begin{defn}\label{defn:strongly-admissible}
	Following \cite[Definition~2.1.6]{KisinZhou2025StronglyCompatible}, a
	Hodge-type Shimura datum with Bruhat--Tits stabilizer scheme $(G,X,\CG)$ is
	called \emph{strongly admissible} if
	\begin{enumerate}
		\item $\CG$ has connected special fiber, equivalently $\CG=\CG^\circ$;
		\item $G_{\Q_p}$ is $R$-smooth in the sense of
		\cite[\S2.1.4]{KisinPappasZhou2024ParahoricII}.
	\end{enumerate}
\end{defn}
\begin{lem}\label{lem:Shimura-data-strongly-admissible}
	Under \Cref{ass:Shimura-data}(\ref{ass:Shimura-data:reductive})--
	(\ref{ass:Shimura-data:R-smooth}), both stabilizer group schemes $\CG_2$ and
	$\CG$ have connected special fiber, and $G_{\Q_p}$ is $R$-smooth.  In
	particular
	$\CG_2=\CG_2^\circ,$ $\CG=\CG^\circ,$
	so $(G,X,\CG)$ is strongly admissible. 
\end{lem}

\begin{proof}
	For $\CG_2$, put $A:=\mathcal O_\FGlobal\otimes_{\Z}\ov{\BF}_p$.
	Under
	\Cref{ass:Shimura-data}(\ref{ass:Shimura-data:reductive}), the geometric
	special fiber of $\CG_2$ is
	$\Res_{A/\ov{\BF}_p}(\mathcal H_A)$.  The morphism induced by
	$A\to A_{\red}$,
	\[
	\Res_{A/\ov{\BF}_p}(\mathcal H_A)
	\longrightarrow
	\Res_{A_{\red}/\ov{\BF}_p}(\mathcal H_{A_{\red}}),
	\]
	has connected unipotent kernel, while the target is a product of connected
	reductive groups because $A_{\red}$ is a finite product of copies of
	$\ov{\BF}_p$.  Hence the geometric special fiber of $\CG_2$ is connected.
	Compare \cite[\S4.2]{HainesRicharz2020}.
	
	For $\CG$, \cite[\S1.1.2]{KisinPappas2018} identifies
	\[
	\CG^\circ(\breve{\Z}_p)
	=\CG(\breve{\Z}_p)\cap\ker(\widetilde\kappa_G),
	\]
	where $\widetilde\kappa_G$ is the Kottwitz homomorphism.  The central exact
	sequence
	\(1\to \kernelalphaG\to G\xrightarrow{\alpha}G_2\to1\)
	gives
	\(0\to X_*(\kernelalphaG)\to\pi_1(G)\to\pi_1(G_2)\to0.\)
	After inertia coinvariants, the possible kernel of
	$X_*(\kernelalphaG)_I\to\pi_1(G)_I$ is torsion.  The torsion-freeness in
	\Cref{ass:Shimura-data}(\ref{ass:Shimura-data:R-smooth}) therefore makes this
	map injective. We thus have an exact
	sequence
	\(0\to X_*(\kernelalphaG)_I\to\pi_1(G)_I\to\pi_1(G_2)_I\to0.\)
	
	Since $\CG_2=\CG_2^\circ$, functoriality of Kottwitz maps shows
	that $\widetilde\kappa_G(\CG(\breve\Z_p))$ lies in the torsion-free kernel of
	$\pi_1(G_{\Q_p})_I\to\pi_1(G_{2,\Q_p})_I$; on the other hand it is finite because it factors
	through $\CG/\CG^\circ$.  Hence it vanishes and $\CG=\CG^\circ$.
	
	Finally, let $T$ be the centralizer of a maximal $\breve\Q_p$-split torus of
	$G$, and let $T_2$ be the corresponding centralizer in $G_2$.  The group
	$T_2$ is quasi-tame under
	\Cref{ass:Shimura-data}(\ref{ass:Shimura-data:reductive}), hence $R$-smooth by
	\cite[Proposition~2.1.5(1)]{KisinPappasZhou2024ParahoricII}.  From
	$1\to \kernelalphaG\to T\to T_2\to1$ and the fact that the $R$-smoothness condition is stable under extensions of tori, $T$ is
	$R$-smooth.  Thus $G$ is $R$-smooth.
\end{proof}

\begin{rmk}\label{rmk:Shimura-data-strongly-admissible}
	Moreover, $G_{2,\Q_p}$ (hence $G_{\Q_p}$) is quasi-split. So the assumptions of
	\cite[Theorem~2.2.7]{KisinZhou2025StronglyCompatible} are satisfied. Thus every isogeny class on $\Sh_{K_pK^p}(G,X)_{\barFp}$, defined in \Cref{defn:Hodge-integral-model}, contains the reduction of a special point of $\ShG_{K_pK^p}(G,X)$. 
\end{rmk}

\begin{rmk}\label{rmk:isogeny}
	The central epimorphism $\alpha:G\to G_2$ is stronger than the minimal
	abelian-type input used in
	\cite[\S7.2]{KisinPappasZhou2024ParahoricII}, where a central isogeny on
	derived groups and an isomorphism on adjoint data suffice.  We retain this
	stronger form because the special representation below is a $G_2$-
	representation and we want to regard it, without further notation, as a
	$G$-representation.  The orthogonal applications in
	\Cref{subsec:orth-parahoric-models} have such a central map explicitly.
\end{rmk}

\begin{rmk}\label{rmk:R-smooth-seq}
	By \cite[Proposition~2.4.13]{KisinZhou2025}, the $R$-smoothness of
	$\kernelalphaG$ gives a short exact sequence of smooth $\Z_p$-group schemes
	\(1\to\kernelalpha\to\CG\to\CG_2\to1,\)
	where $\kernelalpha$ is the smooth Zariski closure of $\kernelalphaG$.
	We glue this sequence with
	\(1\to \kernelalphaG\to G\to G_2\to1\)
	away from $p$ and use the same symbols for the resulting smooth
	$\Z_{(p)}$-models.  Thus, in particular,
	\begin{equation*}
		1\longrightarrow\kernelalpha\longrightarrow\CG
		\longrightarrow\CG_2\longrightarrow1
	\end{equation*}
	is fppf exact over $\Z_{(p)}$ and has the original parahoric sequence as
	its base change to $\Z_p$.  Also, since the integral
	Hodge embedding is a closed immersion, the group $\CG$ can be cut out inside
	$\GL(W_{\Z_{(p)}})$ by a finite collection of tensors
	$(s_\alpha)_{\alpha\in I}$; see
	\cite[Proposition~1.3.2]{Kisin2010}.
\end{rmk}

We have the following variant of \cite[Lemma 2.2.6]{Kisin2010}. It is
slightly weaker because we replace $H(\FGlobal)^+$ by $H(\FGlobal)_+$, the
preimage of the identity component $H^\ad(\FGlobal_\infty)^+$, but it is sufficient for the
applications.
\begin{prop}\label{prop:Kisin2010-2.2.6}
	Let $v_1, \dots, v_s$ be finite places of $\FGlobal$. 
	Let $\CH$ be a connected reductive group over $O_{\FGlobal,(v_1 \dots v_s)}$ with generic fiber $H/\FGlobal$.
	Then
	\[
	\left(\prod_{i=1}^s\CH(O_{\FGlobal,v_i})\right)\cdot H(\FGlobal)_+
	=
	\prod_{i=1}^s H(\FGlobal_{v_i}),
	\]
	where $H(\FGlobal)_+=H(\FGlobal)\cap H(\FGlobal_{\R})_+$ is embedded diagonally in
	$\prod_{i=1}^s H(\FGlobal_{v_i})$, and $H(\FGlobal_{\R})_+$ is the preimage of
	$H^\ad(\FGlobal_{\R})^+$.
\end{prop}

\begin{proof}
	We follow the argument of \cite[Lemma~2.2.6]{Kisin2010}, simultaneously
	at all places \(v_i\).
	
	Let
	\[
	\rho:\widetilde H\longrightarrow H^{\der}\hookrightarrow H
	\]
	be the simply connected cover of the derived group, and put
	\[
	N_v:=
	\rho\left(\prod_{i=1}^s\widetilde H(\FGlobal_{v_i})\right)
	\subset \prod_{i=1}^sH(\FGlobal_{v_i}).
	\]
	By weak approximation for semisimple simply connected groups,
	\(N_v\subset \left(\prod_{i=1}^s\CH(O_{\FGlobal,v_i})\right)H(\FGlobal)_+.\)
	Indeed, given \(x=(x_{v_i})_{i=1}^s\in
	\prod_{i=1}^s\widetilde H(\FGlobal_{v_i})\), choose
	\(y\in\widetilde H(\FGlobal)\) simultaneously sufficiently close to \(x_{v_i}\)
	for every \(i\), so that $\rho(x_{v_i}y^{-1})\in\CH(O_{\FGlobal,v_i})$ for every $i$.
	For every real place \(w\) of \(\FGlobal\), the image of
	\(\widetilde H(\FGlobal_w)\) in \(H^{\ad}(\FGlobal_w)\) is contained in
	\(H^{\ad}(\FGlobal_w)^\circ\); hence \(\rho(y)\in H(\FGlobal)_+\).
	As the quotient \(\left(\prod_{i=1}^sH(\FGlobal_{v_i})\right)/N_v\) is abelian, it suffices to show that $\prod_{i=1}^s\CH(O_{\FGlobal,v_i})$ and $H(\FGlobal)_+$ generate the group $\prod_{i=1}^sH(\FGlobal_{v_i})$. 
	
	Next, we choose a maximal \(\FGlobal\)-torus \(T\subset H\) with suitable
	local properties.  For every \(v=v_i\), choose a maximal
	\(\FGlobal_v\)-torus \(T_v\subset H_{\FGlobal_v}\) which extends to a maximal
	\(O_{\FGlobal,v}\)-torus and contains a maximal \(\FGlobal_v\)-split torus of
	\(H_{\FGlobal_v}\).  For every real place \(w\), choose a fundamental (i.e. maximally compact)  maximal torus \(T_w\subset H_{\FGlobal_w}\); its image in \(H^{\ad}(\FGlobal_w)\) is contained
	in \(H^{\ad}(\FGlobal_w)^\circ\).
	
	The scheme of maximal tori in $H$ is rational over $\FGlobal$ and hence satisfies weak approximation; see the proof of 
	\cite[Lemma 2]{Guo2022SemilocalDedekind}.
	The images of the local conjugation maps are open.  Applying weak
	approximation to the finite set consisting of all \(v_i\) and all
	real places \(\RealPlace\), we may therefore choose a maximal \(\FGlobal\)-torus
	\(T\subset H\) such that  $T\times_{\FGlobal}{\FGlobal_{v_i}}$ is \(\CH(O_{\FGlobal,{v_i}})\)-conjugate to $T_{v_i}$ for each $i$, and $T\times_{\FGlobal}{\FGlobal_\RealPlace}$ is \(H(\FGlobal_\RealPlace)\)-conjugate to $T_\RealPlace$ for each real place $\RealPlace$.

	The schematic closure \(\mathcal T\) of \(T\) in \(\CH\) is then a maximal torus: after completion at each maximal ideal \(v_i\), it is
	\(\CH(O_{\FGlobal,v_i})\)-conjugate to the chosen maximal \(O_{\FGlobal,v_i}\)-torus.
	Moreover, the choice at the real places gives 
	$T(\FGlobal)\subset H(\FGlobal)_+$.
	
	Apply the semilocal decomposition
	\cite[Theorem~3]{Guo2022SemilocalDedekind} to the reductive group \(\mathcal T\); it gives
	\[
	\prod_{i=1}^sT(\FGlobal_{v_i})
	=
	T(\FGlobal)\prod_{i=1}^s\mathcal T(O_{\FGlobal,v_i})
	\subset \left(\prod_{i=1}^s\CH(O_{\FGlobal,v_i})\right)H(\FGlobal)_+.
	\]
	
	For every \(v_i\), the Cartan decomposition for the hyperspecial
	subgroup \(\CH(O_{\FGlobal,v_i})\), together with the choice of \(T_{v_i}\), gives
	\(H(\FGlobal_{v_i})=\CH(O_{\FGlobal,v_i})\,T(\FGlobal_{v_i})\,\CH(O_{\FGlobal,v_i}).\)
	Together with the preceding torus decomposition, this shows that
	$\prod_{i=1}^s\CH(O_{\FGlobal,v_i})$ and $H(\FGlobal)_+$ generate
	$\prod_{i=1}^sH(\FGlobal_{v_i})$, proving the assertion.
\end{proof}

\subsection{Integral models of Shimura varieties}
\label{subsec:int-Shimura}
We recall the canonical semi-global integral models of abelian-type Shimura varieties
with parahoric level from
\cite{KisinPappasZhou2024ParahoricII,DanielsYoucis2025}, beginning with the
Hodge-type model. We use the reflex fields and localizations fixed in
\Cref{ass:Shimura-data}.

\begin{defn}\label{defn:Hodge-integral-model}
	As in \cite[\S7.1.2]{KisinPappasZhou2024ParahoricII}, the symplectic lattice
	$W_{\Z_{(p)}}$ realizes the Siegel Shimura variety as a moduli space of
	polarized abelian varieties with prime-to-$p$ level structure.  We define the
	Hodge-type integral model $\Sh_K(G,X)$ as the normalization of the Zariski
	closure of $\ShG_K(G,X)/\reflexHodge$ in the corresponding Siegel semi-global integral
	model over $O_{\reflexHodge,(v)}$.  We write $\CA$ for the universal abelian
	scheme.  The semi-global models and their bases are specified in
	\Cref{defn:int-model}. Note that by \cite{PappasRapoport2024PadicShtukas} such an integral model is independent of the choice of the Hodge embedding $\CG\to \GSp(W_{\Z_{(p)}})$\footnote{In fact, \cite{PappasRapoport2024PadicShtukas} showed one may replace the closed-immersion condition by the weaker requirement that $\CG(\breve{\Z}_p)=G(\breve{\Q}_p)\cap \GSp(W_{\Z_{(p)}})(\breve{\Z}_p)$.}. 
\end{defn}

We make some preparations before defining the KPZ canonical integral models of Shimura varieties of abelian type.
\begin{defn}
	For a reductive group $H/\Q$, put
	$H(\Q)_+:=H(\Q)\cap H(\R)_+$, where $H(\R)_+$ is the inverse image of
	$H^{\ad}(\R)^+$, and let $H(\Q)_+^-$ (respectively $Z_H(\Q)^-$) denote
	the closure of $H(\Q)_+$ (respectively $Z_H(\Q)$) in $H(\A_f)$
	(respectively $Z_H(\A_f)$).  Following Deligne--Kisin, set
	\begin{equation*}
		\begin{aligned}
			\sA(H)^\circ
			&=
			H(\Q)_+^-/Z_H(\Q)^-
			*_{H(\Q)_+/Z_H(\Q)}
			H^{\ad}(\Q)^+,\\
			\sA(H)
			&=
			H(\A_f)/Z_H(\Q)^-
			*_{H(\Q)_+/Z_H(\Q)}
			H^{\ad}(\Q)^+.
		\end{aligned}
	\end{equation*}
	Here $*$ denotes the quotient of the semidirect product, with the second factor acting by conjugation, by the anti-diagonal image of the group written under $*$; see \cite[\S4.5.6]{KisinPappas2018}.
	The group $\sA(H)^\circ$ depends only on $H^{\der}$.
	
	Now let $\CH/\Z_{(p)}$ be a parahoric model of $H$. 
	Let $\CH^\ad/\Z_{(p)}$ be the parahoric model of $H^\ad$ defined using the same point in the Bruhat--Tits building. 
	Write
	$\CH(\Z_{(p)})_+$ for the inverse image of $H^{\ad}(\R)^+$ under
	$\CH(\Z_{(p)})\to H^{\ad}(\R)$.  If $T\subseteq\CH$, let
	$T(\Z_{(p)})^-$ be the closure of $T(\Z_{(p)})$ in $T(\A_f^p)$.  We use
	the integral analogues from \cite[\S4.5.6]{KisinPappas2018}:
	\begin{equation*}
		\begin{aligned}
			\sA(\CH)^\circ
			&=
			\CH(\Z_{(p)})_+^-/Z(\CH)(\Z_{(p)})^-
			*_{\CH(\Z_{(p)})_+/Z(\CH)(\Z_{(p)})}
			\CH^{\ad}(\Z_{(p)})^+,\\
			\sA(\CH)
			&=
			H(\A_f^p)/Z(\CH)(\Z_{(p)})^-
			*_{\CH(\Z_{(p)})_+/Z(\CH)(\Z_{(p)})}
			\CH^{\ad}(\Z_{(p)})^+.
		\end{aligned}
	\end{equation*}
	The group $\sA(\CH)^\circ$ depends only on the derived group together with its
	parahoric model.  In particular, an ad-isomorphism of parahoric Shimura
	data gives the homomorphisms between these groups used below; cf.
	\cite[Lemmas~3.12--3.13]{DanielsYoucis2025}.
\end{defn}
In the general definition of integral models of abelian-type Shimura varieties (cf. \cite{KisinPappasZhou2024ParahoricII, DanielsYoucis2025}), one needs the $\CH^{\ad}(\Z_{(p)})^+$-action on the Shimura
varieties. Defining this action requires $\CH^{\ad}\cong\CH/Z(\CH)$. This
condition does hold when the center of $H$ is an $R$-smooth torus by
\cite[Proposition 2.4.13]{KisinZhou2025}. But the group $G$ need not have
connected center, especially in the even orthogonal case. To fix this problem, we introduce the following definition.

\begin{defn}
	A \emph{well-adapted Hodge-type datum} (in the sense of \cite[Definition~3.8]{DanielsYoucis2025}) for $(G_2,X_2,\CG_2)$ consists of a Hodge-type parahoric Shimura datum $\fka=(\wt G,\wt X,\wt{\CG})$ together with an isogeny\footnote{A morphism $\wt G\to G_2$ need not exist; the construction of special endomorphisms uses the central epimorphism $G\to G_2$ of \Cref{rmk:isogeny}.}
	\(\beta:\wt G^{\der}\longrightarrow G_2^{\der},\)
	such that
	\begin{enumerate}
		\item $\beta$ induces the prescribed isomorphism of adjoint Shimura data.
		\item Let $\wt {\mathbf{x}}\in\CB(\wt G,\Q_p)$ correspond to a point $\mathbf{x}_2\in \CB(G_2,\Q_p)$ 
		under the building isomorphism induced by $\beta$.  Then
		\[
		\CG_2\simeq(G_2)_{\mathbf{x}_2}^{\circ}
		\quad\Longrightarrow\quad
		\wt{\CG}\simeq(\wt G)_{\wt{\mathbf{x}}}^{\circ}.
		\]
		\item If $\wt E=E(\wt G,\wt X)$ and
		$E':=\reflexAbel\wt E$ is the compositum of reflex fields, then
		every place $v_2\mid p$ of $\reflexAbel$ splits completely in
		$E'/\reflexAbel$; equivalently, for every $w\mid v_2$ the map
		$(\reflexAbel)_{v_2}\to E'_w$ is an isomorphism.
	\end{enumerate}
	From now on, we fix a well-adapted Hodge-type datum $\fka=(\wt G,\wt X,\wt{\CG})$
	and put $\wt K_p:=\wt{\CG}(\Z_p)$.
	We choose the datum furnished by the construction of
	\cite[\S\S7.2.1--7.2.3 and Proposition~7.2.19]{KisinPappasZhou2024ParahoricII}
	with its maximal Hodge-type derived group. The resulting central
	isogeny $\wt G^{\der}\to G^{\der}$ composes with $\alpha$ to give
	$\beta$, and gives the compatible homomorphism
	$\sA(\wt{\CG})^\circ\to\sA(\CG)^\circ$.
	We further assume that $Z(\widetilde G)$ is an $R$-smooth torus and that $\pi_1(\widetilde G^\der)_I\to\pi_1(\widetilde G)_I$ is injective.
	As explained in \cite[proof of Corollary~4.11]{DanielsYoucis2025}, this is possible.
	By \cite[Proposition~7.2.19(1)]{KisinPappasZhou2024ParahoricII} we can even
	assume $\pi_1(\widetilde G^\der)=1$ when $(G^\ad,X^\ad)$ has no factors of type
	$D^\BH$.
\end{defn}

\begin{defn}\label{defn:int-model}
	We now define the semi-global integral model $\Sh_{K_2}=\Sh_{K_2}^\fka(G_2,X_2)$ over
	$O_{\reflexAbel,(v_2)}$. 
	Let $\tilde E=E(\tilde G,\tilde X)$ be the reflex field. 
	Regard $\reflexAbel$ and $\tilde E$ as subfields of a fixed
	algebraic closure $\overline\Q$, and choose an embedding
	$\overline\Q\hookrightarrow\overline\Q_p$ inducing a place $w$ of
	$\reflexAbel\tilde E$ above the fixed place $v_2\mid p$ of
	$\reflexAbel$. Let $(\reflexAbel\tilde E)^p$ be the maximal extension
	of $\reflexAbel\tilde E$ in $\overline\Q$ unramified at every place
	above $p$.  
	Fix a sufficiently small neat $K_2^p\subset G_2(\A_f^p)$, put $K_2=K_{2,p}K_2^p$, and choose a connected component $\wt X^+\subset\wt X$. 
	
	As in \cite[\S3.3, Step~3]{DanielsYoucis2025}, consider the natural map
	\begin{equation*}
		\sA(\wt{\CG})^\circ\backslash\sA(\CG_2)
		\longrightarrow
		\sA(\wt G)^\circ\backslash\sA(G_2)/K_{2,p}.
	\end{equation*}
	We choose a finite set $J_{G_2}\subset G_2(\Q_p)$ containing $1$ such that
	\[
	\sA(\wt G)^\circ\backslash\sA(G_2)/K_{2,p}
	=\coprod_{j\in J_{G_2}}\{[aj]:a\in\sA(\CG_2)\},
	\]
	where brackets mean classes in the displayed double quotient.
	The existence and finiteness of this decomposition are part of the construction in loc.\ cit.
	The component $\wt X^+$ determines a geometric connected component
	$\Sh^+_{\wt K_p}(\wt G,\wt X)$ of $\Sh_{\wt K_p}(\wt G,\wt X)$. Over the base ring $\mathcal O_{(\reflexAbel\tilde E)^p}\otimes_{\mathcal O_{\reflexAbel}}
	O_{\reflexAbel,(v_2)}$ we define 
	\begin{equation*}
		\Sh_{K_{2,p}}(G_2,X_2)
		\simeq
		\coprod_{j\in J_{G_2}}
		\left[
		\Sh^+_{\wt K_p}(\wt G,\wt X)\times\sA(\CG_2)j
		\right]\big/\sA(\wt{\CG})^\circ.
	\end{equation*}
	The semilinear action of
	$\operatorname{Gal}((\reflexAbel\tilde E)^p/\reflexAbel\tilde E)$
	constructed there descends this quotient to
	$\mathcal O_{\reflexAbel\tilde E}\otimes_{\mathcal O_{\reflexAbel}}
	O_{\reflexAbel,(v_2)}$. Localizing at $w$ and taking the quotient by
	$K_2^p$ gives a model over $O_{\reflexAbel\tilde E,(w)}$.
	By well-adaptedness, the completion map induces an isomorphism
	\[
	O_{\reflexAbel,v_2}\xrightarrow{\sim}
	O_{\reflexAbel\tilde E,w}.
	\]
	The base change of this model to $O_{\reflexAbel\tilde E,w}$ therefore
	defines the completed model over $O_{\reflexAbel,v_2}$, as in
	\cite[\S3.3, Step~4]{DanielsYoucis2025}. By
	\cite[Theorem~4.10]{DanielsYoucis2025}, this completed model is
	canonical and independent of $\fka$.
	
	To descend the localized model along
	$O_{\reflexAbel,(v_2)}\to O_{\reflexAbel\tilde E,(w)}$, we combine the  descent datum of $\ShG_{K_2}/\reflexAbel$ and the descent datum of the canonical completed model
	over $O_{\reflexAbel,v_2}$. These data agree over the fraction fields of
	the completions, by uniqueness of the canonical completed model.
	By \cite[Theorem~1.2]{MoretBailly1996Descente}, the two descent data
	glue uniquely over the localized rings; uniqueness gives the cocycle
	condition. Quasi-projectivity makes this descent datum effective. We denote the corresponding scheme over $O_{\reflexAbel,(v_2)}$ by $\Sh_{K_2}=\Sh_{K_2}^\fka(G_2,X_2)$.
\end{defn}
	The central map $\alpha:G\to G_2$ sends $\{\mu\}$ to
$\{\mu_2\}$ and gives $\reflexAbel\subseteq\reflexHodge$.
Applying the same construction to the Hodge type datum $(G,X,\CG,K_p)$,
with the chosen well-adapted datum,
defines a new integral model $\Sh_K^\fka(G,X)$ over
$O_{\reflexHodge,(v)}$.
By \cite[Corollary~4.12(2)]{DanielsYoucis2025}, this integral model coincides with the one $\Sh_K^\emptyset(G,X):=\Sh_K(G,X)$  defined in \Cref{defn:Hodge-integral-model}.

\begin{lem}\label{lem:J2-trivial}
	In the situation of \Cref{ass:Shimura-data}, we have $J_{G_2}=\{1\}$. 
	By contrast, the assumptions made on $(G,X,\CG)$ do not force
	$J_G$ to be trivial.
\end{lem}

\begin{proof}
	Apply \Cref{prop:Kisin2010-2.2.6} to the reductive
	$O_{\FGlobal,(p)}$-group in
	\Cref{ass:Shimura-data}(\ref{ass:Shimura-data:reductive}). It gives
	$G_2(\Q_p)=K_{2,p}G_2(\Q)_+$. By the argument of
	\cite[Lemma~2.2.5]{Kisin2010}, the natural map in
	\Cref{defn:int-model} is therefore surjective. Its image already covers
	the target with $j=1$, so $J_{G_2}=\{1\}$.
	The reductive group used here is over
	$O_{\FGlobal,(p)}$; its Weil restriction need not be reductive over $\Z_p$.
\end{proof}

\begin{prop}\label{prop:Kisin-abelian-(3.3.10)}
	Let	$K^p\subset G(\BA_f^p)$ and
	$K_2^p\subset G_2(\BA_f^p)$ be neat compact open subgroups with
	$\alpha(K^p)\subset K_2^p$; then the central epimorphism $\alpha:G\to G_2$ induces a finite \'etale map
	\[
	q_{K^p,K_2^p}:
	\Sh_{K_pK^p}(G,X)\longrightarrow
	\Sh_{K_{2,p}K_2^p}(G_2,X_2)
	\times_{O_{\reflexAbel,(v)}}O_{\reflexHodge,(v)}.
	\]
	These maps are compatible with the generic Shimura morphism,
	prime-to-$p$ Hecke actions, and changes of level.
\end{prop}

\begin{proof}
	On the canonical local integral models (i.e., over the completed local ring ${O}_{\reflexHodge,v}$), the map is provided by
	\cite[Proposition~3.22]{DanielsYoucis2025}. 
	This local map and the
	generic Shimura map agree over the fraction field of the completion. Since they are quasi-projective, by Beauville--Laszlo gluing \cite[Theorem~1.2]{MoretBailly1996Descente} they glue uniquely to
	a map between the semi-global models of \Cref{defn:int-model}. 
	At the displayed finite levels its completed base change is finite
	\'etale, so the same holds over the localized reflex ring by faithfully
	flat descent. Compatibility with Hecke actions and level changes follows
	from the corresponding local and generic compatibilities and the same
	full faithfulness.
\end{proof}

\subsection{Local models}
Under \Cref{ass:Shimura-data}, the Hodge-type normalization has the
global torsor diagram and tensor compatibility of
\cite[Theorem~7.1.3]{KisinPappasZhou2024ParahoricII}. For the abelian-type
model, the local-model diagram is given by
\cite[Theorem~7.2.21 and Corollary~7.2.24]{KisinPappasZhou2024ParahoricII}.
The resulting local integral models of Shimura varieties are canonical by
\cite[Corollary~4.12]{DanielsYoucis2025}; their semi-global
choice-independence follows from \Cref{defn:int-model}.

Before stating the local-model diagrams, we fix precisely the meaning of the
local model and its functoriality in the groups occurring above.  Let
$(H,\mathcal H,\{\nu\})$ be a local model triple over $\Q_p$, with
$\mathcal H$ parahoric, and let $E(H,\{\nu\})/\Q_p$ be the local reflex
field of $\{\nu\}$.  When the triple is obtained from the global data above, these
local rings are
\[
O_{E(G_{\Q_p},\{\mu\})}=\widehat O_{\reflexHodge},
\qquad
O_{E(G_{2,\Q_p},\{\mu_2\})}=\widehat O_{\reflexAbel}.
\]
Scholze--Weinstein attach to $\mathcal H$ the mixed-characteristic
Beilinson--Drinfeld affine Grassmannian
\(\Gr_{\mathcal H}\longrightarrow \Spd\Z_p;\)
see \cite[\S20.3]{ScholzeWeinstein2020Berkeley}.  If $\{\nu\}$ is minuscule,
the generic fiber over $\Spd E(H,\{\nu\})$ contains the Schubert diamond
\(\Gr_{H,\nu}\simeq \mathcal F_{H,\nu}^{\diamond},\)
where $\mathcal F_{H,\nu}$ is the projective homogeneous variety of parabolics
of type $\{\nu\}$.  
\begin{defn}
	The \emph{$v$-sheaf local model} is defined to be the closure 
	\begin{equation*}
		\BM^v_{\mathcal H,\{\nu\}}
		:=
		\overline{\Gr_{H,\nu}}^{\,v}
		\subset
		\Gr_{\mathcal H}\times_{\Spd\Z_p}\Spd O_{E(H,\{\nu\})}.
	\end{equation*}
	It is the smallest closed
	sub-$v$-sheaf containing the generic Schubert diamond.  This
	is \cite[Definition~2.2 and Definition~4.11]
	{AnschutzGleasonLourencoRicharz2022PadicLocalModels}; compare the original
	formulation in \cite[Conjecture~21.4.1 and Remark~21.4.2]
	{ScholzeWeinstein2020Berkeley}.
	
	The scheme denoted below by
	$\BM^{\loc}_{\mathcal H,\{\nu\}}$ is characterized by
	\begin{equation*}
		\bigl(\BM^{\loc}_{\mathcal H,\{\nu\}}\bigr)^\diamond
		\xrightarrow{\;\sim\;}
		\BM^v_{\mathcal H,\{\nu\}}.
	\end{equation*}
	By the representability
	theorem of
	\cite[Theorem~7.21]{AnschutzGleasonLourencoRicharz2022PadicLocalModels},  for minuscule $\{\nu\}$ there exists a unique, up to unique
	isomorphism, flat projective weakly normal
	$O_{E(H,\{\nu\})}$-scheme $\BM^{\loc}_{\mathcal H,\{\nu\}}$. 
\end{defn}

\begin{prop}\label{prop:local-model-diagram}
	Let $\Sh^{\loc}_{K}
	:=\Sh_{K}(G,X)
	\times_{O_{\reflexHodge}}\widehat O_{\reflexHodge}$ be the local integral model (as opposed to the semi-global integral model $\Sh_{K}(G,X)$). 
	We have a local-model
	diagram over $\widehat O_{\reflexHodge}$
	\[
	\begin{tikzcd}
		\Sh^{\loc}_{K} &  \widetilde{\Sh}^{\loc}_{K}(G,X)
		\arrow[l,"\pi"'] \arrow[r,"q_{\loc}"] & \BM^{\loc}_{\CG,\{\mu\}},
	\end{tikzcd}
	\]
	where $\pi$ is a $\CG$-torsor and $q_{\loc}$ is $\CG$-equivariant and smooth
	of relative dimension $\dim G$.  
\end{prop}
\begin{proof}
	For the Hodge-type Shimura datum $(G, X)$, the hypotheses (A)--(C) of
	\cite[\S7.1.2]{KisinPappasZhou2024ParahoricII} are satisfied by
	\Cref{lem:Shimura-data-strongly-admissible} and
	\Cref{ass:Shimura-data}(\ref{ass:Shimura-data:good}), so the claim follows from \cite[Theorem~7.1.3]{KisinPappasZhou2024ParahoricII}. 
\end{proof}

The torsor $\pi$ is
$\underline{\mathrm{Isom}}^\otimes(W_{\Z_p}\otimes\CO_{\Sh^{\loc}_K},\CW_\DR)$,
where the trivializations preserve the tensors $s_{\alpha,\DR}$ of \Cref{sec:tensor}.
The map $q_\loc$ transports the Hodge filtration on $\CW_\DR$ to
$W_{\Z_p}\otimes\CO$ through such a trivialization. For any $\CG$-representation $U$, the associated de Rham vector bundle $\CU_\DR$ on $\Sh^{\loc}_{K}$ therefore becomes the constant bundle associated to $U$ after pullback by $\pi$.

\begin{prop}\label{eq:local-model-common-adjoint}
	Suppose that $(G, \{\mu\})\to (G_2, \{\mu_2\})$ is a central epimorphism inducing an isomorphism on adjoint groups. Let $\reflexHodgeLocal, \reflexAbelLocal$ be their reflex fields. Then there is a canonical
	isomorphism
	\[
	\BM^{\loc}_{{\CG},\{\mu\}}
	\xrightarrow{\;\sim\;}
	\BM^{\loc}_{\CG_2,\{\mu_2\}}\otimes_{\widehat O_{\reflexAbelLocal}}\widehat O_{\reflexHodgeLocal},
	\]
	obtained by identifying both sides with the base change of the same adjoint
	local model.
\end{prop}
\begin{proof}
	By \cite[Proposition~21.5.1]{ScholzeWeinstein2020Berkeley}, the induced morphism on Beilinson--Drinfeld 
	Grassmannians identifies the $v$-closures of the corresponding minuscule
	Schubert diamonds after extension to $\widehat O_{\reflexHodgeLocal}$. 
	The claim then follows. 
\end{proof}

\subsection{Data used in the definition of special cycles}
We retain the situation in \Cref{ass:Shimura-data}. The
$\FGlobal$-linearity of the Hodge embedding implies that the
$\FGlobal$-action and the $G$-action on $W$ commute. Since
$\CG\hookrightarrow\GL_{\Z_p}(W_{\Z_p})$ is the schematic closure of
$G_{\Q_p}$, the $O_\FGlobal$-action also commutes with the
$\CG$-action on $W_{\Z_p}$.
Moreover, $\psi=\Tr_{\FGlobal/\Q}\circ\Psi$ gives
\(\psi(ax,y)=\psi(x,ay)\), \(a\in \FGlobal,\)
so the $\FGlobal$-action is self-adjoint for the symplectic involution.

The following trace pairing fixes the normalization of the quadratic form
used for special cycles. 

\begin{defn}
	Equip $\End_{O_{\FGlobal,(p)}}(W_{O_{\FGlobal,(p)}})$ with the
	$O_{\FGlobal,(p)}$-valued symmetric bilinear pairing induced by the symplectic
	pairing $\Psi$ on $W_{O_{\FGlobal,(p)}}$.
	Namely, for $f,g\in \End_{O_{\FGlobal,(p)}}(W_{O_{\FGlobal,(p)}})$ we define 
	\[(f,g)_\Psi=\tr_{O_{\FGlobal,(p)}}(f\circ g^\dagger),\]
	where $g^\dagger$ is the integral adjoint characterized by
	$\Psi(ga,b)=\Psi(a,g^\dagger b)$ for all $a,b\in W_{O_{\FGlobal,(p)}}$.
	Its integrality follows from the self-duality hypothesis.
\end{defn}

\begin{ass}\label{ass:special-cycle}
	Let $\VlattiveGlobal$ be a finite projective $O_\FGlobal$-module with a quadratic form $(\ ,\ )$.
	Its localization $\VlattiveGlobal\otimes_{O_\FGlobal}O_{\FGlobal,(p)}$
	carries a $\CG_2$-action preserving this form and admits an
	$O_{\FGlobal,(p)}$-linear $\CG$-equivariant injection into
	$\End_{O_{\FGlobal,(p)}}(W_{O_{\FGlobal,(p)}})$.
	We embed $\VlattiveGlobal$ into $\End_{\Z_{(p)}}(W_{\Z_{(p)}})$. 
	Write $V_{\Z_{(p)}}$ for the underlying $\Z_{(p)}$-lattice of
	$\VlattiveGlobal\otimes_{O_\FGlobal}O_{\FGlobal,(p)}$, and put
	$V_{\Z_p}=V_{\Z_{(p)}}\otimes_{\Z_{(p)}}\Z_p$. We require the
	inclusion at $p$ to be saturated, meaning that
	\[
	Q_{\Z_p}:=\End_{\Z_p}(W_{\Z_p})/V_{\Z_p}
	\]
	is finite free over $\Z_p$. Thus
	$0\to V_{\Z_p}\to\End_{\Z_p}(W_{\Z_p})\to Q_{\Z_p}\to0$
	is an exact sequence of $\CG$-representations. 
	Write $V=V_\Q:=V_{\Z_{(p)}}\otimes_{\Z_{(p)}}\Q$, with its induced
	$\FGlobal$-action.
	
	We assume that $N$ is a positive integer such that the quadratic pairing $(\ ,\ )_\Psi$ on $\End_{O_{\FGlobal,(p)}}(W_{O_{\FGlobal,(p)}})$ restricts to $N(\ ,\ )$ on $\VlattiveGlobal$.  We assume $G_2\to \GL_\Q(V)$ has finite kernel\footnote{This condition  is only used in \Cref{lem:fundamental-emb}.}.
\end{ass}

\begin{defn}\label{defn:d}
	Choose an extension $E'/\reflexHodge$ so that the cocharacter $\mu:\BG_{m, E'}\to G_{E'}$ is defined over $E'$. The $G$-action on $V$ then induces a $\FGlobal\otimes_\Q E'$-linear filtration on $V\otimes_\Q E'$, with weights contained in $\{-1, 0, 1\}$. 		
	Define
	$$d_+:=\dim_{E'} (V\otimes_\Q E')/\Fil^0(V\otimes_\Q E')=\dim_{E'}\gr^{-1}(V\otimes_\Q E').$$
	The deformation-theoretic meaning of $d_+$ as the dimension of the obstruction
	space for special endomorphisms is explained in
	\cite[\S\S4--6]{Madapusi2022Derived}; compare also
	\cite[\S5]{MadapusiPera2016}.  
\end{defn}

\subsection{Tensors on the integral model}\label{sec:tensor}
We use homological realizations of the Hodge representation:
$\CW_\Betti$ is the relative Betti homology with coefficients in
$\Z_{(p)}$, $\CW_p=T_p(\CA)$ on the generic fiber, and
$\CW_{\dR}=(H^1_{\dR}(\CA/\Sh_K))^\vee$. For $\ell\ne p$,
write $\CW_\ell=T_\ell(\CA)\otimes_{\Z_\ell}\Q_\ell$.
The notation $\BD(\CA[p^\infty])$ denotes the dual of the
contravariant crystalline Dieudonn\'e crystal; its evaluation is
$\CW_{\dR}$. For any such realization $\CW$ and integers $i,j\geq0$,
put $\CW^{\otimes(i,j)}=\CW^{\otimes i}\otimes(\CW^\vee)^{\otimes j}$;
$\CW^\otimes$ denotes the tensor constructions containing the
specified defining tensors.

As explained in \cite[(2.2)]{Kisin2010}, the tensors $(s_\alpha)_{\alpha\in I}$ in $W^\otimes$ define Betti tensors
$s_{\alpha,\betti}\in\CW_\betti^\otimes$.
The Betti--\'etale comparison isomorphism gives corresponding tensors
$s_{\alpha,\ell}\in\CW_\ell^\otimes$ for every prime $\ell$.
The Betti--de Rham comparison gives de Rham tensors $s_{\alpha,\dR}\in\CW_\dR^\otimes$.
By \cite[Lemma (2.2.1)]{Kisin2010}, the \'etale and de Rham tensors are defined over $\reflexHodge$. For $\ell\ne p$, the \'etale tensors extend over $O_{\reflexHodge}$, as do the de Rham tensors, by \cite[Proposition 4.2.6]{KisinPappas2018}. The crystalline realization at $p$ is treated in \Cref{prop:cris-tensor}.

By the proof of \cite[Theorem~7.1.3]{KisinPappasZhou2024ParahoricII}, for
any finite field $k$ and $x\in\Sh_K(k)$, one may choose an extension
$O_K/\widehat O_{\reflexHodge}$ of complete discrete valuation rings with
residue field $k$ and a lift
$\widetilde x\in\Sh_K(O_K)$.  The $p$-adic comparison isomorphism and the
\'etale tensors $s_{\alpha,\widetilde x,\et}$ give Frobenius-invariant
tensors
$s_{\alpha,x,0}\in\BD(\CA[p^\infty])^\otimes(W(k))$, independent of the
chosen lift. 
 In fact,
\cite[Proposition~1.3.12]{KisinMadapusiPeraShin2022} shows that
these tensors arise from global crystalline tensors on the special fiber:
\[
s_{\alpha,\cris}'\in
\Gamma\!\left(
(\Sh_K\otimes_{O_{\reflexHodge}} k/W(k))_\cris,
\BD(\CA[p^\infty])^\otimes
\right).
\]
We prove that they extend over the completed local base change when the
semi-global integral model is smooth.  The proof uses prismatic theory and is
contained in \Cref{sec:prism}.
\begin{prop}\label{prop:cris-tensor}
	Consider the formal scheme 
	$	\widehat{\Sh}_{K}:=(\Sh_K)^\wedge_p/\hat{O}_{\reflexHodge}$. There is an
	exact $\Z_p$-linear symmetric monoidal functor
	\[
	\omega_\cris:\Rep_{\Z_p}(\CG)\longrightarrow
	\Vect^\varphi((\widehat{\Sh}_K/\Z_p)_\cris)
	\]
	whose value on $W_{\Z_p}$ is the Dieudonn\'e crystal
	$\CW_\cris=\BD(\CA[p^\infty])$, with the same homological
	normalization as $\CW_p$. Here $\Vect^\varphi$ denotes finite locally free
	crystals $\mathcal E$ equipped with an isomorphism
	$F^*\mathcal E[1/p]\simeq\mathcal E[1/p]$.
	In particular, each defining tensor $s_\alpha$ gives a morphism
	\[
	s_{\alpha,\cris}:\CO_\cris\longrightarrow\CW_\cris^\otimes
	\]
	whose rationalization is Frobenius invariant and whose restriction
	at every geometric special-fiber point $x$ is $s_{\alpha,x,0}$.
	In particular, for every finite extension $L/\reflexHodge_v$, every
	$z\in\Sh_K(O_L)$ and every $U\in\Rep_{\Z_p}(\CG)$, we have the comparison isomorphism
	\[
	\CU_{p,z}\otimes_{\Z_p}\Bcris
	\xrightarrow{\sim}
	\CU_{\cris,z|_{\Spec k_L}}(W(k_L))[1/p]\otimes_{W(k_L)[1/p]}\Bcris.
	\]
	Here $k_L$ is the residue field of $L$, and $\CU_{p,z}$ is the fiber of $z^*\CU_p$ at a geometric generic point of $\Spec O_L$.
\end{prop}

Define
\[
\CV_\cris:=\omega_\cris(V_{\Z_p}),\qquad
\CQ_\cris:=\omega_\cris(Q_{\Z_p}),\qquad
\CH_\cris:=\Endint_{\CO_\cris}(\CW_\cris).
\]
Applying $\omega_\cris$ to the exact sequence in
\Cref{ass:special-cycle} gives an exact sequence of crystals
\begin{equation}\label{eq:special-crystal-quotient}
	0\longrightarrow\CV_\cris
	\xrightarrow{\iota_\cris}\CH_\cris
	\xrightarrow{q_\cris}\CQ_\cris\longrightarrow0.
\end{equation}
Thus both the subcrystal and its quotient are finite locally free;
the sequence stays exact after evaluation on any PD-thickening and
after any base change.

\begin{rmk}
	When we say a
	tensor $s:\CO_\cris\to\BD(\CA[p^\infty])^{\otimes(i,j)}$ is
	$\varphi$-invariant, we always mean that $s[1/p]$ is $\varphi$-invariant.
\end{rmk}
\begin{rmk}\label{rmk-cris-dR}
	For $O_\FGlobal=\Z$, the constructions in \cite{MadapusiPera2016,LZ-orth}
	use a $\CG$-invariant integral idempotent tensor
	$\bpi:\End_{\Z_p}(W_{\Z_p})\to\End_{\Z_p}(W_{\Z_p})$.
	When $O_{\FGlobal,v}/\Z_p$ is ramified for a place $v\mid p$ of $\FGlobal$,
	such an integral projection
	onto $V_{\Z_p}$ may fail to exist.
\end{rmk}
\begin{prop}\label{prop:crystalline-deRham-comparison}
	Let $\Spf R$ be an affine formal scheme over $\widehat{\Sh}_{K}$.
	The evaluation of $s_{\alpha,\cris}$ at
	$(R,R)\in(\widehat{\Sh}_{K}/\Z_p)_\cris$ equals
	$s_{\alpha,\dR}$ under
	$\BD(\CA[p^\infty])(R)\cong\CW_{\dR,R}$.
	The same comparison identifies the evaluation of
	\eqref{eq:special-crystal-quotient} with the exact sequence of
	associated de Rham bundles
	$0\to\CV_{\dR,R}\to\Endint_R(\CW_{\dR,R})\to\CQ_{\dR,R}\to0$.
\end{prop}

\begin{proof}
	We can check this at each completed local ring $\hat{\CO}_{\Sh_K, x}$.
	The assertion for the defining tensors follows from the proof of
	\cite[Proposition~4.2.2]{KisinPappas2018}. The corresponding
	tensor-frame torsors are therefore identified. Applying their
	associated-bundle functors to $V_{\Z_p}$ and $Q_{\Z_p}$ gives the
	assertion for the exact sequence. 
\end{proof}

\subsection{Special endomorphisms for Hodge type}\label{subsec:special-endomorphisms-Hodge-type}
In this section we use the tensors recalled in \Cref{sec:tensor} to define special endomorphisms on the Hodge-type Shimura variety $\Sh_K$. 

For $\diamond\in\{\Betti,\ell,\dR,\cris\}$, let
$\CW_\diamond$ and $\CV_\diamond$ denote the realizations of $W$ and
$V$ and write
\[
\iota_\diamond:\CV_\diamond\lhook\joinrel\longrightarrow
\CH_\diamond:=\Endint(\CW_\diamond)
\]
for the realization of the inclusion $V\hookrightarrow\End(W)$.
We identify each $\CV_\diamond$ with its
image under $\iota_\diamond$, and write $\CQ_\diamond$ for the quotient.
In particular, $f_\cris\in\CV_\cris$ means that the
morphism $\CO_\cris\to\CH_\cris$ defined by $f$ factors through
$\iota_\cris$, or equivalently that $q_\cris(f_\cris)=0$.

\begin{defn}\label{defn:special}
	Let $S$ be a scheme of finite type over the indicated base, and let
	$f\in\End(\CA_S)_{(p)}$.
	\begin{enumerate}
		\item If $S$ is a $\ShG_K(G,X)_\C$-scheme, then $f$ is
		\emph{special} if $f_\Betti\in\Gamma(S(\C),\CV_\Betti)$.
		\item\label{defn:special:generic}  If $S$ is a $\ShG_K(G,X)$-scheme, then $f$ is
		\emph{$\ell$-special} if $f_\ell\in\Gamma(S_{\et},\CV_\ell)$.
		\item\label{defn:special:special} If $S$ is a $\Sh_K\otimes_{O_{\reflexHodge}}k$-scheme,
		then $f$ is \emph{$p$-special} if, for every perfect-field-valued
		point $s\to S$,
		\[
		f_{s,\cris}\in\Gamma((s/\Z_p)_\cris,\CV_{\cris,s}).
		\]
		\item If $S$ is a $\Sh_K$-scheme, then $f$ is
		\emph{special} if $f|_{S\times_{O_{\reflexHodge}}\reflexHodge}$ satisfies (\ref{defn:special:generic}) for
		$\ell=p$ and $f|_{S\times_{O_{\reflexHodge}}k}$ satisfies (\ref{defn:special:special}).
	\end{enumerate}
\end{defn}
This is the fiberwise convention of
\cite[Definitions~5.8 and~5.11]{MadapusiPera2016}.
We call $f$ \emph{crystallinely special} over a base on which $p$
is nilpotent if the entire morphism $f_\cris:\CO_\cris\to\CH_\cris$
factors through $\CV_\cris$. This stronger assertion implies
$p$-specialness. We prove the converse in \Cref{prop:tangent-space-divisor} under \Cref{ass:d=1}.

\begin{lem}[\cite{MadapusiPera2016}, Lemma~5.4]\label{lem:special-char-0-fiberwise}
	Suppose that $T$ is a finite-type $\ShG_K(G, X)$-scheme and $f\in \End(\CA_T)_{(p)}$. Then
	the following are equivalent:
	\begin{enumerate}
		\item $f$ is $\ell$-special for all primes $\ell$.
		\item $f$ is $\ell$-special for some prime $\ell$.
		\item The restriction of $f$ to $T\otimes_{\reflexHodge}\C$ is special. 
	\end{enumerate}
	Moreover, $f$ is $\ell$-special if and only if every connected component
	of $T$ contains a point $s$ for which the fiber $f_s$ is
	$\ell$-special.  In particular,
	$\ell$-specialness is a condition that can be checked at geometric points.
\end{lem}

\begin{lem}[{\cite[Lemma 5.9]{MadapusiPera2016}}]\label{lem:special-char-p-fiberwise}
	Let $T$ be a finite-type $\Sh_K\otimes_{O_{\reflexHodge}}k$-scheme, and let
	$f\in\End(\CA_T)_{(p)}$.
	Then the following are equivalent:
	\begin{enumerate}
		\item For every point $s\to T$ valued in a perfect field, the fiber $f_s\in\End(\CA_s)_{(p)}$ is $p$-special.
		\item In every connected component of $T$, there exists a point $s$ valued in a perfect field such that the fiber $f_s$ is $p$-special.
	\end{enumerate}
\end{lem}

\begin{proof}
	Apply \cite[Lemma~5.10]{MadapusiPera2016} to the image
	of the universal endomorphism in the finite locally free quotient
	crystal $\CQ_\cris$ of \eqref{eq:special-crystal-quotient}, on the special fiber of the 
	locally Noetherian endomorphism scheme $\underline{\End}(\CA)_{(p)}$. 
\end{proof}

\begin{prop}\label{prop:endomorphism-Hilbert-finite}
	Let $A$ be a projective abelian scheme over a Noetherian scheme $S$.
	The fppf sheaf $\Endint(A)$ is represented by an unramified scheme
	locally of finite presentation over $S$.  Each of its connected
	components is finite over $S$.
\end{prop}

\begin{proof}
	Choose a relatively ample invertible sheaf
	\(\mathcal L\) on \(A\times_S A\).
	
	The functor $\Endint(A)$ is the closed subfunctor of
	$\underline{\Mor}_S(A,A)$ cut out by preservation of the identity and
	of the group law.  Here closedness follows from the separatedness of
	$A/S$.
	The graph construction gives a morphism
	\[
	\gamma\colon \underline{\Mor}_S(A,A)
	\longrightarrow
	\Hilb_{A\times_S A/S},
	\qquad
	f\longmapsto \Gamma_f,
	\]
	which identifies 
	\(\underline{\Mor}_S(A,A)\) with an open subscheme of this Hilbert
	scheme; see
	\cite[Chapter~5, Theorem~5.23]{FGA-explained}. 
	For a numerical polynomial \(P\), let
	\[
	E_P
	=
	\Endint(A)
	\times_{\Hilb_{A\times_S A/S}}
	\Hilb^P_{A\times_S A/S}.
	\]
	Here \(\Hilb^P\) denotes the Hilbert scheme of closed subschemes with
	Hilbert polynomial \(P\) with respect to \(\mathcal L\), which is an open and closed subscheme of $\Hilb$. 
	The scheme \(\Hilb^P_{A\times_S A/S}\) is projective over \(S\), and
	\(E_P\) is a locally closed subscheme of it. 
	
	By rigidity of morphisms between abelian varieties, $E_P\to S$ is
	unramified. For a discrete valuation ring $R$ over $S$, every
	endomorphism of $A_{\operatorname{Frac}(R)}$ extends uniquely to an
	endomorphism of $A_R$ by the Néron mapping property
	\cite[Chapter~1, \S2, Proposition~8]{BLR-Neron}.
	Its graph is isomorphic to $A_R$, hence flat over $R$, so its Hilbert
	polynomial with respect to $\mathcal L$ remains $P$. Thus $E_P\to S$
	is proper by the valuative criterion. Since $E_P\to S$ is also of
	finite type and unramified, it is quasi-finite and therefore finite.
	Every connected component of $\Endint(A)$ lies in one of the open
	and closed subschemes $E_P$, so it too is finite over $S$.
\end{proof}

Over $O_{\reflexHodge}$, the same assertions hold for
$\Endint(A)_{(p)}$. Indeed, for every positive integer $n$ prime to $p$,
multiplication by $n$ identifies $\Endint(A)$ with the open and closed
subscheme of endomorphisms $f$ satisfying $f|_{A[n]}=0$.
The condition is open and closed because $A[n]$ is finite \'etale;
such an $f$ factors uniquely through $[n]:A\to A$.
The prime-to-$p$ denominator transitions are therefore open and closed
immersions. Their union represents $\Endint(A)_{(p)}$, and every
connected component occurs at one denominator level, where it is finite.

\begin{lem}[{\cite[Lemma 5.13]{MadapusiPera2016}}]\label{lifts:lem:p-special-consistent}Let $T$ be a finite-type $\Sh_K$-scheme, put
	$T_{\reflexHodge}:=T\times_{O_{\reflexHodge}}\reflexHodge$ and
	$T_k:=T\times_{O_{\reflexHodge}}k$, and assume that $T_{\reflexHodge}$ and $T_k$ are connected and that a point
	of $T_k$ is the specialization of a point of $T_{\reflexHodge}$.  Then
	$f\in\End(\CA_T)_{(p)}$ is $p$-special on $T_{\reflexHodge}$ if and
	only if it is $p$-special on $T_k$.  In particular, if $f$ is
	$p$-special on the special fiber, then it is special.
\end{lem}
\begin{proof}
	First assume that $T=\Spec O_L$, where $L/\reflexHodge_v$ is a finite
	extension of complete discretely valued fields.
	Put $s=\Spec k_L$, choose an algebraic closure $\overline L/L$,
	and write $\overline\eta=\Spec\overline L$.
	The \'etale compatibility in \Cref{prop:cris-tensor} gives
	\[
	\CV_{p,\overline\eta}\otimes_{\Z_p}\Bcris
	\xrightarrow{\sim}
	\CV_{\cris,s}[1/p]\otimes_{W(k_L)[1/p]}\Bcris,
	\]
	compatibly with the ambient endomorphism comparison and the quotient
	in \eqref{eq:special-crystal-quotient}. The comparison isomorphism
	and the finite locally free quotient give the assertion.
	
	For general $T$, by \cite[\href{https://stacks.math.columbia.edu/tag/053U}{Lemma 053U}]{stacks-project} the hypothesis implies that there exist \emph{closed} points $x\in T_{\reflexHodge}$ and $x'\in T_k$ so that $x$ specializes to $x'$. Apply the preceding argument to the normalization of the corresponding trait and apply \Cref{lem:special-char-0-fiberwise,lem:special-char-p-fiberwise}; the claim follows.
\end{proof}

\begin{cor}\label{cor:special-representable}
	The fppf sheaf of special endomorphisms $\sEnd$ on $\Sh_K$ is an open
	and closed subsheaf of $\Endint(\CA)_{(p)}$.  It is therefore represented by a
	scheme unramified and locally of finite presentation over $\Sh_K$, and
	each of its connected components is finite over $\Sh_K$.
\end{cor}

\begin{proof}
	On the generic and closed fibers, specialness is constant on connected
	components by
	\Cref{lem:special-char-0-fiberwise,lem:special-char-p-fiberwise}.
	On a horizontal irreducible component, compare a generic-fiber point
	with a specialization by a trait in its normalization and apply
	\Cref{lifts:lem:p-special-consistent}.
	In each finite-type neighborhood of the endomorphism scheme, any two
	points of a connected component are joined by a finite chain of
	irreducible components. The fiberwise constancy lemmas and the trait
	comparison show that specialness is unchanged along this chain, as in
	\cite[Lemma~6.18 and Proposition~6.19]{Madapusi2022Derived}.
	Thus the fiberwise special locus is a union of connected
	components of $\Endint(\CA)_{(p)}$, and is therefore open and closed. 
\end{proof}

\subsection{Closed points on \texorpdfstring{$\sEnd$}{the special-endomorphism scheme}}
We describe the complex points and finite-field points of $\sEnd$ in this subsection. 
\begin{prop}\label{prop:special-C}
	There is a natural bijection
	\[
	\sEnd(\C)
	\xrightarrow{\ \sim\ }
	G(\Q)\backslash\bigl(Y/K\bigr),
	\qquad
	Y:=\{(v,x, g)|v\in V_\Q,\ x\in X,\ g\in G(\BA_f),\ g_p^{-1}v\in V_{\Z_{p}},\ h_x(i)v=v\}.
	\]
\end{prop} 

\begin{proof}
	Represent a complex point of the Shimura variety by $(x,g)$.  Betti
	realization identifies its special quasi-endomorphisms with the vectors
	$v\in V_\Q$ whose $p$-adic realization lies in $g_pV_{\Z_p}$ and whose
	real realization is fixed by $h_x(i)$.
\end{proof}

\begin{prop}\label{prop:special-F_p}
	Let $k$ be a finite field of
	characteristic $p$.  Then, for every $z\in\Sh_K(k)$, crystalline
	realization and $\ell$-adic realization for each $\ell\neq p$ induce isomorphisms
	\[
	\Z_p\otimes_{\Z_{(p)}}\sEnd(z)
	\xrightarrow{\ \sim\ }
	\CV_{\cris,z}^{\varphi=1}, \qquad \Q_\ell\otimes_{\Z_{(p)}}\sEnd(z)
	\xrightarrow{\ \sim\ }
	\CV_{\ell,z}^{\Frob=1}.
	\]
\end{prop}

\begin{proof}
	We follow the argument of
	\cite[Propositions~6.21 and~6.24]{Madapusi2022Derived}.
	The input used there, \cite[Corollary~2.3.1]{Kisin2017ModP}, treats
	only hyperspecial level. We replace it with
	\cite[Proposition~2.3.3]{KisinZhou2025StronglyCompatible}, which applies
	here by \Cref{rmk:Shimura-data-strongly-admissible}. 
\end{proof}

\subsection{Cycles defined by special endomorphisms}\label{subsec:global-special-cycle}
\begin{lem}\label{lem:special-Z-inv}Let $x \in \Sh_K(T)$ and $f\in \End(\CA_x)_{(p)}$.
	\begin{enumerate}
		\item Assume $T/\Q$. If the $\ell$-adic realization $f_\ell$ is $0$, then $f=0$. \label{lem:special-Z-inv1}

		\item Assume $T/\BF_p$. If the  crystalline realization  $f_\cris$ is $0$, then $f=0$. \label{lem:special-Z-inv2}
		
		\item In general, if both $(f_{\Q})_\ell$ and $(f_{\BF_p})_\cris$ vanish, then $f=0$. \label{lem:special-Z-inv3}
	\end{enumerate}
\end{lem}

\begin{proof}
	For (\ref{lem:special-Z-inv1}), faithfulness of the Tate module gives $f_y=0$ at every geometric point $y$ of $T$.
	For (\ref{lem:special-Z-inv2}), faithfulness of the Dieudonn\'e functor over algebraically closed fields and of $A\mapsto A[p^\infty]$ gives the same conclusion.
	In both cases rigidity implies $f=0$.
	The claim (\ref{lem:special-Z-inv3}) follows from (\ref{lem:special-Z-inv1}) and (\ref{lem:special-Z-inv2}). 
\end{proof}
\begin{defn}
	Let $\lambda:\CA\to\CA^*$ represent the chosen weak polarization.
	For $f\in\End^0(\CA):=\End(\CA)\otimes_\Z\Q$, let $f^*$ be the dual quasi-endomorphism of $\CA^*$, defined by $\Q$-linear extension of duality on endomorphisms. The Rosati involution on $\End^0(\CA)$ is
	\[
	f^\dagger=\lambda^{-1}\circ f^*\circ\lambda,
	\]
	where $\lambda^{-1}$ is the inverse quasi-isogeny.
\end{defn}

\begin{prop}\label{ass:OF-linear-Hodge-realization}
	The $\FGlobal$-linear Hodge embedding determines an action
	$\iota_\FGlobal:O_\FGlobal\longrightarrow\End_{\Sh_K}(\CA)$ with the
	following properties.
	\begin{enumerate}
		\item The action is self-adjoint for the Rosati involution attached to the
		chosen weak polarization:
		\[
		\iota_\FGlobal(a)^\dagger=\iota_\FGlobal(a),\qquad a\in O_\FGlobal.
		\]
		\item Every special endomorphism is $O_\FGlobal$-linear. 
		\item The induced $\FGlobal$-structure is compatible with Betti, $\ell$-adic, de
		Rham, crystalline, and prismatic realizations.
	\end{enumerate}
\end{prop}
\begin{proof}
	Note that multiplication by $O_\FGlobal$ is a family of
	$G$-invariant endomorphism tensors of $W_{O_\FGlobal}$.  Such tensors give endomorphisms of
	the universal abelian variety on the Hodge-type Shimura variety; see
	\cite[\S2.2]{Kisin2010} and
	\cite[Proposition~4.2.6]{KisinPappas2018}. Since the $O_\FGlobal$-action preserves $W_\Z$ and commutes with $G(\A_f)$, these endomorphisms preserve every Tate lattice defined by the level structure. The resulting integral endomorphisms extend over $\Sh_K$ by the extension property; compare
	\cite[Theorem~7.1.3]{KisinPappasZhou2024ParahoricII}.
	
	The self-adjointness follows directly from
	$\psi=\Tr_{\FGlobal/\Q}\circ\Psi$: for $a\in \FGlobal$,
	\[
	\psi(ax,y)=\Tr_{\FGlobal/\Q}(a\Psi(x,y))
	=\Tr_{\FGlobal/\Q}(\Psi(x,ay))=\psi(x,ay).
	\]
	
	Since $V$ is contained in $\End_\FGlobal(W)$ by
	\Cref{ass:special-cycle}, the realization of every special
	endomorphism commutes with $O_\FGlobal$. Faithfulness of realization
	in \Cref{lem:special-Z-inv} implies that the endomorphism itself is
	$O_\FGlobal$-linear. Compatibility across
	realizations follows from the fact that the tensors involved are absolute
	Hodge tensors and from the integral tensor comparison of
	\cites{Kisin2010,KisinPappas2018,KisinPappasZhou2024ParahoricII}.
\end{proof}

For an endomorphism $f$ of an abelian scheme $A$, the usual trace
$\tr(f)$ is computed on each geometric fiber from its rational
$\ell$-adic Tate module, for any prime $\ell$ invertible on that fiber.
It is integer-valued and locally constant;
see \cite[Theorem~10.9 and Proposition~10.20]{milneAV}.
We extend it by $\Q$-linearity to quasi-endomorphisms.

\begin{prop}\label{prop:OF-characteristic-polynomial}
	Let $A$ be an abelian variety over a field $k$, equipped with an action
	$\iota:O_\FGlobal\hookrightarrow\End_k(A)$, and let
	$f\in\End_{O_\FGlobal}(A)$. There is a unique
	$\tr_\FGlobal(f)\in O_{\FGlobal,(p)}$ such that
	\[
	\tr\bigl(\iota(a)f\bigr)
	=\Tr_{\FGlobal/\Q}\bigl(a\,\tr_\FGlobal(f)\bigr),
	\qquad a\in\FGlobal.
	\]
	For every prime $\ell$ invertible in $k$, its image in
	$\FGlobal\otimes_\Q\Q_\ell$ is the relative trace of $f$ on
	$T_\ell(A)\otimes_{\Z_\ell}\Q_\ell$.
\end{prop}

\begin{proof}
	The nondegeneracy of the trace pairing on $\FGlobal/\Q$ gives a unique
	$\tr_\FGlobal(f)\in\FGlobal$ satisfying the displayed identity.
	Restriction of scalars in any $\ell$-adic realization gives the same
	identity for the relative trace, proving the comparison.
	Under each embedding of $\FGlobal$ into $\overline\Q_\ell$, this trace
	is the sum of the eigenvalues of $f$ on the corresponding summand of
	$T_\ell(A)\otimes_{\Z_\ell}\overline\Q_\ell$.
	These eigenvalues are roots of the usual monic characteristic
	polynomial in $\Z[T]$ of $f$. Their sum is therefore an algebraic
	integer, so $\tr_\FGlobal(f)\in O_{\FGlobal,(p)}$.
\end{proof}

For an $O_\FGlobal$-linear endomorphism of an abelian scheme, we define
$\tr_\FGlobal$ on geometric fibers by
\Cref{prop:OF-characteristic-polynomial}. The defining identity, applied
with $a$ in a $\Q$-basis of $\FGlobal$, shows that it is locally constant.
We extend it by $\Q$-linearity to
$\End_{O_\FGlobal}^0(A)=\End_{O_\FGlobal}(A)\otimes_\Z\Q$.

\begin{defn}\label{defn-pairing-special-SO}
	By \Cref{ass:OF-linear-Hodge-realization}, special endomorphisms and
	their Rosati adjoints are $O_\FGlobal$-linear. For each scheme $S/\Sh_{K}$ and $f,g\in \sEnd(S)$ we define the $\frac{1}{N}O_{\FGlobal,(p)}$-valued symmetric
	pairing
	$$h(f,g):=\frac{1}{N}\tr_\FGlobal(f\circ g^\dagger)\in \underline{\frac{1}{N} O_{\FGlobal,(p)}}(S),$$
	where $N$ is defined in \Cref{ass:special-cycle}. 
\end{defn}

\begin{defn}\label{defn:special-cycle}  
	Fix an integer $r>0$.  Let $T\in\Sym_r(\FGlobal)$ and let
	$\Omega\subset V(\BA_f^p)^r$ be a $K^p$-invariant compact subset.  The
	\emph{Kudla--Rapoport cycle} $\CZ(T,\Omega)$ associates to every
	$\Sh_K$-scheme $S$ the set of tuples
	$\bu=(u_1,\ldots,u_r)\in\sEnd(S)^r$ satisfying
	\[
	h(\bu,\bu)=T,
	\qquad
	\realization_{\BA_f^p}(\bu)\in\Omega.
	\]
	It is a separated, unramified scheme over $\Sh_K$.
\end{defn}  
Here $\realization_{\BA_f^p}(\bu)$ denotes the $K^p$-orbit of
$(\varepsilon^p)^{-1}\bu_{\et}\varepsilon^p$ obtained from the prime-to-$p$
level structure at each geometric point. We use the right action
$\mathbf v\cdot k=k^{-1}\mathbf v$. Since $\Omega$ is $K^p$-invariant,
membership in $\Omega$ is independent of the choice of representative of
$\varepsilon^p$.
\begin{prop}\label{prop:fixed-level-special-cycle-finite}
	For every $T$ and $\Omega$ as in \Cref{defn:special-cycle}, the functor
	$\CZ(T,\Omega)$ is represented by a scheme finite and unramified over
	$\Sh_K$.
\end{prop}

\begin{proof}
	By \Cref{cor:special-representable} the functor is represented by a scheme unramified over	$\Sh_K$.
	For finiteness, we use the same
	boundedness argument as in
	\cite[Proposition~6.13]{MadapusiPera2016}; see also
	\cite[\S12.5]{LZ-orth}.
	
	First, $\Omega$ is compact in $V(\BA_f^p)^r$, so there is an
	integer $c\geq1$, prime to $p$, depending only on $\Omega$, such that
	$c\Omega$ is contained in a fixed prime-to-$p$ integral lattice. For
	$\bu=(u_1,\ldots,u_r)\in\CZ(T,\Omega)$, the prime-to-$p$ level condition
	makes $cu_i$ integral away from $p$, and $u_i\in\sEnd(S)$ gives
	integrality at $p$. Thus $v_i:=cu_i$ is an integral special endomorphism.
	
	On every geometric fiber the Rosati form
	\(f\longmapsto\tr(ff^\dagger)\)
	is positive definite.  Moreover,
	\(\tr(v_iv_i^\dagger) = c^2N\Tr_{\FGlobal/\Q}(T_{ii}).\)
	Thus all the $v_i$ lie in a bounded Rosati-norm part of the integral
	endomorphism scheme.
	
	After passing locally on the base, choose a
	symmetric relatively ample line bundle $\mathcal L$ inducing a positive
	multiple of the chosen weak polarization.  For an integral endomorphism
	$f$, use
	\[
	\mathcal M
	=
	\operatorname{pr}_1^*\mathcal L\otimes\operatorname{pr}_2^*\mathcal L
	\]
	on $\CA\times_{\Sh_K}\CA$.  Its restriction to the graph $\Gamma_f$ is
	$\mathcal L\otimes f^*\mathcal L$, whose polarization is
	\(\lambda_{\mathcal L}\circ(1+f^\dagger f).\)
	The eigenvalues of the Rosati-positive endomorphism $f^\dagger f$ are
	non-negative, so a bound for
	$\tr(f^\dagger f)$ gives a bound for
	$\deg(1+f^\dagger f)$, and hence for the Hilbert polynomial of
	$\Gamma_f$ with respect to $\mathcal M$.  Consequently only finitely
	many Hilbert polynomials can occur for the $v_i$.
	
	Applying the Hilbert-scheme argument in the proof of
	\Cref{prop:endomorphism-Hilbert-finite} with this choice of
	$\mathcal M$, the locus $E_P$ corresponding to any fixed Hilbert
	polynomial $P$ is finite and unramified over the base.  Hence the tuples
	$(v_1,\ldots,v_r)$ occurring above lie in a finite union of schemes
	finite and unramified over $\Sh_K$.  On this finite union, specialness,
	the Gram condition
	\(h(\mathbf v,\mathbf v)=c^2T,\)
	and the prime-to-$p$ condition
	$\realization_{\BA_f^p}(\mathbf v)\in c\Omega$ are locally constant.
	They therefore cut out an open and closed subscheme.  
	Thus $\CZ(T,\Omega)$ is represented by a finite $\Sh_K$-scheme.
\end{proof}

\begin{prop}\label{prop:prod-special}
	Let $T_i\in\Sym_{r_i}(\FGlobal)$ and let
	$\Omega_i\subset V(\BA_f^p)^{r_i}$ be right $K^p$-orbits for $i=1,2$.
	Then
	\[
	\CZ(T_1,\Omega_1)\times_{\Sh_K}\CZ(T_2,\Omega_2)
	=
	\coprod_{\substack{
			T=\begin{psmallmatrix}T_1&B\\{}^tB&T_2\end{psmallmatrix}
			\in\Sym_{r_1+r_2}(\FGlobal);\\
			\Omega\in(\Omega_1\times\Omega_2)/K^p}}
	\CZ(T,\Omega),
	\]
	where $K^p$ acts diagonally on the right and the union ranges over the
	nonempty summands.
\end{prop}

\begin{proof}
	An $S$-point of the fiber product is a pair of tuples
	$(\bu_1,\bu_2)$ with Gram matrices $T_1$ and $T_2$.  Concatenating the
	tuples produces a Gram matrix whose off-diagonal block is
	$B=h(\bu_1,\bu_2)$ and a joint diagonal $K^p$-orbit
	$\Omega\in(\Omega_1\times\Omega_2)/K^p$.  Both invariants are locally constant
	on $S$.  They therefore decompose the fiber product into open and closed
	subfunctors, each of which is precisely $\CZ(T,\Omega)$.  The inverse map
	is obtained by projecting a concatenated tuple onto its first $r_1$ and
	last $r_2$ entries.
\end{proof}

\subsection{Deformations of special endomorphisms}\label{subsec:def-special-end}
We show that nonzero special endomorphisms locally define Cartier divisors (cf. \Cref{prop:tangent-space-divisor}) under \Cref{ass:d=1}. The orthogonal calculation in
\Cref{subsec:orth-local-model} verifies this criterion directly.

\begin{defn}\label{defn:special-Hodge-filtration}
	The (homological) de Rham bundle
	$\CW_{\dR}=(H^1_{\dR}(\CA/\Sh_K))^\vee$ has the filtration
	\[
	\Fil^{-1}\CW_{\dR}=\CW_{\dR},\qquad
	\Fil^0\CW_{\dR}
	=\operatorname{Ann}(\Fil^1H^1_{\dR}(\CA/\Sh_K)),\qquad
	\Fil^1\CW_{\dR}=0,
	\]
	where the annihilator is taken under the evaluation pairing.
	On $\BM^{\loc}_{\CG,\{\mu\}}$, set
	\[
	\CW_{\loc}=W_{\Z_p}\otimes\CO_{\BM^{\loc}_{\CG,\{\mu\}}},
	\qquad
	\CV_{\loc}=V_{\Z_p}\otimes\CO_{\BM^{\loc}_{\CG,\{\mu\}}},
	\]
	and give $\CW_{\loc}$ its tautological Hodge filtration.
	For $\diamond\in\{\dR,\loc\}$ and $a\in\Z$, define
	\[
	\begin{aligned}
		\CH_\diamond&=\Endint(\CW_\diamond),\\
		\Fil^a\CH_\diamond
		&=\{u:u(\Fil^b\CW_\diamond)\subseteq
		\Fil^{a+b}\CW_\diamond\text{ for every }b\},
	\end{aligned}
	\]
	and equip $\CV_\diamond$ with the subspace filtration.
\end{defn}

\begin{lem}
	\label{lem:special-filtered-comparison}
	The maps $\pi$ and $q_{\loc}$ of
	\Cref{prop:local-model-diagram} give an isomorphism
	\begin{equation*}
		q_{\loc}^*(\CV_{\loc},\Fil^\bullet)
		\simeq\pi^*(\CV_{\dR},\Fil^\bullet).
	\end{equation*}
\end{lem}
\begin{proof}
	By the definition of the local-model diagram, for each $\CG$-representation $U$ we have a functorial isomorphism between $\pi^*\CU_\DR$ and $U\otimes_{\Z_p}\CO_{\widetilde{\Sh}^{\loc}_K(G,X)}$, and for the representation $W_{\Z_p}$ this isomorphism preserves the filtration on both sides. The claim then follows. 
\end{proof}

Since $\Fil^{-1}\CH_\diamond=\CH_\diamond$, we also have
\[
\gr^{-1}\CV_\diamond
=\CV_\diamond/\Fil^0\CV_\diamond
\lhook\joinrel\longrightarrow
\CH_\diamond/\Fil^0\CH_\diamond,
\qquad \diamond\in\{\dR,\loc\}.
\]
The target is canonically
$\Homint(\Fil^0\CW_\diamond,
\CW_\diamond/\Fil^0\CW_\diamond)$; in particular it is locally free.

\begin{thm}[Grothendieck--Messing]\label{thm:tangent-space}
	Let $R\twoheadrightarrow R_0$ be a surjective local homomorphism
	of complete Noetherian local $O_{\reflexHodge}$-algebras of residue
	characteristic $p$. Assume that its kernel has square zero and carries
	nilpotent divided powers compatible with the canonical divided powers
	on $(p)\subset\Z_p$. Let
	$z\in\Sh_K(R)$ and let $y\in\Sh_K(R_0)$ be its reduction.
	For $f\in\End(\CA_y)$, $p$-adic crystalline evaluation and comparison give
	\[
	f_{\cris,R}\in
	\End_R\bigl(\BD(\CA_y[p^\infty])(R)\bigr)
	=\End_R(\CW_{\dR,R}).
	\]
	Here crystalline evaluation over $R$ is understood in the $p$-adic sense
	of \cite[Remark~12.1.6(5)]{BrinonConrad2009}.
	Then $f$ extends uniquely to $\CA_z$ if and only if
	$f_{\cris,R}\in\Fil^0\CH_{\dR,R}$.
\end{thm}
\begin{proof}
	By \cite[Theorem~12.1.5 and Remarks~12.1.6(2),~(5)]{BrinonConrad2009},
	the morphism $f[p^\infty]$ lifts uniquely to an endomorphism of
	$\CA_z[p^\infty]$ precisely when $f_{\cris,R}$ preserves the Hodge filtration.
	Serre--Tate theory \cite[Chapter~V, Theorem~2.3]{Messing1972} and
	algebraization of morphisms
	\cite[\href{https://stacks.math.columbia.edu/tag/0A42}{Tag 0A42}]{stacks-project}
	identify this with a unique extension of $f$ to $\CA_z$. 
\end{proof}

\begin{ass}\label{ass:d=1}
	Assume that $\gr^{-1} \CV_\loc$ is an invertible sheaf and the cokernel of $\gr^{-1}\CV_\loc\to \CH_\loc/\Fil^0\CH_\loc$ is locally free. 
	By \Cref{lem:special-filtered-comparison} this implies the same assertions for $\CV_\DR, \CH_\DR$. 
	Assume that $\Sh_K$ is proper over $O_\reflexHodge$ and $V^{G^\der}=0$. 
\end{ass}

\begin{rmk}
	The first assertion holds on the generic fiber if and only if $d_+=1$ in \Cref{defn:d}. This generic-rank condition alone leaves open the local freeness of the integral quotient $\CV_\loc/\Fil^0\CV_\loc$ and of the cokernel of its inclusion into $\CH_\loc/\Fil^0\CH_\loc$.
\end{rmk}

\begin{rmk}\label{rmk:d=1}
	Properness and $V^{G^\der}=0$ are used to show that a connected component of $\AbSpecialEnd$, on which the tautology special endomorphism is nonzero, cannot dominant a connected component of the special fiber. These assumptions hold
	in our orthogonal setting when $\FGlobal\ne\Q$. For $\FGlobal=\Q$,
	\cite[Proposition~5.21]{MadapusiPera2016} uses injectivity of the
	Kodaira--Spencer map to show the same result. In the orthogonal setting over $\FGlobal$, the Kodaira--Spencer map
	\[
	\gr^0\CV_{\dR}\longrightarrow
	\gr^{-1}\CV_{\dR}\otimes
	\Omega^1_{\Sh_K/O_{\reflexHodge}}
	\]
	has generic source rank $[\FGlobal:\Q]\dimofV-2$ and target rank
	$\dimofV-2$. Thus this injectivity argument does not extend directly
	when $[\FGlobal:\Q]>1$. We instead use the crystalline
	$H^0$-vanishing proved below.
\end{rmk}

We now prepare to prove \Cref{prop:tangent-space-divisor}. Fix a uniformizer $\varpi$ of $O_{\reflexHodge,(v)}$, so that
$k=O_{\reflexHodge}/(\varpi)$.
Fix a geometric point $x\in\Sh_K(\barFp)$. 
Let $A:=\widehat{\mathscr O_{\Sh_K,x}}$ be the completed local ring at $x$. Write
$
\bar A:=A/\varpi A
\simeq A\otimes_{O_{\reflexHodge}}k$. 
Choose a trivialization of the $\CG$-torsor in the local model diagram \Cref{prop:local-model-diagram} on $A$. We thus have an isomorphism $\CW_{\DR,A}\cong W_{\Z_p}\otimes_{\Z_p} A$. 
 
Let $0\neq f\in\sEnd(x)$. Multiplying $f$ by an integer prime to $p$,
we may assume $f\in\End(\CA_x)$; this multiplication preserves its lifting
ideal because it is invertible in the $\Z_{(p)}$-isogeny category.
Let $I\subset A$ denote the lifting ideal of $f$.
A key point is to show that the universal lift of $f$ over $A/I$ is crystallinely special, i.e., its crystalline realization belongs to $\Gamma(((A/I)/\Z_p)_\cris,\CV_{\cris})$.

\begin{rmk}Even in the case \(F=\mathbb{Q}\), there is a gap at this step in the arguments given in \cite{MadapusiPera2016,Andreatta2018-tt}. 
In \cite[\S5.15]{MadapusiPera2016}, crystalline specialness is used over a PD
thickening, whereas \cite[Definition~5.11]{MadapusiPera2016} imposes it only on
geometric fibers. To deduce $q_\cris(f_\cris)=0$ over a testing ring $R$ from its vanishing at closed points, Dwork's trick below uses the hypothesis that $R/pR$ is a local complete intersection, which gives a $p$-torsion-free PD envelope in \Cref{lem:syntomic-Dwork}. Dwork's trick in
\cite[Lemma~4.3.5]{Andreatta2018-tt} is applied to a hypersurface deformation
ring. The reduction to rings of this kind uses 
\cite[Proposition~4.3.2 and Corollary~4.3.3]{Andreatta2018-tt}, whose proofs
refer back to the deformation results
\cite[Proposition~5.16 and Corollary~5.17]{MadapusiPera2016}. 

If $\FGlobal/\Q$ is unramified at $p$, the fiberwise-to-crystalline implication can be proved by the method of \cite{Madapusi2022Derived}.  
We use \Cref{lem:syntomic-Dwork} to prove crystalline specialness over complete-intersection testing rings and \Cref{prop:complete-intersection} to show that the lifting ideal is principal.
\end{rmk}

By \Cref{ass:d=1}, choose an isomorphism $\CV_{\dR,A}/\Fil^0\CV_{\dR,A}\simlr A$ and an extension
\[
\lambda:\CH_{\dR,A}/\Fil^0\CH_{\dR,A}\longrightarrow A.
\]

The kernel of $A/I^2\twoheadrightarrow A/I$ has its trivial divided
powers, compatible with those on $(p)$ by \Cref{lem:square-zero-PD}.
Let $e_f$ be the image of $f_{\cris,A/I^2}$ in
$\bigl(\CH_{\dR,A}/\Fil^0\CH_{\dR,A}\bigr)\otimes_A A/I^2$.
Since $f$ exists over $A/I$, it preserves the Hodge filtration there,
so
\[
e_f\in\bigl(\CH_{\dR,A}/\Fil^0\CH_{\dR,A}\bigr)\otimes_A I/I^2.
\]
Choose $\fkf\in I$ whose class in $I/I^2$ equals
$(\lambda\otimes1)(e_f)$, and write $\overline{\fkf}$ for
its image in $\bar A$. 

Whenever $f$ deforms over a quotient $B$ of $A$, let $\CQ_{\cris,B}$ be the
pullback of $\CQ_\cris$, let $q_{\cris,B}:\CH_{\cris,B}\to\CQ_{\cris,B}$ be the pullback
of $q_\cris$, and define its crystalline defect by
\[
\Delta_{f,B}:=q_{\cris,B}(f_{\cris,B})
\in\Gamma((B/\Z_p)_\cris,\CQ_{\cris,B}).
\]
Thus $\Delta_{f,B}=0$ if and only if $f_{\cris,B}\in\CV_{\cris,B}$.
Here $\CW_{\cris,B}=\BD(\CA_B[p^\infty])$ is the homological Dieudonn\'e crystal fixed in \Cref{sec:tensor}, and
$\CH_{\cris,B}=\Endint_{\CO_\cris}(\CW_{\cris,B})$.
We consider the Frobenius and Verschiebung operators:
\[
\begin{aligned}
F_{\CW_{\cris,B}}&:=p\varphi_{\CW_{\cris,B}}:F^*\CW_{\cris,B}\longrightarrow\CW_{\cris,B},\\
V_{\CW_{\cris,B}}&:=\varphi_{\CW_{\cris,B}}^{-1}:\CW_{\cris,B}\longrightarrow F^*\CW_{\cris,B}.
\end{aligned}
\]They satisfy
$F_{\CW_{\cris,B}}V_{\CW_{\cris,B}}=p$ and
$V_{\CW_{\cris,B}}F_{\CW_{\cris,B}}=p$. The endomorphism crystal has an
integral Frobenius linearization
\begin{equation}\label{eq:End-crystal-integral-Frobenius}
	\phi_{\CH}:F^*\CH_{\cris,B}\longrightarrow\CH_{\cris,B},\qquad
	u\longmapsto F_{\CW_{\cris,B}}\circ u\circ V_{\CW_{\cris,B}}.
\end{equation}
After inverting $p$, this is $p$ times the Frobenius on the
endomorphism isocrystal. It preserves the subcrystal $\CV_{\cris,B}$
and induces an integral linearization
$\phi_{\CQ}:F^*\CQ_{\cris,B}\to\CQ_{\cris,B}$, which is an isomorphism after
inverting $p$. 

An endomorphism of a $p$-divisible group commutes with its
Dieudonn\'e operators. Consequently,
$
\phi_{\CH}(F^*f_{\cris,B})=pf_{\cris,B}$. 
Passing to the quotient gives
\begin{equation}\label{eq:Delta-F-eigenvalue-p}
	\phi_{\CQ}(F^*\Delta_{f,B})=p\Delta_{f,B}.
\end{equation}
\begin{prop}\label{prop:tangent-space-divisor-special-fiber}
	Assume \Cref{ass:d=1}. If $\bar I$ denotes the image of $I$ in
	$\bar A$, then $
	\bar I=(\overline{\fkf})$. 
	Moreover, 	$f_\cris\in\CV_\cris$ on
	$((\bar A/\bar I)/\Z_p)_\cris$.
\end{prop}

\begin{proof}
	The closed-fiber endomorphism $f_{\barFp}$ is $p$-special, so
	$\Delta_{f,\barFp}=0$.  Let $J\supseteq\bar I$ be a
	complete-intersection ideal of the regular local ring $\bar A$.  Then
	$B:=\bar A/J$ is a complete local complete-intersection
	$\barFp$-algebra.
	The endomorphism $f$ exists over $B$, and
	\eqref{eq:Delta-F-eigenvalue-p} shows that $\Delta_{f,B}$ belongs to the
	$p$-eigenspace of the finite locally free $F$-crystal
	$(\CQ_{\cris,B},\phi_{\CQ})$.
	Therefore \Cref{lem:syntomic-Dwork} gives $\Delta_{f,B}=0$.
	
	For every intermediate ideal
	$J^2\subseteq K\subseteq J,$
	the surjection $\bar A/K\twoheadrightarrow\bar A/J$ is a square-zero
	PD-thickening, with the trivial PD-structure compatible with
	$(\Z_p,(p))$ by \Cref{lem:square-zero-PD}. By specialness, $\lambda$ identifies the
	Grothendieck--Messing obstruction of \Cref{thm:tangent-space} with the image of
	$\overline{\fkf}$ in $J/K$. Hence
	$
	\bar I\subseteq K
	\Leftrightarrow
	\overline{\fkf}\in K
	$ 
	for every intermediate ideal $J^2\subseteq K\subseteq J$.  We conclude that 
	$
	\bar I+J^2=(\overline{\fkf})+J^2.$ 
	By \Cref{prop:complete-intersection}, it follows that
	$\bar I=(\overline{\fkf})$.
	
	The quotient $\bar A/\bar I$ is therefore again a complete local
	complete intersection (with the zero ideal allowed).  Its closed-point
	restriction of $\Delta_{f,\bar A/\bar I}$ vanishes, while
	\eqref{eq:Delta-F-eigenvalue-p} places this section in the $p$-eigenspace.
	A final application of \Cref{lem:syntomic-Dwork} proves the last assertion.
\end{proof}
The possibility $\bar I=0$ will be excluded in
\Cref{prop:tangent-space-divisor-flat}.

\begin{lem}\label{lem:syntomic-Dwork}
	Let $(R,\mathfrak m)$ be a complete Noetherian local $\Z_p$-algebra
	with perfect residue field $\kappa$, and assume that
	$R_0:=R/pR$ is a local complete intersection over $\kappa$.  Let
	$(R/\Z_p)_\cris$ denote the formal crystalline site of
	$\Spf R$ with its $\mathfrak m$-adic topology, and let
	$(\CH,\phi_{\CH})$ be an \emph{effective} $F$-crystal on
	$(R/\Z_p)_\cris$, i.e., $\CH$ is a crystal in finite locally free
	$\CO_\cris$-modules and
	\[
	\phi_{\CH}:F^*\CH\longrightarrow\CH
	\]
	is a Frobenius linearization that becomes an isomorphism after inverting
	$p$.  For every integer $m\geq1$, restriction to the closed point is
	injective on the $p^m$-eigenspace:
	\[
	\left\{
	s\in\Gamma((R/\Z_p)_\cris,\CH):
	\phi_{\CH}(F^*s)=p^m s
	\right\}
	\lhook\joinrel\longrightarrow
	\Gamma((\kappa/\Z_p)_\cris,\CH_\kappa).
	\]
\end{lem}

\begin{proof}
	By \Cref{prop:crys-mod-p} we may replace $R$ by $R_0$.  The resulting
	assertion is proved in \cite[Proposition~2.2.7]{UrbanikYang2026Arithmetic} using Dwork's trick.
\end{proof}

\begin{lem}\label{lem:square-zero-PD}
	Assume $p>2$.  Let $R$ be a $\Z_p$-algebra and let $I\subset R$ satisfy
	$I^2=0$.  Then the trivial divided powers on $I$ are compatible with the
	canonical divided powers on $(p)\subset\Z_p$.
\end{lem}

\begin{proof}
	By the compatibility criterion, it is enough to check that, for every
	$y\in R$ with $py\in I$ and every integer $s\geq 2$, the canonical
	divided power $\delta_s(p)$ acts trivially on $y^s$.  Since
	$\delta_s(p)$ is divisible by $p^2$ for $p>2$, we have
	\[
	\delta_s(p)y^s
	=
	\frac{\delta_s(p)}{p^2}y^{s-2}(py)^2
	\in I^2=(0).
	\]
	Hence the trivial PD-structure on $I$ is compatible with the canonical
	PD-structure on $(p)$.
\end{proof}

\begin{lem}\label{lem:choose-complete-intersection}
	Let $(A,\mathfrak m)$ be a Noetherian regular local ring and
	$0\neq K\subsetneq A$. Then there exists a complete-intersection ideal
	$J\supseteq K$ such that $K\not\subseteq\mathfrak mJ$.
\end{lem}
\begin{proof}
	We argue by induction on $d=\dim A$. If $d=1$, then $A$ is a DVR, so
	$J=K$ works. For $d>1$, choose $0\neq g\in K$ and
	$x\in\mathfrak m\setminus\mathfrak m^2$ with $g\notin(x)$. Such an $x$
	exists because $A$ is a UFD: for a regular system of parameters
	$x_1,\ldots,x_d$, the prime elements $x_1+x_2^N$ are pairwise
	nonassociate, while $g$ has only finitely many prime divisors.
	
	Set $\overline A=A/(x)$ and
	$\overline K=(K+(x))/(x)\neq0$. By induction, there is a
	complete-intersection ideal $\overline J\supseteq\overline K$ such that
	$\overline K\not\subseteq\overline{\mathfrak m}\,\overline J$.
	Lifting a regular sequence generating $\overline J$ and adjoining $x$
	gives a complete-intersection ideal $J\supseteq K$; reducing modulo $x$
	shows that $K\not\subseteq\mathfrak mJ$.
\end{proof}
\begin{lem}\label{prop:complete-intersection}
	Let $(A,\mathfrak m)$ be a Noetherian regular local ring, let
	$I\subsetneq A$ be an ideal, and let $f\in I$. Suppose that
	\[
	I+J^2=(f)+J^2
	\]
	for every complete-intersection ideal $J\supseteq I$. Then $I=(f)$.
\end{lem}
\begin{proof}
	If $I=0$, there is nothing to prove. Assume $I\neq0$ and choose, by
	\Cref{lem:choose-complete-intersection}, a complete-intersection ideal
	$J_0\supseteq I$ such that $I\not\subseteq\mathfrak mJ_0$. Since
	$I+J_0^2=(f)+J_0^2$ and $J_0^2\subseteq\mathfrak mJ_0$, we have
	$f\in J_0\setminus\mathfrak mJ_0$. Hence $f$ can be completed to a regular
	sequence generating $J_0$, say
	\[
	J_0=(f,g_2,\ldots,g_r).
	\]
	If $r=1$, then $I\subseteq J_0=(f)$, and we are done.
	
	For $r\geq2$, successively square the generators $g_2,\ldots,g_r$
	cyclically. This gives complete-intersection ideals
	\[
	J_s=(f,g_2^{a_{2,s}},\ldots,g_r^{a_{r,s}}),\qquad s\geq0,
	\]
	with $a_{i,s}\to\infty$ and $(f)+J_s^2\subseteq J_{s+1}$. Starting from
	$I\subseteq J_0$, the hypothesis therefore implies inductively that
	$I\subseteq J_s$ and $I\subseteq(f)+J_s^2$ for every $s$.
	
	Since $g_2,\ldots,g_r\in\mathfrak m$ and every $a_{i,s}$ tends to
	infinity, for every $N\geq1$ we have
	$J_s\subseteq(f)+\mathfrak m^N$ for all sufficiently large $s$.
	The Krull intersection theorem in $A/(f)$ therefore gives
	\[
	I\subseteq\bigcap_{N\geq1}\bigl((f)+\mathfrak m^N\bigr)=(f).
	\]
	Since $f\in I$, we conclude $I=(f)$.
\end{proof}

\begin{lem}\label{lem:smooth-fiberwise-crystalline-special}
	Let $Y$ be a smooth scheme of finite type over a perfect field
	of characteristic $p$, mapping to the special fiber of $\Sh_K$.
	A $p$-special endomorphism of $\CA_Y$ is crystallinely special.
\end{lem}
\begin{proof}
	Its defect in $\CQ_\cris$ satisfies
	\eqref{eq:Delta-F-eigenvalue-p}. At every geometric closed point,
	the completed local ring of $Y$ is a power series ring over the
	residue field. The defect vanishes on this completion by
	\Cref{lem:syntomic-Dwork}. To check vanishing on $Y$, work
	\'etale locally with a smooth lift to Witt vectors and evaluate
	the crystal there. Its evaluation is finite locally free;
	vanishing can be checked on the faithfully flat completed local
	rings, modulo each power of $p$. Thus the defect is zero.
\end{proof}

\begin{lem}\label{lem:finite-unramified-dominant}
	Let $C$ be a smooth connected scheme of finite type over a field.
	If a finite unramified $C$-scheme $Y$ is connected and has an
	irreducible component dominating $C$, then $Y\to C$ is finite
	\'etale.
\end{lem}
\begin{proof}
	The scheme $C$ is regular and connected, hence integral and
	geometrically unibranch. By 
	\cite[\href{https://stacks.math.columbia.edu/tag/0GS9}{Tag~0GS9}]{stacks-project} $Y\to C$ is \'etale. 
\end{proof}
\begin{prop}
	\label{ass:crystalline-H0-vanishing}
	Under \Cref{ass:d=1}, for every finite extension $k'/k$,
	every connected component $C$ of $\Sh_{K,k'}$, and every
	finite \'etale morphism $Y\to C$, we have
	\[
	H^0_{\cris}(Y/W(k'),\CV_{\cris,Y}[1/p])=0.
	\]
\end{prop}
\begin{proof}
	Let $\reflexHodge'_v/\reflexHodge_v$ be the finite unramified extension
	with residue field $k'$, and write $O_{\reflexHodge',v}$ for its complete
	ring of integers. Base change $\Sh_K$ to $O_{\reflexHodge',v}$.
	By properness and \cite[\href{https://stacks.math.columbia.edu/tag/0A48}{Tag~0A48}]{stacks-project},
	$C$ lifts to a connected component $\bm{C}$ of
	$\Sh_{K,O_{\reflexHodge',v}}$, and $Y\to C$ lifts uniquely to a finite
	\'etale cover $\bm{Y}\to\bm{C}$.
	On the $p$-adic completion of $\bm{Y}$, the crystalline local system
	attached to $\CV_{\cris}$ is the pullback of $\CV_p$ by \Cref{sec:prism}.
	Applying \cite[Corollary~C]{guo2022prismatic} to this completion and using
	\cite[Theorem~3.7.2]{Huber1996} gives
	\[
	H^0_{\et}(\bm{Y}_{\overline{\Q}_p},\CV_p)\otimes_{\Z_p}\Bcris
	\simeq
	H^0_{\cris}(Y/W(k'),\CV_{\cris,Y})
	\otimes_{W(k')}\Bcris.
	\]
	By invariance of geometric finite \'etale covers under algebraically closed extensions of characteristic zero and the Riemann existence theorem, the invariants of $\CV_p$ on each geometric connected component of $\bm{Y}_{\overline{\Q}_p}$ are the invariants of a finite-index subgroup of a congruence subgroup $\Gamma\subseteq G(\Q)$.
	The Zariski closure of this subgroup contains $G^{\der}$ by
	\cite[Theorem~3.19]{MilneIntroductionShimura}, as $G^\ad$ has no compact rational simple factor. 
	Thus its invariants vanish by $V^{G^{\der}}=0$ in
	\Cref{ass:d=1}. 
\end{proof}

\begin{cor}\label{prop:crystalline-no-vertical-component}
	Assume \Cref{ass:d=1}. Let $k'/k$ be a finite extension, and let $C$ be a
	connected component of $\Sh_{K,k'}$. A connected component $Y$ of $\sEnd\times_{O_\reflexHodge}k'$, on which the tautological special endomorphism is nonzero,  cannot dominate $C$.
\end{cor}
\begin{proof}
	Suppose that $Y$ dominates $C$. The morphism $Y\to\Sh_{K,k'}$ is finite
	and unramified by \Cref{prop:endomorphism-Hilbert-finite}, hence
	finite \'etale by \Cref{lem:finite-unramified-dominant}. Thus $Y$ is
	smooth, and \Cref{lem:smooth-fiberwise-crystalline-special} gives a nonzero section
	\[
	u_\cris\in H^0_{\cris}(Y/W(k'),\CV_{\cris,Y}[1/p]).
	\] 
	This contradicts 
	\Cref{ass:crystalline-H0-vanishing}.
\end{proof}

\begin{prop}\label{prop:tangent-space-divisor-flat}
	Under \Cref{ass:d=1}, the element
	$\overline{\fkf}$ is a nonzerodivisor in $\bar A$, and
	$(\fkf,\varpi)$ is an $A$-regular sequence.  Consequently,
	$(\fkf)$ defines a relative Cartier divisor in $\Spec A$.
\end{prop}

\begin{proof}
	Assume $(x,f)$ is defined over a finite extension $k'/k$.
	By \Cref{ass:Shimura-data}\eqref{ass:local-model-smooth} and \Cref{prop:local-model-diagram}, $\Sh_K$ is
	smooth over $O_{\reflexHodge}$.  Thus $\bar A$ is a regular local
	ring. Choose the connected component $Y$ of $\sEnd\times_{O_\reflexHodge}k'$ containing
	$(x,f)$, which is a finite unramified scheme over $\Sh_{K, k'}$ by 
	 \Cref{prop:endomorphism-Hilbert-finite}. 
	By rigidity, completing $Y\to\Sh_{K,k'}$ at $(x,f)$ in the source and at $x$ in the target gives $\Spf \bar A/\bar I\to\Spf \bar A$.  
	If $\bar I=0$, the finite morphism
	$Y\to\Sh_{K,k'}$ dominates the component $C$ through $x$. 
	This contradicts \Cref{prop:crystalline-no-vertical-component}. 
	So $\bar I\ne0$.
	
	By \Cref{prop:tangent-space-divisor-special-fiber},
	$\bar I=(\overline{\fkf})$, so $\overline{\fkf}\ne0$.
	The regular local ring $\bar A$ is a domain, hence
	$\overline{\fkf}$ is a nonzerodivisor. Since $A$ is flat over the
	discrete valuation ring $O_{\reflexHodge}$, the element $\varpi$ is a
	nonzerodivisor in $A$, and $(\varpi,\fkf)$ is an $A$-regular sequence. 
	Thus $A/(\fkf)$ is flat over
	$O_{\reflexHodge}$, and $(\fkf)$ defines a relative Cartier divisor.
\end{proof}

\begin{thm}\label{prop:tangent-space-divisor}
	Assume
	\Cref{ass:d=1}; then $I=(\fkf)$. 
	In particular, the deformation locus of $f$ is an effective Cartier
	divisor, flat over $O_{\reflexHodge}$, and the universal lift of $f$ is
	crystallinely special.
\end{thm}

\begin{proof}
	The special-fiber equality of
	\Cref{prop:tangent-space-divisor-special-fiber} gives
	$
	I\subseteq I_0:=(\varpi,\fkf).
	$ 
	For $n\geq0$, put
	$
	I_n:=(\varpi^{2^n},\fkf).
	$
	We prove inductively that $I\subseteq I_n$.  Suppose this is known for
	$n$, and set $B_n:=A/I_n$.  If $e$ is the ramification index of
	$\widehat O_{\reflexHodge}/\Z_p$, then, up to a unit,
	$p=\varpi^e$, and therefore
	\[
	B_n/pB_n
	\simeq
	A/(\fkf,\varpi^{\min(2^n,e)}).
	\]
	By \Cref{prop:tangent-space-divisor-flat}, the displayed quotient is a
	complete local complete intersection.  Since it is an equicharacteristic
	$p$ complete local ring with residue field $\barFp$, Cohen structure makes
	it a local complete intersection over $\barFp$.  The section
	$\Delta_{f,B_n}$ vanishes at the closed point and satisfies
	$\phi_{\CQ}(F^*\Delta_{f,B_n})=p\Delta_{f,B_n}$ by
	\eqref{eq:Delta-F-eigenvalue-p}; hence
	\Cref{lem:syntomic-Dwork} gives $\Delta_{f,B_n}=0$.
	
	Note that 
	$
	I_{n+1}=(\fkf)+I_n^2
	$, so
	$A/I_{n+1}\twoheadrightarrow A/I_n$ is a square-zero thickening.  Its trivial
	divided powers are compatible with those on $(p)$ by
	\Cref{lem:square-zero-PD}.  Since $\Delta_{f,B_n}=0$, the naturality of the
	obstruction class $e_f$ identifies the Grothendieck--Messing
	obstruction in $I_n/I_{n+1}$ with the image of
	$\fkf$.  That image is zero by the definition of $I_{n+1}$, so $f$
	lifts across the thickening.  By the defining property of $I$, this means
	$I\subseteq I_{n+1}$.  
	
	Since $\fkf\in I$, the Krull intersection theorem gives
	$\bigcap_{n\geq0}(\varpi^{2^n})=0$ in
	$A/(\fkf)$, and hence
	\[
	I\subseteq\bigcap_{n\geq0}I_n=(\fkf)\subseteq I.
	\]
	Thus $I=(\fkf)$.  The flat Cartier assertion follows from
	\Cref{prop:tangent-space-divisor-flat}.  Finally,
	$(A/I)/p(A/I)\simeq A/(\fkf,\varpi^e)$ is, by the same Cohen
	structure argument, a local complete intersection over $\barFp$.
	The section $\Delta_{f,A/I}$ vanishes at the closed point and satisfies
	\eqref{eq:Delta-F-eigenvalue-p}; one last application of
	\Cref{lem:syntomic-Dwork} proves the crystalline specialness assertion.
\end{proof}

\begin{cor}\label{cor:crystalline-special}
	In the terminology of \Cref{defn:special}, $p$-specialness implies
	crystalline specialness under \Cref{ass:d=1}. 
\end{cor}

\begin{proof}
	The zero endomorphism is crystallinely special. For a connected component
	$S$ of $\sEnd$ carrying a nonzero tautological endomorphism,
	\Cref{prop:tangent-space-divisor} gives crystalline specialness on every
	completed local deformation space. Faithfully flat descent of its zero
	image in $\CQ_\cris$ gives the assertion on $S$ and hence on every pullback.
\end{proof}

\subsection{Hecke action on \texorpdfstring{$\sEnd$}{the special-endomorphism scheme}}
\begin{defn}\label{defn-Hecke-equivariant-sheaf}
	A \emph{Hecke-equivariant sheaf} on the Shimura tower $\Sh:=\varprojlim_{K^p\subseteq G(\BA_f^p)}\Sh_{K_pK^p}$ is a collection of fppf sheaves
	\(\{\CF_{K^p}\}_{K^p\subseteq G(\BA_f^p)}\), with $\CF_{K^p}$ on
	$\Sh_{K_pK^p}$, equipped, for every $g\in G(\BA_f^p)$ and every
	$K_1^p\subseteq gK^pg^{-1}$, with identifications $$[\cdot g]^{-1}\CF_{K^p}\xrightarrow{\ \sim\ }\CF_{K_1^p}$$
	satisfying the identity and cocycle conditions.  Here
	$[\cdot g]:\Sh_{K_pK_1^p}\to\Sh_{K_pK^p}$ is induced by right
	multiplication by $g$. 
	For each prime $\ell\ne p$ and each $U\in\mathrm{Rep}_{\Q_\ell}(G)$, the associated $\ell$-adic
	local systems $\CU_\ell=\{\CU_{\ell,K^p}\}$ form a Hecke-equivariant
	sheaf.
\end{defn}

\begin{prop}\label{prop:sEnd-Hecke-equivariant}
	The sheaves $\sEnd$ on the schemes $\Sh_{K_pK^p}$ form a
	Hecke-equivariant sheaf.  For every prime $\ell\neq p$, the
	$\ell$-adic realization map $\realization_\ell: \sEnd\hookrightarrow \CV_\ell$	is Hecke equivariant, and this compatibility determines the Hecke-equivariant structure uniquely.
\end{prop}

\begin{proof}
	Let $g\in G(\BA_f^p)$ and let $K_1^p\subseteq gK^pg^{-1}$ be neat.
	The Hecke morphism
	\[
	[\cdot g]:\Sh_{K_pK_1^p}\longrightarrow\Sh_{K_pK^p}
	\]
	sends $[x,a]$ to $[x,ag]$ on the generic fiber.
	
	Via the fixed Hodge embedding, the prime-to-$p$ isogeny interpretation of
	the Siegel integral moduli problem supplies a tautological prime-to-$p$
	quasi-isogeny
	\begin{equation}\label{eq:Hecke-tautological-quasi-isogeny}
		\nu_g:\CA|_{\Sh_{K_pK_1^p}}\dashrightarrow
		[\cdot g]^*(\CA|_{\Sh_{K_pK^p}})
	\end{equation}
	in the $\Z_{(p)}$-isogeny category. The construction is made at prime-to-$p$ level and descends to
	the indicated neat finite levels.  Since $g$ lies in $G(\BA_f^p)$, it
	preserves the tensors defining the Hodge-type locus, so the Siegel
	quasi-isogeny restricts to the Hodge-type integral models.  See
	\cite[\S\S2.1--2.2]{Kisin2010} and
	\cite[\S\S4.2--4.3]{KisinPappas2018}.
	
	Conjugation by \eqref{eq:Hecke-tautological-quasi-isogeny} gives an
	isomorphism of ambient endomorphism sheaves
	\begin{equation}\label{eq:Hecke-sEnd-conjugation}
		\begin{aligned}
			\Theta_g:
			[\cdot g]^{-1}\Endint(\CA|_{\Sh_{K_pK^p}})_{(p)}
			&\xrightarrow{\ \sim\ }
			\Endint(\CA|_{\Sh_{K_pK_1^p}})_{(p)},\\
			f&\longmapsto \nu_g^{-1}\circ f\circ\nu_g.
		\end{aligned}
	\end{equation}
	The inclusion $V\hookrightarrow\End(W)$ is a morphism of
	$G$-representations. For every $\ell\ne p$, functoriality of the
	associated local systems shows that conjugation by $\nu_{g,\ell}^{-1}$
	carries $[\cdot g]^{-1}\CV_\ell$ to $\CV_\ell$ on $\Sh_{K_pK_1^p}$,
	compatibly with $\realization_\ell$. By
	\Cref{lem:special-char-0-fiberwise}, $\Theta_g$ therefore preserves
	specialness on the generic fiber.
	Since $\nu_g$ is prime-to-$p$, it induces an isomorphism of the
	corresponding $p$-divisible groups. Functoriality and uniqueness of the
	prismatic lift in \Cref{prop:cris-tensor,thm:para-equiv} identify their
	$\CG$-crystals and the subcrystals associated with $V_{\Z_p}$.
	Thus $\Theta_g$ also preserves $p$-specialness on the special fiber.
	Applying the same argument to $\nu_g^{-1}$ shows that $\Theta_g$
	restricts to an isomorphism
	\[
	\theta_g:[\cdot g]^{-1}(\sEnd|_{\Sh_{K_pK^p}})
	\xrightarrow{\ \sim\ }\sEnd|_{\Sh_{K_pK_1^p}}.
	\]
	Compatibility with $\realization_\ell$ and its injectivity give the
	identity and cocycle conditions and the uniqueness of this
	Hecke-equivariant structure.
\end{proof}

\subsection{Special endomorphisms for abelian type}
\label{subsec:special-endomorphisms-abelian-type}

We now construct special endomorphisms on the abelian-type integral model by
descending the special endomorphisms on the Hodge-type integral model defined in \Cref{subsec:special-endomorphisms-Hodge-type}.  The Hodge-type sheaf $\sEnd$ is already compatible 
with $G(\BA_f^p)$-Hecke transport, the only additional compatibility needed in
the quotient presentation is the $\CG_2(\Z_(p))$-action, which can be defined using the twisting on the universal abelian scheme. 
In this subsection we always assume the following.
\begin{ass}\label{ass:J=1}
	The set $J_G$ in the construction of \Cref{defn:int-model} for $(G,X,\CG)$ satisfies $J_G=\{1\}$. 
\end{ass}
\begin{rmk}
	For the group $G^\sharp$ defined in \Cref{defn:G-sharp}, one has
	$J_{G^\sharp}=\{1\}$ by \Cref{defn:orth-proj}. The auxiliary fiber-product groups
	used in \Cref{subsec:indep-hodge} may have several $J$-classes. We only use the condition $J_G=1$ when we descend $\sEnd$ from the Hodge type Shimura variety $\Sh_K$ to $\Sh_{K_2}$. 
\end{rmk}

\subsubsection{Twisting abelian schemes}
We recall the twisting of abelian schemes in \cite[\S4.4]{KisinPappas2018}.
Let $Z$ be a flat affine group scheme of finite type over
$\Spec\Z_{(p)}$, and let $\CP$ be a $Z$-torsor.  Since affineness is fpqc
local on the base, $\CP$ is affine.
We write $O_{Z}$ and $O_{\CP}$ for the coordinate rings of $Z$ and $\CP$, respectively.  
A $Z$-action on a $\Z_{(p)}$-module $M$ amounts to an $O_Z$-comodule structure $\alpha: M\to M\otimes_{\Z_{(p)}}O_Z$. 
We define the $\Z_{(p)}$-submodule $M^Z\subseteq M$ by $M^Z=\{m\in M:\alpha(m)=m\otimes1\}$. 

For a scheme $S$, we define the $\Z_{(p)}$-isogeny category of abelian schemes 
over $S$ to be the category obtained from the category of abelian schemes over $S$ by tensoring the $\Hom$ groups
with $\Z_{(p)}$ over $\mathbb Z$. An object $\CA$ in this category is called an {\it abelian scheme up to $\Z_{(p)}$-isogeny} over $S.$ 
For $T$ an $S$-scheme we set $\CA(T) = \Hom_S(T,\CA)\otimes_{\mathbb Z} \Z_{(p)}.$ 

Let $\CA$ be an abelian scheme up to $\Z_{(p)}$-isogeny over $S.$ 
Denote by $\sAut_{\Z_{(p)}}(\CA)$ the $\Z_{(p)}$-group whose points in a $\Z_{(p)}$-algebra $R$ are given by 
$$\sAut_{\Z_{(p)}}(\CA)(R) = ((\End_S \CA)\otimes_{\Z_{(p)}}R)^\times.$$ 

Let $Z$ and $\CP$ be as above, and suppose that we are given a map of $\Z_{(p)}$-groups 
$Z \rightarrow \sAut_{\Z_{(p)}}(\CA).$ We define a presheaf $\CA^{\CP}$ in the fppf topology of $S$ by setting 
$\CA^{\CP}(T) = (\CA(T)\otimes_{\Z_{(p)}}O_{\CP})^Z.$
By \cite[Lemma 4.4.6]{KisinPappas2018}  $\CA^{\CP}$ is a sheaf, represented by an abelian scheme up to $\Z_{(p)}$-isogeny.

If $Z$ is commutative, then $\CA^{\CP}$ inherits a $Z$-action: $Z$ acts
on $\CA(T)$ in the given way and trivially on $O_\CP$.  This action
commutes with the descent action on
$\CA(T)\otimes_{\Z_{(p)}}O_\CP$ and therefore descends to
$\CA^{\CP}$.  Moreover,
\((\CA^\CP)^{\CP^{-1}}\simeq\CA\) by
\cite[Lemma~4.4.3]{KisinPappas2018}.

If $\CB$ is another abelian scheme with a $Z$-action, every
$Z$-equivariant morphism $f:\CA\to\CB$ twists to a $Z$-equivariant
morphism $f^\CP:\CA^\CP\to\CB^\CP$.  This gives the following
descent statement.
\begin{lem}\label{lem:twisting-bij}
	Assume that $Z$ is commutative, and let $\CA$ and $\CB$ be abelian schemes equipped with $Z$-actions.
	Twisting induces a natural bijection
	\[
	\Hom_{\Z_{(p)},Z}(\CA,\CB)
	\xrightarrow{\ \sim\ }
	\Hom_{\Z_{(p)},Z}(\CA^\CP,\CB^\CP),
	\qquad f\longmapsto f^\CP.
	\]
\end{lem}

\begin{proof}
	Twisting by $\CP^{-1}$ is inverse to twisting by $\CP$, under the
	canonical identifications
	$(\CA^\CP)^{\CP^{-1}}\simeq\CA$ and
	$(\CB^\CP)^{\CP^{-1}}\simeq\CB$.  Functoriality of fpqc descent then
	gives the asserted bijection.
\end{proof}

\begin{prop}
	Let $\CA$ be an abelian scheme equipped with commuting $Z$- and
	$O_\FGlobal$-actions as above, with $Z$ commutative. For the assertion
about $h$, assume that the $O_\FGlobal$-action is self-adjoint for the weak
polarization $\lambda$ and that $Z$ acts by similitudes of $\lambda$.
Equip $\CA^\CP$ with the induced weak polarization and use
$h(f,g)=N^{-1}\tr_\FGlobal(f\circ g^\dagger)$ as in \Cref{defn-pairing-special-SO}.
\begin{enumerate}
		\item Let $f:\CA\to\CA$ be a $Z$-equivariant $O_\FGlobal$-linear
		endomorphism.  Then the twist
		$f^\CP:\CA^\CP\to\CA^\CP$ defined in \Cref{lem:twisting-bij} is
		$O_\FGlobal$-linear and
		$\tr_\FGlobal(f)=\tr_\FGlobal(f^\CP)$.
		
		\item Let $f,g:\CA\to\CA$ be $Z$-equivariant $O_\FGlobal$-linear
		endomorphisms.  Then
		$h(f,g)=h(f^\CP,g^\CP).$
	\end{enumerate}
\end{prop}

\begin{proof}
	\begin{enumerate}
		\item After extending the coefficients of the $\Z_{(p)}$-linear
		isogeny sheaves to $O_\CP$, the twisting isomorphism identifies $f$
		with $f^\CP$ and commutes with $O_\FGlobal$. The induced conjugation
		on any rational Tate realization preserves its relative trace.
		Faithful flatness of $O_\CP/\Z_{(p)}$ gives the asserted equality.
		Here tensoring an abelian scheme in the isogeny category means
		extension of coefficients in its sheaf of $\Z_{(p)}$-modules.

	\item The weak polarization on $\CA^\CP$ is defined by
	\cite[Lemma (3.1.5)]{Kisin2010}. Under the twisting identifications
	\[
	\CA\otimes_{\Z_{(p)}}\CO_\CP\to\CA^\CP\otimes_{\Z_{(p)}}\CO_\CP,\qquad
	\CA^*\otimes_{\Z_{(p)}}\CO_\CP\to(\CA^\CP)^*\otimes_{\Z_{(p)}}\CO_\CP
	\]
	the map $\lambda^\CP\otimes1$ corresponds to
	$\lambda\otimes r$ for $r\in\CO_\CP(\CP)^\times$.  Then
	\begin{align*}
		(f^\CP)^\dagger\otimes1
		&=(\lambda\otimes r)^{-1}
		\circ(f^*\otimes1)\circ(\lambda\otimes r)\\
		&=f^\dagger\otimes1.
	\end{align*}
	Thus $(f^\CP)^\dagger=(f^\dagger)^\CP$. Applying part (1) to
	$f\circ g^\dagger$ gives equality of their relative traces, hence
	the claimed equality of $h$.
	\end{enumerate}
\end{proof}

\subsubsection{\texorpdfstring{$\CG_2(\Z_{(p)})$}{G2(Z(p))}-action}\label{subsubsec:adjoint-act}
We modify the $\CG^\ad(\Z_{(p)})$-action defined in \cite[\S4.5]{KisinPappas2018}. 
By \cite[Lemma~4.5.2]{KisinPappas2018} the universal abelian scheme $\CA$ on $\Sh_K$ admits a natural $Z(\CG)$-action
\(Z(\CG) \rightarrow \sAut_{\Z_{(p)}}(\CA).\)
We restrict this action to $\kernelalpha=\ker(\alpha)\subseteq Z(\CG)$.
Through the fixed Hodge embedding it acts by polarization similitudes.
Hence, given a $\kernelalpha$-torsor $\CP$, we can define the twisted abelian scheme $\CA^\CP$. 

Recall that we have the following exact sequence of parahoric
models from \Cref{rmk:R-smooth-seq},
\[
1\longrightarrow\kernelalpha\longrightarrow\CG
\longrightarrow\CG_2\longrightarrow1.
\]
Let $\gamma \in \CG_2(\Z_{(p)})$, and let $\CP$ be the fiber of $\CG \rightarrow \CG_2$ 
over $\gamma$. From the polarization $\lambda$ and level structure $\varepsilon$ on $\CA$, we can use \cite[Lemma 4.4.8, Lemma 4.5.4]{KisinPappas2018} to produce the polarization $\lambda^\CP$  and level structure $\varepsilon^\CP$ on $\CA^\CP$. 

\begin{prop}[{\cite[Lemma 4.5.7]{KisinPappas2018}}]\label{adjointaction} For $T$ a $\Z_{(p)}$-scheme and $x \in \Sh_K(T),$ the assignment 
	$$ (\CA_x, \lambda_x, \varepsilon^p_x) \mapsto  (\CA_x^\CP, \lambda_x^\CP, \varepsilon^{p,\CP}_x)$$
	induces a map 
	$$ \Sh_{K_p}(G,X) \rightarrow \Sh_{\gamma K_p \gamma^{-1}}(G, \gamma X \gamma^{-1})$$
	whose generic fiber agrees with the map induced by conjugation by $\gamma.$
	When $\gamma\in  \CG_2(\Z_{(p)})_+$, we may identify $\gamma X \gamma^{-1}$ with $X$. 
\end{prop}

\begin{rmk}
	The inherited $\kernelalpha$-action of \Cref{lem:twisting-bij} agrees with
	the central action on the universal abelian scheme at the twisted point
	$\gamma.x$ of \Cref{adjointaction}. After base change to $\CP$, both
	identify with the original $\kernelalpha$-action on $\CA_x$; equality
	follows by fpqc descent.
\end{rmk}

Let $\gamma \in \CG_2(\Z_{(p)}),$ and let $\CP$ be the fiber of
$\CG\rightarrow\CG_2$ over $\gamma$, and let $x\in\Sh_K(T)$.
By \cite[Lemma 4.4.3]{KisinPappas2018} we have an isomorphism
\[\iota: \CA_x^\CP\otimes_{\Z_{(p)}} O_\CP \to \CA_x\otimes_{\Z_{(p)}} O_\CP\]
which induces isomorphisms
\[\iota_\diamond : \CW^\otimes_{\diamond,\gamma.x}\otimes_{\Z_{(p)}} O_\CP \to \CW^\otimes_{\diamond,x} \otimes_{\Z_{(p)}} O_\CP, \diamond\in \{\betti, \ell, \dR, \cris\}.\]
The Betti maps are taken on $T_\C$, the $\ell$-adic maps on the \'etale
site of $T$ (on $T_{\reflexHodge}$ for $\ell=p$), the de Rham maps on $T$,
and the crystalline maps on the crystalline site after $p$-adic completion.
\begin{lem}\label{lem:tensor-inv}
	The isomorphisms $\iota_\diamond$ map $s_{\alpha,\diamond}\otimes 1$ to $s_{\alpha,\diamond}\otimes 1$ for any $\diamond\in \{\betti, \ell, \dR, \cris\}$. 
\end{lem}

\begin{proof}
	First consider the Betti case. It suffices to show that for $x\in \ShG_K(\C)$, any finite Galois extension $L/\Q$ and $\tilde{\gamma}\in \CP(L)\subseteq G(L)$, the map
	\[\iota_{\tilde{\gamma}, \Betti}: \CW^\otimes_{\Betti,\gamma.x}\otimes_{\Z_{(p)}} L \to \CW^\otimes_{\Betti,x} \otimes_{\Z_{(p)}} L\]
	maps $s_{\alpha,\Betti}\otimes 1$ to $s_{\alpha,\Betti}\otimes 1$. 
	Choose a lift $\tilde{x}\in X\times G(\BA_f)$. Applying the $G^{\ad}(\Q)$-action gives a lift $\widetilde{\gamma.x}$ of $\gamma.x$. The two lifts induce trivializations
	$\CW_{\Betti,x}\otimes_{\Z_{(p)}}L\simlr W\otimes_{\Q}L$ and $\CW_{\Betti,\gamma.x}\otimes_{\Z_{(p)}}L\simlr W\otimes_{\Q}L$, under which the tensors $s_{\alpha,\Betti}\otimes 1$ correspond to $s_\alpha\otimes 1\in W^\otimes\otimes_{\Q}L$. As explained in \cite[Lemma (3.2.6)]{Kisin2010}, we therefore have a commutative diagram:
	\[\begin{tikzcd}
		{W^\otimes\otimes_{\Q} L} \arrow[d, "\sim"] \arrow[r, "\tilde{\gamma}^{-1}"]                   & {W^\otimes\otimes_{\Q} L} \arrow[d, "\sim"] \\
		{\CW^\otimes_{\Betti,\gamma.x}\otimes_{\Z_{(p)}} L} \arrow[r, "{\iota_{\tilde{\gamma}, \Betti}}"] & {\CW^\otimes_{\Betti,x} \otimes_{\Z_{(p)}} L}.   
	\end{tikzcd}\]
	The claim then follows from the fact that $\tilde{\gamma}^{-1}s_\alpha=s_\alpha$. 
	
	The claim for $\ell$-adic and de Rham tensors follows because they are realizations of absolute Hodge cycles. The claim for crystalline tensors follows from the crystalline--\'etale comparison map.
\end{proof}

\begin{cor}\label{cor:twisting-special}
	Let $\gamma\in\CG_2(\Z_{(p)})$, let $\CP'$ be its fiber under
	$\CG\to\CG_2$, and let $x\in\Sh_K(T)$. Identify
	$\CA_x^{\CP'}\simeq\CA_{\gamma.x}$ by \Cref{adjointaction}.
	Then the twisting bijection \Cref{lem:twisting-bij}
	\[\Hom_{\Z_{(p)},\kernelalpha}(\CA_x, \CA_x)\simlr \Hom_{\Z_{(p)},\kernelalpha}(\CA_x^{\CP'}, \CA_x^{\CP'}), f\mapsto f^{\CP'}\] 
	preserves special morphisms, i.e., $f$ is special if and only if $f^{\CP'}$ is special. 
\end{cor}

\begin{proof}
	By \Cref{ass:special-cycle}, $\kernelalpha$ acts trivially on $V$.
	It follows from \Cref{lem:special-Z-inv} that every special morphism
	is $\kernelalpha$-invariant. 
	
	Since $G=\CG_\Q$ is reductive, the $\CG_\Q$-subrepresentation $V_\Q \subseteq \End_{\Q}(W_{\Q})$ is a direct summand; we choose a $G$-invariant idempotent $\bpi: \End_{\Q}(W_{\Q})\to \End_{\Q}(W_{\Q})$ with image $V_\Q$. Choose $s\in \BZ_{\geq 0}$ so that $p^s \bpi$ defines a map $\End_{\Z_{(p)}}(W_{\Z_{(p)}})\to \End_{\Z_{(p)}}(W_{\Z_{(p)}})$. We regard $p^s \bpi$ as a $\CG$-invariant tensor and we may assume $p^s \bpi= s_\beta$ for an index $\beta$. 
	
	The locally free quotient in \eqref{eq:special-crystal-quotient} makes
	$\CV_{\cris}$ saturated in $\CH_{\cris}$ after evaluation on $W(\barFp)$.
	Thus \Cref{defn:special} says that $f$ is special if and only if $s_{\beta, \cris}f_\cris = p^s f_\cris$ at each $\barFp$-point and $s_{\beta, \Betti}f_\Betti = p^s f_\Betti$ at each $\C$-point. Then the claim follows from \Cref{lem:tensor-inv}.
\end{proof}

We retain the notation of \Cref{ass:Shimura-data}.  In particular,
$\kernelalpha=\ker(\alpha)$, and the parahoric models attached to $G$, $G_2$, and
an auxiliary Hodge-type datum will always be denoted by $\CG$, $\CG_2$, and
$\wt{\CG}$, respectively.
By \Cref{prop:Kisin-abelian-(3.3.10)}, there is a finite \'etale morphism
\[
q:\Sh_K\longrightarrow
\Sh_{K_2}\times_{O_{\reflexAbel,(v_2)}}O_{\reflexHodge,(v)}.
\]
\begin{prop}\label{prop:main-sEnd-abel}
	{There exists a unique sheaf $\AbSpecialEnd$ on 
		$\Sh_{K_{2,p}}(G_2,X_2)_{O_{\reflexHodge}}$ equipped
		with a Hecke action by $G_2(\BA_f^p)$ such that, for every $\ell\neq p$, there
		is an injective Hecke-equivariant realization map
		$\realization_\ell:\AbSpecialEnd\longrightarrow\CV_\ell,$
		and, via these realization maps, there is a canonical identification
		$q^{-1}\AbSpecialEnd=\sEnd.$  The sheaf is represented by a scheme
		unramified and locally of finite presentation over the base, and each of
		its connected components is finite over the corresponding Shimura model.}
\end{prop}
\begin{proof}
	We use the left-coset presentation obtained from that of
	\Cref{defn:int-model} by $(z,a)\mapsto(a^{-1},z)$.
	Accordingly, $g\cdot y$ denotes $[\cdot g^{-1}](y)$ below.
	For each $T\in\FGlobal$ and $K^p$-invariant compact subset
	$\Omega\subset V(\BA_f^p)$, we construct a finite scheme
	$\CZ_2(T,\Omega)$ whose realization image will enter the
	definition of $\AbSpecialEnd$.
	Fix a sufficiently small neat prime-to-$p$ level $K_2^p$ preserving $\Omega$.
	We work on
	$\Sh_{K_{2,p}K_2^p}(G_2,X_2)_{O_{\reflexHodge}}$,
	with Hecke equivariance understood levelwise as in
	\Cref{defn-Hecke-equivariant-sheaf}.

	We first prove uniqueness.  Fix a prime $\ell\neq p$.  By the proof of
	\cite[Lemma~(2.2.5)]{Kisin2010} and
	\Cref{prop:Kisin2010-2.2.6}, Hecke translates of the image of the
	finite \'etale map $q_{K^p,K_2^p}$ cover this model, with the source level
	chosen compatibly for each translate.  Consequently a
	$G_2(\BA_f^p)$-Hecke-equivariant subsheaf of
	$\CV_\ell$ is determined by its pullback along $q$.  The injective
	$\ell$-adic realization maps identify any two sheaves satisfying the
	proposition with two such subsheaves whose pullbacks are both
	$\realization_\ell(\sEnd)$.  They are therefore equal.

	Since the finite-level model $\Sh_{K_2}$ is of finite type, choose finitely many
	$g_1,\ldots,g_r\in G_2(\BA_f^p)$ and a sufficiently small common neat level $K^p\subseteq G(\BA_f^p)$ at which all
	$g_i\circ q$ are defined such that the images of $g_i\circ q$ cover $\Sh_{K_2}$.

	Set
	$Y:=\coprod_{i=1}^r\Sh_{K_pK^p}(G,X)$ and define
	\[
	\rho:Y\longrightarrow
	\Sh_{K_{2,p}K_2^p}(G_2,X_2)_{O_{\reflexHodge}},
	\qquad
	(i,x)\longmapsto g_i\cdot q(x).
	\]
	By \Cref{prop:Kisin-abelian-(3.3.10)}, each $g_i\circ q$ is finite
	\'etale.  Thus $\rho$ is a finite \'etale surjective morphism between
	schemes of finite type, hence an fppf covering.
	
	By \Cref{prop:fixed-level-special-cycle-finite} each $\CZ(T, g_i^{-1}\Omega)$ is finite over $\Sh_{K_pK^p}(G,X)$. 
	Let
	\[
	\mathscr L:=\coprod_{i=1}^r\CZ(T,g_i^{-1}\Omega)
	\]
	be the finite scheme over $Y$.
	Since the $G$-action on $V$ factors
	through $G_2$ by \Cref{ass:special-cycle}, it defines an automorphic local system $\CV_\ell$ on $\Sh_{K_{2,p}K_2^p}(G_2,X_2)$.  For every
	$\ell\neq p$, there is an injective map
	\(\iota_\ell:\mathscr L\longrightarrow\rho^*\CV_\ell\)
	which on the $i$-th summand is $f\longmapsto g_i\cdot\realization_\ell(f)$.
	
	We now show that the descent datum for $\rho^*\CV_\ell$ restricts to a descent datum for $\mathscr L$. 
	We may check this after the faithfully flat base
	change $O_{\reflexHodge}\to O_{\breve{(\reflexHodge)_v}}$, where $\breve{(\reflexHodge)_v}$ is the completion of the maximal unramified extension of $(\reflexHodge)_v$.  
	Since $J_G=1$ by \Cref{ass:J=1} and
	$J_{G_2}=1$ by \Cref{lem:J2-trivial}, the quotient presentations of
	\Cref{defn:int-model} give
	\begin{equation*}
		\begin{aligned}
			\Sh_{K_{2,p}K_2^p}(G_2,X_2)
			\times_{O_{\reflexAbel}}O_{\breve{(\reflexHodge)_v}}
			&\simeq
			[K_2^p\backslash\sA(\CG_2)\times
			\Sh^+_{\wt{K}_p,O_{\breve{(\reflexHodge)_v}}}]
			/\sA(\wt{\CG})^\circ,\\
			\Sh_{K_pK^p}(G,X)
			\times_{O_{\reflexHodge}}O_{\breve{(\reflexHodge)_v}}
			&\simeq
			[K^p\backslash\sA(\CG)\times
			\Sh^+_{\wt{K}_p,O_{\breve{(\reflexHodge)_v}}}]
			/\sA(\wt{\CG})^\circ.
		\end{aligned}
	\end{equation*}
	On an overlap, pull back to the quotient charts and choose
	representatives, using the $K_2^p$-action and the diagonal
	$\sA(\wt{\CG})^\circ$-action, with the same connected coordinate and
	equal labels in $\sA(\CG_2)$. It suffices to compare the pullback of $\mathscr L$ to the charts
	\[
	\{i\}\times\{(a,b)\}\times
	\Sh^+_{\wt K_p,O_{\breve{(\reflexHodge)_v}}},
	\qquad
	\{j\}\times\{(a',b')\}\times
	\Sh^+_{\wt K_p,O_{\breve{(\reflexHodge)_v}}}.
	\]
	Here $a,a'\in G(\BA_f^p)/\centralHodge(\Z_{(p)})^-$ and
	$b,b'\in\CG^\ad(\Z_{(p)})^+$.
	Recall the homomorphisms
	\[
	\begin{aligned}
	\psi_2:\CG_2(\Z_{(p)})_+/\centralAbel(\Z_{(p)})
	&\longrightarrow G_2(\BA_f^p)/\centralAbel(\Z_{(p)})^-,\\
	\phi_2:\CG_2(\Z_{(p)})_+/\centralAbel(\Z_{(p)})
	&\longrightarrow\CG_2^\ad(\Z_{(p)})^+.
	\end{aligned}
	\]
	The $\ast$-product relation gives $\gamma$ in their common domain with
	\[
	g_i\alpha(a)\ad_b(\psi_2(\gamma))=g_j\alpha(a'),
	\qquad b\phi_2(\gamma^{-1})=b',
	\]
	where the first equality holds in $G_2(\BA_f^p)/\centralAbel(\Z_{(p)})^-$. 
	Choose a representative $\widetilde\gamma\in\CG_2(\Z_{(p)})_+$.
	Apply the twisting morphism of \Cref{adjointaction} for the fiber of
	$\CG\to\CG_2$ over $\ad_b(\widetilde\gamma^{-1})$.
	It changes the labels by
	\[
	(g_i,(a,b))\longmapsto
	\left(g_i\ad_b(\widetilde\gamma),
	\left(\ad_{\ad_b(\widetilde\gamma^{-1})}(a),b'\right)\right).
	\]
	By \Cref{cor:twisting-special}, this twisting preserves special
	endomorphisms and acts on their realizations by
	$v\mapsto\ad_b(\widetilde\gamma^{-1})v$.
	After applying this isomorphism, we may assume $\widetilde\gamma=1$,
	so $g_i\alpha(a)=g_j\alpha(a')$ and $b=b'$.

	Choose lifts of $a,a'$ in $G(\BA_f^p)$ and put
	$h=a'a^{-1}\in G(\BA_f^p)$.
	At a common sufficiently small level, $[\cdot h^{-1}]$ changes the
	labels $(g_i,(a,b))$ further to
	\[
	\left(g_i \alpha(h)^{-1},(a',b')\right).
	\]
	By \Cref{prop:sEnd-Hecke-equivariant}, the Hecke transport by $h$  on
	realizations is $v\mapsto\alpha(h)v$. 
	After this transport, we may assume the two charts have the same labels $a=a'$ and
	$b=b'$, and hence $z:=g_i^{-1}g_j\in\centralAbel(\Z_{(p)})^-
	=\centralAbel(\Z_{(p)})$.
	Here the last equality follows from the Kottwitz condition 
	\Cref{ass:Shimura-data}(\ref{ass:Shimura-data:kottwitz}).

	The two Hecke transports given by $g_i, g_j$ on $\CV_\ell$ now differ by the
	representation action of $z$. We construct the compatible action
	of $\centralAbel(\Z_{(p)})$ on $\sEnd$.
	Take $f\in \sEnd(S)$ where $S$ is a scheme of finite type over $\Sh_K$. 
	By \cite[Lemma~4.5.2]{KisinPappas2018} we have an action
	$\centralHodge\to \Aut(\CA)_{(p)}$. 
	Consider the finite $\Z_{(p)}$-module $\End(\CA_S)_{(p)}$, regarded as a $\Z_{(p)}$-scheme. 
	Then the conjugation on $f\in \End(\CA_S)_{(p)}$ defines a map between $\Z_{(p)}$-schemes
	\[\tau_f: \centralHodge\to \End(\CA_S)_{(p)}, \qquad z\mapsto z\circ f \circ z^{-1}.\]
	The exact sequence of \Cref{rmk:R-smooth-seq} induces the fppf exact sequence
	$1\to\kernelalpha\to\centralHodge\to\centralAbel\to1$.
	Since $f$ is $\kernelalpha$-invariant, $\tau_f$ factors through $\centralAbel$.
	It satisfies $\realization_\ell(\tau_f(z))=z\cdot \realization_\ell(f)$. 
	We define the action map to be
	\[z: \sEnd(S)\to \sEnd(S), \qquad f\mapsto \tau_f(z),\]
	which maps the open and closed subsheaf $\CZ(T,\Omega)$ to
	$\CZ(T,z\Omega)$.

	The chart identifications induce canonical transport in
	$\rho^*\CV_\ell$, so injectivity of $\iota_\ell$ proves their
	independence of representatives and the cocycle condition, and
	descends them from the faithfully flat completed base change.
	Effective fppf descent for finite schemes gives
	$\CZ_2(T,\Omega)$ over $\Sh_{K_2}$ with pullback $\mathscr L$
	and an injective realization into $\CV_\ell$.
	The inclusions for $\Omega\subseteq\Omega'$ commute with this
	datum and descend to
	$\CZ_2(T,\Omega)\subseteq\CZ_2(T,\Omega')$.
	Define $\AbSpecialEnd$ as the union of these realization images
	in the category of subsheaves of the same $\CV_\ell$.
	Its pullback to each summand of $Y$ identifies with $\sEnd$,
	with realization transported by $g_i$; taking $g_1=1$ gives
	$q^{-1}\AbSpecialEnd\simeq\sEnd$. 
	For $g\in G_2(\BA_f^p)$, replacing $g_i,\Omega$ by $gg_i,g\Omega$
	at compatible neat finite levels transports the covering and its
	descent datum by $g$. This gives the $G_2(\BA_f^p)$-Hecke action on
	$\AbSpecialEnd$; the identity and cocycle conditions follow from those
	on $\CV_\ell$ and the injectivity of $\realization_\ell$.
\end{proof}

\begin{defn}\label{defn:special-abelian}
	The sheaf $\AbSpecialEnd$ of \Cref{prop:main-sEnd-abel} is called the sheaf of
	special endomorphisms on the base change of the semi-global abelian-type
	model $\Sh_{K_{2,p}}(G_2,X_2)$ to $O_{\reflexHodge}$.  Its descent to
	$O_{\reflexAbel}$ is carried out in \Cref{rmk:descent-base}.
	For $\ell\neq p$, the injective map
	\(\realization_\ell:\AbSpecialEnd\longrightarrow\CV_\ell\)
	of \Cref{prop:main-sEnd-abel} is called the $\ell$-adic realization of
	abelian-type special endomorphisms.
	
	The Hecke-invariant quadratic pairing
	$\sEnd\times\sEnd\to
	\underline{\frac1N O_{\FGlobal,(p)}}$ of
	\Cref{defn-pairing-special-SO} descends to a quadratic pairing on
	$\AbSpecialEnd$ with the same target.
\end{defn}

\subsection{Independence of the choice of the Hodge embedding}\label{subsec:indep-hodge}

We assume the data $(G,X,\CG)$, $(G_2,X_2,\CG_2)$ and representations $V, W$  satisfy \Cref{ass:Shimura-data}, \Cref{ass:special-cycle}, \Cref{ass:d=1} and \Cref{ass:J=1}. As in \Cref{prop:main-sEnd-abel} we have constructed an unramified scheme $\AbSpecialEnd$ on the Shimura variety 
$\Sh_{K_{2,p}}(G_2,X_2)_{O_{\reflexHodge}}$ equipped
with a Hecke action by $G_2(\BA_f^p)$. Each of its connected components is finite over the corresponding Shimura model. 

In this subsection we prove that $\AbSpecialEnd$ is independent (after passage to a common
localized reflex ring) of the Hodge-type lift and of
the $\FGlobal$-linear Hodge embedding. 
By the $G_2(\BA_f^p)$-Hecke-action, it suffices to check that after restricting to a connected component, the realization map $\realization_\ell$ has the same image. 

We fix a $G_2$-representation $V_\FGlobal$ with a lattice $\VlattiveGlobal$, equipped with a quadratic pairing. 
Let $(G',X')$ be another Hodge-type lifting equipped with a
very good $\FGlobal$-linear Hodge embedding
\(\iota': (G',X')\longrightarrow (\GSp_\Q(W',\psi'),S^{\pm})\)
satisfying \Cref{ass:Shimura-data}, \Cref{ass:special-cycle}, \Cref{ass:d=1} and \Cref{ass:J=1}. 
Put $\reflexHodge':=E(G',X')$.
Write $\CG'$ for the corresponding Bruhat--Tits stabilizer and
$K_p'=\CG'(\Z_p)$.  Choose a sufficiently small prime-to-$p$ level
$K'^p\subset G'(\A_f^p)$ and put $K'=K_p'K'^p$.  
We denote the universal abelian scheme associated to this Hodge embedding by $\CA'$. 
As in \Cref{subsec:special-endomorphisms-Hodge-type}, define the subsheaf
$\sEnd'\subseteq\Endint(\CA')_{(p)}$ of special endomorphisms.  We compare the two
sheaves produced by \Cref{prop:main-sEnd-abel} after base change to a common
localized reflex ring.
\subsubsection{Shimura datum $(\commonG,\commonX)$}
We form a common Shimura datum $(\commonG,\commonX)$ with maps to $(G,X)$ and $(G',X')$.

By real approximation and \Cref{prop:Kisin2010-2.2.6}, we may choose
$\gamma\in\CG_2(\Z_{(p)})$ such that, after replacing $X'$ by
$\gamma X'\gamma^{-1}$, there exist $h\in X$ and $h'\in X'$ with the same
image $h_2:\BS\to G_2$. The twisting bijection of
\Cref{cor:twisting-special} preserves special endomorphisms under this replacement.

Consider the fiber product of $G$ and $G'$ over $G_2\times \Gm$, where each map to the second component is the composite of the corresponding Hodge embedding with the symplectic similitude character. 
Define $\commonG:= G\times_{G_2\times \Gm}G'$. 
Note that the kernel $\commonKernel$ of $G\times_{G_2\times \Gm}G'\to G_2$ fits into a short exact sequence 
\[1\to \kernelalphaG^c\times \kernelalphaG^{',c}\to \commonKernel\to \Gm\to 1.\]
In particular, $\commonKernel$ is connected, so $G\times_{G_2\times \Gm}G'$ is connected. \Cref{ass:Shimura-data}\eqref{ass:Shimura-data:R-smooth} follows. 

Let $\commonh=(h,h'):\mathbb S\to
(G\times_{G_2\times \Gm}G')_\mathbb R$, and let $\commonX$ be its
$(G\times_{G_2\times \Gm}G')(\mathbb R)$-conjugacy class. 
Then the Shimura datum $(\commonG,\commonX)$ admits a Hodge embedding \[
\commoniota: (\commonG,\commonX)
\hookrightarrow
\bigl(\GSp(W\oplus W'),S^{\pm}\bigr).
\]	

\begin{lem}\label{lem:fundamental-emb}
	The Hodge embedding $\commoniota$ is fundamental in the sense of \Cref{defn:fundamental-emb}.
\end{lem}
\begin{proof}
	After base change to $\overline\Q$, each factor of $G^c$
	indexed by an embedding of $\FGlobal$ acts only on the summand of $W$
	with that embedding. The $\FGlobal$-linear equivariant inclusion
	$V\hookrightarrow\End_\FGlobal(W)$ implies that this factor acts trivially
	on every summand of $V$ indexed by a different embedding.
	Thus the actions of $G^c$ and $(G')^c$ on $V$ are obtained by
	restriction of scalars from $\FGlobal$-representations.
	By \cite[Lemma~7.2.5]{KisinPappasZhou2024ParahoricII},
	\[
	(\commonG)^c=(G^c\times_{G_2}(G')^c)^\circ.
	\]
	Since the $G_2$-action on $V$ has finite kernel by
	\Cref{ass:special-cycle}, this neutral component also equals that of
	\[
	G^c\times_{\Res_{\FGlobal/\Q}\GL_\FGlobal(V_\FGlobal)}(G')^c.
	\]
	Indeed, on this fiber product the quotient of the two images in $G_2$
	belongs to the finite kernel and therefore equals $1$ on the neutral component.
	Restriction of scalars commutes with fiber products and neutral components.
	Hence $(\commonG)^c$ satisfies \eqref{eq:KPZ-dagger}, and its
	block-diagonal action on $W\oplus W'$ comes from an $\FGlobal$-representation.
\end{proof}

Its reflex field is $\commonE=\reflexHodge\reflexHodge'$: the common
cocharacter conjugacy class is determined by its two projections, since
both projections induce isomorphisms on adjoint groups. Thus
$\commonE/\reflexAbel$ is unramified at $p$.
The kernel of the simply connected covering of $(\commonG)^{\der}$ is
the intersection of the covering kernels for $G^{\der}$ and
$G'^{\der}$. Its order therefore divides $|\pi_1(G^{\der})|$, so
$p\nmid|\pi_1((\commonG)^{\der})|$.
Let $\commonv$ be the place of $\commonE$ induced by the fixed
embedding into $\overline{\Q}_p$. All integral Shimura varieties in
the following comparison are base changed to the localized ring
$O_{\commonE,(\commonv)}$. For local models and formal deformation arguments,
we further base change to its completion $\widehat O_{\commonE}$ at $\commonv$.

Let $\commonx$ be the point of the building of
$\commonG_{\Q_p}$ corresponding to the points defining $\CG$ and
$\CG'$ under the two adjoint building isomorphisms.  Let
$\commonCG$ be the connected parahoric model attached to
$\commonx$ and put
\(\commonK_p=\commonCG(\Z_p).\)
The two maps to $G_2\times\Gm$ have $R$-smooth kernels
$\kernelalphaG^c$ and $\kernelalphaG^{',c}$. Applying
\cite[Proposition~2.4.13]{KisinZhou2025} to these maps and to
$\commonG\to G$ gives
\[
\commonCG\simeq\CG\times_{\CG_2\times\Gm}\CG'.
\]
The two given integral Hodge embeddings therefore give a closed
immersion $\commonCG\hookrightarrow\GL(W_{\Z_p}\oplus W'_{\Z_p})$.

By \Cref{eq:local-model-common-adjoint} we have isomorphisms between the local models
\[
\begin{aligned}
\BM^\loc_{\commonCG,\{\commonmu\}}&\cong
\BM^\loc_{\CG,\{\mu\}}\times_{\widehat O_\reflexHodge}\widehat O_{\commonE}\\
&\cong \BM^\loc_{\CG',\{\mu'\}}\times_{\widehat O_{\reflexHodge'}}\widehat O_{\commonE}.
\end{aligned}
\]
Here $\{\commonmu\}$ and $\{\mu'\}$ are the $p$-adic conjugacy classes of Hodge cocharacters attached to $\commonX$ and $X'$, respectively.
Hence $\commonCG\to \GL(W_{\Z_p}\oplus W'_{\Z_p})$ induces a closed immersion
on local models, verifying \Cref{ass:Shimura-data}(\ref{ass:Shimura-data:good})
for $\commoniota$.

After shrinking the prime-to-$p$ levels, choose
$\commonKprimeToP\subset \commonG(\A_f^p)$ mapping into the chosen
prime-to-$p$ levels for $G$ and $G'$, and write
$\commonK=\commonK_p\commonKprimeToP$. 
The normalization in \Cref{defn:Hodge-integral-model} and
\cite[Theorem~7.1.3]{KisinPappasZhou2024ParahoricII} give the model
$\Sh_{\commonK}$ and its local-model diagram.
The projections $\commonG\to G,G'$ have $R$-smooth torus kernels.
By \cite[Proposition~3.22]{DanielsYoucis2025} and the gluing argument
in the proof of \Cref{prop:Kisin-abelian-(3.3.10)}, the natural maps
$p_1:\Sh_{\commonK}\longrightarrow\Sh_K$, $p_2:\Sh_{\commonK}\longrightarrow\Sh_{K'}$ are finite \'etale.
\subsubsection{Special cycles on \texorpdfstring{$\Sh_{\commonK}$}{the common cover}}
Equip $\commonW=W\oplus W'$ with
$\commonPsi=\Psi\oplus\Psi'$.  The direct-sum lattice is self-dual
at $p$ for $\commonPsi$ by the standing hypothesis on both summands.
The $\commonG$-action on each summand is $\FGlobal$-linear, so its
block-diagonal action on $\commonW$ commutes with the diagonal
$\FGlobal$-action. Thus $\commoniota$ is an $\FGlobal$-linear Hodge
embedding. The diagonal $O_\FGlobal$-action commutes with $\commonCG$
and is self-adjoint for $\commonPsi$. On the universal abelian scheme
\[
\commonA\simeq p_1^{-1}\CA\times p_2^{-1}\CA',
\]
it is the product of the two given $O_\FGlobal$-actions, and
$\Tr_{\FGlobal/\Q}\circ\commonPsi=\psi\oplus\psi'$.
The two given inclusions induce a diagonal embedding
\[
\VlattiveGlobal\otimes_{O_\FGlobal}O_{\FGlobal,(p)}
\lhook\joinrel\longrightarrow
\End_{O_{\FGlobal,(p)}}
\bigl(W_{O_{\FGlobal,(p)}}\oplus W'_{O_{\FGlobal,(p)}}\bigr).
\]
The $\commonPsi$-trace pairing restricts to
$(N+N')(\ ,\ )$ on this lattice.

Write $\iota_V$ and $\iota'_V$ for the two inclusions of the common
lattice $V_{\Z_p}$. The diagonal inclusion is
\[
\commoniota_V:V_{\Z_p}\longrightarrow
\End_{\Z_p}(W_{\Z_p}\oplus W'_{\Z_p}),\qquad
v\longmapsto\operatorname{diag}(\iota_V(v),\iota'_V(v)).
\]
It is $\commonCG$-equivariant and saturated. Thus it
satisfies \Cref{ass:special-cycle}, with pairing constant $N+N'$.
Applying \Cref{prop:cris-tensor} to this representation defines
$\commonCV_\cris$; its two diagonal components identify with
$p_1^*\CV_\cris$ and $p_2^*\CV'_\cris$, respectively.
The same projections identify the common Betti, de Rham, and rational
\'etale realizations by functoriality. 

\begin{prop}\label{prop:d=1-two-hodge}
	After pullback to $\BM^\loc_{\commonCG,\{\commonmu\}}$, the filtered subbundles $\Fil^{0}\CV_\loc, {\Fil^{0}}'\CV_\loc$ given by $\CV_\loc\subseteq \Endint(\CW_\loc)$ and $\CV_\loc\subseteq \Endint(\CW'_\loc)$ coincide. In particular, the filtered subbundle given by $\CV_\loc\subseteq \Endint(\CW_\loc\oplus \CW'_\loc)$ also equals $\Fil^{0}\CV_\loc$. The tuple $(\commonG,\commonX,V_{\Z_p})$ satisfies \Cref{ass:d=1}.
\end{prop}

\begin{proof}
	The two subbundles coincide on the generic fiber because both
	are defined by the common Hodge cocharacter. By \Cref{ass:d=1},
	both $\Fil^{0}\CV_\loc$ and ${\Fil^{0}}'\CV_\loc$ are locally direct summands of $\CV_\loc$, so equality on the generic fiber gives
	$\Fil^{0}\CV_\loc={\Fil^{0}}'\CV_\loc$. 
	The common $\gr^{-1}\CV_\loc$ is consequently invertible, and the cokernel of the inclusion
	$\gr^{-1}\CV_\loc\hookrightarrow
	\Endint(\CW_\loc\oplus\CW'_\loc)/\Fil^0\Endint(\CW_\loc\oplus\CW'_\loc)$ is locally free. 
	Since $p_1$ is finite \'etale and $\Sh_K$ is proper over the
	localized reflex ring, $\Sh_{\commonK}$ is proper over
	$O_{\commonE,(\commonv)}$.
	Finally, $(\commonG)^{\der}\to G^{\der}$ is surjective, so
	$V^{(\commonG)^{\der}}=V^{G^{\der}}=0$.
\end{proof}

As in \Cref{subsec:special-endomorphisms-Hodge-type} we can define the subsheaf $\commonSpecialEnd\subseteq\Endint(\commonA)_{(p)}$ of special endomorphisms, and by \Cref{lem:special-Z-inv}(\ref{lem:special-Z-inv3}) for any special endomorphism $f$ the components $f_{12}: p_1^{-1}\CA\to p_2^{-1}\CA'$, $f_{21}:  p_2^{-1}\CA'\to p_1^{-1}\CA$ vanish. 
Fix $\ell\ne p$. The sections of $\commonSpecialEnd$ over $T$ are pairs
\[
(f,f')\in p_1^{-1}\sEnd(T)\times p_2^{-1}\sEnd'(T)
\]
such that $f_\ell=f'_\ell$ in the common $\CV_\ell$ at every
characteristic-$0$ geometric point of $T$, and $f_\cris=f'_\cris$
in the common $\CV_\cris$ at every characteristic-$p$ geometric point of $T$.
We thus have two maps of fppf sheaves 
$$\commonSpecialEnd\longrightarrow p_1^{-1}\sEnd, \qquad \commonSpecialEnd\longrightarrow p_2^{-1}\sEnd'$$
which are injective by \Cref{lem:special-Z-inv}(\ref{lem:special-Z-inv3}).
Both maps preserve all these
realizations and the normalized $\FGlobal$-valued pairings: dividing
the trace pairings by $N$, $N'$, and $N+N'$, respectively, restricts
each to the same quadratic pairing on $V$.

Applying \Cref{cor:special-representable} to the three Hodge-type data shows
that $\sEnd$, $\sEnd'$, and $\commonSpecialEnd$ are represented by
unramified schemes locally of finite presentation over their respective
Shimura models, and that their connected components are finite.

\begin{prop}\label{prop:iso-sEnd}
	The two maps between Hecke-equivariant sheaves 
	\[
	\commonSpecialEnd\longrightarrow p_1^{-1}\sEnd,
	\qquad
	\commonSpecialEnd\longrightarrow p_2^{-1}\sEnd'
	\]
	are isomorphisms.
\end{prop}
\begin{proof}
We prove the assertion for the first projection. 
The maps are monomorphisms between schemes locally of finite
presentation by \Cref{lem:special-Z-inv}. 
By \Cref{prop:special-F_p} it is bijective on $\barFp$-points (note that $\Z_p$ is faithfully flat over $\Z_{(p)}$).
We check formal \'etaleness at $\barFp$-points.

Consider a surjection $R\twoheadrightarrow S$ of Artinian
$O_{\commonE,(\commonv)}$-algebras with residue field $\barFp$ whose kernel
is square zero, equipped with the trivial divided power structure.
Fix a map $\Spec R\to\Sh_{\commonK}$ and let
$(f,f')\in\commonSpecialEnd(S)$. Suppose that $f$ lifts to
$p_1^{-1}\sEnd(R)$. We must lift $f'$ to $p_2^{-1}\sEnd'(R)$.
By \Cref{cor:crystalline-special} applied to the common Hodge datum,
evaluation on the PD-thickening $\Spec S\hookrightarrow\Spec R$ gives
$f_{\cris,R}=f'_{\cris,R}$ in the common $\CV_{\dR,R}$.
The given lift of $f$ makes this section belong to $\Fil^0\CV_{\dR,R}$.
By \Cref{prop:d=1-two-hodge} it preserves the Hodge filtration on
$\CW'_{\dR,R}$, so \Cref{thm:tangent-space} gives a unique lift of $f'$.
Its geometric fibers are special and its crystalline realization is the
common section, so the lifted pair belongs to $\commonSpecialEnd(R)$.
Thus $\commonSpecialEnd\longrightarrow p_1^{-1}\sEnd$ is \'etale at every point of the special fiber.

Each connected component of the source and target is proper over
$O_{\commonE,(\commonv)}$ by
\Cref{prop:d=1-two-hodge,cor:special-representable}.
The non-\'etale locus of the first projection is closed and has empty
special fiber, so its intersection with each connected component of the
source is empty. 
Thus the first projection is an \'etale monomorphism, hence an open immersion.
Its closed complement in each connected component of $p_1^{-1}\sEnd$
is proper over $O_{\commonE,(\commonv)}$ and has empty special fiber.
Hence the complement is empty, and the first projection is an isomorphism.
\end{proof}

\begin{lem}\label{lem:etale-map-bij}
	If $\alpha: X\to Y$ is a universally injective \'etale map of finite type between Noetherian schemes and induces a bijection on closed points, then $\alpha$ is an isomorphism.
\end{lem}

\begin{proof}
	A universally injective \'etale morphism is an open immersion.
	Its image contains every closed point of $Y$, so its closed complement
	is empty. Hence $\alpha$ is an isomorphism.
\end{proof}

\begin{cor}\label{cor:same-special}
	Let $(G,X)$ and $(G',X')$ be two Hodge-type lifts of $(G_2,X_2)$
	to which \Cref{prop:iso-sEnd} applies.
	Let $L$ be any number field containing both reflex fields and let $u$ be
	the place selected by the fixed embedding into $\overline\Q_p$.  After
	base change to $\mathcal O_{L,(u)}$, the sheaves furnished by
	\Cref{prop:main-sEnd-abel} are canonically
	isomorphic, compatibly with their Hecke actions, quadratic pairings, and
	$\ell$-adic realization maps.
\end{cor}

\begin{proof}
	As the $\ell$-adic realization map $\realization_\ell$ is injective, we only need to check that the sheaves of special endomorphisms have the same image in $\CV_\ell$. Since $\realization_\ell$ is compatible with the Hecke action, it is enough to check this on one connected component.  The claim then follows from \Cref{prop:iso-sEnd}; compatibility with the pairing follows from \Cref{defn-pairing-special-SO}.
\end{proof}

\begin{cor}\label{rmk:descent-base}
	The fppf sheaf $\AbSpecialEnd$ and the realization maps $\realization_\ell$ on $\Sh_{K_{2,p}}(G_2,X_2)_{O_{\reflexHodge}}$ constructed in \Cref{prop:main-sEnd-abel} descend uniquely to $\Sh_{K_{2,p}}(G_2,X_2)$. The descended sheaf is represented by a union of finite unramified schemes over $\Sh_{K_{2,p}}(G_2,X_2)$.
\end{cor}

\begin{proof}
Injectivity of $\realization_\ell$ makes such descent unique. We show that for each $T\in \FGlobal$ and $K^p$-invariant compact subset
$\Omega\subset V(\BA_f^p)$, the descent datum for $\CV_\ell$ defines a descent datum for $\CZ_2(T,\Omega)$.
On the generic fiber, the descent datum for $\CV_\ell$ preserves $\CZ_2(T,\Omega)$ by
\Cref{lem:special-char-0-fiberwise} and
\cite[II, \S4 and Theorem~5.5; III, Remark~6.1]{Milne1990CanonicalModels}.
Since $\reflexHodge/\reflexAbel$ is unramified above $p$, it remains to compare the two pullbacks
at points above $p$. These overlaps are represented by automorphisms in
$\Gal(\bar\Q/\reflexAbel)$ preserving the place of $\bar\Q$ induced by the fixed embedding into $\overline\Q_p$.
For such a $\sigma$, apply Milne's conjugation construction to a special point of
$X$ and its image in $X_2$.
We obtain a new Hodge-type lifting $(G',X')\to(G_2,X_2)$ with $G'=G$\footnote{The $G_2$-torsor in this case is trivial, so the inner form $G'$ is isomorphic to $G$.}
and $E(G',X')=\sigma(\reflexHodge)$, whose transported Hodge
embedding and representations satisfy \Cref{ass:Shimura-data},
\Cref{ass:special-cycle}, \Cref{ass:d=1}, and \Cref{ass:J=1}. 
For the level $K'$ obtained from $K$ by the finite adelic comparison,
loc.\ cit. gives an isomorphism
\[
\ShG_K(G,X)\times_{\reflexHodge,\sigma}\sigma(\reflexHodge)
\xrightarrow{\sim}\ShG_{K'}(G',X').
\]
By uniqueness of canonical integral models
\cite[Theorem~4.10 and Corollary~4.12]{DanielsYoucis2025}, this extends to
\[
\Sh_K(G,X)\times_{O_{\reflexHodge},\sigma}O_{\sigma(\reflexHodge)}
\xrightarrow{\sim}\Sh_{K'}(G',X').
\]
Thus \Cref{cor:same-special} identifies the two pullbacks
of $\CZ_2(T,\Omega)$ to $O_{\reflexHodge\sigma(\reflexHodge)}$, compatibly with their maps into
the pullback of $\CV_\ell$. The claim follows. 
\end{proof}

\section{\RZname spaces}\label{sec:RZ}

\subsection{\RZname spaces at parahoric level}\label{sec:defn-general-RZ}
Let $G_{\Q_p}$ be a reductive group over $\Q_p$, let $\CG=\CG^\circ/\Z_p$ be the
connected parahoric Bruhat--Tits model attached to a building point,
and let $\{\mu\}$ be a conjugacy class of
minuscule cocharacters of $G_{\Q_p}$. Let
$\reflexHodgeLocal=E(G_{\mathbb Q_p},\{\mu\})$ be the local reflex
field, and let $\breve \reflexHodgeLocal$ be the completion of its maximal unramified
extension. Fix
\(b\in G_{\mathbb Q_p}(\breve{\mathbb Q}_p)\)
such that
\([b]\in B(G_{\mathbb Q_p},\mu^{-1})\)
 is basic. We assume the corresponding local model $\BM^\loc_{\CG, \{\mu\}}$  is smooth over $O_\reflexHodgeLocal$. 
\begin{defn}\label{prop:GSpin-RZ-defn}
The construction of
\cite[Definition~2.4.1]{PappasRapoport2026IntegralLocalShimura}
associates a $v$-sheaf to the local Shimura datum
$ (G_{\mathbb Q_p},b,\{\mu\})$ and the parahoric group scheme $\CG$:
	\[
	\mathcal M^{\mathrm{int}}_{\CG,b,\mu}
	\longrightarrow \operatorname{Spd}(\mathcal O_{\breve \reflexHodgeLocal}).
	\]
	For an affinoid perfectoid space \(S=\Spa(R,R^+)\) of characteristic
	\(p\), choose a topologically nilpotent unit \(\varpi_S\in R^+\).
	For \(r\in\Q_{>0}\), we use the open subspace
	\[
	\CY_{[r,\infty)}=\{|[\varpi_S]|\le|p|^r\ne0\}
	\subset S\dot\times\Spa\Q_p.
	\]
	Here \([\varpi_S]\in W(R^+)\) is the Teichm\"uller lift, and
	the inequalities refer to the valuations in the construction of
	\(S\dot\times\Spa\Z_p\); see
	\cite[\S\S11.2, 22.1]{ScholzeWeinstein2020Berkeley}.
	For rational \(r=a/c>0\), the inequality means
	\(|[\varpi_S]|^c\le|p|^a\ne0\).
	For a general \(S\), framing conditions on these subspaces are
	imposed on affinoid charts. The existence of a framing for sufficiently large \(r\)
	is independent of the choice of \(\varpi_S\).
	For a perfectoid space $S$ over the residue field of
	$\mathcal O_{\breve \reflexHodgeLocal}$, its $S$-points are the isomorphism classes of
	tuples
	\((S^\sharp,\mathscr P,\varphi_{\mathscr P},\iota),\)
	where:
	\begin{enumerate}
		\item $S^\sharp$ is an untilt of $S$ over
		$\operatorname{Spa}(\mathcal O_{\breve \reflexHodgeLocal})$;
		\item $(\mathscr P,\varphi_{\mathscr P})$ is a $\CG$-shtuka
		over $S$ with one leg along $S^\sharp$, bounded by $\{\mu\}$; more precisely, $\sP$ is a $\CG$-torsor on $S\dot{\times} \Spa \Z_p$, and $\varphi_{\mathscr P}$ is an isomorphism
		\[\varphi_{\mathscr P}:( \Frob_S^* \sP)|_{S\dot{\times} \Spa \Z_p\setminus S^\sharp}\simlr \sP|_{S\dot{\times} \Spa \Z_p\setminus S^\sharp}\]
		that is bounded by $\{\mu\}$ at $S^\sharp$. 
		\item $\iota$ is a germ as $r\to\infty$ of isomorphisms
		\(\iota: \mathscr P|_{\CY_{[r,\infty)}}\simlr G\times \CY_{[r,\infty)}\)
		identifying $\varphi_{\mathscr P}$ with $b \times \Frob_S$.
	\end{enumerate}
	By \cite[Proposition~18.4.1]{ScholzeWeinstein2020Berkeley}, the functor
	$\CM\mapsto\CM^\diamondsuit$ from normal flat formal schemes locally formally of finite type over $O_{\breve \reflexHodgeLocal}$
	to the corresponding $v$-sheaves is fully faithful.
	When this $v$-sheaf is represented by a normal flat formal scheme locally formally of finite type over $O_{\breve \reflexHodgeLocal}$, we call the representing formal scheme the parahoric \RZname space associated with
	$(G_{\mathbb Q_p},b,\{\mu\},\CG)$ and denote it by
	$\mathcal M_{\CG,b,\mu}
	/ \Spf\mathcal O_{\breve \reflexHodgeLocal}$,
	with the identifying isomorphism
	\[
	\mathcal M_{\CG,b,\mu}^{\diamondsuit}
	\simeq
	\mathcal M^{\mathrm{int}}_{\CG,b,\mu}.
	\]
\end{defn}
When the local Shimura datum is of abelian type and $p\neq 2$, by
\cite[Theorem~2.5.4]{PappasRapoport2026IntegralLocalShimura} this formal scheme exists and is uniquely determined.  
It carries a natural action of
\[
J_b(\mathbb Q_p)
:=
\bigl\{
g\in G_{\mathbb Q_p}(\breve{\mathbb Q}_p)
\mid g^{-1}b\sigma(g)=b
\bigr\},
\]
given by changing the framing, since this inner form acts on the corresponding $v$-sheaf. 

\begin{examp}[{\cite[Theorem~25.1.2]{ScholzeWeinstein2020Berkeley}}]\label{examp:GL-RZ}
	Let $G=\GL_n$, $\CG=\GL_{n,\Z_p}$, and
	$\mu=\mu_d$, where
	\[
	\mu_d(z)=\operatorname{diag}
	(\underbrace{z,\ldots,z}_{d},\underbrace{1,\ldots,1}_{n-d}),
	\qquad 0\leq d\leq n.
	\]
	Then $\reflexHodgeLocal=\Q_p$ and $\mathcal O_{\breve \reflexHodgeLocal}=W(\barFp)$.
	Fix a $p$-divisible group $\mathbb X_b$ over $\barFp$ of height
	$n$ and dimension $d$ whose rational Dieudonn\'e module, with the
	homological normalization of $\BD$, is isomorphic to
	$(\breve\Q_p^n,b\sigma)$; here $\sigma$ is the Witt-vector Frobenius.
	For a $W(\barFp)$-scheme $S$ on which $p$ is locally nilpotent, let
	$\bar S=S\times_{\Spec W(\barFp)}\Spec\barFp$.
	The \RZname space $\mathcal M_{\CG,b,\mu_d}$ represents the functor
	\[
	S\longmapsto\{(X,\rho)\}/\simeq,
	\]
	where $X$ is a $p$-divisible group over $S$ of height $n$ and
	dimension $d$, and the framing
	\[
	\rho:\mathbb X_b\times_{\barFp}\bar S
	\dashrightarrow X\times_S\bar S
	\]
	is a quasi-isogeny, namely an invertible morphism after inverting
	$p$ in the groups of homomorphisms of $p$-divisible groups.
	Two pairs $(X,\rho)$ and $(X',\rho')$ are isomorphic if there is
	an isomorphism $f:X\xrightarrow{\sim}X'$ over $S$ satisfying
	$(f\times_S\bar S)\circ\rho=\rho'$.
\end{examp}

Suppose that $(G_{\mathbb Q_p},b,\{\mu\})$ is of local Hodge type, i.e., there is a Hodge
embedding
\(\iota: (G_{\mathbb Q_p},\{\mu\}) \hookrightarrow (\GL(W),\{\mu_d\})\), and it extends to a closed immersion $\CG\to \GL(W_{\Z_p})$. 
Let $\mathscr X$ be the pullback of the universal $p$-divisible group to $\mathcal M_{\mathcal G,b,\mu}$. For any scheme $S/\mathcal M_{\mathcal G,b,\mu}$ on which $p$ is locally nilpotent, its framing is $\rho: \mathbb X_b\times_{\barFp}\bar{S} \dashrightarrow \mathscr X\times_{\mathcal M_{\mathcal G,b,\mu}}\bar{S}$.

Note that \cite[Corollary 2.5.6]{PappasRapoport2026IntegralLocalShimura}
implies that for every $x\in\mathcal M_{\mathcal G,b,\mu}(\barFp)$,
there exists $y\in\mathbb M^{\mathrm{loc}}_{\mathcal G,\mu}(\barFp)$
such that
$
\widehat{\mathcal O}_{\mathcal M_{\mathcal G,b,\mu},x}
\cong
\widehat{\mathcal O}_{\mathbb M^{\mathrm{loc}}_{\mathcal G,\mu}\times_{O_\reflexHodgeLocal}O_{\breve{\reflexHodgeLocal}},y}. 
$ 
This implies $\mathcal M_{\mathcal G,b,\mu}/O_{\breve \reflexHodgeLocal}$ is smooth. 

\begin{defn}\label{defn:Weil-descent}
	Let $\sigma_\reflexHodgeLocal$ be the canonical arithmetic Frobenius automorphism of $O_{\breve\reflexHodgeLocal}$ over $O_{\reflexHodgeLocal}$, lifting
	the $p^f$-power Frobenius of its residue field, where $p^f$ is the residue
	cardinality of $\reflexHodgeLocal$. Its restriction to $O_{\breve\Q_p}$ is $\sigma^f$, where $\sigma$ is the Witt-vector Frobenius.
	For a formal scheme $\mathcal X/\Spf O_{\breve\reflexHodgeLocal}$ write
	\[
	\mathcal X^{(\sigma_\reflexHodgeLocal)}:=
	\mathcal X\times_{\Spf O_{\breve\reflexHodgeLocal},\,\sigma_\reflexHodgeLocal}
	\Spf O_{\breve\reflexHodgeLocal}.
	\]
	A \emph{Weil descent datum} on $\mathcal X$ relative to $\reflexHodgeLocal$ is an
	isomorphism
	$\Frob:\mathcal X\xrightarrow{\sim}\mathcal X^{(\sigma_\reflexHodgeLocal)}$.
	The canonical Weil descent datum of
	\cite[(3.1.6)]{PappasRapoport2024PadicShtukas}
	\[
	\Frob^{\mathrm{int}}:
	\mathcal M^{\mathrm{int}}_{\CG,b,\mu}
	\xrightarrow{\sim}
	\left(\mathcal M^{\mathrm{int}}_{\CG,b,\mu}\right)^{(\sigma_\reflexHodgeLocal)}
	\]
	is given on the moduli problem by
	\[
	(S^\sharp,\mathscr P,\varphi_{\mathscr P},\iota)
	\longmapsto
	(S^\sharp,\mathscr P,\varphi_{\mathscr P},
	(b\,\sigma(b)\cdots\sigma^{f-1}(b))^{-1}\circ\iota),
	\]
	where the tuple on the right is regarded as a point of
	$\mathcal M^{\mathrm{int}}_{\CG,\sigma^f(b),\mu}\cong \left(\mathcal M^{\mathrm{int}}_{\CG,b,\mu}\right)^{(\sigma_\reflexHodgeLocal)}$. By the full faithfulness of the functor
	$\mathcal X\mapsto\mathcal X^\diamondsuit$
	in \cite[Proposition~18.4.1]{ScholzeWeinstein2020Berkeley},
	$\Frob^{\mathrm{int}}$ is induced by a unique isomorphism
	of formal schemes
	\[
	\Frob:
	\mathcal M_{\CG,b,\mu}
	\xrightarrow{\sim}
	\mathcal M_{\CG,b,\mu}^{(\sigma_\reflexHodgeLocal)}.
	\]
	We call $\Frob$ the \emph{canonical Weil descent datum}
	on $\mathcal M_{\CG,b,\mu}$ relative to $\reflexHodgeLocal$.
\end{defn}
\subsection{\RZname spaces of abelian type}
\label{subsec:abelian-type-RZ-special-cycles}
Let $\alpha_p:G_{\Q_p}\to G_{2,\Q_p}$ be a central epimorphism inducing an
isomorphism on adjoint groups and carrying $\{\mu\}$ to $\{\mu_2\}$.
Assume that it maps the building point $\mathbf x$ defining $\CG$ to the
point $\mathbf x_2$ defining the connected parahoric model $\CG_2$, and that its kernel $\kernelalphaGp$ 
is $R$-smooth.
Consider the basic element
\(b_2:=\alpha_p(b)\). 
Let
\(\reflexHodgeLocal=E(G_{\mathbb Q_p},\{\mu\})\), \( \reflexAbelLocal=E(G_{2,\Q_p},\{\mu_2\}).\)
Since $\alpha_p$ induces an isomorphism on adjoint groups and
$\{\mu_2\}$ is the image of $\{\mu\}$, one has $\reflexAbelLocal\subset \reflexHodgeLocal$.  
We assume $\reflexHodgeLocal/\reflexAbelLocal$ is unramified, so $O_{\breve \reflexAbelLocal}=O_{\breve \reflexHodgeLocal}$. 
By \Cref{rmk:R-smooth-seq}, $\alpha_p$ extends to a morphism
$\CG\to\CG_2$ with smooth kernel $\kernelalpha$, and hence defines an ad-isomorphism of integral local
shtuka data in the sense of
\cite[\S5.1]{PappasRapoport2026IntegralLocalShimura}:
\begin{equation*}
	\alpha_p:
	(G_{\mathbb Q_p},b,\{\mu\},\CG)
	\longrightarrow
	(G_{2,\Q_p},b_2,\{\mu_2\},\CG_2).
\end{equation*}
The induced morphism of integral local Shimura varieties will be denoted
\(\alpha_{\mathrm{RZ}}:\RZG\longrightarrow\mathcal M_2,\)
where $\mathcal M_2:=\mathcal M_{\CG_2,b_2,\mu_2}$.

\subsubsection{Kottwitz components.}
For a (connected) parahoric group scheme $\mathcal H$ with generic fiber $H$,
Pappas--Rapoport write
\(C_{\mathcal H} = \pi_1(H)_{I}\)
for the target of the Kottwitz map.  
Thus $C_G=\pi_1(G_{\Q_p})_I$ and $C_2=\pi_1(G_{2,\Q_p})_I$ below;
$\varphi$ on these groups denotes the action induced by arithmetic
Frobenius $\sigma$.
The Kottwitz maps are
\begin{equation*}
	\kappa_G:|\RZG|\longrightarrow C_G,
	\qquad
	\kappa_2:|\mathcal M_2|\longrightarrow C_2.
\end{equation*}
By \cite[Lemma 6.1]{HeZhou2020}, the image of $\kappa_G$ is
$c_{b,\mu}+\pi_1(G_{\mathbb Q_p})_{I}^\varphi$, and the image of $\kappa_2$ is
$c_{b_2,\mu_2}+\pi_1(G_{2,\Q_p})_{I}^\varphi$, where $c_{b,\mu}$ is a chosen admissible label and we take
$c_{b_2,\mu_2}=\alpha_*(c_{b,\mu})$.  
For $\tau\in C_G$ and $\tau_2\in C_2$, put
\begin{equation}\label{eq:Kottwitz-component-RZ}
	\RZG^\tau:=\kappa_G^{-1}(\tau),
	\qquad
	\mathcal M_2^{\tau_2}:=\kappa_2^{-1}(\tau_2).
\end{equation}
These are open-and-closed formal subschemes; on the level of integral
local Shimura $v$-sheaves this is exactly
\cite[(4.4.5)]{PappasRapoport2026IntegralLocalShimura}.  We use the term
\emph{Kottwitz component} only for these fibers; no assertion about
topological connectedness is built into the terminology. 

The map $\alpha_p$ induces
\(\alpha_*:c_{b,\mu}+\pi_1(G_{\mathbb Q_p})_{I}^\varphi \longrightarrow c_{b_2,\mu_2}+\pi_1(G_{2,\Q_p})_{I}^\varphi.\)
Its image is the coset
$c_{b_2,\mu_2}+\alpha_*(\pi_1(G_{\mathbb Q_p})_{I}^\varphi)$ of target labels.
For every $\tau\in c_{b,\mu}+\pi_1(G_{\mathbb Q_p})_{I}^\varphi$, with $\tau_2=\alpha_*(\tau)$,
\cite[Theorem~5.1.2]{PappasRapoport2026IntegralLocalShimura}\footnote{There is a typographical error in the statement: $\tau\in C_G$
	should read $\tau\in c_{b,\mu}+\pi_1(G_{\mathbb Q_p})_{I}^\varphi$.} gives an
isomorphism on the corresponding integral local Shimura $v$-sheaves.
By the full faithfulness of
\cite[Proposition~18.4.1]{ScholzeWeinstein2020Berkeley}, the induced morphism
of the representing formal schemes is an isomorphism. We therefore
obtain
\begin{equation}\label{eq:PR-componentwise-ad-isomorphism}
	\alpha_\tau:
	\RZG^\tau
	\xrightarrow{\ \sim\ }
	\mathcal M_2^{\alpha_*(\tau)}.
\end{equation}

\begin{prop}\label{prop:RZ-surj}
	Assume $H^1(\Q_p, \kernelalphaGp)=0$; the map $\alpha_*:c_{b,\mu}+\pi_1(G_{\mathbb Q_p})_{I}^\varphi \longrightarrow c_{b_2,\mu_2}+\pi_1(G_{2,\Q_p})_{I}^\varphi$ is surjective, so the map
	$\alpha_{\mathrm{RZ}}:\RZG\longrightarrow\mathcal M_2$
	is surjective.
\end{prop}

\begin{proof}
	Naturality of the Kottwitz homomorphism gives the commutative diagram
	\[
	\begin{tikzcd}[column sep=large]
		G(\Q_p)\arrow[r,"\alpha_p"]\arrow[d,"\kappa_G"']
		& G_2(\Q_p)\arrow[d,"\kappa_2"]\\
		\pi_1(G_{\mathbb Q_p})_{I}^\varphi\arrow[r,"\alpha_*"']
		&\pi_1(G_{2,\Q_p})_{I}^\varphi.
	\end{tikzcd}
	\]
	Both vertical homomorphisms are surjective. The top arrow is
	surjective by $H^1(\Q_p,\kernelalphaGp)=0$, so the bottom arrow is
	surjective as well. Translating by a compatible pair of admissible
	labels proves the asserted surjectivity on their cosets. Finally,
	\eqref{eq:PR-componentwise-ad-isomorphism} is an isomorphism over each
	label in the image, and hence $\alpha_{\mathrm{RZ}}$ is surjective.
\end{proof}
\subsection{Special cycles on Hodge-type \RZname spaces}\label{subsec:V_b}

Let $V_{\Z_p}\hookrightarrow\End_{\Z_p}(W_{\Z_p})$ be an
inclusion satisfying \Cref{ass:special-cycle}. Let $\sigma$ denote the
arithmetic Frobenius of $\breve\Q_p$. The Frobenius on $W_{\Z_p}\otimes_{\Z_p}\breve\Q_p$ is
$b\sigma$, and the induced Frobenius on
$V_{\breve\Q_p}:=V_{\Z_p}\otimes_{\Z_p}\breve\Q_p\subset\End(W_{\Z_p}\otimes_{\Z_p}\breve\Q_p)$ is
\[
\varphi_{V_{\breve\Q_p}}(x)=b\sigma(x)b^{-1}.
\]
Define the space of special quasi-endomorphisms of the framing
isocrystal by
\[
\mathbb V_b:=V_{\breve\Q_p}^{\varphi_{V_{\breve\Q_p}}=1}
=\{x\in\End_{\breve\Q_p}(W_{\Z_p}\otimes_{\Z_p}\breve\Q_p):
x\circ b\sigma=b\sigma\circ x,\ x\in V_{\breve\Q_p}\}.
\]
It is naturally a $\Q_p$-vector space, stable under $J_b(\Q_p)$.

The Dieudonn\'e correspondence identifies $\mathbb V_b$ with the subspace of
$\End(\mathbb X_b)_\Q$ consisting of special quasi-endomorphisms. The
$\FGlobal$-valued $G$-invariant quadratic pairing on $V$ induces
the pairing
\[
h:\mathbb V_b\times\mathbb V_b\longrightarrow
\FGlobal\otimes_\Q\Q_p.
\]

We have $\dim_{\Q_p}\mathbb V_b\leq \dim_{\breve\Q_p} V_{\breve\Q_p}$, and equality holds if and only if all the Newton slopes of $(V_{\breve\Q_p},\varphi_{V_{\breve\Q_p}})$ equal $0$.

\begin{defn}\label{defn:local-special-cycle}
	Let $x\in\mathbb V_b$. The \emph{local special cycle}
	$\CZ(x)$ is the subfunctor of $\RZG$ consisting of the $S$-points for which
	$\rho\circ x\circ\rho^{-1}$ extends to an element of $\End_S(\mathscr X_S)$.
	It is a closed formal subspace of $\RZG$.
	For a tuple $\mathbf x=(x_1,\ldots,x_r)\in\mathbb V_b^r$, set
	\[
	\CZ(\mathbf x)=
	\CZ(x_1)\times_{\RZG}\cdots\times_{\RZG}\CZ(x_r),
	\qquad
	T(\mathbf x)=\bigl(h(x_i,x_j)\bigr)_{i,j}.
	\]
\end{defn}

\begin{prop}\label{prop:local-special-cycle-Jb}
	For $x\in\mathbb V_b$ and $g\in J_b(\mathbb Q_p)$, $g$ maps $\CZ(x)$ to $\CZ(g\cdot x)$.
\end{prop}

\begin{proof}
	We use the conventions
	\[
	g\cdot(\mathscr X,\rho)=(\mathscr X,\rho\circ g^{-1}),
	\qquad g\cdot x=gxg^{-1}.
	\]
	Let $z=(\mathscr X,\rho)$ be an $S$-point of $\RZG$.  Then
	\[
	(\rho\circ g^{-1})(g x g^{-1})(\rho\circ g^{-1})^{-1}
	=\rho x\rho^{-1}.
	\]
	Thus the special quasi-endomorphism defining $\CZ(x)$ extends over $z$
	if and only if the one defining $\CZ(g\cdot x)$ extends over $g\cdot z$.
	This proves $g\CZ(x)=\CZ(g\cdot x)$ functorially in $S$.
\end{proof}

\begin{prop}\label{prop:special-perfectoid}
	Let $x\in\mathbb V_b$, let $S$ be a characteristic-$p$ perfectoid space and let
	$z=(S^\sharp,\mathscr P,\varphi_{\mathscr P},\iota)\in
	\mathcal M_{\CG,b,\mu}^{\diamondsuit}(S)
	\simeq\mathcal M^{\mathrm{int}}_{\CG,b,\mu}(S)$.
	The representation $V_{\Z_p}$ associates to $\mathscr P$ a vector bundle
	$\mathscr V$ on $Y_S:=S\dot{\times}\Spa\Z_p$, with Frobenius and framing.
	Then $z$ factors through $\CZ(x)^\diamondsuit$ if and only if there is a
	Frobenius-compatible section $x':\CO_{Y_S}\to\mathscr V$ whose restriction
	to $\CY_{[r,\infty)}$ corresponds to $x$ for large enough $r$. This section is unique.
\end{prop}

\begin{proof}
	Work first with $S=\Spa(R,R^+)$ and
	$S^\sharp=\Spa(R^\sharp,R^{\sharp+})$ perfectoid affinoid. 
	Let $\mathscr X_z$ be the pulled-back universal $p$-divisible group and
	$\mathscr W$ the bundle associated with $W_{\Z_p}$. By construction we have an isomorphism
	\[
	 \operatorname{Sht}(\operatorname{BKF}(\mathscr X_z))\simeq\mathscr W,
	\]
	compatible with base change, framing and quasi-endomorphisms.
	We use \cite[Theorem~17.5.2]{ScholzeWeinstein2020Berkeley},
	which identifies the category of $p$-divisible groups over an integral
	perfectoid ring with the category of minuscule Breuil--Kisin--Fargues modules
	in the sense of \cite[Definition~2.2.4]{PappasRapoport2024PadicShtukas};
	see \cite[Remark~2.2.5]{PappasRapoport2024PadicShtukas} for a comparison
	of the two definitions.
	
	If $z\in\CZ(x)^\diamondsuit(S)$, its integral endomorphism gives a
	Frobenius-compatible section $u$ of $\underline{\End}(\mathscr W)$.
	The quotient representation $Q_{\Z_p}$ of \Cref{ass:special-cycle}
	is finite free as a $\Z_p$-module. Its associated bundle
	$\mathscr Q:=\underline{\End}(\mathscr W)/\mathscr V$ is therefore a
	vector shtuka, and the image $\bar u:\CO_{Y_S}\to\mathscr Q$ is a
	morphism of shtukas. It vanishes on $\CY_{[r,\infty)}$. 
	The functor from one-leg shtukas to vector bundles on the relative
	Fargues--Fontaine curve, obtained by restriction to this region and
	Frobenius descent, is faithful
	\cite[\S2.2.2]{PappasRapoport2024PadicShtukas}. Hence $\bar u=0$ on
	$Y_S$. This yields $x'$, and the same faithfulness applied to a difference
	proves uniqueness.
	
	Conversely, restrict a section $x'$ as in the statement to
	$S^\circ:=\Spa(R,R^\circ)$, where $R^\circ$ is the power-bounded subring.
	By \cite[Proposition~2.2.7]{PappasRapoport2024PadicShtukas}, the
	BKF-to-shtuka functor is fully faithful for this power-bounded plus-ring.
	Consequently $x'$ gives an integral endomorphism of $\mathscr X_z$ after
	base change to $R^{\sharp\circ}$, and $z|_{S^\circ}$ factors through
	$\CZ(x)^\diamondsuit$. To recover the original point, work locally in a
	formal chart $\Spf A$ of $\RZG$ with
	$\CZ(x)\cap\Spf A=\Spf(A/I)$.
	The point $z$ gives a continuous map $A\to R^{\sharp+}$; restriction to
	$S^\circ$ composes this same map with the injection
	$R^{\sharp+}\hookrightarrow R^{\sharp\circ}$.
	Since the composite kills $I$, the original map kills $I$ and factors
	continuously through $A/I$. This proves the required factorization.
	The construction and uniqueness commute with base change, so the
	affinoid statements glue for general $S$. 
\end{proof}
\subsection{Special cycles on abelian-type \RZname spaces}
We consider the
representation $V$ and the integral lattice $V_{\Z_p}$ fixed in
\Cref{ass:special-cycle}.  By assumption the $G_{\mathbb Q_p}$-action on $V$ factors
through $G_{2,\Q_p}$.  Moreover, the exact sequence of smooth group schemes
\[
1\longrightarrow\kernelalpha\longrightarrow\CG
\longrightarrow\CG_2\longrightarrow1
\]
of \Cref{rmk:R-smooth-seq} shows that the integral
$\CG$-representation $V_{\Z_p}$ descends uniquely to a
$\CG_2$-representation: indeed, $\kernelalpha$ acts trivially on the generic fiber
$V_{\Q_p}$, hence trivially on $V_{\Z_p}$ because $\kernelalpha$ is flat over
$\Z_p$.  Let
\(\rho_V:G_{2,\Q_p}\longrightarrow\GL(V_{\Q_p})\)
denote the resulting representation.

Using the embedding
$V_{\Z_p}\hookrightarrow\End(W_{\Z_p})$, the space
$\mathbb V_b$ of the preceding subsection can equivalently be written
\begin{equation*}
	\mathbb V_{b_2}
	:=
	\mathbb V_b
	=
	\bigl(V_{\Q_p}\otimes_{\Q_p}\breve{\Q}_p\bigr)^{
		\rho_V(b_2)\sigma=1}.
\end{equation*}
The quadratic pairing on $\mathbb V_{b_2}$ is the restriction
of the same pairing $h$, with values in $\FGlobal\otimes_\Q\Q_p$,
used in \Cref{defn:local-special-cycle}; thus $T(\mathbf x)$ below
has the same meaning as before.

\begin{prop}\label{prop:abelian-local-cycles-independent-Hodge-cover}
	Let $x\in\mathbb V_b=\mathbb V_{b_2}$. Assume $H^1(\Q_p, \kernelalphaGp)=0$ and that
	$X_*(\kernelalphaGp)_I$ is torsion-free. The closed formal subscheme
	$\CZ(x)\subseteq \RZG$ is the pullback of a closed formal subscheme
	$\CZ_2(x)\subseteq \mathcal M_2$ via the surjection  $\alpha_{\mathrm{RZ}}:\RZG\longrightarrow\mathcal M_2$. 
\end{prop}

\begin{proof}
	The torsion-freeness of $X_*(\kernelalphaGp)_I$ makes
	$\kernelalpha$ the connected parahoric model of $\kernelalphaGp$.
	We therefore have an exact sequence
	\[1\to \kernelalpha(\Z_p)=\kernelalphaGp(\Q_p)_1\to  \kernelalphaGp(\Q_p)\to (X_*(\kernelalphaGp)_{I})^\varphi\to 1.\]
	Here $\kernelalphaGp(\Q_p)_1$ is the kernel of the torus Kottwitz map,
	equal to the integral points of its connected parahoric model.
	The central action is compatible with the label maps:
	\[
	\begin{tikzcd}[column sep=large]
		\kernelalphaGp(\Q_p)\times\RZG\arrow[r,"\mathrm{act}"]
		\arrow[d,"\kappa_{\kernelalphaGp}\times\kappa_G"']
		&\RZG\arrow[d,"\kappa_G"]\\
		(X_*(\kernelalphaGp)_{I})^\varphi\times C_G
		\arrow[r,"{(\lambda,\tau)\mapsto\iota_*(\lambda)+\tau}"']&C_G.
	\end{tikzcd}
	\]
	Here $\iota_*$ is induced by $\kernelalphaGp\to G_{\mathbb Q_p}$.
	As in \Cref{lem:Shimura-data-strongly-admissible}, torsion-freeness gives
	$0\to X_*(\kernelalphaGp)_I\to C_G\to C_2\to0$.
	Taking Frobenius invariants yields
	\[
	\ker(\alpha_*:C_G^\varphi\to C_2^\varphi)
	=\iota_*\bigl((X_*(\kernelalphaGp)_I)^\varphi\bigr).
	\]
	The admissible label sets are compatible translates of these invariant
	groups. The displayed torus Kottwitz sequence therefore makes
	$\kernelalphaGp(\Q_p)$ transitive on the labels above a fixed target label.
	On each such label the map is the isomorphism
	\eqref{eq:PR-componentwise-ad-isomorphism}; central elements fixing
	that label consequently act trivially there. This identifies
	$\mathcal M_2$ with the central quotient of $\RZG$. Since $\kernelalphaGp(\Q_p)$ acts trivially on
	$\mathbb V_b$, \Cref{prop:local-special-cycle-Jb} makes $\CZ(x)$
	invariant, and its restrictions descend through these isomorphisms
	to the same closed formal subscheme of $\mathcal M_2$.
\end{proof}

\subsection{Relation with affine Deligne--Lusztig varieties}\label{subsec:BT-RZ-perfection}
We compare the perfection of the \RZname spaces with affine Deligne--Lusztig varieties. These results will be used to show the finiteness of local intersection numbers.

Recall that, if $X$ is a scheme over a perfect field of characteristic $p$,
its perfection is
\[
X^\perf:=\varprojlim\bigl(\cdots\xrightarrow{F_X}X
\xrightarrow{F_X}X\bigr).
\]
For $X=\Spec A$ this is $\Spec A^\perf$, where
\[
A^\perf=\varinjlim\bigl(A\xrightarrow{a\mapsto a^p}A
\xrightarrow{a\mapsto a^p}\cdots\bigr).
\] 
For a formal scheme $X$ we define $X^\perf$ to be the perfection of the underlying reduced closed subscheme $X^\red$. 
\begin{lem}\label{lem:perfection-topology}
	Let $X$ be a scheme over a perfect field of characteristic $p$.  The
	canonical morphism
	$\epsilon_X:X^\perf\longrightarrow X^\red$
	is a universal homeomorphism.  
\end{lem}

\begin{proof}
	Affine locally, $A_\red\to A^\perf$ is injective and every $x\in A^\perf$
	satisfies $x^{p^r}\in A_\red$ for some $r\geq0$. Thus the map is integral
	and radicial, and lying-over gives surjectivity on spectra. Hence
	$\Spec A^\perf\to\Spec A_\red$ is a universal homeomorphism.
\end{proof}
We recall the definition of a generalized affine Deligne--Lusztig variety
over a possibly ramified finite extension of $\Q_p$.
\begin{defn}\label{defn:ADLV}
	Let $K/\Q_p$ be a finite extension with residue field $k$, let
	$\breve K$ be the completion of its maximal unramified extension, and
	let $\CH/O_K$ be a parahoric group scheme with connected reductive
	generic fiber $H/K$. Fix a conjugacy class $\{\mu\}$ of geometric
	cocharacters of $H$ and an element $b\in H(\breve K)$ with
	$[b]\in B(H,\{\mu^{-1}\})$.
	For a perfect $\bar k$-algebra $R$, we denote 
	$W_{O_K}(R):=O_K\otimes_{W(k)}W(R)$, using the canonical inclusion
	$W(k)\hookrightarrow O_K$; see
	\cite[\S\S0.5 and~1.1]{Zhu2017AffineGrassmannians}. 
	The Frobenius $\sigma$ on this ring acts trivially on $O_K$ and as
	$W(\Frob_q)$ on $W(R)$, where $q=\#k$ and $\Frob_q(r)=r^q$.
	Note that 
	$
	O_{\breve K}=O_K\otimes_{W(k)}W(\bar k),
	$ 
	and $\sigma$ induces the arithmetic Frobenius of $\breve K/K$.
	
	Choose an alcove whose closure contains the facet defining $\CH$.
	Let $\operatorname{Adm}(\mu^{-1})$ be the corresponding admissible
	subset of the Iwahori--Weyl group of $H$ over $\breve K$, as in
	\cite[\S2.2]{HeZhou2020}, and choose normalizer representatives
	$\dot w\in H(\breve K)$.
	For each perfect $\bar k$-algebra $R$, define
	\[
	X_{\CH}(b,\mu^{-1})^{\mathrm{pre}}(R)
	\subseteq
	H\bigl((W_{O_K}(R))[1/p]\bigr)
	/\CH\bigl(W_{O_K}(R)\bigr)
	\]
	to consist of the classes represented by $g$ such that, for every
	geometric point $\Spec\Omega\to\Spec R$, its image $g_\Omega$ satisfies
	\[
	g_\Omega^{-1}b\sigma(g_\Omega)\in
	\bigcup_{w\in\operatorname{Adm}(\mu^{-1})}
	\CH\bigl(O_K\otimes_{W(k)}W(\Omega)\bigr)\dot w
	\CH\bigl(O_K\otimes_{W(k)}W(\Omega)\bigr).
	\]
	This is the geometric-point formulation of the admissible relative
	position condition; see
	\cite[\S1 and Remark~1.1]{HamacherViehmann2020Finiteness}.
	It is independent of the representatives of $g$ and $w$.
	The \'etale sheafification of this presheaf is denoted by
	$X_{\CH}(b,\mu^{-1})$. It is represented by a perfect $\bar k$-scheme,
	locally perfectly of finite type, by
	\cite[Theorem~1.2]{HamacherViehmann2020Finiteness}.
	Its $\bar k$-points are
	\[
	X_{\CH}(b,\mu^{-1})(\bar k)
	=\left\{
	\begin{array}{@{}l@{}}
		g\CH(O_{\breve K})\in H(\breve K)/\CH(O_{\breve K}):\\[2pt]
		g^{-1}b\sigma(g)\in
		\displaystyle\bigcup_{w\in\operatorname{Adm}(\mu^{-1})}
		\CH(O_{\breve K})\dot w\CH(O_{\breve K})
	\end{array}
	\right\},
	\]
	where the inverse $\mu^{-1}$ agrees with the convention of
	\Cref{sec:defn-general-RZ}; compare
	\cite[Remark~3.3.2]{PappasRapoport2026IntegralLocalShimura}.
\end{defn}
For the local Shimura datum in \Cref{sec:defn-general-RZ}, we apply this
definition with $K=\Q_p$, $H=G_{\mathbb Q_p}$, and $\CH=\CG$.

\begin{prop}\label{prop:RZ-special-fiber-perfection}
	There is a canonical $J_b(\mathbb Q_p)$-equivariant isomorphism of perfect
	$\bar k$-schemes
	\[
	(\RZG)^\perf\xrightarrow{\ \sim\ }
	X_{\CG}(b,\mu^{-1}).
	\]
\end{prop}

\begin{proof}
	By \cite[Proposition~3.3.1]{PappasRapoport2026IntegralLocalShimura}, the
	reduced locus of the integral local-Shimura $v$-sheaf
	$\mathcal M^{\mathrm{int}}_{\CG,b,\mu}$ is represented by the perfect
	scheme $X_{\CG}(b,\mu^{-1})$.  On the other hand,
	$\RZG^\diamondsuit\simeq\mathcal M^{\mathrm{int}}_{\CG,b,\mu}$ by the
	definition of $\RZG$.  For a formal scheme locally formally of finite type,
	the reduced locus of its diamond is represented by the perfection of its
	reduced underlying closed subscheme: on a perfect $\bar k$-scheme $S$, a morphism
	$S\to\RZG$ factors uniquely through
	$(\RZG)^\perf$ by the universal property of perfection.  Thus both
	perfect schemes represent the same reduced $v$-sheaf.  The construction is
	compatible with changing the framing, hence with the action of
	$J_b(\mathbb Q_p)$.
\end{proof}

\begin{prop}\label{prop:special-ADLV-old}
	For any $x\in\mathbb V_b$, perfect $\bar k$-algebra $R$ and element
	\[g\in X_{\CG}(b,\mu^{-1})^\mathrm{pre}(R),\]
	the corresponding element $a_g\in \RZG(R)$ lies in $\CZ(x)$ if and only if
	\[
	\rho_V(\alpha_p(g)^{-1})(x)=g^{-1}xg
	\in V_{\Z_p}\otimes_{\Z_p}W(R).
	\]
	Here $\rho_V$ is the $G_{2,\Q_p}$-representation defined above,
	and its pullback along $\alpha_p$ acts on
	$V_{\Q_p}\subseteq\End_{\Q_p}(W_{\Q_p})$ by conjugation.
\end{prop}
\begin{proof}
	The lifting condition defines the closed perfect subscheme
	$\CZ(x)^\perf$. The lattice condition also defines a closed perfect
	subscheme: locally on a bounded affine-Grassmannian chart, clear a
	power of $p$ from $g^{-1}xg$. Integrality then says that finitely many
	Witt coordinates of each coefficient vanish.
	On an affinoid perfectoid test object $S=\Spa(R,R^+)$ of
	characteristic $p$, this lattice condition is evaluated in $W(R^+)$.
	Under \Cref{prop:RZ-special-fiber-perfection}, it is the extension
	condition for the section in \Cref{prop:special-perfectoid}.
	Thus the two closed perfect subschemes have the same associated
	$v$-sheaf. The full faithfulness of
	\cite[Proposition~18.3.1]{ScholzeWeinstein2020Berkeley} identifies
	them, proving the assertion for every perfect $\bar k$-algebra $R$.
\end{proof}

We record the restriction-of-scalars form that will be used in the
orthogonal case.  Let $\FLocal/\mathbb Q_p$ have residue degree $f$, let $H/\FLocal$ be
an unramified reductive group with reductive model $\mathcal H/O_\FLocal$, and set
\[
G=\Res_{\FLocal/\mathbb Q_p}H,
\qquad
\mathcal G=\Res_{O_\FLocal/\mathbb Z_p}\mathcal H.
\]
Assume that, under the decomposition
\[
G_{\ov{\Q}_p}\cong
\prod_{\tau:\FLocal\hookrightarrow\ov{\Q}_p}
H\otimes_{\FLocal,\tau}\ov{\Q}_p,
\]
the cocharacter $\mu$ has a single nonzero geometric component, and let
$\mu_\FLocal$ denote its conjugacy class in the distinguished factor.
The restriction of the distinguished embedding to
$W(k_\FLocal)[1/p]$ selects one of the $f$ factors in
\begin{equation}\label{eq:dec-G-breve}
	G_{\breve{\mathbb Q}_p}\cong
	\prod_{\iota:k_\FLocal\hookrightarrow\bar k}
	\bigl(\Res_{\FLocal/W(k_\FLocal)[1/p]}H\bigr)
	\otimes_{W(k_\FLocal)[1/p],\iota}\breve{\mathbb Q}_p.
\end{equation}
Here each residue-field embedding also denotes its Witt-vector lift.
Number these factors by $0,\ldots,f-1$, with the distinguished factor
numbered $0$ and with Witt-vector Frobenius $\sigma$ carrying factor
$i-1$ to factor $i$, cyclically. Write
$\breve\FLocal=\FLocal\breve{\mathbb Q}_p$ using the distinguished
embedding. For $b\in G(\breve{\mathbb Q}_p)$, let $b_i$ denote its
components. The identity
\[
(b\sigma)^f=b_\FLocal\sigma^f
\]
on factor $0$ defines $b_\FLocal\in H(\breve\FLocal)$. The operator
$\sigma^f$ on this factor fixes $\FLocal$ and induces
$r\mapsto r^{p^f}$ on the residue field; it is the relative Frobenius
$\sigma_\FLocal$.

\begin{lem}
	\label{lem:ADLV-Shapiro}
	Projection to the factor containing $\mu_\FLocal$ induces an isomorphism of
	perfect $\bar k$-schemes
	\[
	X_{\mathcal G}(b,\mu^{-1})
	\xrightarrow{\ \sim\ }
	X_{\mathcal H}(b_\FLocal,\mu_\FLocal^{-1}),
	\]
	where the affine Deligne--Lusztig variety on the right is formed over $\FLocal$
	and its Frobenius is $\sigma^f$.
\end{lem}

\begin{proof}
	By \cite[Lemma~6.8(2),(3)]{hove2026cyclessplittingmodelsshimura},
	factoring $\FLocal/\mathbb Q_p$ into its totally ramified and
	unramified parts identifies the bounds with the corresponding
	convolution bounds. The hypothesis on the geometric components of
	$\mu$ leaves exactly one nonzero bound, namely $\mu_\FLocal^{-1}$.
	For comparison, see \cite[Lemma~3.5]{Zhu2017AffineGrassmannians}
	in the unramified case.
	
	We describe the inverse on perfect $\bar k$-algebras $R$.
	Choose a representative $h$ on the distinguished factor and set
	\[
	g_0=h,\qquad g_i=b_i\sigma(g_{i-1})
	\quad(1\leq i\leq f-1).
	\]
	Then $g_i^{-1}b_i\sigma(g_{i-1})=1$ on every factor $i\neq0$,
	so all zero bounds are satisfied. On factor $0$ the relative
	position is determined by
	\[
	g_0^{-1}b_0\sigma(g_{f-1})
	=h^{-1}b_\FLocal\sigma^f(h),
	\]
	which has the required bound precisely when $h$ belongs to the
	affine Deligne--Lusztig variety on the right.
	Changing the representative of $h$ changes each $g_i$ by right
	multiplication by an integral element of the corresponding factor,
	because $\sigma$ carries the integral coefficient ring of factor
	$i-1$ isomorphically onto that of factor $i$.
	Conversely, the zero relative-position bounds determine the coset
	of each $g_i$ recursively from that of $g_0$.
	These constructions are functorial in perfect $R$ and commute
	with \'etale descent. They therefore give mutually inverse maps
	of the \'etale quotient sheaves, and hence of the representing
	perfect schemes.
\end{proof}
\begin{cor}\label{prop:special-ADLV}
	Assume that $V_{\Q_p}$ is obtained by restriction of scalars from an
	$H$-representation over $\FLocal$, and that $V_{\Z_p}$ is an
	$\mathcal H$-stable $O_\FLocal$-lattice in that representation.
	Let $R$ be a perfect $\bar k$-algebra, let
	\[h\in X_{\CH}(b_\FLocal,\mu_\FLocal^{-1})^\mathrm{pre}(R),\]
	and let $x\in\mathbb V_b$. The corresponding element
	$a_h\in\RZG(R)$ lies in $\CZ(x)$ if and only if
	\[
	h^{-1}\cdot x_0\in
	V_{\Z_p}\otimes_{W(k_\FLocal),\iota_0}W(R).
	\]
	Here the dot denotes the action of $H$ on the given
	$\FLocal$-representation, $\iota_0:k_\FLocal\hookrightarrow\bar k$
	is the distinguished residue-field embedding, and $x_0$ is the
	projection of $x$ to the $\iota_0$-summand, after inverting $p$, in
	the integral decomposition
	\[
	V_{\Z_p}\otimes_{\Z_p}W(R)
	\cong\bigoplus_{\iota:k_\FLocal\hookrightarrow\bar k}
	V_{\Z_p}\otimes_{W(k_\FLocal),\iota}W(R).
	\]
\end{cor}
\begin{proof}
	Use the representatives $g_0=h$ and
	$g_i=b_i\sigma(g_{i-1})$ constructed in the proof of
	\Cref{lem:ADLV-Shapiro}, and let $x_i$ be the components of $x$.
	Since $x\in\mathbb V_b$, Frobenius invariance gives
	$x_i=b_i\cdot\sigma(x_{i-1})$ for the special representation.
	Hence
	\[
	g_i^{-1}\cdot x_i
	=\sigma(g_{i-1}^{-1}\cdot x_{i-1})
	=\sigma^i(h^{-1}\cdot x_0)
	\quad(1\leq i\leq f-1).
	\]
	The semilinear map $\sigma$ identifies the integral special
	lattices on consecutive factors. Thus the transported vector
	is integral on all factors if and only if
	$h^{-1}\cdot x_0$ belongs to
	$V_{\Z_p}\otimes_{W(k_\FLocal),\iota_0}W(R)$.
	The assertion follows from
	\Cref{prop:special-ADLV-old}, and the criterion descends along
	\'etale covers of $\Spec R$.
\end{proof}

\subsection{Basic uniformization}\label{subsec:basic-uniformization}
We now assume that the local Shimura datum $(G_{\Q_p}, b, \{\mu\})$ comes from a global Shimura datum $(G,X)$ that satisfies \Cref{ass:Shimura-data}. 
Let $v\mid p$ be a place of the reflex field $\reflexHodge$ of $(G,X)$. Let $\CG$ be the parahoric Bruhat--Tits group scheme over
$\mathbb Z_p$ such that
\(K_p=\CG(\mathbb Z_p)\), \( K=K_pK^p.\)
Let $\{\mu\}$ be the conjugacy class of minuscule cocharacters of
$G_{\mathbb Q_p}$ determined by $X$ and an embedding
$\bar v:\overline\Q\to\overline{\Q}_p$ extending $v$.  The   local reflex field $\reflexHodgeLocal$
can be identified with $\reflexHodge_v$. 

Let $\bar x\in\Sh_K(\ov{\BF}_p)$ be a basic point whose associated local
Shimura datum is $(G_{\mathbb Q_p},b,\{\mu\})$. By \cite[Theorem~1]{KisinMadapusiPeraShin2022} such a point always exists.
Let $A_0$ be the corresponding abelian
variety with  tensors, and let $I=I_{\bar x}$ be the reductive group over
$\Q$ of tensor-preserving self-quasi-isogenies of $A_0$; polarizations are
preserved up to the prescribed similitude. 
Then $I$ is an inner form of $G$ with $I^\ad(\R)$ compact. 
There are natural identifications
$I_{\Q_p}\simeq J_b$, $I(\BA_f^p)\simeq G(\BA_f^p)$, 
the second after fixing a prime-to-$p$ level structure on $A_0$.

\begin{prop}\label{prop:padicunif}
	There is an
	isomorphism of formal schemes over
	$\Spf O_{\breve{(\reflexHodge)_v}}$
	\begin{equation}\label{eq:basic-uniformization}
		\Theta:
		I(\Q)\backslash
		\bigg[\RZG\times G(\BA_f^p)/K^p\bigg]
		\xrightarrow{\ \sim  }
		\widehat{\Sh}_{K,\basic}.
	\end{equation}
\end{prop}
\begin{proof}
	The uniformization near a single basic isogeny class follows from
	\cite[Theorem~4.10.6]{PappasRapoport2024PadicShtukas}.
	By \Cref{rmk:Shimura-data-strongly-admissible} and
	\cite[Theorem~2.2.7]{KisinZhou2025StronglyCompatible}, each basic isogeny
	class contains a point lifting to a special point. The argument of \cite[Proposition~5.2.12]{HeZhouZhu2025}, as applied in
	\cite[proof of Theorem~9.6]{hove2026cyclessplittingmodelsshimura}, implies there is only one basic isogeny class. This gives \eqref{eq:basic-uniformization}.
\end{proof}

\subsection{Compatibility with global special cycles}\label{subsec:unif-cycle}
\begin{defn}
	Let
	$
	\VspecialGlobal:=\sEnd(\bar x)\otimes_{\Z_{(p)}}\Q
	$
	be the space of special quasi-endomorphisms of $A_0$, equipped with its
	quadratic form.  The realization criterion for specialness and the
	faithfulness of the ambient Hom realizations give natural injections
	\[
	\VspecialGlobal(\BA_f^p)\hookrightarrow V(\BA_f^p),
	\qquad
	\VspecialGlobal_{\Q_p}\hookrightarrow\mathbb V_b.
	\]
	Thus, for every $r>0$, each tuple
	$\mathbf x\in (\VspecialGlobal)^r$ defines a local special cycle
	$\CZ(\mathbf x)\subset\RZG$. Let $I_{\mathbf x}\subset I$ be the
	pointwise stabilizer of the tuple $\mathbf x$.
\end{defn}

\begin{prop}\label{prop:uniformization-special-morp}
	The uniformization isomorphism \eqref{eq:basic-uniformization} induces a
	canonical isomorphism of fppf sheaves
	\begin{align*}
		\Psi': \sEnd_{\widehat{\Sh}_{K,\basic}}
		\simeq I(\Q)\backslash\Bigg[
		&\coprod_{x\in \VspecialGlobal}
		\CZ(x)\times G(\BA_f^p)/K^p
		\Bigg].
	\end{align*}
\end{prop}
\begin{proof}
	We first construct the map fppf locally. 
	Let $R$ be an $O_{\breve{(\reflexHodge)_v}}$-algebra in which $p$ is nilpotent,
	and let $z\in\RZG(R)$, $g\in G(\BA_f^p)$. Set $\bar R=R/pR$.
	The universal abelian scheme at $y=\Theta(z,g)\in\widehat{\Sh}_{K,\basic}(R)$
	is equipped with a quasi-isogeny
	$$\iota:A_y\times_{\Spec R}\Spec\bar R\dashrightarrow A_0\times_{\barFp}\Spec\bar R.$$
	Given a special endomorphism $f\in\sEnd(R)$, conjugation by $\iota$ gives a locally constant section
	of the constant sheaf associated with
	$\sEnd(\bar x)\otimes_{\Z_{(p)}}\Q=\VspecialGlobal$ on $\Spec\bar R$.
	After refining the fppf cover so that this section has a fixed value
	$x\in\VspecialGlobal$, the lifting condition gives $z\in\CZ(x)(R)$.
	This gives a section of the quotient sheaf
	\[
	\left[I(\Q)\backslash
	\left(\coprod_{x\in\VspecialGlobal}\CZ(x)\times G(\BA_f^p)/K^p\right)
	\right](R).
	\]
	For $i\in I(\Q)$, changing the representative acts by
	$(z,g,x)\mapsto(i\cdot z,ig,ixi^{-1})$, with
	$i\cdot(\mathscr X,\rho)=(\mathscr X,\rho\circ i^{-1})$.
	The resulting quotient-sheaf section is unchanged.
	The local constructions therefore glue to a morphism $\Psi'$ of formal schemes. 
	The Serre--Tate theorem implies $\Psi'$ is universally injective and formally \'etale. 
	Both sides have adic structure morphisms to $\widehat{\Sh}_{K,\basic}$,
	using \eqref{eq:basic-uniformization} for the right-hand side.
	Since $\Psi'$ respects these structure morphisms, it is adic.
	For a closed point represented by $(z,g,x)$ with $z\in\CZ(x)$ and $y=\Theta(z,g)$, the special quasi-endomorphism $\iota^{-1}x\iota\in\End(A_y)_\Q$ is integral on $A_y[p^\infty]$, hence belongs to $\End(A_y)\otimes_\Z\Z_{(p)}$.
	Thus $\Psi'$ is surjective on closed points.
	\Cref{lem:etale-map-bij}, applied to each scheme-theoretic truncation (or componentwise to obtain finite-type schemes), implies that $\Psi'$ is an isomorphism.
\end{proof}

For $T\in\Sym_r(\FGlobal)$ and a right $K^p$-orbit
$\Omega\subset V(\BA_f^p)^r$, let
\(\widehat{\CZ}_{\basic}(T,\Omega)\)
be the formal completion of the special cycle of
\Cref{defn:special-cycle} along the basic locus.

\begin{prop}\label{prop:uniformization-special-cycles}
	The uniformization isomorphism \eqref{eq:basic-uniformization} induces a
	canonical isomorphism of formal schemes
	\begin{align*}
		\widehat{\CZ}_{\basic}(T,\Omega)
		\simeq I(\Q)\backslash\Bigg[
		&\coprod_{\substack{\mathbf x\in(\VspecialGlobal)^r\\T(\mathbf x)=T}}
		\coprod_{\substack{g\in G(\BA_f^p)/K^p\\g^{-1}\mathbf x\in\Omega}}
		\bigl(\CZ(\mathbf x),g\bigr)
		\Bigg].
	\end{align*}
	Equivalently, after choosing representatives for the $I(\Q)$-orbits,
	the right-hand side is
	\[
	\coprod_{\substack{\mathbf x\in I(\Q)\backslash(\VspecialGlobal)^r\\
			T(\mathbf x)=T}}
	I_{\mathbf x}(\Q)\backslash
	\left[
	\CZ(\mathbf x)\times
	\bigl\{g\in G(\BA_f^p)/K^p:g^{-1}\mathbf x\in\Omega\bigr\}
	\right].
	\]
\end{prop}

\begin{proof}
	Restrict the correspondence of \Cref{prop:uniformization-special-morp} to tuples satisfying $T(\mathbf x)=T$ and $g^{-1}\mathbf x\in\Omega$.
\end{proof}

	\part{Orthogonal case}\label{part2}
	\section{Special cycles on orthogonal Shimura varieties}\label{sec:GSpin}
We now specialize the general construction to orthogonal Shimura varieties. 
We consider orthogonal Shimura varieties associated with quadratic spaces of dimension $\dimofV\geq3$. 
Let $(G_2=\Res_{\FGlobal/\mathbb Q}\SO(\Vcoherent),X_2)$ be the orthogonal Shimura datum defined below in
\Cref{subsec:SO-shimura-datum}.  Fix an even positive integer $\fkd$ and a
sufficiently small neat product level $K_2=\prod_{v<\infty}K_{2,v}\subset G_2(\A_f)$ whose factors away from
$\fkd$ are the parahoric subgroups defined by self-dual lattices. We check that this Shimura datum has Hodge-type covers that satisfy the properties listed in \Cref{ass:Shimura-data,ass:d=1,ass:J=1}. 

We only consider the case $\FGlobal\neq\Q$. 
By the general theory, for every finite place
$v\nmid\fkd$, we can define the canonical
integral model $\sS_{K_2,(v)}$ over $\mathcal O_{\FGlobal,(v)}$ together with a sheaf of special endomorphisms $\sEndint_{2,(v)}/\sS_{K_2,(v)}$. We construct the smooth global model
$\SOIntModel/\mathcal O_\FGlobal[1/\fkd]$ in \Cref{subsec:orth-global-models}.

\subsection{Clifford algebras}
We begin by fixing the Clifford-algebra conventions used in the Hodge-type
cover. 
Let $F$ be a field of characteristic zero, and let $(V,Q)$ be a
nondegenerate quadratic space of dimension $\dimofV$ over $F$. 
Its Clifford
algebra is
\(C(V)=T(V)\big/\langle v\otimes v-Q(v)\cdot 1:v\in V\rangle.\)
It is a free $F$-module of rank $2^\dimofV$. The involution $v\mapsto -v$ induces
its $\mathbb Z/2\mathbb Z$-grading
\(C(V)=C^+(V)\oplus C^-(V),\)
and the natural map $V\hookrightarrow C^-(V)$ identifies $V$ with a
subspace that generates $C(V)$ as an $F$-algebra.

The Clifford algebra also carries the canonical anti-involution
\(*:C(V)\longrightarrow C(V)\), \((v_1\cdots v_r)^*=v_r\cdots v_1.\) 
Write $\tr_{C/F}(a)$ for the trace of left multiplication by $a$
on the $F$-vector space $C(V)$; in particular $\tr_{C/F}(1)=2^{\dimofV}$.
The Clifford trace pairing is the bilinear map $(a,b)\longmapsto\tr_{C/F}(ab^*)$. Every Clifford trace
below uses this left regular trace and its scalar extensions. Traces
on $\End_F(W)$ refer to endomorphisms of the indicated vector space. 

For an $F$-algebra $R$, write $V_R=V\otimes_F R$ and
$C(V)_R=C(V)\otimes_F R$. The spinor similitude group is
\[
\GSpin(V)(R)=
\left\{
g\in C^+(V)_R^\times:
gV_Rg^{-1}=V_R,\quad g^*g\in R^\times
\right\}.
\]
The spinor similitude character is
\(\nu:\GSpin(V)\longrightarrow\Gm\),  \(g\longmapsto g^*g.\) 
Conjugation on $V$ gives a surjective homomorphism
\[
\pi:\GSpin(V)\longrightarrow\SO(V),
\qquad g\longmapsto(v\mapsto gvg^{-1}),
\]
with kernel the central scalar subgroup $\Gm$. Thus
\[
1\longrightarrow\Gm\longrightarrow\GSpin(V)
\xrightarrow{\,\pi\,}\SO(V)\longrightarrow 1.
\]
On the central subgroup, $\nu(z)=z^2$. This subgroup is the identity
component of the center of $\GSpin(V)$, and it is the  center when $\dimofV$
is odd.

\subsection{The Shimura datum and its Hodge-type cover}\label{subsec:SO-shimura-datum}

Assume that $\FGlobal$ is a totally real field. Fix the real place
$\RealPlace_0:\FGlobal\hookrightarrow\mathbb R$.  In accordance with
\eqref{eq:Vincoherent-Vcoherent-bridge}, $\Vcoherent=\Vincoherent^{\RealPlace_0}$ is
the nearby coherent quadratic space attached to the incoherent quadratic space 
$\Vincoherent/\BA_\FGlobal$; in particular, their finite localizations agree.
Suppose that $\Vcoherent$ has signature
$(\dimofV-2,2)$ at $\RealPlace_0$ and is positive definite at every other real place. Set
\[
\widetilde G=\Res_{\FGlobal/\Q}\GSpin(\Vcoherent),
\qquad
G_2=\Res_{\FGlobal/\Q}\SO(\Vcoherent),
\qquad
T_\FGlobal=\Res_{\FGlobal/\Q}\Gm,
\qquad T_\Q=\Gm.
\]
The preceding maps induce
\[
\pi:\widetilde G\longrightarrow G_2,
\qquad
\nu:\widetilde G\longrightarrow T_\FGlobal.
\]
\begin{defn}\label{defn:GSpin-datum}
	Let $\mathbb S=\Res_{\mathbb C/\mathbb R}\Gm$ be the Deligne torus. Choose
	an oriented negative-definite two-plane $z\subset \Vcoherent_{\RealPlace_0}$ and an oriented
	orthogonal basis $e_1,e_2$ with $Q(e_1)=Q(e_2)=-1$. Put
	$J=e_1e_2\in C^+(\Vcoherent_{\RealPlace_0})$, $ J^2=-1$.
	Under the isomorphism
	\(\widetilde G_{\mathbb R} \cong\prod_{v\mid\infty}\GSpin(\Vcoherent_v),\)
	define $\widetilde h_z:\mathbb S\to\widetilde G_{\mathbb R}$ by
	\[
	(\widetilde h_z)_v(x+yi)=
	\begin{cases}
		x+yJ,&v=\RealPlace_0,\\[4pt]
		1,&v\ne \RealPlace_0.
	\end{cases}
	\]
	Let $X$ be the $\widetilde G(\mathbb R)$-conjugacy class of $\widetilde h_z$.
	Let $h_z=\pi_{\mathbb R}\circ\widetilde h_z$, and let $X_2$ be the
	$G_2(\mathbb R)$-conjugacy class of $h_z$. Then $(G_2, X_2)$ is the orthogonal
	Shimura datum attached to $\Vcoherent$, and
	$(\widetilde G,X)\to(G_2,X_2)$ is a central morphism of Shimura data.
\end{defn}

The datum $(G_2, X_2)$ has abelian type, but need not have Hodge type. Choose an imaginary quadratic extension  $E_0/\Q$ and define $E=E_0\FGlobal$. 
Set $T_E=\Res_{E/\Q}\Gm$, and let
\[
N_{E/\FGlobal}:T_E\longrightarrow T_\FGlobal,
\qquad
\Delta:\Gm\longrightarrow T_\FGlobal
\]
be the norm and diagonal maps. 
\begin{defn}\label{defn:G-sharp}
	We now construct a Hodge-type cover. 
	Define
	\[
	G'=
	\left\{
	(\widetilde g,a,c)\in\widetilde G\times T_E\times\Gm:
	\nu(\widetilde g)N_{E/\FGlobal}(a)=\Delta(c)
	\right\}.
	\]
	The central embedding
	$T_\FGlobal\longrightarrow G'$, $t\longmapsto(t,t^{-1},1)$ 
	is well defined because $\nu(t)=t^2$ and $N_{E/\FGlobal}(t^{-1})=t^{-2}$.
	Set
	\(G^\sharp=G'/T_\FGlobal.\)
	Projection to $\widetilde G$, followed by $\pi$, descends to
	\(p:G^\sharp\longrightarrow G_2.\)
\end{defn}
\begin{defn}\label{defn:GSpin-Shimura-datum}
	We fix an embedding $\phi_0: E_0\to \C$. 
	Let $\Phi\subset \Hom_{\Q}(E,\C)$ be the corresponding CM type, and for every real
	place $v$ of $\FGlobal$ let the same symbol $v$ denote the member of
	$\Phi$ above it.  This gives an identification
	\((T_E)_{\mathbb R}\cong \prod_{v\mid\infty}\mathbb S.\)
	For $x+yi\in\mathbb S(\mathbb R)=\mathbb C^\times$, define
	$h'_z:\mathbb S\to G'_{\mathbb R}$ by
	\(h'_z(x+yi)=\bigl((\widetilde h_v(x+yi))_v,(a_v(x+yi))_v,x^2+y^2\bigr),\)
	where
	\[
	(\widetilde h_v(x+yi),a_v(x+yi))=
	\begin{cases}
		(x+yJ,1),&v=\RealPlace_0,\\[4pt]
		(1,x+yi),&v\ne \RealPlace_0.
	\end{cases}
	\]
	The defining relation for $G'$ holds at every real place.  Let $X'$ be
	the $G'(\mathbb R)$-conjugacy class of $h'_z$.  Let $h_z^\sharp$ be the
	composite of $h'_z$ with $G'_{\mathbb R}\to G^\sharp_{\mathbb R}$, and
	let $X^\sharp$ be its
	$G^\sharp(\mathbb R)$-conjugacy class.  The quotient map gives a central
	morphism of Shimura data
	\((G',X')\to(G^\sharp,X^\sharp)\), while $p$ gives a central morphism
	\[
	(G^\sharp,X^\sharp)\longrightarrow(G_2,X_2)
	\]
	inducing an isomorphism on adjoint data.  
	
	Identify $\FGlobal$ with its image under $\RealPlace_0$.  The reflex field of
	$(G_2,X_2)$ is $\FGlobal$.  Write $\varphi_0\in\Phi$ for the embedding
	above $\RealPlace_0$, and put $\Phi^\circ:=\Phi\setminus\{\varphi_0\}$.
	Then the reflex field of $(G^\sharp,X^\sharp)$ is $\reflexField$, characterized by 
	\begin{equation*}
		\operatorname{Gal}(\overline{\mathbb Q}/\reflexField)
		=\{\gamma:\gamma\RealPlace_0=\RealPlace_0,
		                 \ \gamma\Phi^\circ=\Phi^\circ\}.
	\end{equation*}
	In particular, when $\FGlobal\neq \Q$ we have $\reflexField= \phi_0(E_0)\FGlobal\subset \C$. In applications of \Cref{part1} the symbols
	$\reflexHodge$ and $\reflexAbel$ there specialize to $\reflexField$ and
	$\FGlobal$, respectively.
	\end{defn}

\begin{rmk}
	The kernel of $G'\to G_2$ is the torus $\{(t,a,c)\in T_\FGlobal\times T_E\times\Gm|t^2 N_{E/\FGlobal}(a)=\Delta(c)\}$. 
	Quotienting by $T_\FGlobal$ we conclude that the kernel of $p$ is the torus
	\(Z^{\Q}\cong T_E\times_{T_\FGlobal}\Gm\), 
	where the two maps to $T_\FGlobal$ are $N_{E/\FGlobal}$ and $\Delta$. Since
	$H^1(\Q,Z^{\Q})$ need not vanish, the induced finite \'{e}tale morphism
	of Shimura varieties need not be surjective, in contrast to the
	$\FGlobal=\Q$ case; cf. \cite[\S3.2]{MadapusiPera2016}. One must therefore include the
	$G^{\sharp,\mathrm{ad}}(\Q)$-action defined as in \cite[(3.2)]{Kisin2010}. 
\end{rmk}

\subsubsection{Changing the embedding \texorpdfstring{$\FGlobal\to\C$}{F into C}}\label{subsec:changing-iota}
Let $\RealPlace: \FGlobal\to \C$ be another embedding. 
As explained in
\cite[II. \S4]{Milne1990CanonicalModels}, base change of $\ShG_{K_2}(G_2,X_2)/\FGlobal$ along this embedding defines another (complex) Shimura variety $\ShG_{{^\RealPlace K_2}}({^\RealPlace G_2}, {^\RealPlace X_2})$, and ${^\RealPlace G_2}$ is an inner form of $G_2$. 

In our case, we may identify ${^\RealPlace G_2}$ with the Weil restriction of the orthogonal group $\SO(V^{(\RealPlace)})$, where $V^{(\RealPlace)}=\Vincoherent^{\RealPlace}$ is
the nearby coherent quadratic space at $\RealPlace$ attached to the incoherent quadratic space 
$\Vincoherent/\BA_\FGlobal$. In particular, $V^{(\RealPlace_0)}=\Vcoherent$. 

\subsection{Parahoric models}\label{subsec:orth-parahoric-models}
We fix the even integer $\fkd$ and an integral $O_\FGlobal$-lattice
$\VlattiveGlobal\subset\Vcoherent$ such that
$\VlattiveGlobal[1/\fkd]$ is self-dual. Fix the neat product level
$K_2=\prod_{v<\infty}K_{2,v}$, with $K_{2,\fkd}=\prod_{v\mid\fkd}K_{2,v}$,
where for $v\nmid\fkd$
\[
K_{2,v}=\SO(\VlattiveGlobal\otimes_{O_\FGlobal}O_{\FGlobal_v})(O_{\FGlobal_v}). 
\]
Let $p$ be a prime coprime to $\fkd$. Choose the CM pair $(E=E_0\FGlobal,\Phi)$ as above, with $E_0$
split at $p$.
The quadratic lattice $\VlattiveGlobalp=\VlattiveGlobal\otimes_{O_\FGlobal} O_{\FGlobal,(p)}$ is self-dual.

\subsubsection{The Hodge embedding}
We next verify that $(G^\sharp,X^\sharp)$ is of Hodge type.  Put
\(W:=C(\Vcoherent)\otimes_\FGlobal E\), viewed as an
$\FGlobal$-vector space. Let $W_\Q$ denote its underlying
$\Q$-vector space. 
We take
$W_{O_\FGlobal}=C(\VlattiveGlobal)\otimes_{O_\FGlobal}O_E$. 
 Let $\sigma$ be the nontrivial element of
$\Gal(E/\FGlobal)$ and define an anti-involution on $W$ by
\((x\otimes e)^\ddagger=x^*\otimes\sigma(e)\). 
\begin{defn}
	Consider the trace map $\tr_{W,E}:W\to E$, so $\tr_{W,E}(1)=2^{\dimofV}$. 
	For an element $\eta\in W^\times$ satisfying
	$\eta^\ddagger=-\eta$, define $\Psi_\eta: W\times W\to \FGlobal$ by 
	\begin{equation*}
		\Psi_\eta(w_1,w_2)
		:=\Tr_{E/\FGlobal}\!\left(\tr_{W,E}
		(w_1\eta w_2^\ddagger)\right),
		\qquad
		\psi_\eta:=\Tr_{\FGlobal/\Q}\circ\Psi_\eta.
	\end{equation*}
	The identity
	$\tr_{W,E}(b^\ddagger)=\sigma(\tr_{W,E}(b))$ shows that
	$\Psi_\eta$ is alternating.  It is nondegenerate because $\eta$ is
	invertible and the regular trace pairing of the semisimple algebra
	$W$ is nondegenerate.
	
	The formula
	\(\rho([\widetilde g,a,c])w=(\widetilde g\otimes a)w\)
	defines an $\FGlobal$-linear action of $G^\sharp$ on $W$, so the underlying
	$\Q$-representation factors through
	$\Res_{\FGlobal/\Q}\GL_{\FGlobal}(W)$.  Cyclicity of the trace and the defining relation of
	$G'$ give
	\begin{equation*}
		\Psi_\eta(\rho(g)w_1,\rho(g)w_2)
		=c\,\Psi_\eta(w_1,w_2),
		\qquad g=[\widetilde g,a,c]\in G^\sharp,
	\end{equation*}
	because
	$(\widetilde g\otimes a)^\ddagger(\widetilde g\otimes a)
	=\nu(\widetilde g)N_{E/\FGlobal}(a)=\Delta(c)$.
\end{defn}

For each real
place $\RealPlace$ use the member of $\Phi$ above $\RealPlace$ to identify
$E\otimes_{\FGlobal,\RealPlace}\R$ with $\C$, and put
$W_\RealPlace:=W\otimes_{\FGlobal,\RealPlace}\R$.  Let $I_\RealPlace$ be the action of
$h_z^\sharp(i)$ on $W_\RealPlace$.  Thus
\begin{equation*}
	I_{\RealPlace_0}=L_J\otimes 1,
	\qquad
	I_\RealPlace=1\otimes M_i\quad(\RealPlace\ne\RealPlace_0),
\end{equation*}
where $L_J$ denotes left multiplication by $J=e_1e_2$ and $M_i$
denotes multiplication by $i$ on $\C$.  	Let $\mathscr P_\RealPlace$ be the set
of $\eta_\RealPlace\in(W_\RealPlace^\times)^{\ddagger=-1}$ for which
$\Psi_{\eta_\RealPlace,\RealPlace}(w,I_\RealPlace w)>0$ for every $0\ne w\in W_\RealPlace$. This is an open cone containing $J\otimes1$ at $\RealPlace_0$ and $1\otimes i$ at every other real place: the involutions
$w\mapsto-Jw^\ddagger J$ and $w\mapsto w^\ddagger$, respectively, have positive trace forms.
If  $\eta\in W^{\ddagger=-1,\times}$
has image in $\mathscr P_\RealPlace$ for every real place $\RealPlace$, then there is an $\FGlobal$-linear Hodge embedding
\[
(G^\sharp,X^\sharp)\hookrightarrow
\bigl(\GSp_\Q(W_\Q,\psi_\eta),\mathfrak H^\pm\bigr).
\]

At $p$, the two split $E_0$-factors give the skew-adjoint unit $(1,-1)$.
Weak approximation in the
$\FGlobal$-vector space $W^{\ddagger=-1}$ gives a global $\eta\in W^{\ddagger=-1,\times}$ 
whose image belongs to $\mathscr P_v$ for every real place $v$, and at the
places above the fixed odd prime $p$,  
\(
\eta\otimes1\in(C(\VlattiveGlobalpcomplete)\otimes_{O_{\FGlobal,p}}O_{E,p})^{\times,\ddagger=-1}.
\) 
Then $W_{\Z_p}:=W_{O_\FGlobal}\otimes_\Z\Z_p=C(\VlattiveGlobalpcomplete)\otimes_{O_{\FGlobal,p}}O_{E,p}$ is self-dual for $\Psi_\eta$. 

Indeed, the Clifford trace pairing on $C(\VlattiveGlobalpcomplete)$ is perfect because
$\VlattiveGlobalp$ is self-dual and $p\ne2$, while the trace pairing on
$O_{E,p}/O_{\FGlobal,p}$ is perfect because this extension is split.  Their tensor product is therefore a perfect pairing on
$C(\VlattiveGlobalpcomplete)\otimes_{O_{\FGlobal,p}}O_{E,p}$, and right multiplication by
$\eta\otimes1\in\bigl(C(\VlattiveGlobalpcomplete)\otimes_{O_{\FGlobal,p}}O_{E,p}\bigr)^\times$ is an automorphism of $C(\VlattiveGlobalpcomplete)\otimes_{O_{\FGlobal,p}}O_{E,p}$.
\subsubsection{Verifying the hypotheses}
The preceding construction gives the $\FGlobal$-linear Hodge embedding and
the nondegenerate alternating $\FGlobal$-bilinear form $\Psi_\eta$ with
$\psi_\eta=\Tr_{\FGlobal/\Q}\circ\Psi_\eta$ required in \Cref{ass:Shimura-data}.
Under
\[
(G^\sharp)^c\simeq\Res_{\FGlobal/\Q}
\left(\bigl(\Spin(\Vcoherent)\times
\ker(N_{E/\FGlobal}:\Res_{E/\FGlobal}\Gm\to\Gm)\bigr)/\mu_2\right),
\]
with $\mu_2$ embedded diagonally, the restriction of $\rho$ to $(G^\sharp)^c$
comes from the $\FGlobal$-linear Clifford action. Thus the Hodge embedding
is fundamental in the sense of \Cref{defn:fundamental-emb}.
We now define the parahoric models and check that \Cref{ass:Shimura-data}(\ref{ass:Shimura-data:kottwitz})-(\ref{ass:Shimura-data:R-smooth}) and \Cref{ass:special-cycle} hold.  Since $\reflexField=\phi_0(E_0)\FGlobal$ is split over $\FGlobal$
at every place above $p$, we have \Cref{ass:Shimura-data}(\ref{ass:Shimura-data:reflex}). 

Since $G_2$ is semisimple, the Kottwitz condition  \Cref{ass:Shimura-data}(\ref{ass:Shimura-data:kottwitz}) trivially holds. 
The derived group $(G^\sharp)^{\der}=\Res_{\FGlobal/\Q}\Spin(\Vcoherent)$ is simply connected, so \Cref{ass:Shimura-data}(\ref{ass:Shimura-data:pi_1}) holds. 

The self-dual lattice $\VlattiveGlobalpcomplete$ defines the reductive group scheme
$\SO(\VlattiveGlobalpcomplete)$ over $O_{\FGlobal,p}$. Its generic fiber is
$\SO(\Vcoherent\otimes_\FGlobal\FGlobal_p)$ over
$\FGlobal_p:=\FGlobal\otimes_\Q\Q_p$. The product of the corresponding
special vertices in the Bruhat--Tits building at each place $v\mid p$ defines $\mathbf{x}\in\CB(G_2,\Q_p)$. 
We take $\CG^\sharp, \wt{\CG}, \CG_2$ to be the parahoric models for $G^\sharp, \wt{G}, G_2$ associated to $\mathbf{x}$. 
We take $K_p^\sharp=\CG^\sharp(\Z_p)$; the fixed level satisfies $K_{2, p}=\CG_2(\Z_p)$. 
In particular, \Cref{ass:Shimura-data}(\ref{ass:Shimura-data:level}, \ref{ass:Shimura-data:reductive}) is satisfied.

The projection $(a,c)\mapsto c$ gives an exact sequence of tori
\begin{equation}\label{eq:orth-center-exact}
	1\longrightarrow T_E^1
	\longrightarrow Z^\mathbb Q
	\longrightarrow\mathbb G_m
	\longrightarrow1,
	\qquad
	T_E^1:=\ker(N_{E/\FGlobal}:T_E\to T_{\FGlobal}).
\end{equation}
Since $E_0$ splits at $p$, the two $E_0$-factors give
\[
T_{\FGlobal,\mathbb Q_p}\times\mathbb G_{m,\mathbb Q_p}\xrightarrow{\ \sim\ }Z^{\mathbb Q}_{\mathbb Q_p}
=\{(a,b,c)\in T_{\FGlobal,\mathbb Q_p}^2\times\mathbb G_{m,\mathbb Q_p}:
ab=\Delta(c)\}
,
\qquad
(a,c)\longmapsto(a,\Delta(c)a^{-1},c).
\]
This product and $(T_E^1)_{\mathbb Q_p}\simeq T_{\FGlobal,\mathbb Q_p}$
are induced tori, hence quasi-tame and $R$-smooth by
\cite[Proposition~2.1.5(1)]{KisinPappasZhou2024ParahoricII}.
Their cocharacter lattices have permutation bases, whose inertia coinvariants
are free abelian. Thus $X_*(T_E^1)_I$ and $X_*(Z^\mathbb Q)_I$
are torsion-free, and
\Cref{ass:Shimura-data}(\ref{ass:Shimura-data:R-smooth}) is satisfied. 

For a place $\mathfrak p\mid p$ of $E_0$, let $\overline{\mathfrak p}$ be its conjugate.
Writing $a=(a_{\mathfrak p},a_{\overline{\mathfrak p}})$, we fix the splitting
\begin{equation}\label{eq:orth-split-GSpin-projection}
\begin{aligned}
G^\sharp_{\Q_p}&\xrightarrow{\ \sim\ }\widetilde G_{\Q_p}\times\mathbb G_{m,\Q_p},
&[g,a,c]&\longmapsto(g a_{\overline{\mathfrak p}},c),\\
(u,c)&\longmapsto[u,(\Delta(c)\nu(u)^{-1},1),c].
\end{aligned}
\end{equation}
It identifies the chosen parahoric with
$\CG^\sharp\simeq\widetilde{\CG}\times\mathbb G_{m,\Z_p}$. Under this isomorphism, the action of $\CG^\sharp$ on
$W_{\Z_p}\simeq C(\VlattiveGlobalpcomplete)\oplus C(\VlattiveGlobalpcomplete)$ is given by 
\begin{equation}\label{eq:orth-global-repn}
\widetilde{\CG}\times\mathbb G_{m,\Z_p}\longrightarrow \GL_{\Z_p}(C(\VlattiveGlobalpcomplete)\oplus C(\VlattiveGlobalpcomplete)), \qquad
(u,c)\longmapsto
\operatorname{diag}\bigl(c\nu(u)^{-1}L_u,L_u\bigr).
\end{equation}
\begin{prop}\label{prop:group-closed:orth}
	$\CG^\sharp\to \GL_{\Z_p}(W_{O_\FGlobal}\otimes_\Z\Z_p)$ is a closed
	immersion. Hence
	\Cref{ass:Shimura-data}(\ref{ass:Shimura-data:group-closed}) is satisfied. 
\end{prop}

\begin{proof}
	Taking Weil restriction of the closed immersion
	$\GSpin(\VlattiveGlobalpcomplete)\hookrightarrow
	\GL_{O_{\FGlobal,p}}(C(\VlattiveGlobalpcomplete))$
	defined by the Clifford action shows that
	$\widetilde{\CG}\to \GL_{\Z_p}(C(\VlattiveGlobalpcomplete))$
	is a closed immersion.
	The second block of \eqref{eq:orth-global-repn} determines $u$, and multiplication
	of the first block by $L_{u^*}=\nu(u)L_u^{-1}$ gives the scalar $c$.
	The assertion follows.
\end{proof}
\subsubsection{The integral special representation}\label{defn:orth-proj}
We now define the data in \Cref{ass:special-cycle}.  Left multiplication
gives a $\CG^\sharp$-equivariant embedding 
$L_-: \VlattiveGlobalpcomplete\hookrightarrow\End_{O_{\FGlobal,p}}(W_{\Z_p})$.
For $f\in \VlattiveGlobalpcomplete$ and $w_1,w_2\in W_{\Z_p}$, we have 
\(
\Psi_\eta(L_fw_1,w_2)
=\Psi_\eta(w_1,L_{f}w_2)
\), i.e., $L_f^\dagger=L_f$.  Since $\dim_{\FGlobal}\Vcoherent=\dimofV$, we have
$\rank_{O_{\FGlobal,p}}W_{\Z_p}=2^{\dimofV+1}$ and
$L_f^2=((f,f)/2)\id_{W_{\Z_p}}$.  Taking traces and polarizing yields
\[
(L_f,L_g)_\Psi=2^{\dimofV}(f,g),\qquad f,g\in\VlattiveGlobalpcomplete.
\]
Since $\VlattiveGlobalpcomplete$ is self-dual and $2^{\dimofV}\in O_{\FGlobal,p}^\times$, the
induced pairing on $\VlattiveGlobalpcomplete$ is self-dual, so the image of $L_-$ is
an orthogonal direct summand of $\End_{O_{\FGlobal,p}}(W_{\Z_p})$.
Since $\End_{O_{\FGlobal,p}}(W_{\Z_p})$ is saturated in $\End_{\Z_p}(W_{\Z_p})$,
$L_-$ also has saturated image over $\Z_p$.
The required constant in \Cref{ass:special-cycle} is $N=2^{\dimofV}$.

Note that $(G^\sharp)^\der\cong \Res_{\FGlobal/\Q}\Spin(\Vcoherent)$ is simply connected, while
\[
G^\sharp/(G^\sharp)^\der
\cong \Res_{\FGlobal/\Q}(\Res^1_{E/\FGlobal}\Gm)\times\Gm
\]
is $\Q$-rational. The local $H^1$-vanishing and Hasse principle for
simply connected semisimple groups lift sufficiently close rational points
of this quotient to $G^\sharp(\Q)$; weak approximation for
$(G^\sharp)^{\der}$ then adjusts these lifts at $p$ and at the real places.
Thus $G^\sharp$ satisfies weak approximation, and
$G^\sharp(\Q_p)=K_p^\sharp G^\sharp(\Q)_+$.
Consequently $J_{G^\sharp}=\{1\}$ by the same double-coset argument as in
\Cref{lem:J2-trivial}, so \Cref{ass:J=1} holds.
\subsection{Description of local models}\label{subsec:orth-local-model}
Fix a $p$-adic place $v$ of $\FGlobal$ and an embedding
$\wt v:\reflexField\hookrightarrow\overline{\Q}_p$ inducing $v$ on
$\FGlobal$. Let $\mathfrak p$ be the place of $E_0$ induced by
$\wt v\circ\phi_0$, and let $\overline{\mathfrak p}$ be its conjugate. 
The local Clifford model uses
$\Res_{O_{\FGlobal_v}/\Z_p}\GSpin(\VlattiveGlobalv)$ and the cocharacter
class whose inverse is $t\mapsto t^{-1}ef+fe$ at the embedding induced by
$\wt v$ and is trivial at all other embeddings, where $e,f$ are geometric
isotropic vectors with $(e,f)=1$. Its reflex ring is $O_{\FGlobal_v}$.
Since $E_0/\Q$ is split at $p$, both completions of $E_0$ equal $\Q_p$,
and $\reflexField_{\wt v}\simeq\FGlobal_v$. The resulting decomposition is
\[
O_{E,p}\simeq
\bigoplus_{v'\in\Sigma_{\FGlobal,p}}
\left((O_{\FGlobal,v'})_{(v',\mathfrak p)}
\oplus (O_{\FGlobal,v'})_{(v',\overline{\mathfrak p})}\right).
\]
Here the subscripts record the places of $\FGlobal$ and $E_0$ defining
each factor. Fix an embedding $\overline{\Q}\hookrightarrow\overline{\Q}_p$
extending $\wt v$. Through this embedding, the
embeddings in $\Phi$ induce $\mathfrak p$ on $E_0$, and those in
$\overline\Phi$ induce $\overline{\mathfrak p}$.
Let $\QuadHyper$ be the $O_{\FGlobal_v}$-scheme that represents the
rank $1$ isotropic direct summands of $\VlattiveGlobalv$. 
In other words, for an $O_{\FGlobal_v}$-algebra $R$ we have
\[\QuadHyper(R)=\{\mathcal L\subset\VlattiveGlobalv\otimes_{O_{\FGlobal_v}}R:
\mathcal L\text{ is a rank-one direct summand and }Q(\mathcal L)=0\}.\]
It is the quadric hypersurface $Q=0$ in $\mathbb P(\VlattiveGlobalv)$. At a
geometric point represented by a nonzero isotropic vector $u$, the
differential of $Q$ is $w\mapsto(u,w)$.  This functional is nonzero because
$\VlattiveGlobalv$ is self-dual and $2\in O_{\FGlobal_v}^\times$. The
Jacobian criterion therefore shows that $\QuadHyper$ is smooth of relative dimension $\dimofV-2$.  We first describe the map from $\QuadHyper$ to Grassmannians.

Let $\ell\in\VlattiveGlobalv\otimes_{O_{\FGlobal_v}}R$ generate an
isotropic line. Locally on $\Spec R$, self-duality allows us to choose
$s\in\VlattiveGlobalv\otimes_{O_{\FGlobal_v}}R$ such that
$Q(s)=0$ and $(\ell,s)=1$. In the Clifford algebra, this gives
$s^2=\ell^2=0$ and $s\ell+\ell s=1$. These relations give the direct-sum decomposition
\[
C(\VlattiveGlobalv)\otimes_{O_{\FGlobal_v}}R
=\ell\cdot\bigl(C(\VlattiveGlobalv)\otimes_{O_{\FGlobal_v}}R\bigr)
\oplus s\cdot\bigl(C(\VlattiveGlobalv)\otimes_{O_{\FGlobal_v}}R\bigr).
\]
Indeed, $1=\ell s+s\ell$ expresses every element as a sum of the two
summands, and an element in their intersection is killed by both
left multiplications $\ell$ and $s$, hence by $1$.
Left multiplication by $s$ maps the first summand isomorphically
to the second, with inverse left multiplication by $\ell$.
Both summands are therefore locally free of rank $2^{\dimofV-1}$ over $R$. 
\begin{defn}\label{defn:I-wt-local-model}
	Let $R$ be an $O_{\FGlobal_v}$-algebra and let
	$\mathcal L\subset\VlattiveGlobalv\otimes_{O_{\FGlobal_v}}R$ be an
	isotropic line. In the homological convention, the corresponding
	filtration on $C(\VlattiveGlobalv)\otimes_{O_{\FGlobal_v}}R$ has
	$\Fil^{-1}=C(\VlattiveGlobalv)\otimes_{O_{\FGlobal_v}}R$, $\Fil^1=0$, and
	\[
	\Fil^0\bigl(C(\VlattiveGlobalv)\otimes_{O_{\FGlobal_v}}R\bigr)
	=\mathcal L\cdot\bigl(C(\VlattiveGlobalv)\otimes_{O_{\FGlobal_v}}R\bigr).
	\]
	This submodule and its quotient are locally free of rank
	$2^{\dimofV-1}$ over $R$ by the preceding decomposition.
	Under the natural surjection
	\[
	\rho:C(\VlattiveGlobalv)\otimes_{\Z_p}R
	\twoheadrightarrow C(\VlattiveGlobalv)\otimes_{O_{\FGlobal_v}}R,
	\]
	define
	\[
	\Fil^0\bigl(C(\VlattiveGlobalv)\otimes_{\Z_p}R\bigr)
	=\rho^{-1}\!\left(
	\mathcal L\cdot\bigl(C(\VlattiveGlobalv)\otimes_{O_{\FGlobal_v}}R\bigr)
	\right),
	\]
	with $\Fil^{-1}=C(\VlattiveGlobalv)\otimes_{\Z_p}R$ and $\Fil^1=0$.
	Its degree-zero step is an $R$-direct summand of rank
	$2^{\dimofV}[\FGlobal_v:\Q_p]-2^{\dimofV-1}$.
	This defines a map
	\begin{equation*}
		\wt{i}:
		\QuadHyper
		\longrightarrow
		\Gr(C(\VlattiveGlobalv), 2^{\dimofV}[\FGlobal_v:\Q_p]-2^{\dimofV-1}). 
	\end{equation*}
\end{defn}

\begin{defn}\label{defn:I-sharp-local-model}
	For an $O_{\reflexField_{\wt v}}$-algebra $R$ and an isotropic line
	$\mathcal L\subset\VlattiveGlobalv\otimes_{O_{\FGlobal_v}}R$,
	we define the corresponding Hodge filtration on $W_{\Z_p}\otimes_{\Z_p}R$,
	with $\Fil^{-1}=W_{\Z_p}\otimes_{\Z_p}R$ and $\Fil^1=0$.
	The $O_{E,p}$-action gives the decomposition
	\[
	W_{\Z_p}\otimes_{\Z_p}R\simeq
	\bigoplus_{v'\in\Sigma_{\FGlobal,p}}
	\left(
	\bigl(C(\VlattiveGlobalvv)\otimes_{\Z_p}R\bigr)_{\mathfrak p}
	\oplus
	\bigl(C(\VlattiveGlobalvv)\otimes_{\Z_p}R\bigr)_{\overline{\mathfrak p}}
	\right),
	\]
	where the two subscripts specify the $p$-adic place of $E_0$ through
	which it acts. At the chosen place $v$, set
	\begin{equation}\label{eq:orth-integral-Hodge-pair}
	\begin{aligned}
	\Fil^0\bigl(C(\VlattiveGlobalv)\otimes_{\Z_p}R\bigr)_{\overline{\mathfrak p}}
	&=\rho^{-1}\!\left(
	\mathcal L\cdot\bigl(C(\VlattiveGlobalv)\otimes_{O_{\FGlobal_v}}R\bigr)
	\right),\\
	\Fil^0\bigl(C(\VlattiveGlobalv)\otimes_{\Z_p}R\bigr)_{\mathfrak p}
	&=\left[
	\rho^{-1}\!\left(
	\mathcal L\cdot\bigl(C(\VlattiveGlobalv)\otimes_{O_{\FGlobal_v}}R\bigr)
	\right)\right]^\perp.
	\end{aligned}
	\end{equation}
	Here $\rho$ is the surjection in \Cref{defn:I-wt-local-model},
	and $\perp$ is taken for the perfect
	$O_{\FGlobal_v}\otimes_{\Z_p}R$-valued Clifford trace pairing
	$(x,y)\mapsto\tr(xy^*)$ on the two identified copies of
	$C(\VlattiveGlobalv)\otimes_{\Z_p}R$.
	For $v'\ne v$, set
	\[
	\Fil^0\bigl(C(\VlattiveGlobalvv)\otimes_{\Z_p}R\bigr)_{\mathfrak p}=0,
	\qquad
	\Fil^0\bigl(C(\VlattiveGlobalvv)\otimes_{\Z_p}R\bigr)_{\overline{\mathfrak p}}
	=C(\VlattiveGlobalvv)\otimes_{\Z_p}R.
	\]
	To compute orthogonal complements over $R$, compose the Clifford trace pairing
	with the $R$-linear extension of a generator of
	$\Hom_{\Z_p}(O_{\FGlobal_v},\Z_p)$ as an $O_{\FGlobal_v}$-module.
	The resulting $R$-valued pairing is perfect. For an
	$O_{\FGlobal_v}\otimes_{\Z_p}R$-stable submodule its orthogonal complement
	agrees with that for the Clifford trace pairing. Hence the complements
	in \eqref{eq:orth-integral-Hodge-pair} are $R$-direct summands and commute with base change.
	Thus $\Fil^0(W_{\Z_p}\otimes_{\Z_p}R)$ is a direct summand as an $R$-module of rank
	\[(2^{\dimofV}[\FGlobal_v:\Q_p]-2^{\dimofV-1})+
	2^{\dimofV-1}+\sum_{v'\ne v}2^{\dimofV}[\FGlobal_{v'}:\Q_p]=2^{\dimofV}[\FGlobal:\Q].
	\]
	This defines a morphism
	\begin{equation*}
		i^\sharp:
		\QuadHyper\otimes_{O_{\FGlobal_v}}O_{\reflexField_{\wt v}}
		\longrightarrow
		\operatorname{Gr}_{\Z_p}
		\bigl(W_{\Z_p},2^{\dimofV}[\FGlobal:\Q]\bigr)
		\otimes_{\Z_p}O_{\reflexField_{\wt v}}.
	\end{equation*}
Over the fraction field $\reflexField_{\wt v}$ of the completed reflex
ring, extend scalars to $\overline{\Q}_p$ and decompose by the embeddings
$\sigma:E\hookrightarrow\overline{\Q}_p$. Via the fixed embedding
$\overline{\Q}\hookrightarrow\overline{\Q}_p$, the resulting ranks and
filtrations agree with the global Hodge filtration:
 
	\begin{equation*}
		\Fil^0\bigl(W_\Q\otimes_{E,\sigma}\overline{\Q}_p\bigr)=
		\begin{cases}
			0,&\sigma\in\Phi^\circ,\\
			W_\Q\otimes_{E,\sigma}\overline{\Q}_p,&\sigma\in\overline{\Phi^\circ},\\
			\mathcal L\cdot\bigl(W_\Q\otimes_{E,\sigma}\overline{\Q}_p\bigr),
			&\sigma=\varphi_0,\varphi_0^c.
		\end{cases}
	\end{equation*}
\end{defn}

\begin{lem}\label{lem:induced-fil-orth}
	Let $R$ be an $O_{\FGlobal_v}$-algebra and let
	$\mathcal L\subset\VlattiveGlobalv\otimes_{O_{\FGlobal_v}}R$ be
	an isotropic line. The filtration in \Cref{defn:I-wt-local-model}
	induces a filtration on
	$\End_R(C(\VlattiveGlobalv)\otimes_{O_{\FGlobal_v}}R)$,
	whose restriction to $\VlattiveGlobalv\otimes_{O_{\FGlobal_v}}R$ satisfies
	\[
	\Fil^0\bigl(\VlattiveGlobalv\otimes_{O_{\FGlobal_v}}R\bigr)
	=\mathcal L^\perp.
	\]
	The filtration in \Cref{defn:I-sharp-local-model} similarly induces
	a filtration on $\End_R(\CW_\loc\otimes R)$,
	whose restriction to $\CV_\loc\otimes R$ satisfies
	\[
	\Fil^0(\CV_\loc\otimes R)=\kappa^{-1}(\mathcal L^\perp),
	\]
	where
	\[
	\kappa:\CV_\loc\otimes R\twoheadrightarrow
	\VlattiveGlobalv\otimes_{O_{\FGlobal_v}}R
	\]
	is the natural surjection.
	In particular, $\wt{i}$ and $i^\sharp$ are proper monomorphisms,
	hence closed immersions.
\end{lem}

\begin{proof}
	Work locally on $\Spec R$ and choose a generator $\ell$ of
	$\mathcal L$. 	 
	An element $x\in\VlattiveGlobalv\otimes_{O_{\FGlobal_v}}R$ belongs
	to the induced $\Fil^0$ precisely when left multiplication by $x$
	preserves
	$\mathcal L\cdot\bigl(C(\VlattiveGlobalv)\otimes_{O_{\FGlobal_v}}R\bigr)$.
	If $x$ preserves this submodule, then
	$x\ell\in\ell\cdot\bigl(C(\VlattiveGlobalv)\otimes_{O_{\FGlobal_v}}R\bigr)$,
	and hence
	\[
	0=\ell x\ell=(x,\ell)\ell.
	\]
	This implies $(x,\ell)=0$. 
	Conversely, this equality gives $x\ell=-\ell x$, so
	left multiplication by $x$ preserves the submodule.
	This proves the first assertion.

	For $x\in\CV_\loc\otimes R$, preservation of the inverse image
	under $\rho$ in \Cref{eq:orth-integral-Hodge-pair} is equivalent to
	preservation of
	$\mathcal L\cdot\bigl(C(\VlattiveGlobalv)\otimes_{O_{\FGlobal_v}}R\bigr)$
	by $\kappa(x)$, hence to $\kappa(x)\in\mathcal L^\perp$.
	Left multiplication by $x$ is self-adjoint for the Clifford trace
	pairing, so it also preserves the orthogonal complement defining
	the other summand at $v$. At every $v'\ne v$, the degree-zero
	submodules are zero or the entire summand and are preserved by $x$.
	This proves the assertion for $\CV_\loc\otimes R$.
	
	The induced filtrations determine $\mathcal L^\perp$, hence $\mathcal L=(\mathcal L^\perp)^\perp$.
	Thus both maps are monomorphisms. Properness follows from that of $\QuadHyper$.
\end{proof}

\begin{thm}\label{thm:local-model}
Consider the local Hodge embeddings 
$$\wt i:\Res_{O_{\FGlobal,v}/\Z_p}\GSpin(\VlattiveGlobalv)
\longrightarrow\GL_{\Z_p}(C(\VlattiveGlobalv)),$$
$$i^\sharp:\CG^\sharp\longrightarrow \GL_{\Z_p}(W_{O_\FGlobal}\otimes_\Z\Z_p).$$
They induce closed immersions on local models: 
\begin{align*}
\wt i:\quad \QuadHyper&\hookrightarrow
\Gr_{\Z_p}\bigl(C(\VlattiveGlobalv),
2^{\dimofV}[\FGlobal_v:\Q_p]-2^{\dimofV-1}\bigr)
\otimes_{\Z_p}O_{\FGlobal_v},\\
i^\sharp:\quad \QuadHyper\otimes_{O_{\FGlobal_v}}O_{\reflexField_{\wt v}}
&\hookrightarrow
\Gr_{\Z_p}\bigl(W_{\Z_p},2^{\dimofV}[\FGlobal:\Q]\bigr)
\otimes_{\Z_p}O_{\reflexField_{\wt v}}.
\end{align*}
In particular, these local models are smooth over
$O_{\FGlobal_v}$ and $O_{\reflexField_{\wt v}}$, respectively, and
\Cref{ass:Shimura-data}(\ref{ass:local-model-smooth},
\ref{ass:Shimura-data:good}) hold.
\end{thm}

\begin{proof}
First work over the completed local reflex field $\FGlobal_v$ for the
local Clifford representation. Its flag variety is the quadric
$\QuadHyper\otimes_{O_{\FGlobal_v}}\FGlobal_v$, embedded in the Grassmannian of
$C(\VlattiveGlobalv)[1/p]$ by \Cref{defn:I-wt-local-model}.
For $G^\sharp_{\Q_p}$, the flag variety over $\reflexField_{\wt v}$ is
$\QuadHyper\otimes_{O_{\FGlobal_v}}\reflexField_{\wt v}$, embedded in the
Grassmannian of $W_{\Z_p}[1/p]$ by \Cref{defn:I-sharp-local-model}. 

The closed immersions $\wt i$ and $i^\sharp$ have smooth, hence flat,
sources over the indicated completed local reflex rings. Their sources
are therefore the reduced Zariski closures of these generic flag
varieties. The result then follows from
\cite[Lemma~3.4.1]{KisinPappasZhou2024ParahoricII}.
\end{proof}

\begin{prop}\label{prop:d=1-orth}
	The sheaf $\gr^{-1}\CV_\loc$ is invertible, and the natural inclusion
	\[
	\gr^{-1}\CV_\loc\longrightarrow\CH_\loc/\Fil^0\CH_\loc
	\]
	has locally free cokernel. If $\FGlobal\ne\Q$, then
	\Cref{ass:d=1} is satisfied.
\end{prop}
\begin{proof}
	We work on $\QuadHyper\otimes_{O_{\FGlobal_v}}O_{\reflexField_{\wt v}}$
	via \Cref{thm:local-model}, and write $\mathcal L$ for its
	tautological isotropic line. 
	
	Self-duality of $\VlattiveGlobalv$ and the direct-summand property
	of $\mathcal L$ imply that pairing with $\mathcal L$ surjects
	onto $\mathcal L^\vee$. Since $\kappa$ is surjective, \Cref{lem:induced-fil-orth}
	gives an exact sequence
	\[
	0\longrightarrow\Fil^0\CV_\loc\longrightarrow\CV_\loc
	\xrightarrow{\ w\mapsto(\kappa(w),-)\ }
	\mathcal L^\vee\longrightarrow0.
	\]
	Hence $\gr^{-1}\CV_\loc=\CV_\loc/\Fil^0\CV_\loc
	\simeq\mathcal L^\vee$  is invertible.
	Moreover,
	\[
		\CH_\loc/\Fil^0\CH_\loc
		\simeq\Homint\bigl(\Fil^0\CW_\loc,
		\CW_\loc/\Fil^0\CW_\loc\bigr)
	\]
	is locally free. By \Cref{lem:induced-fil-orth} we know that $\gr^{-1}\CV_\loc\longrightarrow\CH_\loc/\Fil^0\CH_\loc$ is injective after arbitrary base change hence has locally free cokernel.
	
	Finally, if $\FGlobal\ne\Q$, then $\Vcoherent$ is positive
	definite at a real place other than $\RealPlace_0$, so it is anisotropic
	over $\FGlobal$. The properness argument of
	\Cref{prop:properness} therefore applies to the Hodge-type integral
	model for $G^\sharp$. The action of $(G^\sharp)^\der$ on $\Vcoherent$
	factors through the standard orthogonal representation, which
	has no fixed vectors since $\dimofV\geq3$. Thus the properness
	and invariant-vanishing conditions in \Cref{ass:d=1} also hold.
\end{proof}

\subsection{Deformation of special cycles}

For $\FGlobal\ne\Q$, \Cref{prop:d=1-orth,prop:tangent-space-divisor}
show that the deformation locus of a nonzero special endomorphism is
an effective Cartier divisor, flat over the local reflex ring. 
Fix $v\nmid\fkd$, and let $p$ be its residue characteristic.
By \Cref{defn:special-abelian,rmk:descent-base,cor:same-special}, the
semi-global model $\sS_{K_2,(v)}$ carries an unramified, locally finitely
presented special-endomorphism scheme $\sEnd_{2,(v)}$.
\begin{lem}\label{prop:cycle-smooth-global}
	Let $v\nmid\fkd$ and $T\in O_{\FGlobal,(v)}^\times$.  The norm-$T$
	locus $\CZ_{T,(v)}\subset\sEnd_{2,(v)}$ is smooth over
	$O_{\FGlobal,(v)}$ of relative dimension $\dimofV-3$.
\end{lem}

\begin{proof}
	After the faithfully flat base change
	$O_{\FGlobal,(v)}\to O_{\FGlobal_v}$, it suffices to prove the assertion
	on the Hodge-type model $\Sh_{K^\sharp}$.
	By \Cref{prop:tangent-space-divisor}, the cycle is a relative Cartier
	divisor and is flat. It remains to show that its geometric special fiber is smooth.
	Fix a geometric point $s$ of the cycle, with underlying Shimura point $x$
	and special endomorphism $f$.
	The tensor-preserving deformation construction in the proof of
	\cite[Theorem~7.1.3]{KisinPappasZhou2024ParahoricII}, using its \S5.1.19,
	identifies the completed deformation ring at $x$ with that of the quadric
	in \Cref{thm:local-model}. On the first-order thickening
	$\barFp\hookrightarrow\barFp[t]/(t^2)$, this identification uses crystalline
	evaluation of the fixed Dieudonn\'e module and its Clifford tensors;
	thus the underlying deformations correspond to lifts of the Hodge isotropic line.
	By \Cref{cor:crystalline-special,lem:induced-fil-orth,thm:tangent-space},
	$f$ extends precisely when this lifted line is perpendicular to its
	crystalline vector $\kappa(f_\cris)$.
	Since $Q(\kappa(f_\cris))=T/2$ is a unit, the orthogonal complement
	is nondegenerate of rank $\dimofV-1$. Its isotropic-line quadric is
	smooth of dimension $\dimofV-3$. The tangent space of the special fiber
	at $s$ therefore has dimension $\dimofV-3$, equal to the dimension
	of the Cartier divisor. This proves the assertion.
\end{proof}

\subsection{Global integral models}\label{subsec:orth-global-models}
As explained in \Cref{subsec:int-Shimura}, we have the canonical
integral model $\sS_{K_2,(v)}$ smooth over $\mathcal O_{\FGlobal,(v)}$ for each place $v$ of $\FGlobal$ not dividing $\fkd$. 
By \cite[Theorem~2.2.1]{Lovering2017}, after inverting finitely many rational primes, including those dividing $\fkd$ and those with non-hyperspecial level, the generic fiber has a smooth integral canonical model. Gluing this model with the semi-global models at the finitely many remaining places $v\nmid\fkd$ gives a smooth scheme
$\SOIntModel/\mathcal O_\FGlobal[1/\fkd]$ 
whose localization at every such $v$ is $\sS_{K_2,(v)}$.
\begin{prop}\label{prop:properness}
	If $\FGlobal\neq \Q$, then the global integral canonical model
	$\SOIntModel$ is proper over $\mathcal O_\FGlobal[1/\fkd]$.
\end{prop}

\begin{proof}
	It is enough to prove that each localization $\sS_{K_2,(v)}$ is proper over
	$\mathcal O_{\FGlobal,(v)}$. By the definition \Cref{defn:int-model}, this reduces to the corresponding
	integral model for an auxiliary well-adapted Hodge-type Shimura datum $(G^\aux,X^\aux)$ whose adjoint datum
	is isomorphic to that of $(G_2,X_2)$. Since $\FGlobal\ne\Q$, choose
a real place $\RealPlace\ne\RealPlace_0$. The quadratic space $\Vcoherent_\RealPlace$
is positive definite, so $\SO(\Vcoherent)$ is anisotropic over
$\FGlobal$. Hence $G_2^{\ad}$, and therefore $(G^\aux)^{\ad}$, is
anisotropic over $\Q$. The properness of the
	Hodge-type integral model follows from the toroidal compactification
	theorem \cite[Theorem 1]{MadapusiPera2012ToroidalCompactifications}, as there exist no cusp labels.
\end{proof}

For $v\nmid\fkd$, let $p$ be its residue characteristic. Define
$\sEndint_{2,(v)}$ to be the open and closed locus on which
$\realization_\ell$ is integral for every $\ell\nmid\fkd p$.  Integrality at
$p$ is already part of the definition of $\sEnd_{2,(v)}$.  

\begin{prop}	
	The schemes
	$\sEndint_{2,(v)}$ over $\Sh_{K_2, (v)}$ for $v\nmid\fkd$, equipped with a quadratic pairing and $\ell$-adic realizations, have a common generic fiber.
\end{prop}

\begin{proof}
	Let $v, v'$ be places of residue characteristics $p, p'$. 
	Choose an imaginary quadratic field $E_0$ split at both $p$ and $p'$,
	and use the cover $(G^\sharp,X^\sharp)$ constructed from $E_0\FGlobal$ at both places.
	On its generic fiber, \Cref{lem:special-char-0-fiberwise} identifies
	rational specialness at $p$ and at $p'$ with Betti specialness, compatibly
	with realizations. For $\sEndint_{2,(v)}$,
	integrality is imposed at $p$ by the endomorphism condition and at every
	$\ell\nmid\fkd p$ by the $\ell$-adic realization; for
	$\sEndint_{2,(v')}$ the same conditions are imposed with $p'$ in place of $p$.
	Over characteristic zero, integrality at $p$ is equivalent to preservation
	of the $p$-adic Tate lattice. Thus the two pullbacks impose integrality
	at exactly the primes outside $\fkd$ and coincide, compatibly with pairing
	and realizations. The $G_2(\BA_f^{p,p'})$-Hecke action is transitive on
	connected components by \Cref{prop:Kisin2010-2.2.6}, so this identification
	extends to the entire generic fiber.
\end{proof}
\begin{defn}
	For
	$x\in\Vadl_{\fkd}/K_{2,\fkd}$ and $T\in\FGlobal$, let
	$\CZ_{T,(v)}(x)$ be the open and closed subscheme of
	$\sEndint_{2,(v)}$ on which the norm is $T$ and the product of the
	$\ell$-adic realizations at $\ell\mid\fkd$ belongs to the orbit $x$.
	
	More generally, for any $K_{2,\fkd}$-invariant compact subset $\Omega_\fkd\subseteq \Vadl_\fkd^r$ and $T\in \Sym_r(\FGlobal)$, we define $\CZ_{T,(v)}(\Omega_\fkd)$ to be the open and closed subscheme of $\sEndint_{2,(v)}^r$ on which the quadratic pairing is $T$ and the $\fkd$-adic realization lies in $\Omega_\fkd$.  
	For a representative $\mathbf x\in\Vadl_\fkd^r$, the notation
	$\CZ_{T,(v)}(\mathbf x)$ means $\CZ_{T,(v)}(K_{2,\fkd}\mathbf x)$;
	we use the same convention for the global cycles.
\end{defn}

\begin{lem}\label{lem:finite-global-gluing-block}
	The morphism
	$\CZ_{T,(v)}(\Omega_\fkd)\to\sS_{K_2,(v)}$ is finite.
\end{lem}

\begin{proof}
	Repeat the argument in the proof of
	\Cref{prop:fixed-level-special-cycle-finite}.
\end{proof}

We now glue the schemes $\sEndint_{2,(v)}$ to a scheme $\sEndint_2/\SOIntModel$. A standard limit argument gives uniqueness, compatibly with the prescribed localizations and their common generic fiber; existence is established in the following proposition.

\begin{prop}\label{prop:global-special-endomorphism-model}
	There exists a scheme $\sEndint_2/\SOIntModel$ whose localization at every such $v$ is $\sEndint_{2,(v)}/\sS_{K_2,(v)}$. 
	Its connected components are finite over $\SOIntModel$, and their localizations at $v$ are finite over $\sS_{K_2,(v)}$. 
	It is equipped with an $\ell$-adic realization for each $\ell\mid \fkd$ and a quadratic pairing. 
\end{prop}
\begin{proof}
	Note that $\sEndint_{2,(v)}=\coprod_{T\in \FGlobal, x\in\Vadl_{\fkd}/K_{2,\fkd}} \CZ_{T,(v)}(x)$. 
	We check this for each open and closed subscheme $\CZ_{T,(v)}(x)$. If $T=0$ then $\CZ_{T,(v)}(0)\cong \sS_{K_2,(v)}$, and all other branches are empty by positive definiteness of the special-endomorphism pairing. If $T\neq 0$, then for all but finitely many $v$, $T\in O_{\FGlobal,(v)}^\times$ so $\CZ_{T,(v)}(x)$ is smooth over $O_{\FGlobal,(v)}$ by \Cref{prop:cycle-smooth-global}. 
	
	For each irreducible component of $\CZ_{T,\FGlobal}(x)$, normalize
	the closure of its image in $\SOIntModel$ in the function field of that component.
	Let $\CZ_T'(x)$ be the disjoint union of these normalizations,
	with $\CZ_T'(x)=\varnothing$ when the generic fiber is empty. Its localization at $v$ coincides with $\CZ_{T,(v)}(x)$ if $\CZ_{T,(v)}(x)$ is normal and flat over $O_{\FGlobal,(v)}$. As a corollary, $\CZ_{T,(v)}(x)\cong \CZ_{T,(v)}'(x)$ for all but finitely many $v$; we then glue the remaining finitely many $\CZ_{T,(v)}(x)$ to define a finite scheme $\CZ_{T}(x)/\SOIntModel$. Its localization $\CZ_{T,(v)}(x)$ is finite over $\sS_{K_2,(v)}$. We define $\sEndint_{2}=\coprod_{T\in \FGlobal, x\in\Vadl_{\fkd}/K_{2,\fkd}} \CZ_{T}(x)$.
	
We now glue the quadratic pairing on each fiber product
\[
\CZ_T(x)\times_{\SOIntModel}\CZ_{T'}(x').
\]
Its generic-fiber pairing spreads out over $O_\FGlobal[1/\fkd']$
for some $\fkd\mid\fkd'$. After enlarging $\fkd'$, this fiber product
has schematically dense generic fiber: its vertical torsion is
supported over finitely many places. For $v\nmid\fkd'$, schematic
density identifies the spread-out pairing with the given local
pairing on
\[
\CZ_{T,(v)}(x)\times_{\sS_{K_2,(v)}}\CZ_{T',(v)}(x').
\]
We then glue in the prescribed pairings at the finitely many places
$v\mid\fkd'$, $v\nmid\fkd$, using their agreement on the generic fiber.

	On each $\CZ_T(x)$, choose a sub $\Z_\ell$-lattice $U\subseteq\CV_\ell$ containing $xK_{2,\fkd}$. Then the image of each realization map  $\realization_\ell(\CZ_{T,(v)}(x))$ lies in $U$, 
	and their reductions modulo $\ell^s$ glue to a section of
	$\Gamma(\CZ_T(x),U/\ell^sU)$ for every $s\geq1$.
	Taking the inverse limit gives the $\ell$-adic realization.
\end{proof}

\section{Special cycles on orthogonal \RZname spaces}\label{sec:GSpin-RZ}

\subsection{GSpin \RZname spaces}\label{subsec:GSpin-RZ}
We first isolate the local datum carried by the distinguished place.
Let $p$ be an odd prime and $\FLocal/\Q_p$ a finite extension.
Let $(\Vlocal,Q)$ be a nondegenerate quadratic space over $\FLocal$ of
dimension $\dimofV\geq3$, and let
$\VlattiveLocaltot\subset\Vlocal$ be a self-dual $O_\FLocal$-lattice.  

Since $p\neq2$ and $\VlattiveLocaltot$ is self-dual,
$\GSpin(\VlattiveLocaltot)$ is reductive over $O_\FLocal$; thus
$\Res_{O_\FLocal/\Z_p}\GSpin(\VlattiveLocaltot)$ is the connected parahoric group scheme
attached to the corresponding point of the Bruhat--Tits building.  

\begin{defn}\label{defn:GSpin-pure-local-Shimura-datum}
	Choose an embedding $\tau_0:\FLocal\hookrightarrow\ov\Q_p$.  Over $\ov\Q_p$ we
	have
	\[
	(\Res_{\FLocal/\Q_p}\GSpin(\Vlocal))_{\ov\Q_p}
	\simeq
	\prod_{\tau:\FLocal\hookrightarrow\ov\Q_p}
	\GSpin(\Vlocal\otimes_{\FLocal,\tau}\ov\Q_p).
	\]
The component of $\mu^{-1}$ at every
	$\tau\neq\tau_0$ is trivial.  On the $\tau_0$-component, choose over
	$\ov\Q_p$ a Witt decomposition
	\(\Vlocal\otimes_{\FLocal,\tau_0}\ov\Q_p
=\ov\Q_p e\oplus V_0\oplus\ov\Q_p f\),  \((e,f)=1.\)
	Then a representative of $\{\mu^{-1}\}$ is the Clifford-algebra
	cocharacter
	\begin{equation}\label{eq:GSpin-standard-cocharacter}
		\mu^{-1}(t)=t^{-1}ef+fe\in C^+(\Vlocal\otimes_{\FLocal,\tau_0}\ov\Q_p)^\times.
	\end{equation}
	Since $ef+fe=1$, this gives a homomorphism
$\mathbb G_{m,\ov\Q_p}\to\GSpin(\Vlocal)_{\ov\Q_p}$
for scalar extension through $\tau_0$; under the
	vector representation it acts by
	\(e\longmapsto t^{-1}e\), \(V_0\longmapsto V_0\), \(f\longmapsto tf,\) 
	and under left multiplication on $C(\Vlocal)$ it has weights $-1$ and $0$.
	This is the standard GSpin minuscule class of
	\cite[\S4.2.1]{HowardPappas2017Spinor}.  Its conjugacy class is defined over
	$\FLocal$.  Therefore the stabilizer of the absolute class $\{\mu\}$ in
	$\Gal(\ov\Q_p/\Q_p)$ is the decomposition subgroup of $\tau_0$, and the
	local reflex field is the embedded copy $\tau_0(\FLocal)$, which we identify with
	$\FLocal$.
	
	Fix a basic element
	\[
	b\in \Res_{\FLocal/\Q_p}\GSpin(\Vlocal)(\breve{\Q}_p),
	\qquad [b]\in B(\Res_{\FLocal/\Q_p}\GSpin(\Vlocal),\mu^{-1}),
	\]
	with the Kottwitz invariant prescribed by $\mu^{-1}$.  Together with the
	parahoric group scheme $\Res_{O_\FLocal/\Z_p}\GSpin(\VlattiveLocaltot)$, this gives the integral
	local Shimura datum
	\((\Res_{\FLocal/\Q_p}\GSpin(\Vlocal),b,\{\mu\},\Res_{O_\FLocal/\Z_p}\GSpin(\VlattiveLocaltot)).\)
	We write
	\[
	\GSpinRZ
	:=\mathcal M_{\Res_{O_\FLocal/\Z_p}\GSpin(\VlattiveLocaltot),b,\mu}
	/\Spf O_{\breve \FLocal}
	\]
	for its integral \RZname space.
	
	The Clifford algebra gives a faithful closed immersion 
	$
	\GSpin(\Vlocal)\hookrightarrow
	\GL_\FLocal(C(\Vlocal))$ 
	by left multiplication.  The lattice $C(\VlattiveLocaltot)$ and the grading,
	right Clifford action, and distinguished subspace
	$\Vlocal\subset\End_{C(\Vlocal)}(C(\Vlocal))$ give the integral Clifford
	tensors cutting out $\GSpin(\VlattiveLocaltot)$. 
\end{defn}

We describe the quadratic space of special quasi-endomorphisms abstractly.  Let
$b_\FLocal\in\GSpin(\Vlocal)(\breve \FLocal)$ be the corresponding relative representative
and let $\sigma_\FLocal$ denote Frobenius on $\breve \FLocal/\FLocal$.  The vector
representation $\rho_V:\GSpin(\Vlocal)\to\SO(\Vlocal)$ defines the semilinear operator $\rho_V(b_\FLocal)\circ\sigma_\FLocal$ on $\Vlocal\otimes_\FLocal\breve \FLocal$. 
\begin{prop}\label{prop:V_b-SO}
	Consider the $\FLocal$-vector space  $\Vspecial	:=	(\Vlocal\otimes_\FLocal\breve \FLocal)^{b_\FLocal \sigma_\FLocal=1}$ as in \Cref{subsec:V_b}. Then $\Vspecial$ is a quadratic space of dimension $\dimofV$ with Hasse invariant $-1$ and satisfies $\chi(\Vspecial)=\chi(\Vlocal)$.
\end{prop}

\begin{proof} Since $b$ is basic, its Newton
	cocharacter factors through the identity component of the center of
	$\GSpin(\Vlocal)$. In this range,
	$Z(\GSpin(\Vlocal))^\circ=\mathbb G_m
	=\ker(\GSpin(\Vlocal)\to\SO(\Vlocal))$.
	Hence the Newton cocharacter acts trivially through the vector
	representation, so the isocrystal
	$(\Vlocal\otimes_\FLocal\breve \FLocal,b_\FLocal \sigma_\FLocal)$ is isoclinic of slope $0$. 
	It follows
	that $\dim_\FLocal\Vspecial=\dim_\FLocal\Vlocal=\dimofV$ and
	$\Vspecial\otimes_\FLocal\breve \FLocal\simeq\Vlocal\otimes_\FLocal\breve \FLocal$.
	
	Let $\nu:\GSpin(\Vlocal)\to\mathbb G_m$ be the spinor similitude character.
	The computation in \cite[Proposition~4.2.5]{HowardPappas2017Spinor}
	shows that the fixed quadratic space of a basic element has the same
	discriminant as $\Vlocal$ and Hasse invariant
	$(-1)^{\valFLocal(\nu(b_\FLocal))}$. 
	The Kottwitz condition and \eqref{eq:GSpin-standard-cocharacter} give
	$\valFLocal(\nu(b_\FLocal))=-1$ so $\Has(\Vspecial)=-1$.
\end{proof}

\subsection{Basic uniformization in the orthogonal case}
We now return to the local datum induced by the global orthogonal Shimura
datum defined in \Cref{subsec:SO-shimura-datum,subsec:orth-parahoric-models}.  
The decomposition $O_{\FGlobal,p}\cong\prod_{v\mid p}O_{\FGlobal,v}$
induces $\VlattiveGlobalpcomplete\cong\operp_{v\mid p}\VlattiveGlobalv$.
Inverting $p$ gives $\VGlobalp\cong\operp_{v\mid p}\VGlobalv$.
For each place $v\mid p$ of the totally real field $\FGlobal$, put
\[
\widetilde{G}_v:=\Res_{\FGlobal_v/\Q_p}\GSpin(\VGlobalv),
\qquad
 \widetilde{\mathcal G}_v
:=\Res_{O_{\FGlobal_v}/\Z_p}\GSpin(\VlattiveGlobalv).
\]
Then
\[
\widetilde G_p\simeq\prod_{v\mid p}\widetilde{G}_v,
\qquad
\widetilde{\mathcal G}_p\simeq\prod_{v\mid p} \widetilde{\mathcal G}_v.
\]

\begin{defn}\label{defn:GSpin-local-Shimura-datum}
	Fix the distinguished place $v_0\mid p$.  Let
	$\{\mu^{(v_0)}\}$ be the cocharacter class of $\widetilde G_p$ which is
	trivial on $\widetilde{G}_v$ for $v\neq v_0$ and whose $v_0$-component, denoted
	$\{\mu_{v_0}\}$, is the class of
	\Cref{defn:GSpin-pure-local-Shimura-datum}.  The local reflex field of
	$\{\mu^{(v_0)}\}$ is $\FGlobal_{v_0}$, which is identified with the
	fixed local field $\FLocal$.  Let
	$b=(b_v)_{v\mid p}$ be the image of the global basic element under the
	first projection of \eqref{eq:orth-split-GSpin-projection}, using the
	place $\mathfrak p$ induced by the selected reflex place above $v_0$. Then
	\([b]\in B\bigl(\widetilde G_p,(\mu^{(v_0)})^{-1}\bigr).\)
	Then $[b_v]$ is the neutral class for every $v\neq v_0$, while
	$(\widetilde G_{v_0},b_{v_0},\{\mu_{v_0}\},\widetilde{\mathcal G}_{v_0})$ is the pure local GSpin
	datum of \Cref{defn:GSpin-pure-local-Shimura-datum} under the identification
	$\FGlobal_{v_0}\simeq\FLocal$.
	We call
	\((\widetilde G_p,b,\{\mu^{(v_0)}\},\widetilde{\mathcal G}_p)\)
	the GSpin local Shimura datum attached to $v_0$.
\end{defn}

The product decomposition of the group gives a product decomposition of integral local Shimura
spaces over $\Spf O_{\breve{\FGlobal}_{v_0}}$:
\[
\mathcal M_{\widetilde{\mathcal G}_p,b,\mu^{(v_0)}}
\simeq
\GSpinRZ^{(v_0)}
\times
\prod_{v\mid p, v\neq v_0}
\underline{
	\GSpin(\VGlobalv)(\FGlobal_v)/
	\GSpin(\VlattiveGlobalv)(O_{\FGlobal_v})
},
\]
where an underlined coset set denotes the disjoint union of copies of
$\Spf O_{\breve{\FGlobal}_{v_0}}$ indexed by that set, or its associated
constant $v$-sheaf, and
$\GSpinRZ^{(v_0)}:=\mathcal M_{\widetilde{\mathcal G}_{v_0},b_{v_0},\mu_{v_0}}$.
Indeed, for $v\neq v_0$ the local datum is $(\widetilde{G}_v,1,\{0\})$, so its 
integral local Shimura space is the constant coset sheaf
$\underline{\widetilde{G}_v(\Q_p)/ \widetilde{\mathcal G}_v(\Z_p)}$ by \cite[(6.1.1)]{PappasRapoport2026IntegralLocalShimura}.

We use the local  Hodge embedding $\widetilde \CG_p\to \GL_{\Z_p}(C(\VlattiveGlobalpcomplete))$ and the left multiplication map  $\VlattiveGlobalpcomplete\subseteq \End_{\Z_{p}}(C(\VlattiveGlobalpcomplete))$ to define the special cycles $\CZ_{\mathrm{GSpin}}(x)$ on $\mathcal M_{\widetilde{\mathcal G}_p,b,\mu^{(v_0)}}$.

\begin{prop}\label{prop:special-orth-global}
	As in \Cref{subsec:unif-cycle}, let $\VspecialGlobal$ denote the
	$\FGlobal$-quadratic space of special quasi-endomorphisms of the abelian variety at a basic point. Then $\VspecialGlobal$ is the nearby quadratic space of $\Vincoherent$ at the finite place $v_0$, i.e., $\VspecialGlobal$
has Hasse invariant $-1$ at $v_0$ and
$\VspecialGlobal_w\cong\Vincoherent_w$ for every place $w\ne v_0$.
\end{prop}

\begin{proof}
	Note that $\Vadl_p\cong \operp_{v\mid p}\VGlobalv$ and we may choose $b$ so that it acts trivially on $\Vadl_{v}$ for $v\neq v_0$. 
	By \Cref{prop:special-F_p} and \Cref{prop:V_b-SO}, we have a decomposition $\VspecialGlobal_p\cong \VspecialGlobal_{v_0}\operp \operp_{v\mid p, v\neq v_0}\VGlobalv$ where $\VspecialGlobal_{v_0}$ has Hasse invariant $-1$, and for every non-$p$-adic finite place $v$ we have $\VspecialGlobal_v\cong \VGlobalv$.  At every real
	place the special-endomorphism pairing is positive definite by
	Rosati positivity. 
\end{proof}
\begin{prop}\label{prop:dec-special-cycle}
	For $x=(x_v)\in\mathop{\operp}_{v\mid p}(\VspecialGlobal)_v$, the special
	cycle $\CZ_{\mathrm{GSpin}}(x)$ is
	\[\CZ_{\mathrm{GSpin}}(x_{v_0}) \times 
	\prod_{v\mid p, v\neq v_0} \underline{\{g\in \GSpin(\VGlobalv)(\FGlobal_v)| g^{-1}x_v\in \VlattiveGlobalv\}/
		\GSpin(\VlattiveGlobalv)(O_{\FGlobal_v})}.\]
\end{prop}
\begin{proof}
	The universal $p$-divisible group on $\underline{
		\GSpin(\VGlobalv)(\FGlobal_v)/
		\GSpin(\VlattiveGlobalv)(O_{\FGlobal_v})
	}$ associated to the Hodge embedding $\Res_{\FGlobal_v/\Q_p}\GSpin(\VGlobalv)\to \GL_{\Q_p}(C(\VGlobalv))$ is the constant $p$-divisible group associated to the $O_{\FGlobal,v}$-lattice $C(\VlattiveGlobalv)$ with framing given by elements in $\GSpin(\VGlobalv)(\FGlobal_v)/
	\GSpin(\VlattiveGlobalv)(O_{\FGlobal_v})$. The special endomorphisms are
	elements of $\VGlobalv$ acting by left multiplication. The claim follows.
\end{proof}
For each $r\geq1$, we use
$\phi^{v_0}_p=\prod_{v\mid p, v\neq v_0}\charfun_{\VlattiveGlobalv^r}$
on $\prod_{v\mid p,\,v\neq v_0}\VGlobalv^r$. We also write
\[\tilde{G}^{v_0}(\Q_p)=\prod_{v\mid p, v\neq v_0}\GSpin(\VGlobalv)(\FGlobal_v), \qquad \tilde{K}^{v_0}_p =\prod_{v\mid p, v\neq v_0}\GSpin(\VlattiveGlobalv)(O_{\FGlobal_v}).\]

Now we consider the \RZname space arising from the global Shimura datum $(G^\sharp, X^\sharp)$. 
\begin{prop}\label{prop:uniformization-sharp-Shimura}
	We have a canonical uniformization isomorphism of formal schemes
	over $O_{\breve{\reflexField},{\wt v}}$
	\begin{align*}
		\widehat{\Sh_{K^\sharp}}_{\basic}
		\simeq I^\sharp(\Q)\backslash\Bigg[
		\GSpinRZ^{(v_0)}
		\times
		\underline{\tilde G^{v_0}(\Q_p)/\tilde K^{v_0}_p}
		\times
		\underline{\Q_p^\times/\Z_p^\times}
		\times
		G^\sharp(\BA_f^p)/K^{\sharp,p}
		\Bigg].
	\end{align*}
	For $r>0$, $T\in\Sym_r(\FGlobal)$, and a right $K^{\sharp,p}$-orbit
	$\Omega\subset\Vadl(\BA_f^p)^r$, the special cycles correspond as follows:
	\begin{align*}
		\widehat{\CZ^{\sharp}(T,\Omega)}_{\basic}
		\simeq I^\sharp(\Q)\backslash\Bigg[
		\coprod_{\substack{\mathbf x\in(\VspecialGlobal)^r\\T(\mathbf x)=T}}
		\coprod_{g}
		\bigl(\CZ_{\mathrm{GSpin}}(\mathbf x_{v_0}),g\bigr)
		\times
		\underline{\Q_p^\times/\Z_p^\times}
		\Bigg],
	\end{align*}
	where $g=(g^p,g_p)\in G^\sharp(\BA_f^p)/K^{\sharp,p}
	\times \tilde G^{v_0}(\Q_p)/\tilde K^{v_0}_p$ satisfies
	$(g^p)^{-1}\mathbf x^p\in\Omega$ and
	$\phi_p^{v_0}(g_p^{-1}\mathbf x_p^{v_0})=1$. 
The inner form $I^\sharp$ fits into a short exact sequence
	\[1\to Z^\Q\to I^\sharp\to I\cong \Res_{\FGlobal/\Q}\SO(\VspecialGlobal)\to 1.\]
\end{prop}
\begin{proof}
	Under $\CG^\sharp_{\Q_p}\simlr\widetilde \CG_{\Q_p}\times\mathbb G_{m,\Q_p}$ from \eqref{eq:orth-split-GSpin-projection}, the Hodge cocharacter is $\mu=(\mu^{(v_0)},\operatorname{id}_{\mathbb G_m})$, with $\mu^{(v_0)}$ as in \Cref{defn:GSpin-local-Shimura-datum}. By \cite[Prop. 25.2.1]{ScholzeWeinstein2020Berkeley} the associated \RZname space is isomorphic to 
	\[\mathcal M_{\widetilde{\mathcal G}_p,b,\mu^{(v_0)}}\times	\underline{\Q_p^\times/\Z_p^\times}\cong \GSpinRZ^{(v_0)}
	\times
	\prod_{v\mid p, v\neq v_0}
	\underline{
		\GSpin(\VGlobalv)(\FGlobal_v)/
		\GSpin(\VlattiveGlobalv)(O_{\FGlobal_v})}\times 
		\underline{\Q_p^\times/\Z_p^\times}.\]
	The first isomorphism then follows from \Cref{prop:padicunif}. 
	
	In \eqref{eq:orth-global-repn}, the two representation blocks are $c\nu(u)^{-1}L_u$ and $L_u$; special endomorphisms act by left multiplication on both blocks.
	Let $\mathscr{X}$ be the universal $p$-divisible group on $\mathcal M_{\widetilde{\mathcal G}_p,b,\mu^{(v_0)}}$. The universal $p$-divisible group on $\mathcal M_{{\CG^\sharp}_p,b,\mu}$ decomposes as $\mathscr{X}\oplus \mathscr{X}^*$, where $\mathscr{X}^*$ is the Cartier dual of $\mathscr{X}$.
	Hence, over any test scheme, each $x\in \Vspecialp$ defines a special endomorphism on $\mathscr{X}$ if and only if it defines a special endomorphism on $\mathscr{X}\oplus \mathscr{X}^*$. We conclude that the special cycle $\CZ^\sharp(x)$ is isomorphic to $\CZ_{\mathrm{GSpin}}(x)\times 
	\underline{\Q_p^\times/\Z_p^\times}.$ 
	The second isomorphism then follows from \Cref{prop:uniformization-special-cycles,prop:dec-special-cycle}.
	
	Let $I$ be the inner form of $G_2$ given by the same cohomology class $[I^\sharp]\in H^1(\Q, G^{\sharp,\ad})$. 
	By \Cref{prop:special-orth-global}, conjugation on special quasi-endomorphisms gives
	$I\to\Res_{\FGlobal/\Q}\SO(\VspecialGlobal)$.
	For every prime $\ell\ne p$, the $\ell$-adic realization identifies
	this map with the standard orthogonal representation of
	$I_{\Q_\ell}\cong G_{2,\Q_\ell}$, so it is faithful.
	Both groups are connected and have the same dimension, giving
	$I\cong\Res_{\FGlobal/\Q}\SO(\VspecialGlobal)$.
\end{proof}

\begin{prop}\label{prop:local-special-cycle-properties}
	For any $x\in\Vspecial\setminus\{0\}$, the closed formal subscheme
	$\CZ_{\mathrm{GSpin}}(x)\subseteq\GSpinRZ$ is a relative effective Cartier divisor over
	$\Spf O_{\breve\FLocal}$. If moreover $(x,x)\in O_\FLocal^\times$,
	then $\CZ_{\mathrm{GSpin}}(x)$ is formally smooth over $O_{\breve\FLocal}$, and hence over
	$O_\FLocal$, of relative dimension $\dimofV-3$.
\end{prop}
\begin{proof}
	Let $\varpi$ be a uniformizer of
	$O_\FLocal$. The framing action identifies
	$J_b(\Q_p)$ with $\GSpin(\Vspecial)(\FLocal)$ and satisfies
	\[
	g\CZ_{\mathrm{GSpin}}(v)=\CZ_{\mathrm{GSpin}}(\bar g v),\qquad
	g\in\GSpin(\Vspecial)(\FLocal),\quad v\in\Vspecial,
	\]
	where $\bar g$ is the image in $\SO(\Vspecial)(\FLocal)$.
	Also $\CZ_{\mathrm{GSpin}}(av)=\CZ_{\mathrm{GSpin}}(v)$ for $a\in O_\FLocal^\times$.
	
	Choose an auxiliary totally real field $\FGlobal_\aux\ne\Q$ with a place
	$v_\aux\mid p$ and an identification
	$(\FGlobal_\aux)_{v_\aux}\simeq\FLocal$, together with a global
	orthogonal datum as in \Cref{subsec:SO-shimura-datum} having the
	prescribed local quadratic space and self-dual lattices at every place above $p$.
	These choices use only the field and Gram-matrix approximation in the
	first two paragraphs of the proof of \Cref{prop:globalization}.
	Take $E_\aux=E_{0,\aux}\FGlobal_\aux$, where $E_{0,\aux}/\Q$ is imaginary quadratic
	and split at $p$, with the induced CM type of \Cref{defn:GSpin-Shimura-datum}.

	Let $\Vspecial_\aux$ be the space of
	special quasi-endomorphisms at a basic point, so
	$(\Vspecial_\aux)_{v_\aux}\simeq\Vspecial$ by \Cref{prop:special-orth-global}.
	Weak approximation permits us to approximate any vector of
	$\Vspecial$ by vectors of $\Vspecial_\aux$ which are integral at all the
	other places above $p$. For every such nonzero vector $v$, clearing
	denominators away from $p$ and applying
	\Cref{prop:tangent-space-divisor,prop:uniformization-special-morp}
	shows that $\CZ_{\mathrm{GSpin}}(v)$ is a relative effective Cartier divisor.
	The other local factors in the uniformization are discrete, and the
	integrality just imposed makes the chosen factors nonempty.
	If $(v,v)$ is a unit at $v_\aux$, then
	\Cref{prop:cycle-smooth-global} gives formal smoothness of relative
	dimension $\dimofV-3$ as well.
	
	Suppose first that $(x,x)\ne0$. Choose such a global vector $x'$
	sufficiently close to $x$ that
	\[
	\frac{(x',x')}{(x,x)}\in1+\varpi O_\FLocal.
	\]
	Since $p\ne2$, this ratio is $a^2$ for some $a\in O_\FLocal^\times$.
	Witt's extension theorem gives an isometry carrying $x'$ to $ax$.
	If its determinant is $-1$, compose it with a reflection in a
	nonisotropic vector of $(ax)^\perp$; thus we may choose the isometry
	in $\SO(\Vspecial)(\FLocal)$. It lifts to
	$\GSpin(\Vspecial)(\FLocal)$, since the kernel of
	$\GSpin\to\SO$ is $\Gm$ and $H^1(\FLocal,\Gm)=0$.
	The resulting action gives an isomorphism
	\[
	\CZ_{\mathrm{GSpin}}(x')\simeq\CZ_{\mathrm{GSpin}}(ax)=\CZ_{\mathrm{GSpin}}(x).
	\]
	This proves the Cartier and flatness assertions when $(x,x)\ne0$.
	It also proves the smoothness assertion, since $(x',x')$ is a unit
	whenever $(x,x)$ is a unit.
	
	We next prove local principality for an arbitrary nonzero $x$.
	Fix a geometric point $z\in\CZ_{\mathrm{GSpin}}(x)(\kb)$, let
	\[
	A=\widehat{\mathcal O}_{\GSpinRZ,z}
	\simeq\OFb[[t_1,\ldots,t_{\dimofV-2}]],\qquad
	\mathfrak m=(\varpi,t_1,\ldots,t_{\dimofV-2}),
	\]
	and let $I_x\subset A$ be the defining ideal of the completed cycle.
	For $n\ge1$ the scheme $U_n=\Spec(A/\mathfrak m^n)$ is a
	finite infinitesimal neighborhood of $z$ in $\GSpinRZ$.
	The denominator-clearing argument of
	\cite[Lemma~4.4]{Mihatsch2022LocalConstancy} gives an
	$O_\FLocal$-lattice $\Lambda_n\subset\Vspecial$ all of whose vectors
	extend over $U_n$. Indeed, transport an $\FLocal$-basis of
	$\Vspecial$ by the framing, lift the resulting quasi-endomorphisms
	over the nilpotent thickening $U_n$, and clear their finitely many
	denominators.
	Choose a global vector $x'_n$ as above with
	$x-x'_n\in\Lambda_n$. The two lifting conditions are then equivalent
	on every $U_n$-scheme. Since $\CZ_{\mathrm{GSpin}}(x'_n)$ has a local Cartier equation
	$f_n\in A$, we obtain
	\begin{equation}\label{eq:local-ideal-congruence}
		I_x+\mathfrak m^n=(f_n)+\mathfrak m^n.
	\end{equation}
	By Artin--Rees, choose $n$ sufficiently large that
	$I_x\cap\mathfrak m^n\subset\mathfrak m I_x$. Choose $h_n\in I_x$
	with $h_n\equiv f_n\pmod{\mathfrak m^n}$.
	Equation \eqref{eq:local-ideal-congruence} implies
	\[
	I_x=(h_n)+(I_x\cap\mathfrak m^n)
	=(h_n)+\mathfrak m I_x.
	\]
	Nakayama's lemma gives $I_x=(h_n)$. Faithful flatness of completed
	local rings therefore shows that the coherent ideal of $\CZ_{\mathrm{GSpin}}(x)$
	in $\GSpinRZ$ is locally principal. At this stage the generator is
	allowed to be zero.
	
	It remains to prove flatness and the nonzerodivisor assertion when
	$(x,x)=0$. 
	By smoothness it suffices to show that $\CZ_{\mathrm{GSpin}}(x)$ cannot contain a connected component of the special fiber of $\GSpinRZ$. 
	Assume that a connected component $\mathfrak C$ of the special fiber
	of $\GSpinRZ$ is contained in $\CZ_{\mathrm{GSpin}}(x)$.
	By \Cref{prop:SO-RZ-connected-components}, $\mathfrak C$ is a fiber
	of the Kottwitz map. An element of $\Spin(\Vspecial)(\FLocal)$ has
	trivial Kottwitz invariant and hence preserves $\mathfrak C$.

	Choose $y\in\Vspecial$ with
	$(y,y)=0$ and $(x,y)=1$; such a vector exists because $x$ is nonzero
	and isotropic. In the Clifford algebra $C(\Vspecial)$ one has
	$xy+yx=1$, and the element
	\[
	s=\varpi xy+\varpi^{-1}yx
	\]
	satisfies $\nu(s)=1$ and $sxs^{-1}=\varpi^2x$. In particular
	$s\in\Spin(\Vspecial)(\FLocal)$ and $s$ preserves $\mathfrak C$.
	Equivariance and $\mathfrak C\subseteq\CZ_{\mathrm{GSpin}}(x)$ give
	\[
	\mathfrak C_{\kb}\subseteq\CZ_{\mathrm{GSpin}}(\varpi^{-2n}x)
	\qquad\text{for every }n\ge0.
	\]
	Fix $z_0\in\mathfrak C(\kb)$. Inside the rational Dieudonn\'e module of
	the framing object, the endomorphisms preserving the Dieudonn\'e
	lattice at $z_0$ form an integral lattice. Its intersection
	\[
	L_{z_0}:=\{v\in\Vspecial:\ v\text{ extends to an endomorphism of }
	\mathscr X_{z_0}\}
	\]
	is an $O_\FLocal$-lattice in $\Vspecial$.
	The displayed inclusions would put $\varpi^{-2n}x$ in $L_{z_0}$
	for every $n$, contradicting $x\ne0$. 
\end{proof}

\subsection{The \texorpdfstring{$\dimofV=2$}{dimension-two} case}
\label{subsec:dim-two}

Assume that $\dimofV=2$.
Fix an embedding $\tau_0:\FLocal\hookrightarrow\ov\Q_p$ and an ordered
isotropic basis $e,f$ of $\Vlocal\otimes_{\FLocal,\tau_0}\ov\Q_p$
with $(e,f)=1$. Define $\{\mu\}$ by
$\mu^{-1}(t)=t^{-1}ef+fe$ on the $\tau_0$-component and by the trivial
cocharacter on every other component. Choose a basic class
$[b]\in B(\Res_{\FLocal/\Q_p}\GSpin(\Vlocal),\mu^{-1})$, and let
$b_\FLocal\in\GSpin(\Vlocal)(\breve\FLocal)$ be its relative representative
under Shapiro's identification.
We have the following analogue of \Cref{prop:V_b-SO}. 
\begin{prop}\label{prop:V_b-SO-2}
	Consider the $\FLocal$-vector space  $\Vspecial	:=	(\Vlocal\otimes_\FLocal\breve \FLocal)^{b_\FLocal \sigma_\FLocal=1}$ as in \Cref{subsec:V_b}. 
	If $\chi(\Vlocal)\neq 1$, then $\Vspecial$ is a quadratic space of dimension $\dimofV$ with Hasse invariant $-1$ and satisfies $\chi(\Vspecial)=\chi(\Vlocal)$. 
	 If $\chi(\Vlocal)= 1$, then $\Vspecial=0$. 
\end{prop}

\begin{proof}Since $\VlattiveLocaltot$ is
	self-dual, the discriminant of $\Vlocal$ is a unit, and hence
	$\chi(\Vlocal)\in\{\pm1\}$.
	
	Assume first that $\chi(\Vlocal)=1$. Then $\Vlocal$ is the split
	hyperbolic plane. Choose isotropic vectors $e,f$ with $(e,f)=1$. The
	idempotents $ef$ and $fe$ identify
	$C^+(\Vlocal)\simeq \FLocal\times \FLocal$ and
	$\GSpin(\Vlocal)\simeq\mathbb G_m\times\mathbb G_m$. Under this
	identification, $\mu^{-1}(t)=(t^{-1},1)$, and the basic class prescribed
	by $\mu^{-1}$ may, after $\sigma_\FLocal$-conjugation, be represented by
	$b_\FLocal=(\varpi^{-1},1)$.
	
	If $(a,c)\in\mathbb G_m\times\mathbb G_m$, then its action through the
	vector representation is given by
	$e\mapsto (a/c)e$ and $f\mapsto(c/a)f$. Since $\sigma_\FLocal(e)=e$ and
	$\sigma_\FLocal(f)=f$ in the split case, we have
	\(b_\FLocal \sigma_\FLocal(xe+yf)=\varpi^{-1}\sigma_\FLocal(x)e+\varpi\sigma_\FLocal(y)f.\)
	Thus a fixed vector satisfies
	$x=\varpi^{-1}\sigma_\FLocal(x)$ and $y=\varpi\sigma_\FLocal(y)$. Taking valuations
	shows that $x=y=0$. Hence $\Vspecial=0$.
	
	Finally, assume that $\chi(\Vlocal)=-1$. Let $K=\FLocal(\sqrt d)$, where
	$d\in O_\FLocal^\times$ is a nonsquare representing $\disc(\Vlocal)$. Then
	$K/\FLocal$ is the unramified quadratic extension,
	$C^+(\Vlocal)\simeq K$, and
	$\GSpin(\Vlocal)\simeq\Res_{K/\FLocal}\mathbb G_m$. After base change to
	$\breve \FLocal$, we identify
	$(\Res_{K/\FLocal}\mathbb G_m)_{\breve \FLocal}\simeq\mathbb G_m\times\mathbb G_m$
	so that $\sigma_\FLocal$ interchanges the two factors. We may choose isotropic
	vectors $e,f$ with $(e,f)=1$, $\sigma_\FLocal(e)=f$, and
	$\sigma_\FLocal(f)=e$. As above, after $\sigma_\FLocal$-conjugation we may take
	$b_\FLocal=(\varpi^{-1},1)$.
	
	It follows that
	\(b_\FLocal \sigma_\FLocal(xe+yf) =\varpi^{-1}\sigma_\FLocal(y)e+\varpi\sigma_\FLocal(x)f.\)
	Hence $b_\FLocal \sigma_\FLocal(xe+yf)=xe+yf$ is equivalent to
	$y=\varpi\sigma_\FLocal(x)$ and $x=\sigma_\FLocal^2(x)$. Since the fixed field of
	$\sigma_\FLocal^2$ in $\breve \FLocal$ is $K$, the map
	$x\mapsto xe+\varpi\tau(x)f$, where
	$\tau\in\Gal(K/\FLocal)$ is the nontrivial element, identifies $K$ with
	$\Vspecial$. Moreover, using $ef+fe=1$, the quadratic form is
	$Q(xe+\varpi\tau(x)f)=\varpi x\tau(x)=\varpi N_{K/\FLocal}(x)$. Therefore
	$\Vspecial\simeq(K,\varpi N_{K/\FLocal})$.
	
	Writing $K=\FLocal(\sqrt d)$, we obtain
	$\Vspecial\simeq\langle\varpi,-\varpi d\rangle$. Its discriminant is
	$-\varpi(-\varpi d)=\varpi^2d\equiv d\pmod{(\FLocal^\times)^2}$, so
	$\chi(\Vspecial)=-1=\chi(\Vlocal)$. Furthermore,
	$\Has(\Vspecial)=(\varpi,-\varpi d)_\FLocal=(\varpi,d)_\FLocal=-1$, because
	$K/\FLocal$ is unramified and every norm from $K^\times$ has even valuation.
	This proves the proposition.
\end{proof}
In the split case, $\Vspecial=0$ by \Cref{prop:V_b-SO-2}. 
For the rest of this subsection, assume that $\chi(\Vlocal)=-1$. 
Let $K/\FLocal$ be the unramified quadratic extension associated with
the common discriminant of $\Vlocal$ and $\Vspecial$, as computed in
\Cref{prop:V_b-SO-2}. Then
\[
C^+(\Vlocal)\simeq K,
\qquad
C^+(\Vspecial)\simeq K,
\qquad
\GSpin(\Vlocal)\simeq\Res_{K/\FLocal}\mathbb G_m,
\qquad \SO(\Vlocal)\simeq\Res^1_{K/\FLocal}\mathbb G_m.
\]
Their parahoric models are $\Res_{O_K/O_\FLocal}\mathbb G_m$ and $\Res^1_{O_K/O_\FLocal}\mathbb G_m$ respectively. 
\begin{rmk}\label{rmk:GSpin-RZ-rank-two-reflex}
	Let $\iota:K\hookrightarrow\breve\FLocal$ be the $\FLocal$-embedding
	determining the fixed cocharacter $\mu$, and write $\tau$ for the
	nontrivial element of $\Gal(K/\FLocal)$.  The two embeddings $\iota$
	and $\iota\circ\tau$ determine two distinct cocharacters
	$\mu_\iota$ and $\mu_{\iota\circ\tau}$.  Since the group is a torus,
	these are distinct geometric conjugacy classes, each with reflex
	field $K$.
\end{rmk}

For the fixed class $\{\mu\}=\{\mu_\iota\}$, let $\GSpinRZ$ be the
Hodge-type space of \Cref{prop:GSpin-RZ-defn}, and let
$b_{\mathrm{SO}}$ and $\{\mu_{\mathrm{SO}}\}$ be the images of $b$ and
$\{\mu\}$ under the central quotient. Define
\[
\SORZ:=\mathcal M_{\Res_{O_\FLocal/\Z_p}\SO(\VlattiveLocaltot),
 b_{\mathrm{SO}},\mu_{\mathrm{SO}}}/\Spf O_{\breve\FLocal}.
\]
The reflex field $K$ is embedded in $\breve\FLocal$ through $\iota$.
The central quotient induces $\pi_{\mathrm{RZ}}:\GSpinRZ\to\SORZ$.
For $\mathbf x\in\Vspecial^r$, let $\CZ_{\mathrm{SO}}(\mathbf x)$ be the
closed formal subscheme obtained by descending the GSpin lifting locus
along this map, as in \Cref{prop:abelian-local-cycles-independent-Hodge-cover}.
The other embedding $\iota\circ\tau$ gives a separate datum and the same
constructions for its cocharacter class.

\begin{lem}\label{lem:n2-RZ-component}
	The SO space is given by
	\[
	\SORZ\simeq\Spf O_{\breve\FLocal}.
	\]
	The map induced by central quotient $\GSpinRZ\to\SORZ$ restricts to an isomorphism
	on every connected component, and $\FLocal^\times/O_\FLocal^\times$ acts simply transitively on these
	components.
\end{lem}

\begin{proof}
	Since $\SO(\Vlocal)$ and $\GSpin(\Vlocal)$ are tori, the \RZname spaces are given by
	\[\SORZ\cong \underline{\SO(\Vlocal)(\FLocal)/\SO(\VlattiveLocaltot)(O_\FLocal)}\cong \Spf O_{\breve\FLocal}, \qquad \GSpinRZ\cong \underline{\GSpin(\Vlocal)(\FLocal)/\GSpin(\VlattiveLocaltot)(O_\FLocal)}\cong \underline{\FLocal^\times/O_\FLocal^\times}\]
	where both constant formal schemes are taken over $\Spf O_{\breve\FLocal}$,
by \cite[Proposition 25.2.1]{ScholzeWeinstein2020Berkeley}.
\end{proof}
Since $\GSpin(\VlattiveLocaltot)$ is the group of $O_K$-linear
automorphisms of $O_K$, $\GSpinRZ$ is the rank-one EL \RZname space for
$O_K$, with $\mu$ of weight one at $\iota$ and zero at the other
embeddings of $K$; see \cite[Corollary~25.1.3 and Proposition~18.4.1]
{ScholzeWeinstein2020Berkeley}.
Let $\mathbb X$ be the even factor of the Clifford framing object over
the residue field of $O_{\breve\FLocal}$. 
For an $O_{\breve\FLocal}$-scheme $S$ on which $\varpi$ is locally
nilpotent, the moduli problem parametrizes isomorphism classes of pairs
$(\mathscr X,\rho)$, where $\mathscr X$ is a one-dimensional formal $O_K$-module of
$O_K$-height one, $O_K$ acts on $\Lie \mathscr X$ through $\iota$, and
\[
\rho:\mathbb X_{S/\varpi S}\dashrightarrow \mathscr X_{S/\varpi S}
\]
is an $O_K$-linear quasi-isogeny.

On each connected component, the universal EL group $\mathscr X$ is
the canonical $O_K$-linear lift of $\mathbb X$ over $O_{\breve\FLocal}$;
see \cite[Theorem~2.1.2(iii) and Remark~2.1.3(iv)]{CarterStrauch2021}.
Choose $e\in\VlattiveLocaltot$ such that
\[
C(\VlattiveLocaltot)=O_K\oplus O_Ke,\qquad
e^2=1,\qquad ec=\tau(c)e\quad(c\in O_K).
\]
The universal Clifford $p$-divisible group, with its even and odd factors,
is
\[
\mathscr X_{\mathrm{Cl}}=\mathscr X^+\oplus\mathscr X^-
\simeq\mathscr X\oplus\mathscr X^\tau,
\]
where $\mathscr X^\tau$ has the same underlying formal
$O_\FLocal$-module as $\mathscr X$ and the $O_K$-action twisted by
$\tau$. The displayed $O_K$-actions are the right Clifford actions.
Right multiplication by $e$ exchanges the factors and gives an integral
isomorphism $j:\mathscr X^+\xrightarrow{\sim}\mathscr X^-$.
The group $\mathscr X_{\mathrm{Cl}}$ has dimension two and $O_\FLocal$-height four;
the central action changes only the framing.

The algebra $D:=\End_{O_\FLocal}(\mathbb X)_\Q$ is the quaternion division
algebra over $\FLocal$ of invariant $1/2$. Using the given $O_K$-action,
choose $\Pi\in D$ such that 
\[
D=K\oplus K\Pi,\qquad O_D=O_K\oplus O_K\Pi,\qquad
\Pi c=\tau(c)\Pi\ (c\in K),\qquad \Pi^2=\varpi.
\]
For a vector $y\in \Vspecial$, let $L_y: \mathbb X\oplus \mathbb X^\tau\to \mathbb X\oplus \mathbb X^\tau$ be its action
on the full Clifford framing object $\mathbb X\oplus\mathbb X^\tau$ and set $f_y=j^{-1}L_y|_{\mathbb X}\in D$.
Since $L_y$ commutes with the right actions of $e$ and $O_K$,
identifying the two factors by $j$ gives
\[
L_y=\begin{pmatrix}0&f_y\\ f_y&0\end{pmatrix},\qquad
f_y c=\tau(c)f_y\ (c\in K),\qquad
f_y^2=Q(y)\operatorname{id}=\tfrac12(y,y)\operatorname{id}.
\]
As $j$ is integral, $L_y$ lifts to $\mathscr X_{\mathrm{Cl}}$
if and only if $f_y$ lifts to $\mathscr X$. 
Thus the special cycle is exactly the lifting locus of $f_y$ on the
canonical deformation. Via $y\mapsto f_y$, we identify
\[
\Vspecial
=
\{y\in D:yc=\tau(c)y\text{ for every }c\in K\}
=K\Pi.
\]
Thus every $0\ne x\in\Vspecial$ is uniquely of the form $x=b\Pi$ with
$b\in K^\times$. The relation $\Pi^2=\varpi$ gives directly
\[
(x,x)=2\varpi N_{K/\FLocal}(b).
\]
Let $\val_K:K^\times\to\Z$ be the discrete valuation normalized by
$\val_K(\varpi)=1$.
Since $p\ne2$ and $K/\FLocal$ is unramified,
we conclude that $\valFLocal((x,x))=1+2\val_K(b).$ 
In particular, the valuation of the norm of every nonzero special
quasi-endomorphism is odd.

\begin{prop}
	\label{prop:n-two-special}
	Let $0\ne x=b\Pi\in\Vspecial$ and put
	\(r=\val_K(b)=\frac{\valFLocal((x,x))-1}{2}.\)
	Then the special cycle is given by
	\[
	\CZ_{\mathrm{SO}}(x)\simeq
	\begin{cases}
		\varnothing,& r<0,\\[4pt]
		\Spf\bigl(O_{\breve \FLocal}/\varpi^{r+1}\bigr),& r\geq0.
	\end{cases}
	\]
	Consequently, $\CZ_{\mathrm{SO}}(x)\ne\varnothing$ if and only if $\valFLocal((x,x))\geq1$, and in this case
	\[
	\length_{O_{\breve \FLocal}}O_{\CZ_{\mathrm{SO}}(x)}
	=r+1
	=\frac{\valFLocal((x,x))+1}{2}.
	\]
\end{prop}

\begin{proof}
	On each GSpin component $\Spf O_{\breve\FLocal}$, the even factor of
	the universal group is the canonical $O_K$-linear deformation
	$\mathscr X$ of $\mathbb X$, and special cycles are its lifting loci,
	as explained above. Using the framing
	to regard endomorphisms as elements of $D$, Gross's calculation
	\cite[Proposition~3.3]{Gross1986Canonical} gives, for $j\geq1$,
	\[
	\End_{O_\FLocal}(\mathscr X\bmod\varpi^j)
	=O_K+\varpi^{j-1}O_D
	=O_K\oplus\varpi^{j-1}O_K\Pi.
	\]
	Hence $x=b\Pi$ extends modulo $\varpi^j$ if and only if
	$b\in\varpi^{j-1}O_K$, equivalently $j\leq r+1$.
	
	The defining ideal of the special cycle is therefore contained in
	$(\varpi^j)$ exactly when $j\leq r+1$.  Since $O_{\breve\FLocal}$ is
	a discrete valuation ring, this ideal is the unit ideal for $r<0$ and
	$(\varpi^{r+1})$ for $r\geq0$.  This proves the description of the
	cycle.
\end{proof}

\subsection{\texorpdfstring{$\mathrm{SO}$}{SO} \RZname space and its components}
\label{subsec:SO-RZ-components}
Assume $\dimofV\geq3$ and retain the GSpin datum $(b,\{\mu\})$, the
self-dual lattice $\VlattiveLocaltot$, and the space $\Vspecial$ of
\Cref{subsec:GSpin-RZ}.  Write $\varpi$ for a uniformizer of $\FLocal$,
$k_\FLocal$ for its residue field, and $\bar k$ for an algebraic closure.

\begin{defn}\label{defn:SO-RZ-space}
	The central quotient
	\begin{equation*}
		1\longrightarrow\Gm\longrightarrow\GSpin(\Vlocal)
		\xrightarrow{\,\pi\,}\SO(\Vlocal)\longrightarrow1
	\end{equation*}
	induces, after restriction of scalars to $\Q_p$, the images
	$b_{\mathrm{SO}}$ and $\{\mu_{\mathrm{SO}}\}$ of $b$ and $\{\mu\}$.
	The associated abelian-type \RZname space is
	\[
	\SORZ:=
	\mathcal M_{\Res_{O_\FLocal/\Z_p}\SO(\VlattiveLocaltot),
	             b_{\mathrm{SO}},\mu_{\mathrm{SO}}}
	/\Spf O_{\breve\FLocal}.
	\]
	The morphism of integral local Shimura data induces
	$\pi_{\mathrm{RZ}}:\GSpinRZ\to\SORZ$, where $\GSpinRZ$ is the
	Hodge-type space of \Cref{defn:GSpin-pure-local-Shimura-datum}.
\end{defn}

\begin{defn}\label{defn:SO-Kottwitz-components}
	Consider the  Kottwitz maps 
	\[
	\kappa_{\mathrm{GSpin}}:|\GSpinRZ|\longrightarrow\Z,
	\qquad
	\kappa_{\mathrm{SO}}:|\SORZ|\longrightarrow\Z/2\Z.
	\]
	Normalize the first by
	$\kappa_{\mathrm{GSpin}}([g])=\valFLocal(\nu(g))$ on geometric
	special-fiber points, where $g\in\GSpin(\Vlocal)(\breve\FLocal)$ is the
	relative representative under \Cref{lem:ADLV-Shapiro}. Since $\nu(z)=z^2$ on the central $\Gm$,
	$\pi_{\mathrm{RZ}}$ induces
	\begin{equation}\label{eq:pi1-GSpin-SO-parity}
		\Z\longrightarrow\Z/2\Z,\qquad \ell\longmapsto\bar\ell.
	\end{equation}
	Define the open-and-closed formal subschemes
	\[
	\GSpinRZ^{(\ell)}:=\kappa_{\mathrm{GSpin}}^{-1}(\ell),
	\qquad
	\mathcal N_{\bar\kappa}:=\kappa_{\mathrm{SO}}^{-1}(\bar\kappa).
	\]
	Fix $\bar\kappa\in\Z/2\Z$ and set
	\begin{equation*}
		\SORZcomp:=\mathcal N_{\bar\kappa}\subset\SORZ.
	\end{equation*}
	Write $\SORZcomp^\red$ for its reduced underlying closed subscheme.
\end{defn}

\begin{prop}\label{prop:SO-RZ-connected-components}
	The connected-component decompositions are
	\begin{equation*}
		\GSpinRZ=\coprod_{\ell\in\Z}\GSpinRZ^{(\ell)},
		\qquad
		\SORZ=\mathcal N_{\bar0}\coprod\mathcal N_{\bar1}.
	\end{equation*}
	For every $\ell\in\Z$, the central quotient restricts to an isomorphism
	\begin{equation}\label{eq:GSpin-SO-component-isomorphism}
		\pi_\ell:\GSpinRZ^{(\ell)}
		\xrightarrow{\ \sim\ }\mathcal N_{\bar\ell}.
	\end{equation}
\end{prop}

\begin{proof}
	By \eqref{eq:dec-G-breve},
	\[
	\begin{aligned}
		\pi_1(\Res_{\FLocal/\Q_p}\GSpin(\Vlocal))_I
		&\simeq\Z^{[k_\FLocal:\mathbb F_p]},\\
		\pi_1(\Res_{\FLocal/\Q_p}\SO(\Vlocal))_I
		&\simeq(\Z/2\Z)^{[k_\FLocal:\mathbb F_p]}.
	\end{aligned}
	\]
	Frobenius cyclically permutes the factors, so the Kottwitz maps in
	\eqref{eq:Kottwitz-component-RZ} have images
	\[
	c_{b,\mu}+\{(\ell,\ldots,\ell):\ell\in\Z\},
	\qquad
	\pi_*(c_{b,\mu})
	+\{(\bar\ell,\ldots,\bar\ell):\bar\ell\in\Z/2\Z\}.
	\]
	Projection $(a_i)_i\mapsto a_0$ identifies these cosets with $\Z$
	and $\Z/2\Z$, and $\pi_*:(a_i)_i\mapsto(a_i\bmod2)_i$ induces
	$\ell\mapsto\bar\ell$. Thus
	\eqref{eq:PR-componentwise-ad-isomorphism} gives
	\eqref{eq:GSpin-SO-component-isomorphism}.

	By \Cref{prop:RZ-special-fiber-perfection,lem:ADLV-Shapiro},
	\[
	\SORZ^\perf\simeq
	X_{\SO(\VlattiveLocaltot)}
	\bigl(\pi(b_\FLocal),(\pi\circ\mu_\FLocal)^{-1}\bigr).
	\]
	The element $\pi(b_\FLocal)$ is basic, and
	$\pi\circ\mu_\FLocal$ has nontrivial image in every
	geometric simple factor of $\SO(\Vlocal)^\ad$, including both $A_1$
	factors when $\dimofV=4$. Hence
	\cite[Example~4.6 and Theorem~1.2]{GleasonLimXu2026Connected} and
	\Cref{lem:perfection-topology} give
	\[
	\kappa_{\mathrm{GSpin}}: \pi_0(\GSpinRZ)\simlr \Z,
	\qquad
	\kappa_{\mathrm{SO}}: \pi_0(\SORZ)\simlr \Z/2\Z.
	\]
\end{proof}

\begin{rmk}\label{rmk:different-on}
	The framing action of $g\in\GSpin(\Vspecial)(\FLocal)$ sends
	$\GSpinRZ^{(\ell)}$ to
	$\GSpinRZ^{(\ell+\valFLocal(\nu(g)))}$.
	The vanishing $H^1(\FLocal,\Spin(\Vspecial))=1$ makes $\nu$ surjective
	on $\FLocal$-points, so the two SO components are noncanonically
	isomorphic.  Such an isomorphism can move special cycles; the
	Frobenius comparison in \Cref{lem:SO-Frob-components} preserves them.
\end{rmk}

\begin{defn}\label{defn:SO-local-special-cycle}
	For $\mathbf x=(x_1,\ldots,x_r)\in\Vspecial^r$, let
	$\CZ_{\mathrm{SO}}(\mathbf x)\subset\SORZ$ be the closed formal
	subscheme obtained by descent from the GSpin special cycle as in
	\Cref{prop:abelian-local-cycles-independent-Hodge-cover}.
	The central kernel is $\Res_{\FLocal/\Q_p}\Gm$; its cocharacter
	lattice is a permutation lattice, so its inertia coinvariants are free,
	and its first Galois cohomology over $\Q_p$ vanishes by Hilbert~90.
	The cycle is characterized by
	\begin{equation*}
		\pi_\ell^{-1}\bigl(\CZ_{\mathrm{SO}}(\mathbf x)
		\cap\mathcal N_{\bar\ell}\bigr)
		=\CZ_{\mathrm{GSpin}}(\mathbf x)\cap\GSpinRZ^{(\ell)}
		\qquad(\ell\in\Z).
	\end{equation*}
	For an $O_\FLocal$-lattice $L\subset\Vspecial$ with basis
	$x_1,\ldots,x_r$, define
	\[
	\CZ_{\mathrm{SO}}(L):=\CZ_{\mathrm{SO}}(x_1,\ldots,x_r),
	\qquad
	\CZ_{\SORZcomp}(L):=\CZ_{\mathrm{SO}}(L)\times_{\SORZ}\SORZcomp.
	\]
	These cycles are independent of the basis: their pullbacks to
	$\GSpinRZ$ are the loci where every element of $L$ lifts to an
	endomorphism.
\end{defn}
\subsection{Bruhat-Tits stratification}
We generalize \cite[Definition~5.2.1]{HowardPappas2017Spinor} from the case $\FLocal=\Q_p$ to the general case.
\begin{defn}\label{defn:SO-special-lattice}
	Let $R$ be a perfect $\bark$-algebra. Write
	$q_\FLocal=\#k_\FLocal=p^{[k_\FLocal:\mathbb F_p]}$.
	The relative Frobenius on the coefficient ring is
	\[
	\sigma_\FLocal:=1\otimes W(\operatorname{Frob}_{q_\FLocal})
	\quad\text{on}\quad
	O_\FLocal\otimes_{W(k_\FLocal)}W(R),
	\]
	where $\operatorname{Frob}_{q_\FLocal}(r)=r^{q_\FLocal}$.
	Since this power map fixes $k_\FLocal$, its Witt-vector lift fixes
	$W(k_\FLocal)$, so the displayed tensor-product automorphism is
	well defined. It fixes $O_\FLocal$ and acts on $W(R)$ by the
	$[k_\FLocal:\mathbb F_p]$-th power of Witt-vector Frobenius.
	For $R=\bar k$, it agrees with the previously defined
	$\sigma_\FLocal$ on $O_{\breve\FLocal}$.
	Let $\Phi$ act on
	$\Vspecial\otimes_{W(k_\FLocal)[1/p]}W(R)[1/p]$
	by fixing $\Vspecial$ and applying this relative Frobenius to
	the coefficients. A \emph{special lattice over $R$}
	is an $O_\FLocal\otimes_{W(k_\FLocal)}W(R)$-submodule
	\[
	L\subset\Vspecial\otimes_{W(k_\FLocal)[1/p]}W(R)[1/p]
	\]
	that is finite locally free of rank $\dimofV$ over
	$O_\FLocal\otimes_{W(k_\FLocal)}W(R)$, satisfies
	$L[1/p]=\Vspecial\otimes_{W(k_\FLocal)[1/p]}W(R)[1/p]$, and has the properties
	\begin{equation*}
		L=L^\vee,\qquad
		\bigl(L+\Phi L\bigr)/L
		\text{ is locally free of rank one over }R.
	\end{equation*}
	Here $L^\vee$ is the dual for the extended $\FLocal$-valued pairing,
	with integrality measured in $O_\FLocal\otimes_{W(k_\FLocal)}W(R)$.
	The second condition includes the requirement that the action of
	$O_\FLocal\otimes_{W(k_\FLocal)}W(R)$ on the quotient factor through
	\[
	O_\FLocal\otimes_{W(k_\FLocal)}W(R)\longrightarrow R,
	\qquad a\otimes w\longmapsto\bar a\,(w\bmod p),
	\]
	using the fixed embedding $k_\FLocal\hookrightarrow\bark$.

	By \Cref{prop:RZ-special-fiber-perfection,lem:ADLV-Shapiro}, the
	perfect scheme $(\SORZ)^\perf$ is canonically identified with
	the \'etale sheafification of the presheaf
	\[
	R\longmapsto
	\left\{L\subset\Vspecial\otimes_{W(k_\FLocal)[1/p]}W(R)[1/p]:
	L\text{ is a special lattice over }R\right\}.
	\]
	For $y\in(\SORZ)^\perf(R)$, denote by $L_y$ the corresponding
	special lattice, with descent understood in the \'etale topology.
	Its Kottwitz invariant $\kappa_{\mathrm{SO}}(L_y)$ is a locally
	constant function on $\Spec R$ with values in $\Z/2\Z$.
	We also use $\kappa_{\mathrm{SO}}$ on all self-dual lattices in
	$\Vspecial\otimes_\FLocal\breve\FLocal$: transport the lattice to
	$\Vlocal\otimes_\FLocal\breve\FLocal$, choose an element of
	$\SO(\Vlocal)(\breve\FLocal)$ carrying
	$\VlattiveLocaltot\otimes_{O_\FLocal}O_{\breve\FLocal}$ to it, and take
	its group Kottwitz invariant. This is independent of the choice because
	the self-dual-lattice stabilizer has trivial Kottwitz invariant,
	and restricts to the preceding invariant on special lattices.
\end{defn}

\begin{defn}\label{defn:SO-vertex-variety}
	Let $\Lambda\subset\Vspecial$ be a vertex lattice, with the convention
	\[
	\Lambda\subseteq\Lambda^\vee\subseteq\varpi^{-1}\Lambda,
	\qquad t(\Lambda)=\dim_{k_\FLocal}(\Lambda^\vee/\Lambda).
	\]
	Define the finite quadratic space
	\[
	W_\Lambda:=\Lambda^\vee/\Lambda,
	\qquad \bar Q_\Lambda(\bar x):=\varpi Q(x)\bmod\varpi
	\quad(x\in\Lambda^\vee).
	\]
	Use $\Phi$ also for the induced relative Frobenius
	$1\otimes\operatorname{Frob}_{q_\FLocal}$ on
	$W_\Lambda\otimes_{k_\FLocal}\bar k$.
	Let
	$\operatorname{OGr}^{\epsilon_{\SORZcomp}}
	(t(\Lambda)/2,W_\Lambda\otimes_{k_\FLocal}\bar k)$
	be the component of the maximal isotropic Grassmannian consisting of
	subspaces whose inverse-image lattices in
	$\Lambda^\vee\otimes_{O_\FLocal}O_{\breve\FLocal}$ have SO Kottwitz
	invariant $\bar\kappa$.  Define the reduced closed subvariety
	\begin{equation}\label{eq:YLambda-closed}
		Y_\Lambda:=\left\{
		U\in\operatorname{OGr}^{\epsilon_{\SORZcomp}}
		(t(\Lambda)/2,W_\Lambda\otimes_{k_\FLocal}\bar k):
		\dim(U+\Phi U)=t(\Lambda)/2+1
		\right\}.
	\end{equation}
\end{defn}

Here $t(\Lambda)$ is positive and even: the Gram determinant of
$\Lambda$, computed in any $O_\FLocal$-basis, has even valuation,
and $\Has(\Vspecial)=-1$ excludes a self-dual lattice.
The quadratic space $W_\Lambda$ is nondegenerate and nonsplit, since
a rational Lagrangian would lift to a self-dual lattice in $\Vspecial$.
The Grassmannian component in the definition is determined by
\[
\kappa_{\mathrm{SO}}(L)-\kappa_{\mathrm{SO}}(L')
=\length_{O_{\breve\FLocal}}(L/(L\cap L'))\pmod2:
\]
for lattices between $\Lambda$ and $\Lambda^\vee$, this is precisely the
parity distinguishing the two maximal isotropic Grassmannian components.
Frobenius exchanges those components, so $U\ne\Phi U$ and
\[
\dim(U+\Phi U)\geq t(\Lambda)/2+1.
\]
Consequently the equality in \eqref{eq:YLambda-closed} is equivalent
to the determinantal closed condition
$\dim(U+\Phi U)\leq t(\Lambda)/2+1$.
By the finite-field calculation of
\cite[\S4.7]{LZ-orth} and
\cite[Proposition~5.3.2]{HowardPappas2017Spinor},
$Y_\Lambda$ is smooth, projective, and connected of dimension
$t(\Lambda)/2-1$.

\begin{defn}\label{defn:SO-BT-strata}
	For a vertex lattice $\Lambda\subset\Vspecial$, define its \emph{closed
	Bruhat--Tits stratum} in $\SORZcomp^\red$ by
	\[
	\BTstratum{\Lambda}:=\CZ_{\SORZcomp}(\Lambda)_{\bar k}^\red.
	\]
\end{defn}

\begin{prop}\label{lem:orthogonal-BT-dictionary}
	The Bruhat--Tits strata have the following properties.
	\begin{enumerate}
		\item \label{lem:orthogonal-BT-dictionary:1}For $y\in\SORZcomp(\bar k)$,
		$
		y\in\BTstratum{\Lambda}(\bar k)
		$ if and only if $		\Lambda_{O_{\breve\FLocal}}\subset L_y
		\subset\Lambda^\vee_{O_{\breve\FLocal}}.
		$
		\item \label{lem:orthogonal-BT-dictionary:2}The closed strata cover $\SORZcomp^\red$. 
		\item \label{lem:orthogonal-BT-dictionary:3}There are canonical isomorphisms of perfect $\bar k$-schemes
		$
		\BTstratum{\Lambda}^\perf\simeq Y_\Lambda^\perf$.  
		In particular, $\CZ_{\SORZcomp}(\Lambda)_{\bar k}^\red$ is proper
		over $\bar k$.
	\end{enumerate}
\label{prop:irreducible-components-SO}
\end{prop}

\begin{proof}
	By \Cref{prop:special-ADLV}, the vertex cycle is the locus where
	$\Lambda\subset L_y$. Since $L_y$ is self-dual,  \eqref{lem:orthogonal-BT-dictionary:1} follows. In this case the quotient map
	\[
	L_y\longmapsto L_y/\Lambda_{O_{\breve\FLocal}}
	\subset W_\Lambda\otimes_{k_\FLocal}\bar k
	\]
	identifies self-duality with maximal isotropy and the length-one
	condition with $\dim(U+\Phi U)=t(\Lambda)/2+1$.
	Its inverse takes the inverse image of $U$ in
	$\Lambda^\vee_{O_{\breve\FLocal}}$. 
	It therefore gives $\BTstratum{\Lambda}^\perf\simeq Y_\Lambda^\perf$. Note that  $\BTstratum{\Lambda}$ is separated and locally of finite
	type.  The projectivity of $Y_\Lambda$ and
	\Cref{lem:perfection-topology} imply that it is quasi-compact, hence
	of finite type.  Since $Y_\Lambda^\perf\to Y_\Lambda$ is integral,
	$\BTstratum{\Lambda}^\perf$ is universally closed over $\bar k$.
	Universal closedness descends along the universal homeomorphism
	$\BTstratum{\Lambda}^\perf\to\BTstratum{\Lambda}$, proving properness. Thus \eqref{lem:orthogonal-BT-dictionary:3} follows.

	For \eqref{lem:orthogonal-BT-dictionary:2}, we spell out the part of
	\cite[Proposition~5.2.2]{HowardPappas2017Spinor} needed for the covering.
	Let $L$ be a special lattice over $\bar k$, with $\Phi$ as in
	\Cref{defn:SO-special-lattice}.
		The increasing sums $\sum_{i=0}^j\Phi^iL$ lie in the scalar extension
	of a sufficiently large $O_\FLocal$-lattice in $\Vspecial$, so they
	stabilize.  Let $d\geq1$ be the first index of stabilization, i.e., $\sum_{i=0}^d\Phi^iL=\sum_{i=0}^{d+1}\Phi^iL$. 
	Then $\sum_{i=0}^d\Phi^iL$ is $\Phi$-stable.  Each successive quotient
	is a quotient of a Frobenius translate of $(L+\Phi L)/L$, so every
	strict inclusion has colength one.  
	
	For every $1\leq j<d$, the inclusion $\sum_{i=1}^{j}\Phi^iL\subsetneq \sum_{i=1}^{j+1}\Phi^iL$ has colength one. 
	Thus the intersection $\left(\sum_{i=0}^j\Phi^iL\right)
	\cap\left(\sum_{i=1}^{j+1}\Phi^iL\right)$ is either
	$\sum_{i=1}^j\Phi^iL$ or $\sum_{i=1}^{j+1}\Phi^iL$.
	In the second case,
	$\Phi\bigl(\sum_{i=0}^j\Phi^iL\bigr)\subseteq\sum_{i=0}^j\Phi^iL$;
	the same colength comparison makes this an equality, contradicting
	the minimality of $d$. Thus $\left(\sum_{i=0}^j\Phi^iL\right)
	\cap\left(\sum_{i=1}^{j+1}\Phi^iL\right)=\sum_{i=1}^j\Phi^iL$. 
	
	Now we use descending induction to show that \[
	\varpi\sum_{i=0}^d\Phi^iL\subseteq\sum_{i=0}^{j}\Phi^iL
	\]
	for every $0\leq j<d$. The quotient
	$(\sum_{i=0}^d\Phi^iL)/(\sum_{i=0}^{d-1}\Phi^iL)$ has length one,
	so multiplication by $\varpi$ kills it; this proves the case $j=d-1$.
	Suppose the displayed inclusion holds for $1\leq j<d$.
	Since $\varpi\sum_{i=0}^d\Phi^iL$ is $\Phi$-stable,
	applying $\Phi$ also gives its containment in
	$\sum_{i=1}^{j+1}\Phi^iL$. Thus the
	preceding intersection identity gives
	\[
	\varpi\sum_{i=0}^d\Phi^iL
	\subseteq\left(\sum_{i=0}^j\Phi^iL\right)
	\cap\left(\sum_{i=1}^{j+1}\Phi^iL\right)
	=\sum_{i=1}^j\Phi^iL.
	\]
	Applying $\Phi^{-1}$ leaves the left-hand side unchanged and yields
	$\varpi\sum_{i=0}^d\Phi^iL\subseteq\sum_{i=0}^{j-1}\Phi^iL$.
	
	We conclude that $\varpi\sum_{i=0}^d\Phi^iL\subseteq L$. 
	Applying powers of $\Phi$ and using $L=L^\vee$, we obtain
	\[
	\varpi\sum_{i=0}^d\Phi^iL
	\subseteq\bigcap_{i=0}^d\Phi^iL
	=\left(\sum_{i=0}^d\Phi^iL\right)^\vee
	\subseteq L\subseteq\sum_{i=0}^d\Phi^iL.
	\] 
	Both $\sum_{i=0}^d\Phi^iL$ and its dual are $\Phi$-stable, and clearly $\left(\bigcap_{i=0}^d\Phi^iL\right)^{\Phi=1}
	=L\cap\Vspecial$. 
	Frobenius descent and duality give
	\[
	\bigcap_{i=0}^d\Phi^iL=(L\cap\Vspecial)_{O_{\breve\FLocal}},
	\qquad
	\sum_{i=0}^d\Phi^iL=(L\cap\Vspecial)^\vee_{O_{\breve\FLocal}}.
	\]
	The preceding inclusions descend to
	$\varpi(L\cap\Vspecial)^\vee\subseteq L\cap\Vspecial
	\subseteq(L\cap\Vspecial)^\vee$.
	Thus $L\cap\Vspecial$ satisfies the stated vertex convention, and
	\[
	(L\cap\Vspecial)_{O_{\breve\FLocal}}\subseteq L
	\subseteq(L\cap\Vspecial)^\vee_{O_{\breve\FLocal}}.
	\]
	By \eqref{lem:orthogonal-BT-dictionary:1}, the corresponding point
	lies in $\BTstratum{L\cap\Vspecial}$.  The same argument over every
	algebraically closed extension of $\bar k$ proves the covering. So \eqref{lem:orthogonal-BT-dictionary:2} follows. 
\end{proof}

\begin{cor}\label{cor:cycle-noetherian}
	For every nondegenerate lattice $L\subset\Vspecial$ of rank
	$\dimofV-1$, the cycle $\CZ_{\mathrm{SO}}(L)\cap\SORZcomp$ is Noetherian.
\end{cor}

\begin{proof}
	The cycle is locally noetherian.  By the Bruhat--Tits covering, each
	geometric point lies on a vertex cycle for some $\Lambda_0$.
	At that point $L+\Lambda_0$ is integral, and
	\[
	L+\Lambda_0\subseteq(L+\Lambda_0)^\vee
	\subseteq\Lambda_0^\vee\subseteq\varpi^{-1}(L+\Lambda_0).
	\]
	Thus $L+\Lambda_0$ is a vertex lattice containing $L$.  Consequently
	$\CZ_{\mathrm{SO}}(L)\cap\SORZcomp$ is covered, on underlying
	topological spaces, by vertex cycles with $\Lambda\supseteq L$.

	There are only finitely many such vertex lattices, by the lattice
	bound of \cite[Lemma~4.12.2(i)]{LZ-orth}.  Indeed, the projection of
	$\Lambda$ to $L_\FLocal$ lies in $L^\vee$, and integrality bounds its
	projection to the anisotropic line $L_\FLocal^\perp$ in a fixed lattice.
	This gives a common upper bound for $\Lambda$; dualizing and using
	$\varpi\Lambda^\vee\subseteq\Lambda$ gives a common full-rank lower
	bound.  There are finitely many lattices between these bounds.
	Their vertex cycles are quasi-compact by
	\Cref{prop:irreducible-components-SO}, so the original cycle is
	quasi-compact and hence Noetherian.
\end{proof}

\begin{lem}\label{lem:SO-Frob-components}
	The Weil descent datum (cf. \Cref{defn:Weil-descent}) acts on components by
	\begin{equation}\label{eq:Frob-on-GSpin-SO-components}
		\begin{aligned}
		\Frob\bigl(\GSpinRZ^{(\ell)}\bigr)
		&=\bigl(\GSpinRZ^{(\ell+1)}\bigr)^{(\sigma_\FLocal)},\\
		\Frob\bigl(\mathcal N_{\bar\kappa}\bigr)
		&=\bigl(\mathcal N_{\overline{\kappa+1}}\bigr)^{(\sigma_\FLocal)}.
		\end{aligned}
	\end{equation}
	For every $O_\FLocal$-lattice $L\subset\Vspecial$, it restricts to
	\begin{equation}\label{eq:Frob-on-SO-special-cycles}
		\Frob:\CZ_{\mathrm{SO}}(L)\cap\mathcal N_{\bar\kappa}
		\xrightarrow{\ \sim\ }
		\bigl(\CZ_{\mathrm{SO}}(L)\cap
		\mathcal N_{\overline{\kappa+1}}\bigr)^{(\sigma_\FLocal)}.
	\end{equation}
\end{lem}

\begin{proof}
	By \Cref{prop:RZ-special-fiber-perfection,lem:ADLV-Shapiro}, we may
	check \eqref{eq:Frob-on-GSpin-SO-components} on
	$X_{\GSpin(\VlattiveLocaltot)}(b_\FLocal,\mu_\FLocal^{-1})$.
	Here $b_\FLocal\in\GSpin(\Vlocal)(\breve\FLocal)$ is determined on
	the distinguished factor of \eqref{eq:dec-G-breve} by
	\[
	(b\sigma)^{[k_\FLocal:\mathbb F_p]}
	=b_\FLocal\sigma_\FLocal,
	\]
	and $\mu_\FLocal$ is the cocharacter class on that factor.
	Since $\GSpin(\VlattiveLocaltot)$ and $\{\mu_\FLocal\}$ are defined
	over $O_\FLocal$ and $\FLocal$, respectively, the
	$\sigma_\FLocal$-twist of this variety is
	$X_{\GSpin(\VlattiveLocaltot)}
	(\sigma_\FLocal(b_\FLocal),\mu_\FLocal^{-1})$.
	The change of framing in \Cref{defn:Weil-descent} therefore induces
	\[
	\begin{gathered}
	\Frob:
	X_{\GSpin(\VlattiveLocaltot)}(b_\FLocal,\mu_\FLocal^{-1})
	\xrightarrow{\ \sim\ }
	X_{\GSpin(\VlattiveLocaltot)}
	(\sigma_\FLocal(b_\FLocal),\mu_\FLocal^{-1}),\\
	[g]\longmapsto[b_\FLocal^{-1}g].
	\end{gathered}
	\]
	This preserves the admissible relative-position condition because
	\[
	(b_\FLocal^{-1}g)^{-1}\sigma_\FLocal(b_\FLocal)
	\sigma_\FLocal(b_\FLocal^{-1}g)
	=g^{-1}b_\FLocal\sigma_\FLocal(g).
	\]
	By \eqref{eq:GSpin-standard-cocharacter} and
	$[b]\in B(\Res_{\FLocal/\Q_p}\GSpin(\Vlocal),\mu^{-1})$, one has
	$\valFLocal(\nu(b_\FLocal))=-1$.  Hence
	\[
	\kappa_{\mathrm{GSpin}}([b_\FLocal^{-1}g])
	=\valFLocal(\nu(g))+1.
	\]
	The central quotient is compatible with this change of framing,
	and its map on component labels is reduction modulo $2$ by
	\eqref{eq:pi1-GSpin-SO-parity}.  This proves
	\eqref{eq:Frob-on-GSpin-SO-components}.

	For every $x\in\Vspecial$, the identity
	$\rho_V(b_\FLocal)\sigma_\FLocal(x)=x$ implies
	\[
	\rho_V(b_\FLocal^{-1}g)^{-1}\sigma_\FLocal(x)
	=\rho_V(g)^{-1}x.
	\]
	Thus the change of framing identifies the endomorphism lifting
	condition for $x$ with its $\sigma_\FLocal$-twist.
	Applying this to every $x\in L$ and descending through the central
	quotient proves \eqref{eq:Frob-on-SO-special-cycles}.
\end{proof}

\section{Intersection number}\label{sec:Int-number}
\subsection{Derived special cycles and intersection numbers}
We retain the notation of \Cref{sec:GSpin}, especially \Cref{subsec:orth-global-models}. 

\begin{defn}\label{defn:derived-special-cycles-Sh}
	Recall that each connected component of $\sEndint_2$ is finite and unramified over $\SOIntModel$ by \Cref{prop:global-special-endomorphism-model}.
	We define the derived structure sheaf $\CO^{\Lder}$ on the $r$-fold derived fiber product
	\[\sEndint_2\times_{\SOIntModel}^{\Lder} \cdots \times_{\SOIntModel}^{\Lder} \sEndint_2,\]
	regarded as a bounded complex on $\sEndint_2^r$, hence an element of the $K$-group $K'_0(\sEndint_2^r)$  (see \Cref{subsec:nota-K-gp}). Here we define $K'_0(\sEndint_2^r)$ to be the direct product of the $K'_0$-groups of the connected components of $\sEndint_2^r$, which are of finite type by \Cref{prop:global-special-endomorphism-model}. 
	Let $U\subseteq \sEndint_2^r$ be the open and closed subscheme on which the Gram matrix of the tautological tuple is nondegenerate.
	By \Cref{prop:tangent-space-divisor} and \cite[Proposition 23.58(1)]{GoertzWedhorn2023}, $[\CO^{\Lder}|_U]$ lies in $\Fil_{\dimofV-1-r} K'_0(\sEndint_2^r)_\Q$.  
	For each nondegenerate $T\in \Sym_{r}(\FGlobal)$ and each $x\in \Vadl_{\fkd}^{r}$, we define 
	$\CZ^\Lder_T(x)$ to be the pushforward of $\CO^\Lder|_{\CZ_T(x)}$ to $\SOIntModel$, regarded as an element of $\Fil_{\dimofV-1-r} K'_0(\im(\CZ_T(x)\to \SOIntModel))_\Q\cong \Fil^{r} K_0^{\CZ_T(x)}(\SOIntModel)_\Q$. 	
	More generally, for each nondegenerate $T\in \Sym_{r}(\FGlobal)$ and each $\ov{\Q}$-valued \emph{admissible} Schwartz function, meaning a $K_2$-invariant function of the form $\phi=\phi_{\fkd}\charfun_{(\VlattiveGlobal\otimes \widehat{\Z}^\fkd)^{r}}$ on $(\Vadl\otimes_\FGlobal \BA_{\FGlobal,f})^{r}$, we define
	$$\CZ^\Lder_T(\phi)=\sum_{x\in \Vadl_{\fkd}^{r}/K_{2,\fkd}, T(x)=T}\phi_{\fkd}(x) \cdot \CZ^\Lder_T(x)\in \Fil^{r} K_0^{\CZ_T(\supp \phi)}(\SOIntModel)_{\ov{\Q}}.$$
	This is a finite sum and the scheme $\CZ_T(\supp \phi)$ is finite over $\Sh_{K_2}$ by \Cref{lem:finite-global-gluing-block}. For each $v\nmid\fkd$ and the chosen reflex place $\wt v\mid v$, base change to $O_{\reflexField,(\wt v)}$ and define
	$\CZ^{\sharp,\Lder}_T(\phi)\in \Fil^{r} K_0^{\CZ^{\sharp}_T(\phi)}(\Sh_{K^\sharp})_{\ov{\Q}}$ by pullback along the Hodge-type cover $\Sh_{K^\sharp}\to \Sh_{K_2}\otimes_{O_{\FGlobal,(v)}}O_{\reflexField,(\wt v)}$.
	
\end{defn}

\begin{defn}\label{defn:derived-special-cycles-RZ}
	For a tuple $\mathbf x=(x_1,\ldots,x_r)$ of special
	quasi-endomorphisms, consider
	\(\CZ_{\mathrm{SO}}(x_i)\subset\SORZ\).
	For $\dimofV\geq3$, this is the descended cycle of
	\Cref{defn:SO-local-special-cycle}; in rank two it is the cycle on the
	 SO space for the chosen cocharacter, as in \Cref{subsec:dim-two}.
	Define the 
	derived cycle by
	\[
	\CZ_{\SORZ}^{\Lder}(\mathbf x)
	:=
	\mathcal O_{\CZ_{\mathrm{SO}}(x_1)}
	\otimes^{\Lder}_{\mathcal O_{\SORZ}}\cdots
	\otimes^{\Lder}_{\mathcal O_{\SORZ}}
	\mathcal O_{\CZ_{\mathrm{SO}}(x_r)}.
	\]
	Its pullback to $\GSpinRZ$ is isomorphic to the derived tensor complex
	$\CZ^{\Lder}_{\GSpin}(\mathbf x)$ by construction. 
	
	For a nonsingular tuple of length $\dimofV-1$, define
	\[
	\Int^\epsilon(\mathbf x)
	:=
	\begin{cases}
	\,\chi\!\left(R\Gamma\!\left(\SORZ,
	\CZ_{\SORZ}^{\Lder}(\mathbf x)\right)\right),&\dimofV=2,\\
	\tfrac12 \,\chi\!\left(R\Gamma\!\left(\SORZ,
	\CZ_{\SORZ}^{\Lder}(\mathbf x)\right)\right),&\dimofV\ge3.
	\end{cases}
	\]
	In the nonsplit two-dimensional case, $\SORZ$ is connected for each chosen cocharacter;
	for $\dimofV\ge3$, $\SORZ$ has two connected components with equal derived
	degrees by \Cref{lem:SO-Frob-components}. The properness and finiteness
	needed for these Euler characteristics are proved in \Cref{cor:int-finite}.
	
	Assume first that $\dimofV\geq3$, and let $\bx'$ have the same Gram
	matrix as $\bx$.  By Witt's theorem there is an element
	\(g\in J_b(\Q_p)\cong\GSpin(\Vspecial)(\FLocal)\)
	with \(g\bx=\bx'\).  We conclude
	\[
	\Int^\epsilon(\mathbf x)=\Int^\epsilon(\mathbf x').
	\]
	For $\dimofV=2$, Gram invariance follows directly from
	\Cref{prop:n-two-special}, where the branchwise length depends only on the
	valuation of $(x,x)$.  Thus $\Int^\epsilon$ may be regarded throughout as
	a function of the Gram matrix.
\end{defn}

\begin{thm}\label{thm:int-in-basic}
	For any nonsingular $T\in\Sym_{\dimofV-1}(\FGlobal)$ and admissible Schwartz function $\phi\in \sS((\Vadl\otimes_\FGlobal\BA_{\FGlobal,f})^{\dimofV-1})$, $\CZ_T(\phi)$ has empty generic fiber and its fiber at every $v\nmid\fkd$ is supported on the basic locus of $\sS_{K_2,(v)}$.
\end{thm}

\begin{proof}
	
At the distinguished real place, \(\Vcoherent\) has signature
	\((\dimofV-2,2)\).  At a complex point of a special cycle, the
	realizations of special endomorphisms span a positive-definite
	subspace for the special-endomorphism quadratic form.  Hence a
	nonsingular tuple of length \(\dimofV-1\) would give a
	positive-definite \((\dimofV-1)\)-dimensional subspace of
	\(\Vcoherent_{\RealPlace_0}\), which is impossible because its positive
	index is only \(\dimofV-2\).  Thus \(\CZ_T\) has empty generic fiber.
	
	It remains to prove basicness in characteristic \(p\).  Let \(k\) be a
	finite field and let \((x,f)\in\CZ_T(k)\), with
	\(x\in\SOIntModel(k)\) and
	\(f\in\sEndint_{2,x}^{\dimofV-1}(k)\).  
	After extending $k$, choose a Hodge-type lift of $(x,f)$.
	Its crystalline structure and the orthogonal quotient define
	$b_x\in B(G_2)$. By \Cref{prop:cris-tensor}, the vector isocrystal
	\(\CV_x\) has a
	nondegenerate symmetric pairing with values in
	\((\FGlobal\otimes_\Q W(\barFp)[1/p],1\otimes\sigma)\).
	The crystalline realizations of the special
	endomorphisms span a \(\phi\)-fixed nondegenerate subisocrystal of rank
	\(\dimofV-1\) over \(\FGlobal\otimes_\Q W(\barFp)[1/p]\),
	and hence all its Newton slopes are zero.
	Its orthogonal complement \(M\) has rank one over \(\FGlobal\otimes_\Q W(\barFp)[1/p]\)  and inherits a nondegenerate Frobenius-compatible pairing.
	For every $v\mid p$, the factor $M_v$ is a rank-one
	$\FGlobal_v$-isocrystal. Its pairing gives $M_v\simeq M_v^\vee$.
	If $\lambda_v$ is its Newton slope, duality gives
	$\lambda_v=-\lambda_v$, hence $\lambda_v=0$.
	Thus all slopes of $M$, and therefore of $\CV_x$, vanish.
	The standard representation of $G_2$ is faithful, so the Newton
	cocharacter of $b_x\in B(G_2)$ is trivial. Hence $x$ is basic. 
\end{proof}

\begin{cor}\label{cor:int-finite}   
	Let $L\subset\Vspecial$ be a nondegenerate $O_\FLocal$-lattice of
	rank $\dimofV-1$. Then
	$\CZ_{\mathrm{SO}}(L)$ is a proper scheme over
	$\Spec O_{\breve\FLocal}$ with empty generic fiber.
	In particular, $\Int^\epsilon(\mathbf x)$ is finite for every
	$O_\FLocal$-basis $\mathbf x$ of $L$.
\end{cor}
\begin{proof}
	If $\dimofV=2$, the description in \Cref{prop:n-two-special}
	shows directly that $\CZ_{\mathrm{SO}}(L)$ is a finite scheme with empty
	generic fiber. We may therefore assume
	$\dimofV\geq3$. By \Cref{prop:SO-RZ-connected-components}, it suffices to work with
	$\CZ_{\mathrm{GSpin}}(L)\cap\GSpinRZ^{(\ell)}$ for each $\ell$.
	This formal scheme is Noetherian by \Cref{cor:cycle-noetherian}.
	We first rule out its horizontal parts.
	
	 Suppose such a
	horizontal part exists.  Since $\GSpinRZ$ is locally formally of finite
	type over $\Spf O_{\breve\FLocal}$, the non-archimedean Nullstellensatz
	shows that its nonempty rigid generic fiber contains a classical point
	whose residue field is finite over $\breve\FLocal$.
	Thus for a finite extension
	$K/\breve\FLocal$ there is a point
	$z\in\bigl(\CZ_{\mathrm{GSpin}}(\Lattice)\cap
	\GSpinRZ^{(\ell)}\bigr)(O_K)$
	on the horizontal part. Let $\mathscr X_z/O_K$ be the corresponding
	$p$-divisible group, and let $\bar\eta$ be a geometric generic point of
	$\Spec O_K$. By \cite[\S2.5.1 and Proposition~2.5.3]
	{PappasRapoport2024PadicShtukas}, the vector representation of the
	generic-fiber shtuka gives a quadratic $O_\FLocal$-lattice
	\[
	V_{\et,z}\subset\End_{\Z_p}(T_p\mathscr X_{z,\bar\eta}).
	\]
	Trivializing the associated pro-\'{e}tale torsor at $\bar\eta$ identifies
	$V_{\et,z}$ with $\VlattiveLocaltot$, hence
	$V_{\et,z}[1/p]\simeq\Vlocal$ as quadratic $\FLocal$-spaces.
	
	Choose an $O_\FLocal$-basis $x_1,\ldots,x_{\dimofV-1}$ of $\Lattice$.
	The lifting condition at $z$ extends the quasi-endomorphism corresponding
	to each $x_i$ under the framing quasi-isogeny to an endomorphism
	$x_{i,z}$ of $\mathscr X_z$.
	By \Cref{prop:special-perfectoid}, its Tate realization lies in
	$V_{\et,z}$; here we use the compatibility of the \'{e}tale realization
	of the shtuka with the Tate module, as in
	\cite[Proposition~22.6.1]{ScholzeWeinstein2020Berkeley}.
	This construction is $O_\FLocal$-linear and preserves the quadratic
	pairing. Since $\Lattice$ is nondegenerate, it gives an isometric embedding
	\begin{equation}\label{eq:Tate-embedding-special-subspace}
		\Lattice\otimes_{O_\FLocal}\FLocal\lhook\joinrel\longrightarrow
		V_{\et,z}[1/p]\simeq\Vlocal,
		\qquad x_i\longmapsto T_p(x_{i,z}).
	\end{equation}
	
	Such an embedding is impossible.  The spaces $\Vspecial$ and $\Vlocal$
	have dimension $\dimofV$ and the same discriminant, but opposite Hasse
	invariants: $\Has(\Vspecial)=-1$, whereas
	$\Has(\Vlocal)=+1$ because $\VlattiveLocaltot$ is self-dual.  Since
	$\Lattice$ is nondegenerate of codimension one, its orthogonal complements
	in $\Vspecial$ and in $\Vlocal$ are one-dimensional.  Equality of the
	two ambient discriminants forces these one-dimensional complements to
	be isometric.  Together with
	\eqref{eq:Tate-embedding-special-subspace}, this would give an isometry
	$\Vspecial\simeq\Vlocal$, contradicting their opposite Hasse
	invariants.  

	The formal-curve argument of \cite[Proposition~6.1.1]{LZ-orth}
	now shows that the ideal of the reduced RZ locus acts nilpotently
	on this cycle. Indeed, otherwise that argument produces an affine
	formal curve whose closed point lies in the reduced RZ locus and whose
	algebraic generic point lies outside that locus.
	The algebraic generic point of the curve carries the same nondegenerate
	rank-$(\dimofV-1)$ lattice of special endomorphisms. The slope argument
	of \Cref{thm:int-in-basic} makes this generic point basic, contradicting
	the choice of the formal curve. Thus the cycle is a scheme.
	Its reduced subscheme is contained in finitely many proper vertex strata
	by \Cref{cor:cycle-noetherian,prop:irreducible-components-SO}.
	It is therefore proper, as in \cite[Lemma~4.13.1]{LZ-orth}.
\end{proof}
\begin{thm}[{\cite{Mihatsch2022LocalConstancy}}]\label{thm:local-constancy}
	On the open subset $\{\bu\in\Vspecial^\nInLZorth:
	\det T(\bu)\ne0\}$
	the arithmetic intersection number \(\Int^\epsilon(\bu)\) is locally
	constant.
\end{thm}

\begin{proof}
	For $\dimofV=2$ this follows from \Cref{prop:n-two-special}, so we assume $\dimofV\geq 3$. 
	We spell out the application of the local constancy theorem of
	\cite{Mihatsch2022LocalConstancy}.  
	
	Let $S$ be a profinite set with a continuous map $u:S\to\Vspecial\setminus\{0\}$. We define a closed formal subscheme $\CZ(u)\subseteq (\SORZ)_S$ as in \cite{Mihatsch2022LocalConstancy}*{\S4(8)}. Each fiber is a Cartier divisor by \Cref{prop:local-special-cycle-properties}.
	For any $O_\FLocal$-lattice $\Lambda\subset\Vspecial$ and any
affine scheme $U$ with a map $U\to\SORZ$, this map factors through
$\CZ_{\mathrm{SO}}(p^n\Lambda)$ for sufficiently large $n$; see
\cite[Lemma~4.4]{Mihatsch2022LocalConstancy}. 
	By the arguments in \cite[Proposition~4.2]{Mihatsch2022LocalConstancy}, $\CZ(u)$ is a Cartier divisor on $(\SORZ)_S$ flat over $S$. 
	
	Now take continuous maps $u_1,\ldots,u_{\nInLZorth}:S\to\Vspecial$
	such that $T(u_1(s),\ldots,u_{\nInLZorth}(s))$ is nonsingular for
	every $s\in S$. We must show that the intersection-number function
	$s\in S\mapsto\Int^\epsilon(u_1(s),\ldots,u_{\nInLZorth}(s))$ is locally constant.
	By \cite{Mihatsch2022LocalConstancy}*{Theorem 3.4}, it remains to check \cite{Mihatsch2022LocalConstancy}*{\S1(Z3)}: each $\CZ_{\mathrm{SO}}(u_1(s),\ldots,u_{\nInLZorth}(s))$ is a proper scheme, and these subschemes of $\SORZ$ are locally constant in $s$.
	
	The first assertion follows from \Cref{cor:int-finite}. 
	By the arguments in \cite{Mihatsch2022LocalConstancy}*{Lemma 4.5}, it suffices to check that $\bigcap_{i=1}^{\nInLZorth}\CZ(u_i)\subset(\SORZ)_S$ is quasi-compact.
	This follows from  \Cref{lem:vertex-lattices-locally-constant} below and compactness of $S$. 
\end{proof}
\begin{lem}	\label{lem:vertex-lattices-locally-constant}
	Let $S$ be a profinite set and let
	$u_i:S\to\Vspecial$, $1\le i\le\dimofV-1$, be continuous maps such
	that $T(u_1(s),\ldots,u_{\dimofV-1}(s))$ is nonsingular for every
	$s\in S$.  Put
	\[
	L_s:=\sum_iO_\FLocal u_i(s),\qquad
	\operatorname{Vert}(L_s)
	:=\{\Lambda\subset\Vspecial:\ \Lambda\text{ is a vertex lattice and }
	L_s\subset\Lambda\}.
	\]
	Then $s\mapsto\operatorname{Vert}(L_s)$ is locally constant.  
\end{lem}

\begin{proof}
	Fix $s_0$.  Integrality of the finitely many entries of the Gram matrix
	is a clopen condition. If $(L_{s_0},L_{s_0})\nsubseteq O_\FLocal$, then no integral
	vertex lattice contains it, and after shrinking about $s_0$ all the sets
	$\operatorname{Vert}(L_s)$ are empty.  We may therefore assume that the
	nearby $L_s$ are integral.
	
	Choose $z\in\Vspecial$ whose projection
	to $(L_{s_0}\otimes\FLocal)^\perp$ is nonzero.  After shrinking to a
	clopen neighborhood on which
	$T(u_1(s),\ldots,u_{\dimofV-1}(s))$ remains invertible,
	\[
	d(s):=z-\sum_{i,j}(z,u_i(s))
	\bigl(T(u_1(s),\ldots,u_{\dimofV-1}(s))^{-1}\bigr)_{ij}u_j(s)
	\]
	is a continuous nonzero generator of the one-dimensional anisotropic
	space $(L_s\otimes\FLocal)^\perp$.  Shrinking once more, the valuation of $(d(s),d(s))$
	is constant.
	
	The lattice argument in the proof of
	\Cref{cor:cycle-noetherian} is uniform here.  To see the bound explicitly,
	write an element of a vertex lattice $\Lambda$ containing $L_s$ as
	$x=w+t d(s)$ with $w\in L_s\otimes\FLocal$.  Since
	$(x,L_s)\subset O_\FLocal$, one has $w\in L_s^\vee$, where this dual is
	taken inside $L_s\otimes\FLocal$.  The entries of
	$T(u_1(s),\ldots,u_{\dimofV-1}(s))^{-1}$ are uniformly bounded, so there is one $c$ such that
	$(w,w)\in\varpi^{-c}O_\FLocal$ for all nearby $s$.  On the other hand
	$(x,x)\in O_\FLocal$ because $\Lambda$ is integral.  Hence
	\[
	(t d(s),t d(s))=(x,x)-(w,w)\in\varpi^{-c}O_\FLocal.
	\]
	The valuation of $(d(s),d(s))$ is fixed and the line is anisotropic, so
	$\valFLocal(t)$ has a uniform lower bound, say $t\in\varpi^{-b}O_\FLocal$.
	The coordinates of $u_i(s)$ and $d(s)$ are bounded on the chosen compact
	clopen neighborhood, as are the entries of the inverse Gram matrix.
	Hence there is one full $O_\FLocal$-lattice $B\subset\Vspecial$ such that
	$L_s^\vee+\varpi^{-b}O_\FLocal d(s)\subset B$ for every nearby $s$.
	Enlarge $B$ so that $\varpi B^\vee\subset B$.
	Every vertex lattice containing $L_s$ then satisfies
	\begin{equation}\label{eq:uniform-vertex-sandwich}
	\varpi B^\vee\subset\varpi\Lambda^\vee\subset\Lambda\subset B.
	\end{equation}
	The quotient $B/\varpi B^\vee$ is finite, so only finitely many vertex
	lattices occur in this sandwich. For each of them, the condition
	\[
	L_s\subset\Lambda
	\quad\Longleftrightarrow\quad
	u_i(s)\in\Lambda\ \text{for every }i
	\]
	is clopen by continuity of the $u_i$ and openness of $\Lambda$.
	Intersecting the finitely many neighborhoods of $s_0$ on which these
	conditions have their value at $s_0$ proves
	$\operatorname{Vert}(L_s)=\operatorname{Vert}(L_{s_0})$ there.
\end{proof}
\subsection{Special cycles on the generic fiber}

We recall the construction of special cycles on the generic fiber; see 
\cite[\S12.2]{LZ-orth} for the GSpin case. In this subsection, all Shimura
varieties and morphisms are taken over $\C$, using the embeddings of the
reflex fields into $\C$ specified by the Shimura data.
Let \(1\le r\le\dimofV-2\), and let
\[
\mathbf x=(x_1,\ldots,x_r)\in(\Vcoherent)^r
\]
have a totally positive definite Gram matrix
\[
T(\mathbf x)=((x_i,x_j))_{i,j}\in\Sym_r(\FGlobal)_{>0}.
\]
The span of $\mathbf x$ is then totally positive definite.  Let
$\Vcoherent_\mathbf x$ be its orthogonal complement.  Its signatures are
\[
\operatorname{sgn}(\Vcoherent_{\mathbf x,v})=
\begin{cases}
	(\dimofV-r-2,2),&v=\RealPlace_0,\\
	(\dimofV-r,0),&v\ne\RealPlace_0.
\end{cases}
\]
Put 
$
G_\mathbf x:=\Res_{\FGlobal/\Q}\SO(\Vcoherent_\mathbf x),
$ 
and let $\mathcal D_\mathbf x$ be the Hermitian symmetric domain of
oriented negative $2$-planes in
$\Vcoherent_{\mathbf x,\RealPlace_0}$.  
When $\dimofV-r\geq3$, there is a natural embedding of Shimura data
\[
(G_\mathbf x,\mathcal D_\mathbf x)\hookrightarrow(G_2,X_2).
\]
When $\dimofV-r=2$, each of the two orientations of
$\Vcoherent_{\mathbf x,\RealPlace_0}$ defines a one-point Shimura datum
for the torus $G_\mathbf x$ and an embedding into $(G_2,X_2)$.
In this case, $\ShG(G_\mathbf x,\mathcal D_\mathbf x)$ denotes the disjoint
union of the two canonical models after base change to $\C$, and the following
morphisms are taken on each summand.

For $g\in G_2(\BA_f)$, let $Z(\mathbf x,g)_{K_2}$ denote the source
of the finite morphism
\[
\ShG_{gK_2g^{-1}\cap G_\mathbf x(\BA_f)}
(G_\mathbf x,\mathcal D_\mathbf x)_\C
\longrightarrow
\ShG_{gK_2g^{-1}}(G_2,X_2)_\C
\xrightarrow{\cdot g}
\SOShim(G_2,X_2)_\C.
\]
It defines a class in $\Fil^rK_0(\SOShim(G_2,X_2)_\C)$. 
The complex uniformization gives an isomorphism
\[
Z(\mathbf x,g)_{K_2}(\mathbb C)
\simeq
G_\mathbf x(\Q)\backslash
\bigl[\mathcal D_\mathbf x\times G_\mathbf x(\BA_f)gK_2/K_2\bigr].
\]
Let $\varphi_{K_2}\in\sS((\Vadl\otimes_\FGlobal\BA_{\FGlobal,f})^r)$ be  $K_2$-invariant and let
\(T\in\Sym_r(\FGlobal)_{>0}\).  The weighted special cycle is
\[
Z_T(\varphi_{K_2})
:=
\sum_{\substack{\mathbf x\in G_2(\Q)\backslash(\Vcoherent)^r\\
		T(\mathbf x)=T}}
\ \sum_{g\in G_\mathbf x(\BA_f)\backslash G_2(\BA_f)/K_2}
\varphi_{K_2}(g^{-1}\mathbf x)\,[Z(\mathbf x,g)_{K_2}]_{K_0}
\in\Fil^{r}K_0(\SOShim(G_2,X_2)_\C)_{\ov{\Q}}.
\]

\begin{prop}\label{prop:generic-special-cycle-comparison}
	Assume $\varphi_{K_2}=\varphi_{\fkd}\charfun_{(\VlattiveGlobal\otimes \widehat{\Z}^\fkd)^{r}}$ is admissible, so we can define the special cycles $\CZ_T(\varphi_{K_2})$ on the integral model, using \Cref{prop:global-special-endomorphism-model}. 
	Then we have the equality
	\[
	Z_T(\varphi_{K_2})=\CZ_T(\varphi_{K_2})_\C= \CZ^\Lder_T(\varphi_{K_2})_\C
	\quad\text{in }\Fil^rK_0(\SOShim(G_2,X_2)_\C)_{\ov{\Q}}.
	\]
\end{prop}

\begin{proof}
	Fix a $p$-adic place $v\nmid \fkd$. We show that for any $T\in\Sym_r(\FGlobal)_{>0}$ and any $K_{2,\fkd}$-invariant compact subset $\Omega_\fkd\subseteq \Vadl_\fkd^r$, under the complex uniformization we have an identification
	\[\CZ_{T, (v)}(\Omega_\fkd)_\C\simeq \coprod_{\substack{\mathbf x\in G_2(\Q)\backslash(\Vcoherent)^r\\T(\mathbf x)=T}}\coprod_{\substack{g\in G_\mathbf x(\BA_f)\backslash G_2(\BA_f)/K_2\\g^{-1}\mathbf x\in \Omega_\fkd(\VlattiveGlobal\otimes \widehat{\Z}^\fkd)^{r}}} Z(\mathbf x,g)_{K_2},\]
	where the $\ell$-adic realization map $\realization_\ell$ maps $[y,hg]\in Z(\mathbf x,g)_{K_2}(\C)$ to $$[y,hg,g^{-1}\mathbf x]\in G_2(\Q)\backslash
	\bigl[X_2\times G_2(\BA_f) \times V_\ell^r\bigr/K_2]\cong \CV_\ell^r.$$
	The comparison is compatible with prime-to-$p$ Hecke morphisms as the level $K_2$ and the condition $\Omega_\fkd$ vary. By \Cref{prop:Kisin2010-2.2.6}, it suffices to check the comparison on a connected component of $\SOShim(G_2,X_2)_\C$, hence after pullback to the Hodge-type Shimura variety $\ShG(G^\sharp, X^\sharp)_\C$.
	
	We denote the right hand side by $\CZ_{T, (v)}(\Omega_\fkd)_\C'$.
	On the Hodge-type cover its tautological tuple of special endomorphisms
	defines a comparison map
	$\alpha: \CZ_{T, (v)}(\Omega_\fkd)_\C'\to \CZ_{T, (v)}(\Omega_\fkd)_\C$,
	which descends by \Cref{cor:same-special}.
	Since both schemes are finite over $\ShG_{K_2,\C}$, the map $\alpha$ is finite.
	By \Cref{prop:special-C} it is bijective on $\C$-points.
	At every complex point, the period description identifies the local deformation locus of
	the tuple with the quadric of isotropic lines in the orthogonal complement
	of its span. Since $T$ is nondegenerate, this complement is nondegenerate,
	so the quadric is smooth of dimension $\dimofV-r-2$.
	Thus both schemes are smooth. In characteristic zero, the finite map
	$\alpha$ is birational on each irreducible component. Normality of the target
	then makes $\alpha$ an isomorphism.

	The scheme $Z(\mathbf x,g)_{K_2}$ has dimension $\dimofV-2-r$.
	Thus its image has codimension $r$ in the smooth Shimura variety
	$\SOShim(G_2,X_2)_\C$. The $r$ special-divisor equations therefore
	form a regular sequence, and their derived tensor
	product has no higher homology, which proves $\CO^\BL|_{Z(\mathbf x,g)_{K_2}}\cong \CO|_{Z(\mathbf x,g)_{K_2}}$.
\end{proof}

\subsection{Local--global intersection formula}\label{subsec:local-global-int}

\begin{prop}\label{prop:uniformization-special-cycles-derived}
	Consider a nondegenerate $T\in\Sym_\nInLZorth(\FGlobal)$ and a $K^{\sharp,p}$-invariant compact subset $\Omega\subseteq(\Vadl\otimes_\FGlobal\BA_{\FGlobal,f}^p)^{\dimofV-1}$. Under the uniformization of \Cref{prop:uniformization-sharp-Shimura}, we have the equality of classes in the $K$-group with support on $\CZ^\sharp(T,\Omega)$:
	\begin{align*}
		\CZ^{\sharp, \BL}(T,\Omega)_{O_{\breve{\reflexField},{\wt v}}}
		={}&\sum_{\substack{\mathbf x\in I^\sharp(\Q)\backslash(\VspecialGlobal)^\nInLZorth\\
				T(\mathbf x)=T}}\sum_{g}\\
		&\Bigg[
		(Z^\Q(\Q)\cap K^{\sharp,p} \tilde{K}^{v_0}_p) \backslash
		\left(
		(\CZ^\BL_{\GSpin}(\mathbf x_{v_0}),g)\boxtimes \CO_{\underline{\Q_p^\times/ \Z_p^\times}}
		\right)\Bigg]_{K_0},
	\end{align*}
	where $g=(g^p,g_p)$ ranges over representatives of
	$Z^\Q(\Q)\backslash (G^\sharp(\BA_f^p)\times \tilde{G}^{v_0}(\Q_p))/(K^{\sharp,p} \times \tilde{K}^{v_0}_p)$
	satisfying $(g^p)^{-1}\mathbf x^p\in\Omega$ and
	$\phi_p^{v_0}(g_p^{-1}\mathbf x_p^{v_0})=1$.
	Each summand denotes the class obtained by pulling back the local complex
	to the product with the displayed constant formal scheme, descending through the central quotient,
	and pushing forward to $\Sh_{K^\sharp}\otimes O_{\breve{\reflexField},{\wt v}}$.
\end{prop}
\begin{proof}
	By \Cref{prop:uniformization-sharp-Shimura}, completion of the derived
	structure complex pulls back on each $\mathbf x$-summand to
	$\CZ^{\BL}_{\GSpin}(\mathbf x_{v_0})\boxtimes\CO_{\underline{\Q_p^\times/\Z_p^\times}}$,
	since completion commutes with derived tensor products.
	Its cohomology sheaves are supported on the basic locus by
	\Cref{thm:int-in-basic}, so this computes the original class in the
	$K$-group with support. A nondegenerate $(\dimofV-1)$-tuple has trivial
	pointwise stabilizer in $\SO(\VspecialGlobal)$, so its stabilizer in $I^\sharp(\Q)$
	is therefore $Z^\Q(\Q)$. Taking $I^\sharp(\Q)$-orbits of the tuples
	and then $Z^\Q(\Q)$-orbits of $g$ gives the stated decomposition with stabilizer
	$Z^\Q(\Q)\cap gK^{\sharp,p}\widetilde K_p^{v_0}g^{-1}$.
	Centrality of $Z^\Q$ gives
	\[
	Z^\Q(\Q)\cap gK^{\sharp,p}\widetilde K_p^{v_0}g^{-1}
	=Z^\Q(\Q)\cap K^{\sharp,p}\widetilde K_p^{v_0},
	\]
	which proves the stated description.
\end{proof}

Let $K^\sharp=K^\sharp_pK^{\sharp, p}\subseteq G^\sharp(\BA_f)$ be neat and put
$L:=K^\sharp\cap Z^\Q(\BA_f).$
Since the central torus $Z^\Q$ satisfies the Kottwitz condition
$\rank_\Q Z^\Q=\rank_\R Z^\Q$, the group $Z^\Q(\Q)\cap L$ has rank zero;
neatness therefore gives $Z^\Q(\Q)\cap L=1$.  
Define $L'\subseteq Z^\Q(\BA_f^p)\times \prod_{v\mid p, v\neq v_0}\FGlobal_{v}^\times/O_{\FGlobal,v}^\times $ to be the projection of $L$ under
\[Z^\Q(\BA_f)\cong Z^\Q(\BA_f^p)\times \prod_{v\mid p}\FGlobal_{v}^\times \times \Q_p^\times\twoheadrightarrow Z^\Q(\BA_f^p)\times \prod_{v\mid p, v\neq v_0}\FGlobal_{v}^\times/O_{\FGlobal,v}^\times.\]
This projection is an open continuous map so $L'$ is also compact open. We have a short exact sequence
\begin{align}\label{eq:ses-class-number}
	1\longrightarrow& 
	(Z^\Q(\Q)\cap K^{\sharp, p} \tilde{K}^{v_0}_p)\backslash
	(\FGlobal_{v_0}^\times\times\Q_p^\times)/
	(O_{\FGlobal,v_0}^\times\times\Z_p^\times)\longrightarrow Z^\Q(\Q)\backslash Z^\Q(\BA_f)/L\\
	\longrightarrow& Z^\Q(\Q) \backslash \left(Z^\Q(\BA_f^p)\times \prod_{v\mid p, v\neq v_0}\FGlobal_{v}^\times/O_{\FGlobal,v}^\times\right)/ L'\longrightarrow 1.\notag
\end{align}
We denote
\begin{align*}
\tau_L(Z^\Q)&:=
\#\bigl(Z^\Q(\Q)\backslash Z^\Q(\BA_f)/L\bigr),\\
N_{v_0}&:=\#\Bigl(
(Z^\Q(\Q)\cap K^{\sharp, p} \tilde{K}^{v_0}_p)\backslash
(\FGlobal_{v_0}^\times\times\Q_p^\times)/
(O_{\FGlobal,v_0}^\times\times\Z_p^\times)
\Bigr).
\end{align*}

\begin{prop}\label{thm:global-int-Hodge}
	For an admissible Schwartz function
	$\phi\in \sS((\Vadl\otimes_\FGlobal \BA_{\FGlobal,f})^{\nInLZorth})$
	and a nondegenerate $T\in\Sym_{\nInLZorth}(\FGlobal)$, we have
	\[
	\deg_{v_0}\CZ_T^{\sharp,\BL}(\phi)
	=2N_{v_0}\sum_{\substack{\mathbf x\in I^\sharp(\Q)\backslash(\VspecialGlobal)^{\nInLZorth}\\
			T(\mathbf x)=T}}
	\sum_{g}
	\phi^{v_0}(g^{-1}\mathbf x)\Int^{\epsilon_{v_0}}_{v_0}(\mathbf x_{v_0})
	\]
	where $g\in Z^\Q(\Q) \backslash (G^\sharp(\BA_f^p)\times \tilde{G}^{v_0}(\Q_p))/(K^{\sharp, p} \times \tilde{K}^{v_0}_p)$ and the degree $\deg_{v_0}$ is defined in \eqref{eq:deg-v-convention}, well-defined by \Cref{thm:int-in-basic}. 
\end{prop}

\begin{proof}
	Recall that $\SORZ\cong\GSpinRZ/
\underline{\FGlobal_{v_0}^\times/O_{\FGlobal,v_0}^\times}$. 
	The complex descending from
	$\CZ^\BL_{\GSpin}(\mathbf x_{v_0})$ on the quotient
	\[(Z^\Q(\Q)\cap K^{\sharp, p} \tilde{K}^{v_0}_p)\backslash
	\left[
	(\CZ_{\mathrm{GSpin}}(\mathbf x_{v_0}),g)\times \underline{\Q_p^\times/ \Z_p^\times}
	\right]\]
	has Euler characteristic $2N_{v_0}\Int^{\epsilon_{v_0}}_{v_0}(\mathbf x_{v_0})$
	by \Cref{defn:derived-special-cycles-RZ} and the definition of $N_{v_0}$.
	We then apply \Cref{prop:uniformization-special-cycles-derived}.
\end{proof}

We now descend from $G^\sharp$ to $G_2\cong \Res_{\FGlobal/\Q}\SO(\Vcoherent)$.  Let
\[
1\longrightarrow Z^\Q\longrightarrow G^\sharp
\xlongrightarrow{p}G_2\longrightarrow1
\]
be the central extension and assume $K_2=p(K^\sharp)\subseteq G_2(\BA_f)$. For any sufficiently small  $K_2$ we can find such $K^\sharp$. 
\begin{defn}\label{defn:obst}
	For every place $w$ of $\Q$, let
	$\delta_w:G_2(\Q_w)\longrightarrow H^1(\Q_w,Z^\Q)$
	be the connecting homomorphism, and let
	$\delta:G_2(\Q)\to H^1(\Q,Z^\Q)$ be the global one.  At the finite places set
	\[
	\gamma:=(\delta_\ell)_\ell:
	G_2(\BA_f)\longrightarrow\bigoplus_\ell H^1(\Q_\ell,Z^\Q).
	\]
	Recall that
	$Z^\Q_{\Q_p}\simeq T_{\FGlobal,\Q_p}\times\mathbb G_{m,\Q_p}$, so $H^1(\Q_p, Z^\Q)=0$. 
	Let $X_{2,0}:=p_\R(X^\sharp)\subset X_2$, $G_2(\R)_0:=\{a\in G_2(\R):aX_{2,0}=X_{2,0}\}$. 
	The morphism of Shimura data identifies $X^\sharp$ with $X_{2,0}$, and we have identifications
	\begin{equation*}G_2(\R)_0=p(G^\sharp(\R))=\ker(\delta_\infty).
	\end{equation*}
	Put $G_2(\Q)_0:=G_2(\Q)\cap G_2(\R)_0$, $\beta_0:=\gamma|_{G_2(\Q)_0}$, $\delta_0:=\delta|_{G_2(\Q)_0}$, $$\Gamma_0:=\ker(\beta_0)=
	\{a\in G_2(\Q):\delta_w(a)=1\text{ for every place }w\text{ of }\Q\}.$$
	Define
	\begin{equation*}
		\Sha^1(\Q,Z^\Q):=
		\ker\!\left(
		H^1(\Q,Z^\Q)\longrightarrow
		\prod_w H^1(\Q_w,Z^\Q)
		\right).
	\end{equation*}
\end{defn}

\begin{defn}\label{defn:obst-inner-form}
	Similarly, we can consider the inner forms $I\cong \Res_{\FGlobal/\Q}\SO(\VspecialGlobal)$, $I^\sharp$. We have a central extension
	\[
	1\longrightarrow Z^\Q\longrightarrow I^\sharp
	\xlongrightarrow{p}I\longrightarrow1.
	\]
	For every place $w$ of $\Q$, let
	$\delta_w':I(\Q_w)\longrightarrow H^1(\Q_w,Z^\Q)$
	be the connecting homomorphism, and let
	$\delta':I(\Q)\to H^1(\Q,Z^\Q)$ be the global one. 
	Since $I_\R$ is compact, $I_\R(\R)$ is connected. Hence the map $\delta'_\infty$ is trivial.  Define the finite-adelic connecting map
	\[
	\gamma':=(\delta_\ell')_\ell:
	I(\BA_f)\longrightarrow\bigoplus_\ell H^1(\Q_\ell,Z^\Q),
	\]
	and set
	\[
	\Gamma_0':=\ker\bigl(\gamma'|_{I(\Q)}\bigr)
	=\{a\in I(\Q):\delta_w'(a)=1\text{ for every place }w\text{ of }\Q\}.
	\]
\end{defn}

\begin{lem}\label{lem:global-obstruction-inner-same}
	We have the identifications
	\[
	\im(\gamma')=\im(\gamma),\qquad \im(\delta')=\im(\delta_0).
	\]
	Moreover, the equal quotients
	$\im(\gamma)/\im(\beta_0)=\im(\gamma')/\gamma'(I(\Q))$ are finite.
\end{lem}

\begin{proof}
	Note that $H^1(\Q_p,Z^\Q)$ is trivial and that $G_2$ and $I$ are isomorphic at the places $w\neq p,\infty$. We thus have $\im(\gamma')=\im(\gamma)$. 
	For the second identity, it suffices to verify that
	\[\ker(H^1(\Q,Z^\Q)\to H^1(\Q, G^\sharp)\times H^1(\R,Z^\Q))\simeq \ker(H^1(\Q,Z^\Q)\to H^1(\Q, I^\sharp)\times H^1(\R,Z^\Q)).\]
	Note that
	$(G^\sharp)^{\der}\simeq\Res_{\FGlobal/\Q}\Spin(\Vcoherent)$ is simply
	connected; its inner form $(I^\sharp)^{\der}$ is also simply
	connected. By \cite[Theorem~(3.4.6)]{Kisin2017ModP} we may identify $$\ker(H^1(\Q,G^\sharp)\to H^1(\R,G^\sharp))\cong \ker(H^1(\Q,G^\sharp/(G^\sharp)^{\der})\to H^1(\R,G^\sharp/(G^\sharp)^{\der})),$$ $$\ker(H^1(\Q,I^\sharp)\to H^1(\R,I^\sharp))\cong \ker(H^1(\Q,I^\sharp/(I^\sharp)^{\der})\to H^1(\R,I^\sharp/(I^\sharp)^{\der})).$$
	Since $G^\sharp/(G^\sharp)^{\der}\simeq I^\sharp/(I^\sharp)^{\der}$, we may identify both $\im(\delta')$ and $\im(\delta_0)$ with the kernel of 
	\[H^1(\Q,Z^\Q)\to H^1(\Q,G^\sharp/(G^\sharp)^{\der})\times H^1(\R,Z^\Q).\]
	
	For finiteness, the double-coset set
	$G_2(\Q)\backslash G_2(\BA_f)/K_2$ is finite and $\gamma(K_2)=1$, so
	$\im(\gamma)$ has only finitely many cosets modulo the image of
	$G_2(\Q)$.  Since $X_{2,0}$ is a connected component of $X_2$, the subgroup
	$G_2(\Q)_0$ has finite index in $G_2(\Q)$, and hence the same is true modulo
	$\im(\beta_0)$.
\end{proof}

\begin{lem}\label{lem:local-global-central-obstruction}
	The global connecting maps induce isomorphisms
	\begin{equation*}
		\delta: \Gamma_0/p(G^\sharp(\Q))
		\xrightarrow{\ \sim\ }
		\Sha^1(\Q,Z^\Q),
		\qquad \delta': \Gamma_0'/p(I^\sharp(\Q))
		\xrightarrow{\ \sim\ }\Sha^1(\Q,Z^\Q).
	\end{equation*}
\end{lem}

\begin{proof}
	Injectivity is immediate from exactness. We need to show that every class in $\Sha^1(\Q,Z^\Q)$ maps to the trivial class in $H^1(\Q,G^\sharp)$. We know that its image is a locally trivial class, so it remains to check the Hasse principle for $G^\sharp$. 
	Note that
	\[
	G^\sharp/(G^\sharp)^{\der}
	\simeq
	\Res_{\FGlobal/\Q}(\Res^1_{E/\FGlobal}\Gm)\times\Gm
	\]
	satisfies the Hasse principle by the Hasse norm theorem for the quadratic extension $E/\FGlobal$. The claim then follows from \cite[Theorem~(3.4.6)]{Kisin2017ModP}. The same argument holds for the inner form. 
\end{proof}

Choose $g_1,\ldots,g_s\in G_2(\BA_f^p)$ such that
$\im(\gamma)=
\coprod_{i=1}^s\gamma(g_i)\im(\beta_0),$
and let
$K_i^\sharp:=\ad(g_i)(K^\sharp)$, using
$(G^\sharp)^{\ad}(\BA_f)\simeq G_2^{\ad}(\BA_f)$.
For each $i$, let
$q_i:\Sh_{K_i^\sharp}\to
\Sh_{g_iK_2g_i^{-1}}\times_{O_\FGlobal}O_\reflexField$
be the morphism induced by $p:G^\sharp\to G_2$.

\begin{prop}\label{lem:finite-equal-fiber-KP}
	The morphism
	\[
	p_{K_2}:=\coprod_{i=1}^s(\cdot g_i)\circ q_i:
	\coprod_{i=1}^s\Sh_{K_i^\sharp}\longrightarrow\Sh_{K_2}\times_{O_\FGlobal}O_\reflexField
	\]
	is finite \'etale and surjective, of constant degree
	\begin{equation*}
		\deg(p_{K_2})=\tau_L(Z^\Q)\,\#\Sha^1(\Q,Z^\Q).
	\end{equation*}
\end{prop}

\begin{proof}
	Each summand is finite \'etale by
	\Cref{prop:Kisin-abelian-(3.3.10)}, so it remains to describe its fiber on complex points.  Every point of $\Sh_{K_2}(\C)$ has a representative
	$[x,h]$ with $x\in X_{2,0}$; the two connected components of $X_2$ are exchanged
	by an element of $G_2(\Q)$ by real approximation.  For such a representative define
	\[
	\kappa([x,h]):=\gamma(h)\im(\beta_0)\in
	\bigoplus_\ell H^1(\Q_\ell,Z^\Q)/\im(\beta_0).
	\]
	This does not depend on the choice of representative, because changing a representative while keeping the
	first coordinate in $X_{2,0}$ changes $h$ by an element of $G_2(\Q)_0$ and by an
	element of $K_2$, and $\gamma(K_2)=1$.  Moreover we have
	\[
	\im\bigl((\cdot g_i)\circ q_i\bigr)
	=
	\kappa^{-1}\!\left(
	\gamma(g_i)\im(\beta_0)
	\right).
	\]
	Hence the images of the summands are disjoint and exhaust $\Sh_{K_2}(\C)$. 
	It remains to compute the fibers.  
	
	Fix $y$ in the $i$-th image and write
	$y=[x_0,p(\widetilde h_0)g_i],$ $x_0\in X_{2,0}=X^\sharp.$
	Then the fiber of $(\cdot g_i)\circ q_i$ over $y$ can be identified with the set
	\[G^\sharp (\Q)\backslash \{[x_0', g_f \widetilde h_0]|x_0'\in X_{2,0}, g_f\in G^\sharp(\BA_f), \exists a\in G_2(\Q), ap(g_f)\in \ad_{p(\widetilde h_0)g_i}K_2, ax_0'=x_0\}/K_i^\sharp.\]
	Since the level $K_2$ is neat, for such a pair $[x_0', g_f \widetilde h_0]$ the element $a$ is unique. 
	Also it necessarily satisfies $\delta_\infty(a)=1$, $\gamma(a)=1$ so lies in $\Gamma_0$. 
	Given $a\in \Gamma_0$, we can find a lift $g_{a,f}\in G^\sharp(\BA_f)$ of $a^{-1}\in G_2(\BA_f)$.  
	We can rewrite the set as
	\[\coprod_{a\in \Gamma_0/p(G^\sharp(\Q))} Z^\Q(\Q)\backslash Z^\Q(\BA_f)/ (Z^\Q(\BA_f)\cap g_{a,f}K_i^\sharp g_{a,f}^{-1})\]
	Again, since $Z^\Q$ is central, we may identify $Z^\Q(\BA_f)\cap g_{a,f}K_i^\sharp g_{a,f}^{-1}$ with $Z^\Q(\BA_f)\cap K_i^\sharp=L$. The claim then follows from \Cref{lem:local-global-central-obstruction}.
\end{proof}

\begin{thm}\label{thm:global-int}
	For an admissible Schwartz function
	$
	\phi\in
	\sS((\Vadl\otimes_\FGlobal \BA_{\FGlobal,f})^{\nInLZorth})
	$
	and a nonsingular \(T\in\Sym_{\nInLZorth}(\FGlobal)\), let
	\(\mathbf x\in(\VspecialGlobal)^{\nInLZorth}\)
	have Gram matrix \(T\).  At the distinguished finite place \(v_0\),
	\[
	\deg_{v_0}\CZ_T^{\BL}(\phi)
	=
	2\Int^{\epsilon_{v_0}}_{v_0}(\mathbf x_{v_0})
	\sum_{g\in
		\SO(\Vadl)(\BA_{\FGlobal,f}^{v_0})/K_2^{v_0}}
	\phi^{v_0}(g^{-1}\mathbf x).
	\]
	Equivalently, for any Haar measure on $\SO(\VspecialGlobal)(\BA_{\FGlobal,f}^{v_0})$, we have
	\[
	\deg_{v_0}\CZ_T^{\BL}(\phi)
	=
	\frac{2}{\vol(K_2^{v_0})}
	\Int^{\epsilon_{v_0}}_{v_0}(\mathbf x_{v_0})
	\operatorname{Orb}_{\mathbf x}(\phi^{v_0}),
	\]
	where the stabilizer of \(\mathbf x\) is trivial and
	\[
	\operatorname{Orb}_{\mathbf x}(\phi^{v_0})
	=
	\int_{\SO(\VspecialGlobal)(\BA_{\FGlobal,f}^{v_0})}
	\phi^{v_0}(h^{-1}\mathbf x)\,dh,
	\qquad
	\vol(K_2^{v_0})=\int_{K_2^{v_0}}dh.
	\] 
\end{thm}

\begin{proof}
	For the $i$-th cover, denote 
	$K_i^{\sharp,p}=\ad(g_i)(K^{\sharp,p})$. 
	The local intersection number depends only on the nonsingular Gram
	matrix $T$ by \Cref{defn:derived-special-cycles-RZ}.
	Applying \Cref{lem:finite-equal-fiber-KP,thm:global-int-Hodge} and
	factoring out $\Int^{\epsilon_{v_0}}_{v_0}(\mathbf x_{v_0})$ gives
	\begin{align}
		\deg_{v_0}\CZ_T^{\BL}(\phi)
		={}&
		\Int^{\epsilon_{v_0}}_{v_0}(\mathbf x_{v_0})\,
		\frac{2N_{v_0}}{\tau_L(Z^\Q)\,\#\Sha^1(\Q,Z^\Q)}
		\sum_{i=1}^s
		\sum_{\substack{\mathbf y\in I^\sharp(\Q)\backslash(\VspecialGlobal)^{\nInLZorth}\\
				T(\mathbf y)=T}}
		\sum_g\phi^{v_0}(g_i^{-1}\bar g^{-1}\mathbf y).
		\label{eq:global-int-unfold-1}
	\end{align}
	Here $g=(g^p,g_p^{v_0})$ ranges over
	\[
	Z^\Q(\Q)\backslash
	\bigl(G^\sharp(\BA_f^p)\times\widetilde G^{v_0}(\Q_p)\bigr)/
	\bigl(K_i^{\sharp,p}\times\widetilde K_p^{v_0}\bigr),
	\]
	and $\bar g\in\SO(\Vadl)(\BA_{\FGlobal,f}^{v_0})$ is its image
	under the orthogonal quotient maps.
	
	For each $\mathbf y$ there is a unique
	$a_{\mathbf y}\in\SO(\VspecialGlobal)(\FGlobal)$ with
	$\mathbf y=a_{\mathbf y}\mathbf x$. Using this fact, we can rewrite the sum as $\sum_{a,g}\phi^{v_0}(g_i^{-1}\bar g^{-1}a\mathbf x)$, where 
	\[(a,g^p,g_p^{v_0})\in I^\sharp(\Q)\backslash \left(\SO(\VspecialGlobal)(\FGlobal) \times  G^\sharp(\BA_f^p)\times \tilde{G}^{v_0}(\Q_p)\right)/(K_i^{\sharp, p} \times \tilde{K}^{v_0}_p).\]
	
	Consider the following map
	\[\Xi=(\Xi_i): \coprod_{i=1}^s I^\sharp(\Q)\backslash \left(\SO(\VspecialGlobal)(\FGlobal) \times  G^\sharp(\BA_f^p)\times \tilde{G}^{v_0}(\Q_p)\right)/(K_i^{\sharp, p} \times \tilde{K}^{v_0}_p)\to \SO(\Vadl)(\BA_{\FGlobal,f}^{v_0})/K_2^{v_0} \]
	defined by $\Xi_i(a,g)=[a^{-1}\bar g g_i]$. The diagonal
	action of $I^\sharp(\Q)$ leaves this class unchanged. The image of the
	right level lies in $g_iK_2^{v_0}g_i^{-1}$, so right multiplication by
	the level also preserves this class. Thus $\Xi_i$ is well defined.
	We claim that every fiber has
	the same cardinality, namely
	\[
	\frac{\tau_L(Z^\Q)\,\#\Sha^1(\Q,Z^\Q)}{N_{v_0}}.
	\]
	
	Since $H^1(\Q_p,Z^\Q)=0$, the map $\gamma'$ descends to
	$\SO(\Vadl)(\BA_{\FGlobal,f}^{v_0})/K_2^{v_0}$ through the
	identification of $I$ with $G_2$ away from $p$ and $\infty$;
	the factors at places above $p$ contribute the trivial class. 
	Given $b\in \SO(\Vadl)(\BA_{\FGlobal,f}^{v_0})/K_2^{v_0}$,
	$b$ lies in the image of $\Xi_i$ precisely when
	$\gamma'(b)\in\gamma(g_i)\gamma'(I(\Q))$.
	These are the distinct cosets $\gamma(g_i)\im(\beta_0)$ by
	\Cref{lem:global-obstruction-inner-same}, so the index $i$ is unique.
	
	A class $a\in p(I^\sharp(\Q))\backslash I(\Q)$ admits $g$ with
	$\Xi_i(a,g)=b$ precisely when
	$\gamma'(a)\gamma'(b)\gamma(g_i)^{-1}=1$,
	since $\bar g$ must lie in the right coset $abg_i^{-1}$,
	whose connecting class is this product.
	By \Cref{lem:local-global-central-obstruction}, these classes form
	a torsor under
	$\Gamma_0'/p(I^\sharp(\Q))\simeq\Sha^1(\Q,Z^\Q)$. 
	
	Choose a representative $a\in I(\Q)$ of such a class. With $g$ ranging
	over $G^\sharp(\BA_f^p)\times\widetilde G^{v_0}(\Q_p)$, define
	\[
	S_a:=Z^\Q(\Q)\backslash
	\left\{
	g(K_i^{\sharp,p}\times\widetilde K_p^{v_0}):
	[a^{-1}\bar g g_i]=b
	\right\}.
	\]
	The residual stabilizer of $a$ under the diagonal action is
	$Z^\Q(\Q)$. Central multiplication makes $S_a$ a torsor under
	\[Z^\Q(\Q) \backslash \left(Z^\Q(\BA_f^p)\times \prod_{v\mid p, v\neq v_0}\FGlobal_{v}^\times/O_{\FGlobal,v}^\times\right)/ L'.\]
	Here the rational center is projected as in \eqref{eq:ses-class-number},
	and conjugation of the level by $g_i$ preserves $L'$.
	Combining this description with \eqref{eq:ses-class-number} gives
	$\#S_a=\tau_L(Z^\Q)/N_{v_0}$.  Hence
	\[\#\Xi^{-1}(b)=\sum_{a}\# S_a=\frac{\tau_L(Z^\Q)\,\#\Sha^1(\Q,Z^\Q)}{N_{v_0}}.\]
	Substituting this fiber cardinality into
	\eqref{eq:global-int-unfold-1} gives
	\[\deg_{v_0}\CZ_T^{\BL}(\phi)
	=
	2\Int^{\epsilon_{v_0}}_{v_0}(\mathbf x_{v_0})\sum_{g\in \SO(\Vadl)(\BA_{\FGlobal,f}^{v_0})/K_2^{v_0}}\phi^{v_0}(g^{-1}\mathbf x).\]
	The claim for orbit integrals follows immediately.
\end{proof}

\subsection{A cancellation law}
Throughout this subsection, $\dimofV\geq3$, and $\GSpinRZ$ and $\SORZ$
denote the  GSpin and SO \RZname spaces.  We identify the special
cycle of a self-dual line with the  spaces attached to its
orthogonal complement.
With the factor numbering of \eqref{eq:dec-G-breve}, we choose the
framing representative $b=(b_\FLocal,1,\ldots,1)$, where
$f=[k_\FLocal:\BF_p]$ is the number of factors.

\begin{lem}
	\label{lem:self-dual-line-nonempty}
	Let $L\subset\Vspecial$ be a self-dual quadratic line.  Then
	$\CZ_{\mathrm{GSpin}}(L)$ and $\CZ_{\mathrm{SO}}(L)$ are nonempty.
\end{lem}
\begin{proof}
	Adjoin to $L$ a sufficiently small  lattice in $L_\FLocal^\perp$
	and extend their orthogonal sum to a maximal integral lattice
	$\Lambda'\subset\Vspecial$.  Maximal-lattice theory over the
	non-dyadic discrete valuation ring $O_\FLocal$ gives
	$\varpi(\Lambda')^\vee\subset\Lambda'\subset(\Lambda')^\vee$;
	see \cite[Sections~91--93]{OMeara1973}.  Thus $\Lambda'$ is a vertex
	lattice.  The vertex-stratum description in
	\Cref{prop:irreducible-components-SO} gives a nonempty reduced
	locus of $\CZ_{\mathrm{SO}}(\Lambda')$.  Since $L\subset\Lambda'$,
	this cycle is contained in $\CZ_{\mathrm{SO}}(L)$.
	The central projection $\GSpinRZ\to\SORZ$ is surjective, and the
	GSpin special cycle is the inverse image of the SO special cycle
	by \Cref{defn:SO-local-special-cycle}.  This proves both assertions.
\end{proof}

Let $\Lattice=O_\FLocal x\subset\Vspecial$ be a self-dual quadratic
line, so $(x,x)\in O_\FLocal^\times$, and put
$\Vspecial^\flat:=x^\perp$.  Choose a geometric point
$z_0\in\CZ_{\mathrm{GSpin}}(\Lattice)$ using
\Cref{lem:self-dual-line-nonempty}.  We use this point  to 
choose a framing for the lower-dimensional datum.

\begin{lem}\label{lem:unit-special-vector-gauge}
	There exist $x_0\in\VlattiveLocaltot$ and
	$g_0\in\GSpin(\Vlocal)(\breve\FLocal)$ such that
	$g_0\GSpin(\VlattiveLocaltot)(O_{\breve\FLocal})$ represents $z_0$
	and $\rho_V(g_0)x_0=x$. Define
	\[
	\Vlocal^\flat:=x_0^\perp\subset\Vlocal,\qquad
	\VlattiveLocaltot^\flat:=\VlattiveLocaltot\cap x_0^\perp,\qquad
	b^\flat:=g_0^{-1}b_\FLocal\sigma_\FLocal(g_0).
	\]
	Then $\VlattiveLocaltot=O_\FLocal x_0\operp\VlattiveLocaltot^\flat$
	is a decomposition into self-dual lattices,
	$b^\flat\in\GSpin(\Vlocal^\flat)(\breve\FLocal)$ is basic, and
	$\rho_V(g_0)$ induces an isometry
	\[
	(\Vlocal^\flat\otimes_\FLocal\breve\FLocal)^{
		\rho_V(b^\flat)\sigma_\FLocal=1}
	\xrightarrow{\ \sim\ }\Vspecial^\flat.
	\]
\end{lem}
\begin{proof}
	The reduction of $\VlattiveLocaltot$ contains a hyperbolic plane,
	since it is nondegenerate of rank at least three over a finite
	field of odd characteristic. Hensel lifting on the smooth quadric
	$(x_0,x_0)=(x,x)$ gives $x_0\in\VlattiveLocaltot$.
	Choose an affine Deligne--Lusztig representative
	$g\in\GSpin(\Vlocal)(\breve\FLocal)$ of $z_0$ using
	\Cref{lem:ADLV-Shapiro}. By \Cref{prop:special-ADLV},
	$\rho_V(g)^{-1}x$ is integral and has the same unit norm as $x_0$.
	Their unimodular orthogonal complements over $O_{\breve\FLocal}$
	are isometric. Extending the line isometry and correcting its
	determinant by a reflection in the complement gives an element of
	$\SO(\VlattiveLocaltot)(O_{\breve\FLocal})$ sending $x_0$ to
	$\rho_V(g)^{-1}x$. Since $\operatorname{Pic}(O_{\breve\FLocal})=0$,
	it lifts to $k\in\GSpin(\VlattiveLocaltot)(O_{\breve\FLocal})$.
	Take $g_0=gk$.

	The unit norm gives the self-dual orthogonal decomposition.
	The equality $\rho_V(b_\FLocal)\sigma_\FLocal(x)=x$ implies that
	$b^\flat$ fixes $x_0$, so it lies in $\GSpin(\Vlocal^\flat)$.
	Conjugation by $\rho_V(g_0)$ intertwines the Frobenius operators
	and gives the fixed-space isometry. The vector Newton slopes
	are zero, so the Newton cocharacter of $b^\flat$ lies in the
	scalar center; hence $b^\flat$ is basic.
\end{proof}
\begin{lem}\label{lem:lower-cocharacter-admissibility}
	Let $\mathscr U^\flat$ be the set of standard minuscule classes of
	$\GSpin(\Vlocal^\flat)$ whose images belong to $\{\mu_\FLocal\}$.
	It has one element when $\dimofV\geq4$ and two Galois-conjugate
	elements when $\dimofV=3$. We have
	\[
	[b^\flat]\in B(\GSpin(\Vlocal^\flat),\lambda^{-1})
	\qquad\text{for every }\lambda\in\mathscr U^\flat.
	\]
	In particular, we have a local Shimura datum $(\GSpin(\Vlocal^\flat),b^\flat,\lambda)$ for each $\lambda\in \mathscr U^\flat$. 
\end{lem}
\begin{proof}
	The inclusion $\GSpin(\Vlocal^\flat)\hookrightarrow\GSpin(\Vlocal)$
	preserves the spinor similitude character $\nu$. Since $b^\flat$
	is $\sigma_\FLocal$-conjugate to $b_\FLocal$,
	\eqref{eq:GSpin-standard-cocharacter} gives
	\[
	\valFLocal\bigl(\nu(b^\flat)\bigr)=-1,
	\qquad \nu\circ\lambda^{-1}(t)=t^{-1}
	\quad(\lambda\in\mathscr U^\flat).
	\]
	For $\dimofV\geq4$, the geometric isotropic lines of
	$\Vlocal^\flat$ form one $\SO(\Vlocal^\flat)$-orbit. Together
	with the prescribed similitude weight, this gives one standard
	minuscule class. Moreover, $\ker\nu=\Spin(\Vlocal^\flat)$ is
	simply connected, so $\nu$ induces a Galois-equivariant
	isomorphism $\pi_1(\GSpin(\Vlocal^\flat))\simeq\Z$.

	Consider the case $\dimofV=3$. 
	If $\Vlocal^\flat$ were split, then
	$\GSpin(\Vlocal^\flat)\simeq\Gm^2$ with
	$\nu(t_1,t_2)=t_1t_2$. The zero vector Newton slopes of
	$b^\flat$ from \Cref{lem:unit-special-vector-gauge} would force
	its two coordinate valuations to agree, making
	$\valFLocal(\nu(b^\flat))$ even. Thus $\Vlocal^\flat$ is anisotropic. Since $\VlattiveLocaltot^\flat$ is self-dual, $C^+(\Vlocal^\flat)$ is
	an unramified quadratic extension of $\FLocal$. Frobenius exchanges the two
	geometric isotropic lines, giving two Galois-conjugate classes
	in $\mathscr U^\flat$. The two embeddings of $C^+(\Vlocal^\flat)$
	identify the cocharacter lattice of $\GSpin(\Vlocal^\flat)$
	with $\Z^2$. Frobenius exchanges the coordinates, and $\nu_*$
	is their sum. Hence $\nu_*$ identifies the Galois coinvariants
	with $\Z$.

	In both cases, $[b^\flat]$ and every $\lambda^{-1}$ have the
	same Kottwitz invariant, namely $-1$ under the identification
	induced by $\nu$. Since $[b^\flat]$ is basic by
	\Cref{lem:unit-special-vector-gauge},
	\cite[\S6.4]{Kottwitz1997IsocrystalsII} gives
	$[b^\flat]\in B(\GSpin(\Vlocal^\flat),\lambda^{-1})$
	for every $\lambda\in\mathscr U^\flat$.
\end{proof}

Identify the lower special space with
$\Vspecial^\flat$ via $\rho_V(g_0)$, and put
$\epsilon':=\chi(\Vlocal^\flat)$. For each
$\lambda\in\mathscr U^\flat$, let
$\mathcal M_{\mathrm{GSpin},\lambda}^\flat$ and
$\mathcal M_{\mathrm{SO},\lambda}^\flat$ be the \RZname spaces over
$O_{\breve\FLocal}$
for the lower GSpin datum and its SO quotient, respectively.
The absolute framing element is $(b^\flat,1,\ldots,1)$, and the
absolute cocharacter is obtained from $\lambda$ as in
\Cref{defn:GSpin-pure-local-Shimura-datum}.
In rank two, the unramified quadratic reflex field has completed
maximal unramified extension $\breve\FLocal$. Define
\[
\RZGS:=\coprod_{\lambda\in\mathscr U^\flat}
\mathcal M_{\mathrm{GSpin},\lambda}^\flat,
\qquad
\RZSOflat:=\coprod_{\lambda\in\mathscr U^\flat}
\mathcal M_{\mathrm{SO},\lambda}^\flat.
\]
For every $\lambda\in\mathscr U^\flat$, restriction of scalars of the
inclusion $\GSpin(\Vlocal^\flat)\hookrightarrow\GSpin(\Vlocal)$ and
the change of framing by $(\sigma^i(g_0))_{0\le i<f}$ give a morphism
to $\GSpinRZ$.
Their disjoint union is
$\delta_\Lattice:\RZGS\to\GSpinRZ$. Passing to the central quotients
gives $\delta_\Lattice^{\mathrm{SO}}:\RZSOflat\to\SORZ$.
The integral decomposition
$\VlattiveLocaltot=O_\FLocal x_0\operp\VlattiveLocaltot^\flat$ gives
\begin{equation*}
 C(\VlattiveLocaltot)\simeq
 C(\VlattiveLocaltot^\flat)\otimes_{O_\FLocal}^{\mathrm{gr}}C(O_\FLocal x_0).
\end{equation*}
The graded tensor product uses Clifford parity and multiplication
$(c\otimes b)(a\otimes d)=(-1)^{\deg a\deg b}ca\otimes bd$
for homogeneous $a,b,c,d$. As a representation of
$\GSpin(\VlattiveLocaltot^\flat)$, the right-hand side is the direct sum
of two copies of $C(\VlattiveLocaltot^\flat)$. Thus the pulled-back
universal $p$-divisible group on $\RZGS$ is $\sX^\flat\oplus\sX^\flat$,
where $\sX^\flat$ is the universal $p$-divisible group on $\RZGS$ for the lower GSpin spaces.
The vector $x$ acts integrally through the Clifford grading and multiplication
by $x_0$. For an integral special endomorphism $y$ over $S/\RZGS$,
the integral special pairing is $O_\FLocal$-valued.
Since $(x,x)$ is a unit, both
$\frac{(x,y)}{(x,x)}x$ and $y-\frac{(x,y)}{(x,x)}x$ are integral.
The latter acts diagonally on the two copies of $\sX_S^\flat$, through
its lower Clifford action on each copy. Its integrality is therefore
equivalent to integrality of the lower special endomorphism. This identifies
$\sEndint(\sX_S)$ with $\sEndint(\sX_S^\flat)\oplus\Lattice$.
In particular, $\delta_\Lattice$ factors through
\[
\gamma_\Lattice:\RZGS\longrightarrow\CZ_{\mathrm{GSpin}}(\Lattice).
\] 
\begin{prop}\label{prop:cancellation-pullback}
	For $y^\flat\in\Vspecial^\flat$, write
	$\CZ_{\mathrm{GSpin}}^\flat(y^\flat)\subset\RZGS$ and
	$\CZ_{\mathrm{SO}}^\flat(y^\flat)\subset\RZSOflat$ for the  lower
	special cycles, including both cocharacters when $\dimofV=3$.
	For $y\in\Vspecial\setminus O_\FLocal x$ with $(x,y)\in O_\FLocal$, define
	\[
	y_\Lattice:=\frac{(x,y)}{(x,x)}x\in\Lattice,
	\qquad y^\flat:=y-y_\Lattice\in\Vspecial^\flat\setminus\{0\}.
	\]
	Then
	\[
	\delta_\Lattice^{-1}\CZ_{\mathrm{GSpin}}(y)
	=\CZ_{\mathrm{GSpin}}^\flat(y^\flat),\qquad
	(\delta_\Lattice^{\mathrm{SO}})^{-1}\CZ_{\mathrm{SO}}(y)
	=\CZ_{\mathrm{SO}}^\flat(y^\flat).
	\]
	The same identities hold for the derived pullbacks of their
	structure sheaves.
\end{prop}
\begin{proof}
	The GSpin identity follows from the preceding decomposition of $\sEndint(\sX_S)\cong \sEndint(\sX_S^\flat)\oplus\Lattice$.
	Descent under the central action of $\FLocal^\times/O_\FLocal^\times$ gives the SO identity.
	
	The lower cycles $\CZ_{\mathrm{GSpin}}^\flat(y^\flat)$ and
	$\CZ_{\mathrm{SO}}^\flat(y^\flat)$ are Cartier divisors by
	\Cref{prop:local-special-cycle-properties} (in rank two use
	\Cref{prop:n-two-special}). Hence the pulled-back local equations
	are nonzerodivisors, so the higher Tor sheaves vanish.
\end{proof}
\begin{prop}\label{prop:gamma-iso-GSpin}
	The morphism $\gamma_\Lattice$ is an isomorphism of formal schemes
	over $O_{\breve\FLocal}$, equivariant for the central framing action
	of $\FLocal^\times/O_\FLocal^\times$.
\end{prop}
\begin{proof}
	The source is formally smooth by the self-dual local model, using
	\Cref{lem:n2-RZ-component} in lower dimension two. Since $(x,x)$ is a
	unit, the target is formally smooth by
	\Cref{prop:local-special-cycle-properties}. Both are therefore normal
	and flat, so it suffices to identify their diamonds by
	\cite[Proposition~18.4.1]{ScholzeWeinstein2020Berkeley}.

	Let $S$ be a characteristic-$p$ perfectoid space over the residue field
	of $O_{\breve\FLocal}$. For an object of
	$\CZ_{\mathrm{GSpin}}(\Lattice)^\diamondsuit(S)$, let $\sP$ be the $\Res_{O_\FLocal/\Z_p}\GSpin(\VlattiveLocaltot)$ 
	torsor of \Cref{prop:GSpin-RZ-defn} and let
	$x'\in\sV$ be the section supplied by \Cref{prop:special-perfectoid}.
	It has norm $(x_0,x_0)\in O_\FLocal^\times$. The sheaf
	\[
	\sP^\flat:=\{a\in\sP:a(x_0)=x'\}
	\]
	is a torsor under
	$\Res_{O_\FLocal/\Z_p}\GSpin(\VlattiveLocaltot^\flat)$.
	Indeed, the smooth quadric of vectors of norm $(x_0,x_0)$ is the
	homogeneous space of $\GSpin(\VlattiveLocaltot)$ with stabilizer
	$\GSpin(\VlattiveLocaltot^\flat)$; this follows from the unimodular
	orthogonal decomposition and $p\ne2$.

	Frobenius preserves $\sP^\flat$ because $x'$ is Frobenius compatible.
	Composing the framing with $(\sigma^i(g_0^{-1}))_{0\le i<f}$ gives
	the framing element $(b^\flat,1,\ldots,1)$. The ambient minuscule local model parametrizes
	isotropic lines, and the extension condition for $x'$ cuts out
	the quadric of isotropic lines in $(x')^\perp$.
	By \Cref{lem:lower-cocharacter-admissibility}, the lower modification
	is therefore bounded by a class $\lambda\in\mathscr U^\flat$ satisfying $\nu\circ\lambda^{-1}(t)=t^{-1}$. 
	When \(\dimofV=3\), the quadric of isotropic lines in \((x')^\perp\) consists of two geometric points, corresponding to the two classes in \(\mathscr U^\flat\).
	Thus $\sP^\flat$ defines an object of $\RZGS^\diamondsuit(S)$.

	Conversely, extension of structure group from
	$\Res_{O_\FLocal/\Z_p}\GSpin(\VlattiveLocaltot^\flat)$ to
	$\Res_{O_\FLocal/\Z_p}\GSpin(\VlattiveLocaltot)$,
	with framing changed by $(\sigma^i(g_0))_{0\le i<f}$, gives an ambient object together with
	the integral section $x_0$, framed by $x$. These constructions are
	mutually inverse and commute with base change in $S$. This proves the
	isomorphism of diamonds. Since central scalars fix $x_0$, it is
	equivariant for the central framing action.
\end{proof}

\begin{cor}\label{prop:gamma-iso}
	The induced morphism
	$\gamma_\Lattice^{\mathrm{SO}}:
	\RZSOflat\to\CZ_{\mathrm{SO}}(\Lattice)$
	is an isomorphism of formal schemes over $O_{\breve\FLocal}$.
\end{cor}
\begin{proof}
	The central projection realizes $\SORZ$ as the $v$-sheaf quotient
	$\GSpinRZ/(\FLocal^\times/O_\FLocal^\times)$; this is the $v$-sheaf form of the central comparison
	\cite[Theorem~5.1.2]{PappasRapoport2026IntegralLocalShimura}
	used in \Cref{subsec:SO-RZ-components}.
	The same assertion gives $\RZGS/(\FLocal^\times/O_\FLocal^\times)\simeq\RZSOflat$, using
	\Cref{lem:n2-RZ-component} for each rank-two datum.
	Central scalars commute with special quasi-endomorphisms, and
	\Cref{defn:SO-local-special-cycle} identifies
	$\CZ_{\mathrm{GSpin}}(\Lattice)/(\FLocal^\times/O_\FLocal^\times)$ with
	$\CZ_{\mathrm{SO}}(\Lattice)$.
	Taking the quotients of the equivariant isomorphism in
	\Cref{prop:gamma-iso-GSpin} proves the claim on $v$-sheaves.
	Normality and flatness then give the formal-scheme assertion.
\end{proof}

\begin{cor}[Cancellation law]\label{cor:cancellation}
	Let $x\in\Vspecial$ span the self-dual line $M=O_\FLocal x$.
	Let $\bu=(x,u_2,\ldots,u_{\dimofV-1})$ have integral nonsingular
	Gram matrix, and let
	$\bu^\flat=(u_2^\flat,\ldots,u_{\dimofV-1}^\flat)$ be the tuple
	of orthogonal projections to $\Vspecial^\flat=x^\perp$.  Then
	\[
	\Int^\epsilon(\bu)=\Int^{\epsilon'}(\bu^\flat),
	\qquad \epsilon'=\chi(\Vlocal^\flat).
	\]
	When $\dimofV=3$, the right-hand side may be computed on the 
	rank-two SO space for either cocharacter.
\end{cor}
\begin{proof}
	The projected Gram matrix is the Schur complement of the unit
	entry $(x,x)$, hence is integral and nonsingular.  In particular,
	all $u_i^\flat$ are nonzero.  By \Cref{prop:gamma-iso}, the 
	divisor $\CZ_{\mathrm{SO}}(M)$ is $\RZSOflat$.
	Using \Cref{prop:cancellation-pullback} and the projection formula
	for its closed immersion into $\SORZ$, we obtain
	\[
	R\Gamma\!\left(\SORZ,\CZ_{\SORZ}^{\Lder}(\bu)\right)
	\simeq
	R\Gamma\!\left(\RZSOflat,
	\bigotimes_{2\leq i\leq\dimofV-1}^{\Lder}
	\CO_{\CZ_{\mathrm{SO}}^\flat(u_i^\flat)}\right),
	\]
	where the derived tensor product on the right is over
	$\CO_{\RZSOflat}$.

	If $\dimofV\geq4$, both intersection numbers are one half of the
	respective total derived degrees, by
	\Cref{defn:derived-special-cycles-RZ}.  The displayed isomorphism
	therefore proves their equality.
	If $\dimofV=3$, the lower space is the disjoint union of the 
	SO spaces indexed by $\mathscr U^\flat$.
	Their derived degrees agree by \Cref{prop:n-two-special}, and each
	equals $\Int^{\epsilon'}(\bu^\flat)$.  
\end{proof}

	\part{Proof of the arithmetic Siegel--Weil formula}\label{part3}
	
\section{Local Weil representations and Whittaker functions}\label{sec:Weil}
In this section we fix a local field $\FLocal$ of characteristic zero and a nontrivial unitary additive character $\psi:\FLocal\to\BC^\times$. We recall the definition of local Weil representations and Whittaker functions, cf. \cite{Kudla1997CentralDerivatives}.

Let \(V\) be a nondegenerate quadratic \(\FLocal\)-space of
dimension \(\dimofV=\nInLZorth+1\), with symmetric bilinear form \((\ ,\ )\)
and quadratic form \(Q(x)=(x,x)/2\).
\(\Or(V)\) is the group of linear isometries of \(V\).

\subsection{Definition of Gelfand--Leray Haar measures}\label{conv:haar-local}
In this subsection we recall the definition of Gelfand--Leray measure (cf. 
\cite[Notation~2.3.6 and \S4.2]{AizenbudGourevitchKazhdanSayag2026Integrability}). 
Consider the Lie algebra
\[
\mathfrak{so}(V)=\{A\in\operatorname{End}_{\FLocal}(V):
(Au,w)+(u,Aw)=0\text{ for all }u,w\in V\}
\]
with the nondegenerate symmetric bilinear pairing
\[
(A_1,A_2)\longmapsto-\frac12\operatorname{tr}(A_1A_2),
\qquad A_1,A_2\in\mathfrak{so}(V),
\]
where the trace is taken on \(\operatorname{End}_{\FLocal}(V)\).
\begin{defn}
	For an integer \(r\ge1\), let \(d\mathbf y^{\mathrm{sd}}\),
	\(dS^{\mathrm{sd}}\), and \(dA^{\mathrm{sd}}\) be the Haar measures
	on \(V^r\), \(\Sym_{\nInLZorth}(\FLocal)\), and
	\(\mathfrak{so}(V)\), respectively, that are self-dual with respect
	to \(\psi\) and the respective pairings
	\[
	(\mathbf x,\mathbf y)\longmapsto\sum_{i=1}^r(x_i,y_i),\qquad
	(b_1,b_2)\longmapsto\tfrac12\operatorname{tr}(b_1b_2),\qquad
	(A_1,A_2)\longmapsto-\tfrac12\operatorname{tr}(A_1A_2).
	\]
	We use the self-duality convention of \textup{(H\ref{defn:Haar-self-dual})}.
\end{defn}

\begin{defn}
	Recall the Gram matrix map of \eqref{eq:Gram-matrix-convention}:
	\[
	T:V^{\nInLZorth}\longrightarrow\Sym_{\nInLZorth}(\FLocal),
	\qquad
	T(\mathbf y)=((y_i,y_j))_{1\le i,j\le\nInLZorth}.
	\]
	Its differential at \(\mathbf y\) sends \(\mathbf z\) to
	\(((z_i,y_j)+(y_i,z_j))_{i,j}\), and is surjective whenever
	\(T(\mathbf y)\) is nonsingular.  For a nonsingular represented
	matrix \(S\), the quotient of the source density
	\(d\mathbf y^{\mathrm{sd}}\) by the target density \(dS^{\mathrm{sd}}\),
	formed using the tangent exact sequence of this submersion, defines
	a measure \(d\mu_S^{\mathrm{GL}}\) on the fiber \(T^{-1}(S)\).
	When \(\FLocal=\FGlobal_v\), we write \(d\mu_{S,v}^{\mathrm{GL}}\)
	for this measure.
	For every continuous compactly supported function
	\(f:V^{\nInLZorth}\to\mathbb C\),
	the following disintegration formula holds\footnote{For the non-archimedean fiber-integration formula, see
		\cite[Theorem~7.6.1]{Igusa2000LocalZeta}.}
	\begin{equation}\label{eq:GL-disintegration}
		\int_{V^{\nInLZorth}}f(\mathbf y)\,d\mathbf y^{\mathrm{sd}}
		=
		\int_{\Sym_{\nInLZorth}(\FLocal)}
		\left(\int_{T^{-1}(S)}
		f(\mathbf y)\,d\mu_S^{\mathrm{GL}}(\mathbf y)\right)
		dS^{\mathrm{sd}}.
	\end{equation}
	Here the inner integral is zero for empty fibers and is extended by zero
	on singular matrices, whose inverse image has measure zero.
\end{defn}
\begin{defn}\label{defn:Gelfand--Leray}
	For a nonsingular \(S\in\Sym_{\nInLZorth}(\FLocal)\) represented
	by \(V\), choose \(\mathbf x\in V^{\nInLZorth}\) with
	\(T(\mathbf x)=S\). By Witt's theorem and the one-dimensional
	orthogonal complement of \(\sum_i\FLocal x_i\), every isometry
	\(x_i\mapsto y_i\), with \(\mathbf y\in T^{-1}(S)\), has a unique
	extension of determinant \(1\) to \(V\). Thus the map
	\begin{equation}\label{eq:GL-torsor-map}
		\SO(V)(\FLocal)\longrightarrow T^{-1}(S),
		\qquad h\longmapsto h^{-1}\mathbf x,
	\end{equation}
	is an isomorphism of \(\FLocal\)-analytic manifolds.  Let \(dh\) be the pullback of
	\(d\mu_S^{\mathrm{GL}}\) under this map, which is Haar measure on \(\SO(V)(\FLocal)\). We call it the
	\emph{Gelfand--Leray} measure. 
\end{defn}

\begin{lem}\label{lem:GL-Whittaker-measure}
	Identify the tangent space
	\(T_1\SO(V)\) at the identity with \(\mathfrak{so}(V)\), and let
	\((dh)_1\) denote the density of \(dh\) at $1$. Then we have
	\begin{equation}\label{eq:GL-trace-density}
		(dh)_1=\lvert2\rvert^{-\nInLZorth/2}
		dA^{\mathrm{sd}}. 
	\end{equation}
	In particular, \(dh\) depends only on \(V\) and \(\psi\).
\end{lem}

\begin{proof}
	A column change \(\mathbf y\mapsto\mathbf y a\), with
	\(a\in\GL_{\nInLZorth}(\FLocal)\), changes the source and target measures
	by \(\lvert\det a\rvert^{\dimofV}\) and
	\(\lvert\det a\rvert^{\nInLZorth+1}\), respectively.
	These factors agree, so we may diagonalize \(S\) and take an orthogonal
	basis \(e_1,\ldots,e_{\dimofV}\) of \(V\) with
	\(\mathbf x=(e_1,\ldots,e_{\nInLZorth})\).
	For \(1\le i<j\le\dimofV\), the endomorphisms defined by
	\[
	A_{ij}e_i=(e_i,e_i)e_j,\qquad
	A_{ij}e_j=-(e_j,e_j)e_i,\qquad
	A_{ij}e_k=0\quad(k\ne i,j)
	\]
	form a basis of \(\mathfrak{so}(V)\).  Define coordinates by
	\[
	y_i=\sum_{j=1}^{\dimofV}y_{ji}e_j,\qquad
	A=\sum_{1\le i<j\le\dimofV}a_{ij}A_{ij},
	\]
	and use the independent entries \(s_{ij}\),
	\(1\le i\le j\le\nInLZorth\), on
	\(\Sym_{\nInLZorth}(\FLocal)\).
	Absolute values of coordinate top forms below use the self-dual
	measure for \(\psi(uv)\) in each coordinate; wedges are ordered
	lexicographically.  The three pairings on \(V^{\nInLZorth}\),
	\(\Sym_{\nInLZorth}(\FLocal)\), and \(\mathfrak{so}(V)\) are diagonal in these
	coordinates, with respective diagonal entries
	\[
	(e_j,e_j)\quad(1\le i\le\nInLZorth,\ 1\le j\le\dimofV),
	\qquad
	\begin{cases}\tfrac12&i=j,\\1&i<j,\end{cases}
	\qquad
	(e_i,e_i)(e_j,e_j)\quad(i<j).
	\]
	Hence
	\begin{align*}
		d\mathbf y^{\mathrm{sd}}
		&=\left(\prod_{j=1}^{\dimofV}|(e_j,e_j)|^{\nInLZorth/2}\right)
		\left|\bigwedge_{i=1}^{\nInLZorth}
		\bigwedge_{j=1}^{\dimofV}dy_{ji}\right|,\qquad
		dS^{\mathrm{sd}}
		=|2|^{-\nInLZorth/2}
		\left|\bigwedge_{1\le i\le j\le\nInLZorth}ds_{ij}\right|,\\
		dA^{\mathrm{sd}}
		&=\left(\prod_{1\le i<j\le\dimofV}
		|(e_i,e_i)(e_j,e_j)|^{1/2}\right)
		\left|\bigwedge_{1\le i<j\le\dimofV}da_{ij}\right|
		=\left(\prod_{j=1}^{\dimofV}|(e_j,e_j)|^{\nInLZorth/2}\right)
		\left|\bigwedge_{1\le i<j\le\dimofV}da_{ij}\right|.
	\end{align*}
	At \(\mathbf x\), the Gram differential is
	\[
	ds_{ii}=2(e_i,e_i)\,dy_{ii},\qquad
	ds_{ij}=(e_i,e_i)\,dy_{ij}+(e_j,e_j)\,dy_{ji}
	\quad(i<j).
	\]
	Use the following splitting of the tangent exact sequence:
	\[
	\begin{aligned}
		\mathfrak{so}(V)\oplus\Sym_{\nInLZorth}(\FLocal)
		&\longrightarrow T_{\mathbf x}V^{\nInLZorth},\\
		\bigl(A,(s_{ij})\bigr)
		&\longmapsto
		\left(-Ae_i+\frac{s_{ii}}{2(e_i,e_i)}e_i
		+\sum_{j<i}\frac{s_{ji}}{(e_j,e_j)}e_j
		\right)_{i=1}^{\nInLZorth}.
	\end{aligned}
	\]
	Its coordinate blocks are
	\[
	\begin{pmatrix}dy_{ji}\\dy_{ij}\end{pmatrix}
	=
	\begin{pmatrix}-(e_i,e_i)&0\\(e_j,e_j)&(e_i,e_i)^{-1}\end{pmatrix}
	\begin{pmatrix}da_{ij}\\ds_{ij}\end{pmatrix}
	\quad(i<j\le\nInLZorth),
	\]
	\[
	\begin{pmatrix}dy_{ii}\\dy_{\dimofV,i}\end{pmatrix}
	=
	\begin{pmatrix}0&(2(e_i,e_i))^{-1}\\-(e_i,e_i)&0\end{pmatrix}
	\begin{pmatrix}da_{i,\dimofV}\\ds_{ii}\end{pmatrix}
	\quad(1\le i\le\nInLZorth).
	\]
	Their determinants are \(-1\) and \(1/2\), respectively, so the
	pullback of the source top form is
	\[
	\bigwedge_{i=1}^{\nInLZorth}\bigwedge_{j=1}^{\dimofV}dy_{ji}
	=\pm2^{-\nInLZorth}
	\left(\bigwedge_{1\le i<j\le\dimofV}da_{ij}\right)
	\wedge\left(\bigwedge_{1\le i\le j\le\nInLZorth}ds_{ij}\right).
	\]
	Taking the quotient by the target density gives
	\[
	\begin{aligned}
		(dh)_1
		&=\frac{|2|^{-\nInLZorth}}{|2|^{-\nInLZorth/2}}
		\left(\prod_{j=1}^{\dimofV}|(e_j,e_j)|^{\nInLZorth/2}\right)
		\left|\bigwedge_{1\le i<j\le\dimofV}da_{ij}\right|=|2|^{-\nInLZorth/2}dA^{\mathrm{sd}}.
	\end{aligned}\qedhere
	\]
\end{proof}

\begin{lem}\label{lem:GL-hyperspecial-volume}
	Assume that \(\FLocal\) is non-archimedean of odd residue characteristic,
	\(\psi\) has conductor \(O_{\FLocal}\), and \(\Lambda\subset V\)
	is a self-dual lattice for \((\ ,\ )\).
	The group scheme \(\SO(\Lambda)\) of determinant-one isometries
	of \(\Lambda\) is smooth and reductive over \(O_{\FLocal}\), and
	\(K^\circ=\SO(\Lambda)(O_{\FLocal})\) is a hyperspecial subgroup
	of \(\SO(V)(\FLocal)\).  Its \(dh\)-measure is
	\begin{equation}\label{eq:GL-hyperspecial-volume}
		\vol(K^\circ)
		=q^{-\nInLZorth\dimofV/2}
		\#\SO(\Lambda\otimes_{O_{\FLocal}}k_\FLocal)(k_\FLocal).
	\end{equation}
	More explicitly,
	\begin{equation}\label{eq:GL-hyperspecial-volume-product}
		\vol(K^\circ)=
		\begin{cases}
			\displaystyle\prod_{i=1}^{(\dimofV-1)/2}(1-q^{-2i}),
			&\dimofV\text{ odd},\\[6pt]
			\displaystyle(1-q^{-\dimofV/2})
			\prod_{i=1}^{\dimofV/2-1}(1-q^{-2i}),
			&\dimofV\text{ even, split reduction},\\[6pt]
			\displaystyle(1+q^{-\dimofV/2})
			\prod_{i=1}^{\dimofV/2-1}(1-q^{-2i}),
			&\dimofV\text{ even, nonsplit reduction}.
		\end{cases}
	\end{equation}
	In the even-dimensional cases, split and nonsplit mean that the
	quadratic space \(\Lambda\otimes_{O_{\FLocal}}k_\FLocal\) has Witt index
	\(\dimofV/2\) and \(\dimofV/2-1\), respectively. Empty products have value \(1\).
	If \(\epsilon=\chi(\Lambda)\), with \(\chi\) as in
	\Cref{sec:notat-quadr-latt-1}, then
	\[
	\vol(K^\circ)=\operatorname{Nor}^{\epsilon}(1,L)
	\]
	for every nondegenerate integral quadratic \(O_{\FLocal}\)-lattice
	\(L\) of rank \(\nInLZorth\); the normalizing polynomial
	\(\operatorname{Nor}^{\epsilon}(X,L)\) is defined in
	\Cref{defn:local-density}.
\end{lem}
\begin{proof}
	Choose an orthogonal \(O_{\FLocal}\)-basis
	\(e_1,\ldots,e_{\dimofV}\) of \(\Lambda\).  Self-duality implies
	\((e_i,e_i)\in O_{\FLocal}^{\times}\) for every \(i\).
	The endomorphisms \(A_{ij}\) in the preceding proof form an
	\(O_{\FLocal}\)-basis of
	\[
	\mathfrak{so}(\Lambda)=
	\{A\in\operatorname{End}_{O_{\FLocal}}(\Lambda):
	(Au,w)+(u,Aw)=0\text{ for all }u,w\in\Lambda\}.
	\]
	The trace pairing in this basis has diagonal unit entries
	\((e_i,e_i)(e_j,e_j)\).  The conductor assumption therefore gives
	\(\mathfrak{so}(\Lambda)\) volume \(1\) for
	\(dA^{\mathrm{sd}}\).  Since \(\lvert2\rvert=1\),
	\Cref{lem:GL-Whittaker-measure} shows that the invariant group density
	has coefficient \(1\) in this integral tangent basis.
	
	Since \(\SO(\Lambda)\) is smooth of relative dimension
	\(\nInLZorth\dimofV/2\), the reduction map
	\(K^\circ\to\SO(\Lambda\otimes_{O_{\FLocal}}k_\FLocal)(k_\FLocal)\)
	is surjective, with each fiber of \(dh\)-measure
	\(q^{-\nInLZorth\dimofV/2}\). This proves
	\eqref{eq:GL-hyperspecial-volume}.
	The orders of the finite special orthogonal groups give
	\eqref{eq:GL-hyperspecial-volume-product}.  
\end{proof}

\subsection{Local Weil representations}\label{subsec:genus-one-conventions}

\begin{defn}\label{defn:Weil-repn}
	For an integer \(r\ge1\), equip
	\(\FLocal^r\oplus\FLocal^r\) with the standard symplectic form
	\[
	((x,y),(x',y'))\longmapsto{}^txy'-{}^tyx',
	\]
	where \(x,y,x',y'\) are columns.  The Siegel parabolic
	\(P_r=M_rN_r\subset\Sp_{2r}\) stabilizes
	\(\FLocal^r\oplus\{0\}\), and
	\[
	m(a)=\begin{pmatrix}a&0\\0&{}^ta^{-1}\end{pmatrix},\qquad
	n(b)=\begin{pmatrix}1_r&b\\0&1_r\end{pmatrix},\qquad
	w_r=\begin{pmatrix}0&1_r\\-1_r&0\end{pmatrix},
	\]
	with \(a\in\GL_r(\FLocal)\) and \(b\in\Sym_r(\FLocal)\).
	The local metaplectic double cover is
	\[
	1\longrightarrow\{\pm1\}\longrightarrow\Mp_{2r}(\FLocal)
	\longrightarrow\Sp_{2r}(\FLocal)\longrightarrow1.
	\]
	We use the canonical lifts of \(n(b)\) and the Schr\"odinger lifts
	of \(m(a)\) and \(w_r\). Thus \((m(a),t)\),
	\(t\in\{\pm1\}\), denotes an element of the inverse image of \(M_r\).

	Define
	\begin{equation}\label{eq:weil-fourier-transform}
		\widehat\varphi^{\,\psi}(\mathbf x)
		=\int_{V^r}\varphi(\mathbf y)
		\psi\!\left(\sum_{i=1}^r(x_i,y_i)\right)\,d\mathbf y^{\mathrm{sd}}.
	\end{equation}
	Let \(\gamma(\psi)\) be the Weil index of
	\(x\mapsto\psi(x^2/2)\), and put
	\(\gamma(a,\psi):=\gamma(\psi(a\,\cdot))/\gamma(\psi)\) for
	\(a\in\FLocal^\times\); see \cite[\S I.4 and \S II.4]{Kudla1996LocalTheta}.
	The Weil index of \((V,Q)\) is
	\[
	\gamma_{V,\psi}
	=\gamma(\det V,\psi)\gamma(\psi)^{\dimofV}\Has(V),
	\]
	where \(\Has(V)\) is defined in \Cref{sec:notat-quadr-latt-1}.
	The quadratic character is
	\(\chi_{V}(a)=(a,\disc(V))_{\FLocal}\), where
	\((\ ,\ )_{\FLocal}\) is the Hilbert symbol.  Its extension to
	the metaplectic Levi is
	\[
	\chi(m(a),t)=
	\begin{cases}
		t^{\dimofV}\chi_{V}(\det a),&\dimofV\text{ even},\\
		t^{\dimofV}\chi_{V}(\det a)\gamma(\det a,\psi)^{-1},&\dimofV\text{ odd}.
	\end{cases}
	\]
	Extend this character to the inverse image of \(P_r\), trivially on
	the canonically lifted \(N_r\).
	The Weil representation \(\omega=\omega_{V,\psi}\) of
	\(\Mp_{2r}(\FLocal)\times\Or(V)\) on \(\sS(V^r)\) is given by
	\begin{align}
		\omega(n(b))\varphi(\mathbf x)
		&=\psi\!\left(\frac12\operatorname{tr}(bT(\mathbf x))\right)
		\varphi(\mathbf x),\label{eq:weil-n}\\
		\omega(m(a),t)\varphi(\mathbf x)
		&=\chi(m(a),t)|\det a|^{\dimofV/2}\varphi(\mathbf x a),\notag\\
		\omega(w_r)\varphi(\mathbf x)
		&=\gamma_{V,\psi}^{\,r}\widehat\varphi^{\,\psi}(\mathbf x),
		\label{eq:weil-w}\\
		\omega(h)\varphi(\mathbf x)
		&=\varphi(h^{-1}\mathbf x)\qquad(h\in\Or(V)).\notag
	\end{align}
	The covering-kernel element \((1,t)\) acts by \(t^{\dimofV}\).  Thus the
	representation is genuine for odd \(\dimofV\) and descends to
	\(\Sp_{2r}(\FLocal)\) for even \(\dimofV\).
\end{defn}
When \(\FLocal\) is non-archimedean, the Weil representation preserves
\(\ov{\BQ}\)-valued Schwartz functions; see
\cite[\S\S I.4, II.4]{Kudla1996LocalTheta}.
\subsection{Local principal series and Siegel--Weil sections}

\begin{defn}\label{defn:local-principal-series}
	Let
	\(\widetilde P_{\nInLZorth}(\FLocal)\) denote the inverse image of the
	Siegel parabolic in \(\Mp_{2\nInLZorth}(\FLocal)\).  Define
	\begin{equation*}
		I_{\nInLZorth}(s,\chi)
		=\operatorname{Ind}_{\widetilde P_{\nInLZorth}(\FLocal)}
		^{\Mp_{2\nInLZorth}(\FLocal)}
		\bigl(\chi|\det|^{s+(\nInLZorth+1)/2}\bigr).
	\end{equation*}
	This is smooth unnormalized induction, using locally constant functions
	when \(\FLocal\) is non-archimedean and \(C^\infty\) functions when
	\(\FLocal\) is archimedean.  Its functions satisfy
	\begin{equation}\label{eq:local-induced-transformation}
		f((m(a),t)n(c)g)
		=\chi(m(a),t)|\det a|^{s+(\nInLZorth+1)/2}f(g),
	\end{equation}
	with \(a\in\GL_{\nInLZorth}(\FLocal)\),
	\(c\in\Sym_{\nInLZorth}(\FLocal)\), and \(t\in\{\pm1\}\).
	The group acts by \((R(h)f)(g)=f(gh)\).
	
\end{defn}
\begin{defn}\label{defn:standard-section}
	For each integer \(r\ge1\), let \(K^{\mathrm{std}}_r\) be the
	inverse image in \(\Mp_{2r}(\FLocal)\) of
	\[
	\begin{cases}
		\Sp_{2r}(O_{\FLocal}),&\FLocal\text{ non-archimedean},\\
		\Sp_{2r}(\mathbb R)\cap\operatorname O(2r),&\FLocal=\mathbb R,\\
		\Sp_{2r}(\mathbb C)\cap\operatorname U(2r),&\FLocal=\mathbb C.
	\end{cases}
	\]
	Here \(\operatorname O(2r)\) preserves the standard positive definite
	real inner product, and \(\operatorname U(2r)\) preserves the standard
	positive definite Hermitian form.  In the real case,
	\(\Sp_{2r}(\mathbb R)\cap\operatorname O(2r)\simeq\operatorname U(r)\)
	is the stabilizer of \(i1_r\) under the action
	\[
	\begin{pmatrix}A&B\\C&D\end{pmatrix}\cdot Z
	=(AZ+B)(CZ+D)^{-1},
	\qquad Z\in\Sym_r(\mathbb C),\quad\operatorname{Im}Z>0.
	\]
	A \emph{standard section} is a family
	\(f_s\in I_{\nInLZorth}(s,\chi)\), \(s\in\mathbb C\), such that
	\(f_s|_{K^{\mathrm{std}}_{\nInLZorth}}\) is independent of \(s\).
	The Iwasawa decomposition determines the family uniquely from any
	one of its members, and \(s\mapsto f_s(g)\) is entire for every fixed
	\(g\).  At an archimedean place, a section is called
	\(K^{\mathrm{std}}_{\nInLZorth}\)-finite if the complex span of its
	right translates by \(K^{\mathrm{std}}_{\nInLZorth}\) is finite-dimensional.
\end{defn}

\begin{defn}\label{defn:SW-section}
	The \emph{Siegel--Weil map} $\Phi_-: \sS(V^{\nInLZorth})\to I_{\nInLZorth}(0,\chi)$ associates to each \(\varphi\in\sS(V^{\nInLZorth})\) the section 
	\[
	\Phi_{\varphi}(-,0)\in I_{\nInLZorth}(0,\chi),\qquad
	\Phi_{\varphi}(g,0)=(\omega(g)\varphi)(0).
	\]
	This determines a standard section $\Phi_{\varphi}(-,s)\in I_{\nInLZorth}(s,\chi)$.  
	For an Iwasawa decomposition \(g=(m(a),t)n(c)k\), with
	\(a\in\GL_{\nInLZorth}(\FLocal)\),
	\(c\in\Sym_{\nInLZorth}(\FLocal)\), \(t\in\{\pm1\}\), and
	\(k\in K^{\mathrm{std}}_{\nInLZorth}\), it satisfies
	\begin{equation}\label{eq:local-standard-section-Iwasawa}
		\Phi_{\varphi}(g,s)
		=\chi(m(a),t)|\det a|^{s+(\nInLZorth+1)/2}
		(\omega(k)\varphi)(0)=|\det a|^s(\omega(g)\varphi)(0). 
	\end{equation}
	The factor \(|\det a|\) is independent of the Iwasawa decomposition. 
	For every
	\(h\in \Mp_{2\nInLZorth}(\FLocal)\),
	\[
	R(h)\Phi_{\varphi}(\,\cdot\,,0)
	=\Phi_{\omega(h)\varphi}(\,\cdot\,,0).
	\]
	For \(h\in K^{\mathrm{std}}_{\nInLZorth}\), this identity holds
	for every \(s\).
\end{defn}
\begin{defn}\label{defn:local-Whitt}
	For a standard section
	\(\Phi(-,s)\in I_{\nInLZorth}(s,\chi)\), and a nonsingular
	\(T\in\Sym_{\nInLZorth}(\FLocal)\), define
	\begin{equation}\label{eq:local-Whittaker-integral}
		W_{T/2}^{\psi}(g,s,\Phi)
		=\int_{\Sym_{\nInLZorth}(\FLocal)}
		\Phi(w_{\nInLZorth}^{-1}n(b)g,s)
		\psi\!\left(-\frac12\operatorname{tr}(Tb)\right)\,db.
	\end{equation}
	The measure \(db\) is self-dual for the pairing
	\((b_1,b_2)\mapsto\psi(\tfrac12\operatorname{tr}(b_1b_2))\)
	on \(\Sym_{\nInLZorth}(\FLocal)\).
	The integral converges for \(\Re(s)\gg0\); elsewhere it denotes its
	analytic continuation.  For a Siegel--Weil section, write
	\(W_{T/2}^{\psi}(g,s,\varphi)
	=W_{T/2}^{\psi}(g,s,\Phi_{\varphi})\).
\end{defn}

\subsection{Archimedean local Whittaker functions}
In this subsection, assume \(\FLocal=\mathbb R\), that \(V\) is positive
definite of dimension \(\dimofV=\nInLZorth+1\), and that
\[
\psi(u)=\exp(2\pi i\kappa u),\qquad u\in\mathbb R,\quad\kappa>0.
\]
We use the same exponential formula to extend \(\psi\) to \(\mathbb C\).
\begin{defn}\label{defn:Gaussian}
	For every integer \(r\ge1\), every
	\(\sz=\sx+i\sy\in\Sym_r(\mathbb C)\) with
	\(\sx,\sy\in\Sym_r(\mathbb R)\) and \(\sy\) positive definite,
	and every \(T\in\Sym_r(\mathbb R)\), define
	\begin{equation*}
		q^{T/2}(\sz)
		=\psi\!\left(\frac12\operatorname{tr}(T\sz)\right)
		=\exp\!\left(\pi i\kappa\operatorname{tr}(T\sz)\right).
	\end{equation*}
	The Gaussian \(\varphi_\psi\in\sS(V^{\nInLZorth})\) is
	\begin{equation*}
		\varphi_\psi(\mathbf x)
		=\exp\!\left(-\pi\kappa\operatorname{tr}T(\mathbf x)\right),
		\qquad \mathbf x\in V^{\nInLZorth},
	\end{equation*}
	where \(T(\mathbf x)\) is the Gram matrix defined in
	\Cref{sec:notat-quadr-latt-1}.
	For \(\sz=\sx+i\sy\) of genus \(\nInLZorth\), let
	\(a=\sy^{1/2}\) be the positive definite square root, and set
	\(g_\sz=n(\sx)m(a)\), using the Schr\"odinger lift of
	\(m(a)\) fixed in \Cref{defn:Weil-repn}.
	The Siegel gamma factor is
	\[
	\Gamma_{\nInLZorth}(s)
	=\pi^{\nInLZorth(\nInLZorth-1)/4}
	\prod_{j=1}^{\nInLZorth}\Gamma(s-(j-1)/2).
	\]
\end{defn}

\begin{lem}\label{lem:matrix-Laplace}
	Let \(b\in\Sym_{\nInLZorth}(\mathbb R)\).  The notation \(Y>0\)
	means that \(Y\in\Sym_{\nInLZorth}(\mathbb R)\) is positive definite,
	and \(dY=\prod_{1\le j\le k\le\nInLZorth}dY_{jk}\) is product
	Lebesgue measure in its independent entries.  Then
	\[
	\det(1_{\nInLZorth}-ib)^{-(\nInLZorth+1)/2}
	=
	\Gamma_{\nInLZorth}((\nInLZorth+1)/2)^{-1}
	\int_{Y>0}e^{-\operatorname{tr}Y}e^{i\operatorname{tr}(bY)}\,dY.
	\]
	Here \(\Gamma_{\nInLZorth}\) is defined in
	\Cref{defn:Gaussian}.  If
	\(\lambda_1,\ldots,\lambda_{\nInLZorth}\) are the real eigenvalues
	of \(b\), counted with multiplicity, the determinant power means
	\[
	\exp\!\left(-\frac{\nInLZorth+1}{2}
	\sum_{j=1}^{\nInLZorth}\operatorname{Log}(1-i\lambda_j)\right),
	\]
	where \(\operatorname{Log}\) is the principal logarithm on the right
	half-plane.
\end{lem}

\begin{proof}
	Apply \cite[equation~(1.4.6), p.~19]{GuptaNagar2000}
	with matrix argument \(1_{\nInLZorth}-ib\) and parameter
	\((\nInLZorth+1)/2\).  The real part of the matrix argument is
	\(1_{\nInLZorth}>0\), and
	\((\nInLZorth+1)/2>(\nInLZorth-1)/2\), so the integral converges
	absolutely.  The power of \(\det Y\) in that formula is zero.
	Theorem~1.4.1 on the same page gives the normalization
	\(\Gamma_{\nInLZorth}\) used here.  The determinant power in
	the cited formula is defined by analytic continuation from positive
	definite real matrices; along \(1_{\nInLZorth}-itb\),
	\(0\le t\le1\), it agrees with the principal-logarithm expression above.
\end{proof}

\begin{prop}\label{prop:archimedean-whittaker-factor-0}
	For every positive definite
	\(T\in\Sym_{\nInLZorth}(\mathbb R)\), one has
	\begin{equation}\label{eq:archimedean-whittaker-local-explicit-0}
		\begin{aligned}
			W_{T/2}^{\psi}(g_\sz,0,\varphi_\psi)
			&=\gamma_{V,\psi}^{-\nInLZorth}
			\kappa^{\nInLZorth(\nInLZorth+1)/4}
			\frac{
				2^{\nInLZorth/2}
				\pi^{\nInLZorth(\nInLZorth+1)/2}}
			{\Gamma_{\nInLZorth}((\nInLZorth+1)/2)}(\det\sy)^{\dimofV/4}
			q^{T/2}(\sz).
		\end{aligned}
	\end{equation}
\end{prop}

\begin{proof}
	Using \eqref{eq:local-standard-section-Iwasawa}
	and the Gaussian of \Cref{defn:Gaussian}, the Fourier-transform formula in
	\Cref{subsec:genus-one-conventions} gives
	\[
	\Phi_{\varphi_\psi}(w_{\nInLZorth}^{-1}n(b),0)
	=\gamma_{V,\psi}^{-\nInLZorth}
	\det(1_{\nInLZorth}-ib)^{-(\nInLZorth+1)/2}.
	\]
	The determinant power uses the branch specified in
	\Cref{lem:matrix-Laplace}.
	Hence
	\begin{equation*}
		W_{T/2}^{\psi}(1,0,\varphi_\psi)
		=\gamma_{V,\psi}^{-\nInLZorth}
		\int_{\Sym_{\nInLZorth}(\mathbb R)}
		\det(1_{\nInLZorth}-ib)^{-(\nInLZorth+1)/2}
		e^{-\pi i\kappa\operatorname{tr}(T b)}\,db.
	\end{equation*}
	At \(s=0\), the integral is interpreted as a Fourier transform of
	tempered distributions, consistently with the analytic continuation
	in \Cref{defn:local-Whitt}.
	
	Let \(db_{\mathrm{Leb}}\) be product Lebesgue measure in the independent
	entries of \(b\).  Self-duality for
	\(\psi(\tfrac12\operatorname{tr}(b_1b_2))\) gives
	\begin{equation}\label{eq:sym-trace-self-dual-measure}
		db
		=\kappa^{\nInLZorth(\nInLZorth+1)/4}
		2^{-\nInLZorth/2}\,db_{\mathrm{Leb}}.
	\end{equation}
	Indeed, the diagonal coordinates have coefficient \(\kappa/2\) in the
	bicharacter and the off-diagonal coordinates have coefficient \(\kappa\).
	Fourier inversion in the independent entries gives the distributional identity
	\[
	\int_{\Sym_{\nInLZorth}(\mathbb R)}
	e^{i\operatorname{tr}(b(Y-\pi\kappa T))}\,db_{\mathrm{Leb}}
	=2^{\nInLZorth}\pi^{\nInLZorth(\nInLZorth+1)/2}
	\delta_{\pi\kappa T}(Y),
	\]
	where \(\delta_{\pi\kappa T}\) is the Dirac distribution with respect
	to \(dY\).  The integrable density
	\(e^{-\operatorname{tr}Y}\) on the positive definite cone, extended
	by zero, is smooth near the positive definite matrix \(\pi\kappa T\).  Applying distributional
	Fourier inversion to the Laplace formula therefore evaluates this density
	at \(Y=\pi\kappa T\).  Combining this with
	\eqref{eq:sym-trace-self-dual-measure} yields
	\begin{equation*}
		W_{T/2}^{\psi}(1,0,\varphi_\psi)
		=\gamma_{V,\psi}^{-\nInLZorth}
		\kappa^{\nInLZorth(\nInLZorth+1)/4}
		\frac{
			2^{\nInLZorth/2}
			\pi^{\nInLZorth(\nInLZorth+1)/2}}
		{\Gamma_{\nInLZorth}((\nInLZorth+1)/2)}
		e^{-\pi\kappa\operatorname{tr}T}.
	\end{equation*}
	The Weil representation formulas and the change of variables
	\(b\mapsto ab{}^ta\) in the Whittaker integral give
	\[
	W_{T/2}^{\psi}(g_\sz,0,\varphi_\psi)
	=\chi(m(a))|\det a|^{\dimofV/2}
	e^{\pi i\kappa\operatorname{tr}(T\sx)}
	W_{{}^taTa/2}^{\psi}(1,0,\varphi_\psi),
	\]
	where we use \(\dimofV=\nInLZorth+1\).
	Since \(a=\sy^{1/2}\) and \(\chi(m(a))=1\), the definition of
	\(q^{T/2}(\sz)\) gives
	\eqref{eq:archimedean-whittaker-local-explicit-0}.
\end{proof}

\subsection{Non-archimedean local Whittaker functions}
In the rest of this section we assume $\FLocal$ is a non-archimedean local field with residue cardinality \(q\) and
uniformizer \(\varpi\).  
For Whittaker functions on \(\Sp_{2\nInLZorth}(\FLocal)\), see
\cite[Theorem~3.6, Corollary~3.6.1]{Karel1979Whittaker}.
The proofs of
\Cref{lem:bruhat-boundary-coinvariants,lem:bad-finite-whittaker-rational}
adapt the arguments of
\cite[Theorems~2.7 and~3.6]{Karel1979Whittaker} to the metaplectic group.

Consider the families \(F(-,s)\) satisfying
\eqref{eq:local-induced-transformation}, taking values in
\(\ov{\Q}[q^{-s},q^s]\), and fixed under right translation by some
compact open subgroup independent of \(s\). Denote their
\(\ov{\Q}[q^{-s},q^s]\)-module by \(I_{\nInLZorth}(\chi)\), with
\(\Mp_{2\nInLZorth}(\FLocal)\) acting by right translation.
For \(T\in\Sym_{\nInLZorth}(\FLocal)\) and a smooth representation
\(U\) of this group over \(\ov{\Q}[q^{-s},q^s]\), its twisted coinvariant module
\(U_{N_{\nInLZorth}(\FLocal),\psi,T}\) is the quotient by the submodule generated by
\[
n(a)u-\psi\!\left(\tfrac12\operatorname{tr}(Ta)\right)u,
\qquad u\in U,\quad a\in\Sym_{\nInLZorth}(\FLocal).
\]
Here \(N_{\nInLZorth}(\FLocal)\) acts through its canonical lift.
This coinvariant functor is exact; see
\cite[\S1.3.3]{Karel1979Whittaker}.

\begin{lem}\label{lem:bruhat-boundary-coinvariants}
	Assume that \(T\in\Sym_{\nInLZorth}(\FLocal)\) is nonsingular.
	Let
	\(U\subset \widetilde P_{\nInLZorth}(\FLocal)\backslash
	\Mp_{2\nInLZorth}(\FLocal)\) consist of the cosets
	\(\widetilde P_{\nInLZorth}(\FLocal)h\) satisfying
	\[
	h^{-1}(\FLocal^{\nInLZorth}\oplus\{0\})
	\cap(\FLocal^{\nInLZorth}\oplus\{0\})= 0,
	\]
	where \(h\) acts through its image in the symplectic group.
	The locus \(U\) is open and stable under right translation by
	\(N_{\nInLZorth}(\FLocal)\).  Let \(I_{\nInLZorth}'(\chi)\) be the submodule
	whose support modulo \(\widetilde P_{\nInLZorth}(\FLocal)\) is contained in
	\(U\).  Then $I_{\nInLZorth}(\chi)/I_{\nInLZorth}'(\chi)$ has zero \((N_{\nInLZorth}(\FLocal),\psi,T)\)-coinvariants.
\end{lem}
\begin{proof}
	Identify the quotient $\widetilde P_{\nInLZorth}(\FLocal)\backslash
	\Mp_{2\nInLZorth}(\FLocal)$ with the Lagrangian Grassmannian by
	\(\widetilde P_{\nInLZorth}(\FLocal)h\mapsto
	h^{-1}(\FLocal^{\nInLZorth}\oplus\{0\})\). 
	Restriction identifies
	\(I_{\nInLZorth}(\chi)/I_{\nInLZorth}'(\chi)\) with locally constant
	functions on the inverse image in \(\Mp_{2\nInLZorth}(\FLocal)\)
	of the complement of \(U\). These functions take values in
	\(\ov{\Q}[q^{-s},q^s]\), satisfy
	\eqref{eq:local-induced-transformation}, and have compact support
	modulo \(\widetilde P_{\nInLZorth}(\FLocal)\); see
	\cite[Lemma~1.4.3]{Karel1979Whittaker}.
	Let \(F\) be any such function. It has a compact open right stabilizer in
	\(N_{\nInLZorth}(\FLocal)\) independent of \(s\).
	
	Fix \(c\in\FLocal\) with \(\psi(c)\ne1\), and choose an integer
	\(r\ge0\) such that
	\(cT^{-1}O_{\FLocal}^{\nInLZorth}
	\subset\varpi^{-r}O_{\FLocal}^{\nInLZorth}\).
	For a Lagrangian \(L\) in $\widetilde P_{\nInLZorth}(\FLocal)\backslash
	\Mp_{2\nInLZorth}(\FLocal)-U$, choose
	\(v\in O_{\FLocal}^{\nInLZorth}\) with \((v,0)\in L\) and a unit
	coordinate \(v_k\).  Such a vector exists because
	\(L\cap(\FLocal^{\nInLZorth}\oplus\{0\})\ne0\).
	Let \(e_k\) be the \(k\)-th standard column vector, and set
	\(u=cv_k^{-1}T^{-1}e_k\in \varpi^{-r}O_{\FLocal}^{\nInLZorth}\) and \(a_0=v\,{}^tu+u\,{}^tv\in\Sym_{\nInLZorth}(\varpi^{-r}O_{\FLocal})\).
	Note $\tfrac12\operatorname{tr}(Ta_0)={}^tuTv=c$. 
	Isotropy gives \({}^tvy=0\) for every \((x,y)\in L\), so
	\[
	n(ta_0)(x,y)=(x+t v\,{}^tuy,y)\in L
	\qquad(t\in\FLocal).
	\]
	For \(h\in \Mp_{2\nInLZorth}(\FLocal)\) satisfying
	\(L=h^{-1}(\FLocal^{\nInLZorth}\oplus\{0\})\), the element
	\(hn(ta_0)h^{-1}\) belongs to \(\widetilde P_{\nInLZorth}(\FLocal)\).
	The action of its symplectic image on
	\(\FLocal^{\nInLZorth}\oplus\{0\}\) is conjugate to
	\(n(ta_0)|_L\).  The preceding formula gives
	\((n(ta_0)|_L-\operatorname{id}_L)^2=0\), so its Levi component
	has determinant \(1\).  Since
	\(\chi_V(1)=\gamma(1,\psi)=1\), the formula for \(\chi\) in
	\Cref{defn:Weil-repn} shows that
	\(\chi(hn(ta_0)h^{-1})\in\{\pm1\}\).
	The canonical lifts satisfy
	\(n(ta_0)=n((t/2)a_0)^2\).  As \(\chi\) is a character of
	\(\widetilde P_{\nInLZorth}(\FLocal)\), it follows that
	\[
	\chi(hn(ta_0)h^{-1})
	=\chi(hn((t/2)a_0)h^{-1})^2=1
	\qquad(t\in\FLocal).
	\]
	Thus both factors of the inducing character
	\(\chi|\det|^{s+(\nInLZorth+1)/2}\) are \(1\) on these elements.
	Using \eqref{eq:local-induced-transformation} and the commutativity
	of \(N_{\nInLZorth}(\FLocal)\), we obtain, for every
	\(a\in\Sym_{\nInLZorth}(\FLocal)\),
	\[
	F(hn(a+a_0),s)
	=F\bigl((hn(a_0)h^{-1})hn(a),s\bigr)
	=F(hn(a),s).
	\]
	Translation by \(a_0\) preserves
	\(\Sym_{\nInLZorth}(\varpi^{-r}O_{\FLocal})\) and multiplies the
	twisted integral by \(\psi(c)\).  Since \(1-\psi(c)\) is invertible
	in \(\ov{\Q}[q^{-s},q^s]\), this proves
	\[
	\int_{\Sym_{\nInLZorth}(\varpi^{-r}O_{\FLocal})}
	F(hn(a),s)\psi\!\left(-\tfrac12\operatorname{tr}(Ta)\right)\,da=0.
	\]
	In the twisted coinvariant module, \(F\) has the same image as the
	function obtained by dividing the displayed integral by
	\(\operatorname{vol}(\Sym_{\nInLZorth}(\varpi^{-r}O_{\FLocal}))\).
	The displayed identity therefore shows that the image of \(F\) is zero.
\end{proof}

\begin{lem}\label{lem:bad-finite-whittaker-rational}
	Let \(\varphi\in\sS(V^{\nInLZorth},\ov{\Q})\).
	For every nonsingular \(T\in\Sym_{\nInLZorth}(\FLocal)\) and
	\(g\in\Mp_{2\nInLZorth}(\FLocal)\), the integral
	\eqref{eq:local-Whittaker-integral} agrees in its half-plane of
	absolute convergence with a Laurent polynomial in \(q^{-s}\)
	with coefficients in \(\ov{\Q}\).  Thus it has an entire
	continuation satisfying
	\begin{equation*}
		W_{T/2}^{\psi}(g,s,\varphi)
		\in\ov{\Q}[q^{-s},q^s]
	\end{equation*} 
\end{lem}

\begin{proof}
	Fix \(g\). The family
	\(h\mapsto\Phi_{\varphi}(hg,s)\) belongs to
	\(I_{\nInLZorth}(\chi)\) by
	\eqref{eq:local-standard-section-Iwasawa} and the smoothness of the
	Weil representation.
	
	By \Cref{lem:bruhat-boundary-coinvariants},  the quotient \(I_{\nInLZorth}(\chi)/I_{\nInLZorth}'(\chi)\) has zero \((N_{\nInLZorth}(\FLocal),\psi,T)\)-coinvariants. 
	By exactness, we can find 
	\(\Phi^{\circ}(-,s)\in I_{\nInLZorth}'(\chi)\) such that the difference
	\(R(g)\Phi_{\varphi}-\Phi^{\circ}\) has trivial image in $I_{\nInLZorth}(\chi)_{N_{\nInLZorth}(\FLocal),\psi,T}$. By the definition of this quotient, the difference
	is a finite sum of families
	\[
	h\longmapsto F(hn(a),s)
	-\psi\!\left(\tfrac12\operatorname{tr}(Ta)\right)F(h,s),
	\]
	with fixed \(a\in\Sym_{\nInLZorth}(\FLocal)\) and fixed families
	\(F(-,s)\in I_{\nInLZorth}(\chi)\).
	
	Choose an integer \(R\ge0\) so that
	\(\Sym_{\nInLZorth}(\varpi^{-R}O_{\FLocal})\) contains the support
	of \(b\mapsto\Phi^{\circ}(w_{\nInLZorth}^{-1}n(b),s)\) for every
	\(s\), and contains every \(a\) occurring in that finite sum.
	The Whittaker integral of \(\Phi^{\circ}\) can therefore be taken
	over this lattice, for every \(s\).  For each summand and every
	\(R'\ge R\), translation by \(a\) preserves
	\(\Sym_{\nInLZorth}(\varpi^{-R'}O_{\FLocal})\).  The substitution
	\(b\mapsto b+a\) thus gives
	\[
	\begin{aligned}
		&\int_{\Sym_{\nInLZorth}(\varpi^{-R'}O_{\FLocal})}
		\Bigl(F(w_{\nInLZorth}^{-1}n(b+a),s)
		-\psi\!\left(\tfrac12\operatorname{tr}(Ta)\right)
		F(w_{\nInLZorth}^{-1}n(b),s)\Bigr)
		\psi\!\left(-\tfrac12\operatorname{tr}(Tb)\right)\,db=0.
	\end{aligned}
	\]
	Hence the truncated integrals of
	\(h\mapsto\Phi_{\varphi}(hg,s)\) and \(\Phi^{\circ}(-,s)\) agree for every
	\(R'\ge R\).  This implies 
	\begin{equation}\label{eq:bad-finite-whittaker-lattice}
		W_{T/2}^{\psi}(g,s,\varphi)
		=\int_{\Sym_{\nInLZorth}(\varpi^{-R}O_{\FLocal})}
		\Phi_{\varphi}(w_{\nInLZorth}^{-1}n(b)g,s)
		\psi\!\left(-\tfrac12\operatorname{tr}(Tb)\right)\,db.
	\end{equation}
	The integral over this fixed compact lattice is a finite sum of
	Laurent polynomials in \(q^{-s}\) with algebraic coefficients.
\end{proof}

The local group measure of \Cref{defn:Gelfand--Leray} gives the
following comparison with orbital integrals.
\begin{lem}\label{lem:local-Whittaker-GL-exact}
	Let
	\(\varphi\in\sS(V^{\nInLZorth})\).
	For a nonsingular \(S\in\Sym_{\nInLZorth}(\FLocal)\) represented by \(V\),
	choose \(\mathbf x\) with \(T(\mathbf x)=S\).  Then
	\begin{equation*}
		W_{S/2}^{\psi}(1,0,\varphi)
		=\gamma_{V,\psi}^{-\nInLZorth}
		\int_{\SO(V)(\FLocal)}\varphi(h^{-1}\mathbf x)\,dh.
	\end{equation*}
\end{lem}
\begin{proof}
	By complex linearity and a finite decomposition of \(\varphi\)
	into characteristic functions of compact open subsets, we may assume
	that \(\varphi\) takes values in \(\ov{\Q}\).
	By \eqref{eq:GL-disintegration} applied to \(\lvert\varphi\rvert\),
	the function
	\[
	U\longmapsto\int_{T^{-1}(U)}
	\varphi(\mathbf y)\,d\mu_U^{\mathrm{GL}}(\mathbf y),
	\qquad U\in\Sym_{\nInLZorth}(\FLocal),
	\]
	is absolutely integrable and, by the submersion property of \(T\)
	and the compact support of \(\varphi\), locally constant on a
	neighborhood of \(S\) consisting of nonsingular matrices.
	For an integer \(R\ge0\), the superscript \(\perp\) below denotes
	the annihilator of \(\Sym_{\nInLZorth}(\varpi^{-R}O_{\FLocal})\)
	under \((U,b)\mapsto\psi(\tfrac12\operatorname{tr}(Ub))\).
	By \eqref{eq:weil-w}, \eqref{eq:GL-disintegration}, and character
	orthogonality, local constancy and the reciprocal volumes for the
	self-dual measures give, for all sufficiently large \(R\),
	\begin{align*}
		&\int_{\Sym_{\nInLZorth}(\varpi^{-R}O_{\FLocal})}
		\Phi_{\varphi}(w_{\nInLZorth}^{-1}n(b),0)
		\psi\!\left(-\tfrac12\operatorname{tr}(Sb)\right)\,db\\
		&\quad=\gamma_{V,\psi}^{-\nInLZorth}
		\operatorname{vol}\bigl(\Sym_{\nInLZorth}(\varpi^{-R}O_{\FLocal})\bigr)
		\int_{S+\Sym_{\nInLZorth}(\varpi^{-R}O_{\FLocal})^{\perp}}
		\left(\int_{T^{-1}(U)}
		\varphi(\mathbf y)\,d\mu_U^{\mathrm{GL}}(\mathbf y)\right)
		\,dU^{\mathrm{sd}}\\
		&\quad=\gamma_{V,\psi}^{-\nInLZorth}
		\int_{T^{-1}(S)}
		\varphi(\mathbf y)\,d\mu_S^{\mathrm{GL}}(\mathbf y).
	\end{align*}
	By \eqref{eq:bad-finite-whittaker-lattice}, the left side equals
	\(W_{S/2}^{\psi}(1,0,\varphi)\); \eqref{eq:GL-torsor-map} and the
	definition of \(dh\) identify the last fiber integral with the
	stated orbital integral.
\end{proof}
\subsection{Local densities}\label{subsec:local-density}
In this subsection we assume that \(\FLocal\) has odd residue characteristic and that
\(\psi\) is unramified, i.e., has conductor \(O_\FLocal\). 
\begin{defn}
	Let \(L\) and \(M\) be \emph{integral} (i.e., the symmetric pairing takes values in $O_\FLocal$) non-degenerate quadratic
	\(O_\FLocal\)-lattices, with \(\operatorname{rk}L=\nInLZorth\) and
	\(\operatorname{rk}M\ge\nInLZorth\).  Following
	\cite[Definition~3.1.1]{LZ-orth}, let \(\Rep_{M,L}\) be the
	\(O_\FLocal\)-scheme characterized by
	\begin{equation*}
		\Rep_{M,L}(R)
		=\operatorname{QHom}(L\otimes_{O_\FLocal}R,M\otimes_{O_\FLocal}R)
	\end{equation*}
	for every \(O_\FLocal\)-algebra \(R\), where \(\operatorname{QHom}\) denotes
	the \(R\)-linear maps preserving the quadratic forms.  Its local
	representation density is
	\begin{equation*}
		\operatorname{Den}(M,L)
		:=\lim_{N\to\infty}
		q^{-N\nInLZorth(2\operatorname{rk}M-\nInLZorth-1)/2}
		\#\Rep_{M,L}(O_\FLocal/\varpi^N).
	\end{equation*}
\end{defn}

\begin{defn}\label{defn:local-density}
	For \(\epsilon\in\{\pm1\}\), let \(H_m^\epsilon\) be the self-dual lattice
	fixed in \Cref{sec:notat-quadr-latt-1}. Following
	\cite[Definitions~3.4.1--3.4.2]{LZ-orth}, define the normalizing
	polynomial in the odd-corank case by
	\begin{equation*}
		\operatorname{Nor}^\epsilon(X,L)
		=
		\left(1-\operatorname{sgn}(H_{\nInLZorth+1}^\epsilon)
		q^{-(\nInLZorth+1)/2}X\right)
		\prod_{1\le i<(\nInLZorth+1)/2}(1-q^{-2i}X^2).
	\end{equation*}
	The normalized local Siegel series is the polynomial
	\(\operatorname{Den}^\epsilon(X,L)\in\mathbb Z[X]\) determined by
	\begin{equation}\label{eq:normalized-density-interpolation}
		\operatorname{Den}^\epsilon(q^{-k},L)
		=
		\frac{\operatorname{Den}(H_{\nInLZorth+1+2k}^\epsilon,L)}
		{\operatorname{Nor}^\epsilon(q^{-k},L)},\qquad k\ge0.
	\end{equation}
	The derived local density is
	\(\pDen^\epsilon(L):=-\left.\frac{d}{dX}\operatorname{Den}^\epsilon(X,L)\right|_{X=1}\).
	If \(L\) is nonintegral, we define \(\operatorname{Den}^\epsilon(X,L)\) and hence
	\(\pDen^\epsilon(L)\) to be zero.
	
	For a nonsingular tuple \(\bu=(u_1,\ldots,u_{\nInLZorth})\) in a fixed
	quadratic \(\FLocal\)-space, write
	\(L(\bu):=\sum_i O_\FLocal u_i\) and set
	\(\pDen^\epsilon(\bu):=\pDen^\epsilon(L(\bu))\).
\end{defn}

\begin{prop}\label{prop:analytic-locally-constant}
	On the open subset of tuples \(\bu\) with \(\det T(\bu)\ne0\), the function
	\(\bu\mapsto\pDen^\epsilon(\bu)\) is locally constant.
\end{prop}
\begin{proof}
	For integral \(L(\bu)\), sufficiently small perturbations of its
	nondegenerate Gram matrix preserve its integral isometry class, hence its
	density polynomial and derivative; cf.\ \cite[Remark~3.4.3]{LZ-orth}.
	The nonintegral locus is open, and the derivative is zero there by definition.
\end{proof}

\begin{prop}[{\cite[Remark~3.4.9]{LZ-orth}}]\label{prop:pden-whittaker-ratio}
	Let $\Lambda$ be a self-dual lattice in $V$, and set
	\(\epsilon=\chi(\Lambda)\).
	Let \(L\) be a non-degenerate lattice of rank
	\(\nInLZorth=\dimofV-1\), and let \(T\) be the Gram matrix of a
	basis of \(L\). If $L$ is integral, then for every \(s\in\C\) we have
	\[
	W_{T/2}^{\psi}(1,s,\mathbf1_{\Lambda^{\nInLZorth}})
	=\operatorname{Den}^\epsilon(q^{-s},L)
	\operatorname{Nor}^\epsilon(q^{-s},L).
	\]
	If $V$ does not represent $T$,
	differentiation at \(s=0\) gives
	\begin{equation*}
		(W_{T/2}^{\psi})'(1,0,\mathbf1_{\Lambda^{\nInLZorth}})
		=\pDen^\epsilon(L)\operatorname{Nor}^\epsilon(1,L)\log q.
	\end{equation*}
	For nonintegral $L$, one has $W_{T/2}^{\psi}(1,s,\mathbf1_{\Lambda^{\nInLZorth}})=0$.
\end{prop}

\section{Global Weil representations and Siegel Eisenstein series}\label{sec:analytic}
\begin{defn}\label{defn:Weil-data}
	Retain the totally real field \(\FGlobal\), the dimension
	\(\dimofV=\nInLZorth+1\), and the additive character
	\(\psi:\FGlobal\backslash\BA_{\FGlobal}\to\BC^\times\)
	of \Cref{sec:notat-quadr-spac,subsec:nota-hilbert}.
	Let \(\Sigma_F\) be the set of places of \(\FGlobal\), with finite
	and real subsets \(\Sigma_{F,f}\) and \(\Sigma_{F,\infty}\).
	Set \(\Sigma_{F,f}^{\fkd}=\{v\in\Sigma_{F,f}:v\nmid\fkd\}\).
	Let \(\psi_{\BQ}\) be the character of \(\BQ\backslash\BA_{\BQ}\)
	with real component \(x\mapsto e^{2\pi i x}\).
	Define the global character
	\(\psi_{\mathrm{std}}=\psi_{\BQ}\circ\operatorname{Tr}_{\FGlobal/\BQ}\).
	There is a unique \(a_\psi\in\FGlobal^\times\) such that
	\begin{equation}\label{eq:additive-character-scalar}
		\psi(x)=\psi_{\mathrm{std}}(a_\psi x),
		\qquad \RealPlace(a_\psi)=\kappa_\RealPlace>0\quad(\RealPlace\mid\infty).
	\end{equation}
	For finite \(v\), the conductor of \(\psi_v\) is the annihilator of
	\(O_{\FGlobal_v}\) under \((x,y)\mapsto\psi_v(xy)\).  Define
	\[
	S_\psi=\{v\in\Sigma_{F,f}:\operatorname{cond}(\psi_v)\ne O_{\FGlobal_v}\}.
	\]
	The even integer \(\fkd\) fixed in \Cref{sec:GSpin} is divisible by
	every rational prime below a place of \(S_\psi\), and satisfies
	\((\fkd,p)=1\).
\end{defn}
\subsection{Global Weil representations}

\begin{defn}\label{defn:global-metaplectic-extension}
	For \(r\ge1\), let \(W_r=\FGlobal^r\oplus\FGlobal^r\) have the
	standard symplectic form of \Cref{defn:Weil-repn}.
	The adelic metaplectic double cover \(\Mp_{2r}(\BA_{\FGlobal})\)
	is constructed by
	\eqref{eq:adelic-metaplectic-cover} in genus \(r\), with
	\(\Mp_2\) and \(\SL_2\) replaced by \(\Mp_{2r}\) and \(\Sp_{2r}\).
	In the global Rao cocycle model its elements are \((g,t)\), with
	\(g\in\Sp_{2r}(\BA_{\FGlobal})\) and \(t\in\{\pm1\}\).
	Multiplication is given by the global Rao cocycle.
	The extension splits canonically over \(\Sp_{2r}(\FGlobal)\) and
	\(N_r(\BA_{\FGlobal})\); see
	\cite[p.~549]{Kudla1997CentralDerivatives}.
	Let \(\widetilde P_r(\BA_{\FGlobal})\) and
	\(\widetilde M_r(\BA_{\FGlobal})\) denote the inverse images of
	\(P_r(\BA_{\FGlobal})\) and \(M_r(\BA_{\FGlobal})\).
	Then \(\widetilde P_r(\BA_{\FGlobal})=
	\widetilde M_r(\BA_{\FGlobal})N_r(\BA_{\FGlobal})\), using the
	canonical splitting of the unipotent subgroup.
	We write \(m(a)=(m(a),1)\) for the chosen Levi lift.
	Let \(K^{\mathrm{std}}_{r,f}\) be the inverse image of
	\(\prod_{v<\infty}\Sp_{2r}(O_{\FGlobal_v})\) in
	\(\Mp_{2r}(\BA_{\FGlobal,f})\).
\end{defn}

\begin{defn}\label{defn:global-Weil-representation}
	Let \(\Vincoherent\) be the totally positive definite incoherent
	quadratic space over \(\BA_{\FGlobal}\) fixed in
	\Cref{sec:notat-quadr-spac}, of dimension \(\dimofV\).
	The Schwartz space and Weil actions on the double cover are
	\[
	\begin{gathered}
		\sS(\Vincoherent^r)=\bigotimes_v'\sS(\Vincoherent_v^r),\qquad
		\omega=\bigotimes_v'\omega_{\Vincoherent_v,\psi_v},\\
		\omega_f=\bigotimes_{v<\infty}'\omega_{\Vincoherent_v,\psi_v},\qquad
		\omega_\infty=\bigotimes_{\RealPlace\mid\infty}\omega_{\Vincoherent_\RealPlace,\psi_\RealPlace}.
	\end{gathered}
	\]
	At almost all finite places the reference vector is
	\(\mathbf1_{\Lambda_v^r}\), with \(\Lambda_v\subset\Vincoherent_v\)
	self-dual.  Write \(\Vincoherent_f\) for the finite adelic quadratic
	space and \(\Or(\Vincoherent)=\prod_v'\Or(\Vincoherent_v)\) for its
	adelic isometry group. By the local covering-kernel action in
	\Cref{defn:Weil-repn} and the quotient in
	\eqref{eq:adelic-metaplectic-cover}, \(\omega\) is a representation of
	\(\Mp_{2r}(\BA_{\FGlobal})\times\Or(\Vincoherent)\), genuine for odd
	\(\dimofV\) and descending for even \(\dimofV\) to
	\(\Sp_{2r}(\BA_{\FGlobal})\times\Or(\Vincoherent)\).
	
	Let \(\chi=\chi_{\Vincoherent}\) be the character of
	\(\BA_{\FGlobal}^\times/\FGlobal^\times\) defined in
	\Cref{sec:notat-quadr-spac}.
	Its extension to the Levi is
	\[
	\chi(m(a),t)=
	\begin{cases}
		t^{\dimofV}\chi(\det a),&\dimofV\text{ even},\\
		t^{\dimofV}\chi(\det a)\gamma(\det a,\psi)^{-1},
		&\dimofV\text{ odd},
	\end{cases}
	\]
	where \(t\in\{\pm1\}\); it extends trivially on
	\(N_r(\BA_{\FGlobal})\).
	The global indices are products of the local indices in \Cref{defn:Weil-repn}:
	\[
	\gamma(a,\psi)=\prod_v\gamma(a_v,\psi_v),\qquad
	\gamma(\psi)=\prod_v\gamma(\psi_v),\qquad
	\chi_\infty=\prod_{\RealPlace\mid\infty}\chi_\RealPlace.
	\]
	The Weil constant is $
	\gamma_{\Vincoherent,\psi}
	=\gamma(\det\Vincoherent,\psi)
	\gamma(\psi)^{\dimofV}\Has(\Vincoherent)$.  
	For \(\varphi\in\sS(\Vincoherent^r)\),
	\begin{align*}
		\omega(m(a))\varphi(\mathbf x)
		&=\chi(m(a))|\det a|_{\FGlobal}^{\dimofV/2}\varphi(\mathbf x a),\notag\\
		\omega(n(b))\varphi(\mathbf x)
		&=\psi\!\left(\frac12\operatorname{tr}(bT(\mathbf x))\right)\varphi(\mathbf x),
		\\
		\omega(w_r)\varphi&=\gamma_{\Vincoherent,\psi}^{\,r}\widehat\varphi,\notag\\
		\omega(h)\varphi(\mathbf x)&=\varphi(h^{-1}\mathbf x)
		\qquad(h\in\Or(\Vincoherent)).\notag
	\end{align*}
	The norm \(|\cdot|_{\FGlobal}\) is the product of the local normalized
	absolute values, and \(\widehat\varphi\) is the product Fourier transform
	\eqref{eq:weil-fourier-transform}.
\end{defn}

\begin{defn}\label{defn:Siegel-Fourier-Gaussian}
	For every \(r\ge1\), every \(\sz=\sx+i\sy\) in
	\[
	\mathbb H_r
	=\prod_{\RealPlace\mid\infty}
	\{z\in\Sym_r(\BC):\operatorname{Im}z>0\},
	\]
	and \(T\in\Sym_r(\FGlobal)\), define
	\begin{equation*}
		q^{T/2}(\sz)
		=\psi_\infty\!\left(\frac12\operatorname{tr}(T\sz)\right)
		=\prod_{\RealPlace\mid\infty}
		\exp\!\left(\pi i\kappa_\RealPlace\operatorname{tr}(T_\RealPlace\sz_\RealPlace)\right).
	\end{equation*}
	In genus \(\nInLZorth\), the archimedean Gaussian is
	\begin{equation}\label{eq:psi-adapted-gaussian}
		\varphi_{\infty,\psi}(\mathbf x)
		=\prod_{\RealPlace\mid\infty}
		\exp\!\left(-\pi\kappa_\RealPlace\operatorname{tr}T(\mathbf x_\RealPlace)\right).
	\end{equation}
	When only finite Schwartz data are specified, the archimedean component is $\varphi_{\infty,\psi}$.
\end{defn}

\subsection{Siegel Eisenstein series}\label{subsec:Siegel-Eisenstein}

\begin{defn}\label{defn:Siegel-Eisenstein}
	In genus \(\nInLZorth=\dimofV-1\), define
	\begin{equation}\label{eq:degenerate-principal-series}
		I_{\nInLZorth}(s,\chi)
		=\operatorname{Ind}_{\widetilde P_{\nInLZorth}(\BA_{\FGlobal})}
		^{\Mp_{2\nInLZorth}(\BA_{\FGlobal})}
		\bigl(\chi|\det|_{\FGlobal}^{s+(\nInLZorth+1)/2}\bigr).
	\end{equation}
	The induction is unnormalized.  A standard section has restriction
	to the standard maximal compact subgroup independent of \(s\).
	If \(\chi_v\) denotes the local inducing character of
	\Cref{defn:Weil-repn}, this space is the restricted tensor product
	\(\bigotimes_v'I_{\nInLZorth}(s,\chi_v)\) of the local principal
	series in \Cref{defn:local-principal-series}.
	Its covering-kernel character is \(t\mapsto t^{\dimofV}\).
	For a standard section \(\Phi(-,s)\), set
	\[
	E(g,s,\Phi)
	=\sum_{\gamma\in P_{\nInLZorth}(\FGlobal)\backslash
		\Sp_{2\nInLZorth}(\FGlobal)}\Phi(\gamma g,s),
	\qquad g\in\Mp_{2\nInLZorth}(\BA_{\FGlobal}),
	\]
	where the rational groups use the canonical splitting.
	The sum converges for \(\Re(s)\gg0\) and admits meromorphic
	continuation to \(\BC\); see \cite[\S11.1]{LZ-orth}.
	
	For \(\varphi\in\sS(\Vincoherent^{\nInLZorth})\), define
	\(\Phi_\varphi(g,0)=(\omega(g)\varphi)(0)\).
	Its standard extension is the restricted tensor product of the local
	sections in \Cref{defn:SW-section}, extended linearly in \(\varphi\).
	For global and finite Schwartz data, respectively, define
	\[
	E(g,s,\varphi)=E(g,s,\Phi_\varphi),\qquad
	E(g,s,\varphi_f)=E(g,s,\Phi_{\varphi_f\otimes\varphi_{\infty,\psi}}),
	\quad\varphi_f\in\sS(\Vincoherent_f^{\nInLZorth}).
	\]
	The Siegel--Weil Eisenstein series is regular at \(s=0\); see
	\cite[Appendix~II]{Langlands1976FunctionalEquations} and
	\cite[proof of Proposition~7.2]{Ikeda1996Residue}.
	By incoherence \cite[\S11.3]{LZ-orth}, \(E(g,0,\varphi)=0\). We set
	\[
	\pEis(g,\varphi)=\left.\frac{d}{ds}E(g,s,\varphi)\right|_{s=0}.
	\]
	Eisenstein summation gives
	\(\pEis(\gamma g,\varphi)=\pEis(g,\varphi)\) for
	\(\gamma\in\Sp_{2\nInLZorth}(\FGlobal)\).
	
\end{defn}

\subsection{Fourier coefficients and central derivatives}\label{sec:four-coeff-deriv}

We have a Fourier expansion \[
E(g,s,\Phi)=\sum_{T\in\Sym_{\nInLZorth}(\FGlobal)}E_{T/2}(g,s,\Phi),
\]
where 
\begin{equation}\label{eq:eisenstein-fourier-coefficient}
	E_{T/2}(g,s,\Phi)
	=\int_{\Sym_{\nInLZorth}(\FGlobal)\backslash
		\Sym_{\nInLZorth}(\BA_{\FGlobal})}
	E(n(b)g,s,\Phi)\psi\!\left(-\frac12\operatorname{tr}(Tb)\right)\,db.
\end{equation}
The measure \(db\) is induced on the quotient by the product, over all
places \(v\) of \(\FGlobal\), of the local self-dual measures for the pairings
\((b_1,b_2)\mapsto\psi_v(\tfrac12\operatorname{tr}(b_1b_2))\)
on \(\Sym_{\nInLZorth}(\FGlobal_v)\), as in \Cref{defn:local-Whitt}.
For nonsingular \(T\) and factorizable \(\Phi=\bigotimes_v\Phi_v\) we have
\[
E_{T/2}(g,s,\Phi)=\prod_v W_{T/2,v}(g_v,s,\Phi_v).
\]
Here \(W_{T/2,v}\) denotes \(W_{T_v/2}^{\psi_v}\) from
\Cref{defn:local-Whitt}, and \((g_v)_v\) is a representative of
\(g\) in the restricted product defining the adelic double cover.
At almost all finite places, \(g_v\) fixes the unramified section and
\(T_v\) is the Gram matrix of a self-dual lattice. By
\Cref{prop:pden-whittaker-ratio,defn:local-density} and
\cite[Theorem~3.4.5]{LZ-orth},
\[
W_{T/2,v}(g_v,s,\Phi_v)
=\left(1-\operatorname{sgn}(\Lambda_v)q_v^{-\dimofV/2-s}\right)
\prod_{1\le i<\dimofV/2}(1-q_v^{-2i-2s}).
\]
Since \(\dimofV\ge3\), these factors are \(1+O(q_v^{-3/2})\), with
derivatives \(O((\log q_v)q_v^{-3/2})\), uniformly for
\(\lvert\operatorname{Re}(s)\rvert\le1/4\). Thus the unramified Euler
product and its derivative converge locally uniformly near zero.
Differentiating the factorization gives
\begin{equation}\label{eq:eissum}
	E'_{T/2}(g,0,\Phi)=\sum_{v\in\Sigma_F}E'_{T/2,v}(g,0,\Phi),
\end{equation}
where \begin{equation}\label{eq:eisenfactor}
	E'_{T/2,v}(g,s,\Phi)
	=W'_{T/2,v}(g_v,s,\Phi_v)
	\prod_{v'\ne v}W_{T/2,v'}(g_{v'},s,\Phi_{v'}).
\end{equation}              
For a Siegel--Weil section set
\[
\begin{aligned}
	\pEis_{T/2}(g,\varphi)=E'_{T/2}(g,0,\Phi_\varphi),\qquad 
	\pEis_{T/2,v}(g,\varphi)=E'_{T/2,v}(g,0,\Phi_\varphi).
\end{aligned}
\]
The second definition applies to factorizable data and extends linearly.
\begin{defn}\label{defn:classic-Eis}
	For \(\sz=\sx+i\sy\in\mathbb H_{\nInLZorth}\), set
	\(g_\sz=n(\sx)m(\sy^{1/2})1_f\), using the positive definite
	square root at each real place; \(1_f\) is the finite identity.
	For \(T\in\Sym_{\nInLZorth}(\FGlobal)\), define
	\begin{equation}\label{eq:classical-eisenstein-normalization}
		\begin{aligned}
			\pEis_{T/2}(\sz,\varphi)
			&=\left(\prod_{\RealPlace\mid\infty}(\det\sy_\RealPlace)^{-\dimofV/4}\right)
			\pEis_{T/2}(g_\sz,\varphi),\\
			\pEis_{T/2,v}(\sz,\varphi)
			&=\left(\prod_{\RealPlace\mid\infty}(\det\sy_\RealPlace)^{-\dimofV/4}\right)
			\pEis_{T/2,v}(g_\sz,\varphi).
		\end{aligned}
	\end{equation}
	The second definition applies to nonsingular \(T\) and \(v\in\Sigma_F\). 
\end{defn}

\subsection{Nonsingular indices and the difference set}\label{rem:DiffTV}

\begin{defn}\label{defn:Eisenstein-difference-set}
	For nonsingular \(T\in\Sym_{\nInLZorth}(\FGlobal)\), let
	\(\Has(T_v)\) be the Hasse invariant of \(\FGlobal_v^{\nInLZorth}\)
	with Gram matrix \(T_v\).  Recall
	\[
	\Diff(T,\Vincoherent)
	=\{v:\text{there is no }\mathbf x\in\Vincoherent_v^{\nInLZorth}
	\text{ with }T(\mathbf x)=T_v\}.
	\]
\end{defn}

At a finite place \(v\), \cite[Lemma~3.7.1]{LZ-orth} gives
\begin{equation}\label{eq:diff}
	v\in\Diff(T,\Vincoherent)
	\quad\Longleftrightarrow\quad
	\Has(\Vincoherent_v)
	\ne\Has(T_v)\bigl(\det(T)_v,
	\det\Vincoherent_v\det(T)_v\bigr)_{\FGlobal_v}.
\end{equation}
At a real place \(\RealPlace\), the space \(\Vincoherent_\RealPlace\) is positive definite, so
\(\RealPlace\in\Diff(T,\Vincoherent)\) if and only if \(T_\RealPlace\) fails to be positive
definite.  
By \cite[Proposition~1.4]{Kudla1997CentralDerivatives},
\(W_{T/2,v}(g_v,0,\varphi_v)\ne0\) only when
\(v\notin\Diff(T,\Vincoherent)\).  Therefore
\(\pEis_{T/2,v}\ne0\) only if \(\Diff(T,\Vincoherent)=\{v\}\). In particular, if the difference
set has at least two places, then the nonsingular central derivative vanishes.

\begin{lem}\label{lem:block-positive-Diff-parity}
	For every totally positive definite
	\(T\in\Sym_{\nInLZorth}(\FGlobal)\), the set
	\(\Diff(T,\Vincoherent)\) has odd cardinality.
\end{lem}
\begin{proof}
	At every real place, \(T\) is represented by \(\Vincoherent\).
	At finite places, multiply the Hasse invariant criterion
	\eqref{eq:diff}.  The product formulas give
	\(\prod_v\Has(T_v)=1\) and
	\(\prod_v(\det T_v,\det\Vincoherent_v\det T_v)_{\FGlobal_v}=1\),
	whereas \(\prod_v\Has(\Vincoherent_v)=-1\).
	Consequently \((-1)^{\#\Diff(T,\Vincoherent)}=-1\).
\end{proof}

\subsection{Singular coefficients}
\label{subsec:singular-coefficients}
\begin{thm}
	\label{thm:two-place-singular}
	Assume \(\dimofV\ge3\).
	Let \(\varphi=\bigotimes_v\varphi_v\in
	\sS(\Vincoherent^{\nInLZorth})\) be factorizable, with
	standard Siegel--Weil section \(\Phi_\varphi(-,s)\) in
	\eqref{eq:degenerate-principal-series}, and take
	\(\varphi_\infty=\varphi_{\infty,\psi}\) as in
	\eqref{eq:psi-adapted-gaussian}.  Suppose that two distinct finite
	places \(v_1,v_2\) satisfy
	\[
	\supp(\varphi_{v_i})\subset
	\{\mathbf x\in\Vincoherent_{v_i}^{\nInLZorth}:\det T(\mathbf x)\ne0\},
	\qquad i=1,2.
	\]
	Then, for every singular matrix
	\(S\in\Sym_{\nInLZorth}(\FGlobal)\)
	and every \(g_\infty\in \Mp_{2\nInLZorth}(\FGlobal_\infty)\), one has
	\begin{equation}\label{eq:orthogonal-Reg}
		\operatorname{ord}_{s=0}E_{S/2}(g_\infty1_f,s,\varphi)\ge2.
	\end{equation}
	In particular, \(E'_{S/2}(g_\infty1_f,0,\varphi)=0\).
\end{thm}

\begin{proof}
	Let \(e_s\in I_{\nInLZorth}(s,\chi)\) be a holomorphic section family
	with fixed finite level and finitely many archimedean compact types.
	Then \(E(g,s,e_s)\) is holomorphic at \(s=0\); see
	\cite[Appendix~II]{Langlands1976FunctionalEquations} and
	\cite[proof of Proposition~7.2]{Ikeda1996Residue}.
	The Fourier coefficients of \(E(g,s,e_s)\) are likewise holomorphic, since
	\(N(\FGlobal)\backslash N(\BA_{\FGlobal})\) is compact. 
	
	Fix \(S\) and \(g_\infty\) as in the statement.  For a global section
	\(e_s\), a finite place \(v\), and \(B\in\Sym_{\nInLZorth}(\FGlobal_v)\),
	translation in \eqref{eq:eisenstein-fourier-coefficient} gives
	\begin{equation}\label{eq:singular-Fourier-equivariance}
		E_{S/2}\bigl(g_\infty1_f,s,e_s(\,\cdot\,n(B))\bigr)
		=\psi_v\!\left(\tfrac12\operatorname{tr}(SB)\right)
		E_{S/2}(g_\infty1_f,s,e_s).
	\end{equation}
	Here \(n(B)\) acts at \(v\) and has identity components at all other
	places.
	The identity holds first in the convergence region and then
	meromorphically in \(s\).
	
	For \(v\in\{v_1,v_2\}\), the compact set
	\(T(\supp\varphi_v)\) does not contain \(S\).
	Choose a compact open additive subgroup
	\(H_v\subset\Sym_{\nInLZorth}(\FGlobal_v)\) whose annihilator
	\[
	H_v^\perp:=\left\{Y\in\Sym_{\nInLZorth}(\FGlobal_v):
	\psi_v\!\left(\tfrac12\operatorname{tr}(YB)\right)=1
	\text{ for every }B\in H_v\right\}
	\]
	satisfies \((S+H_v^\perp)\cap T(\supp\varphi_v)=\varnothing\).
	Give \(H_v\) probability Haar measure \(dB\), and define the averaging
	operator on sections by
	\[
	\mathcal P_{v,s}e_s:=\int_{H_v}
	\psi_v\!\left(-\tfrac12\operatorname{tr}(SB)\right)
	e_s(\,\cdot\,n(B))\,dB.
	\]
	By \eqref{eq:singular-Fourier-equivariance},
	\begin{equation}\label{eq:singular-average-preserves-coefficient}
		E_{S/2}(g_\infty1_f,s,\mathcal P_{v,s}e_s)
		=E_{S/2}(g_\infty1_f,s,e_s).
	\end{equation}
	
	At \(s=0\), the Siegel--Weil map is equivariant.  Thus character
	orthogonality gives
	\begin{align*}
		\mathcal P_{v,0}\Phi_{\varphi_v}(\,\cdot\,,0)
		&=\Phi_{\int_{H_v}
			\psi_v(-\tfrac12\operatorname{tr}(SB))
			\omega_v(n(B))\varphi_v\,dB}(\,\cdot\,,0)\\
		&=\Phi_{\mathbf1_{S+H_v^\perp}(T(\,\cdot\,))\varphi_v}
		(\,\cdot\,,0)=0.
	\end{align*}
	At a sufficiently small fixed level, the averaged family
	\(\mathcal P_{v,s}\Phi_{\varphi_v}(\,\cdot\,,s)\) is a finite sum of
	right translates in the compact picture, hence is holomorphic.
	Since \(\mathcal P_{v,0}\Phi_{\varphi_v}(\,\cdot\,,0)=0\), we can write
	\[
	\mathcal P_{v,s}\Phi_{\varphi_v}(\,\cdot\,,s)=s\,h_{v,s}
	\]
	for a holomorphic local section \(h_{v,s}\) at the same level.
	
	The two averaging operators act on different tensor factors, so
	\eqref{eq:singular-average-preserves-coefficient} gives
	\begin{align*}
		E_{S/2}(g_\infty1_f,s,\varphi)
		&=E_{S/2}\!\left(g_\infty1_f,s,
		\mathcal P_{v_1,s}\mathcal P_{v_2,s}\Phi_\varphi(\,\cdot\,,s)\right)\\
		&\quad=s^2E_{S/2}\!\left(g_\infty1_f,s,
		h_{v_1,s}\otimes h_{v_2,s}
		\otimes\bigotimes_{v\ne v_1,v_2}\Phi_{\varphi_v}(\,\cdot\,,s)\right).
	\end{align*}
	The last coefficient is holomorphic at zero by the central
	holomorphy recalled above, proving \eqref{eq:orthogonal-Reg}.
\end{proof}

\subsection{Convention about global Haar measures}\label{conv:haar-global-constant}
We fix our notation about the Haar measure on $\prod_{w<\infty}'\SO(\Vincoherent_w)(\FGlobal_w)$ in this subsection. 
At every finite place, let \(dh_v\) be the measure of
\Cref{defn:Gelfand--Leray} for \(V=\Vincoherent_v\), and define
\[
\vol^{\mathrm{GL}}_v(K_{2,v}):=\int_{K_{2,v}}dh_v,\qquad
\vol^{\mathrm{GL}}(K_2)=\int_{K_2}\prod_{v<\infty}dh_v.
\]
The product converges: \eqref{eq:GL-hyperspecial-volume-product} gives compact
volume \(1+O(q_v^{-2})\) at almost every finite place, since
\(\dimofV\ge3\) in the global setting.  
\begin{defn}\label{defn:Haar-measure}
	For a finite set \(S\) of finite places, the Haar measure
	\(dg_f^S\) on the restricted product
	\(\prod_{w<\infty,\,w\notin S}'
	\SO(\Vincoherent_w)(\FGlobal_w)\) is defined by
	$
	dg_f^S=(a_{\infty,\psi}^{\operatorname{tr}})^{-1}
	\prod_{\substack{w<\infty\\w\notin S}}dh_w$,
	where the positive constant \(a_{\infty,\psi}^{\operatorname{tr}}\)
	is defined in \Cref{lem:global-Haar-phase-a}.
	For a compact subset \(K_2^S\subseteq \prod_{w<\infty,\,w\notin S}'
	\SO(\Vincoherent_w)(\FGlobal_w)\) we denote its volume by 
	\(\vol(K_2^S):=\int_{K_2^S}dg_f^S\).  In particular,
	\[
	\vol(K_2)=(a_{\infty,\psi}^{\operatorname{tr}})^{-1}
	\prod_{v<\infty}\vol^{\mathrm{GL}}_v(K_{2,v}).
	\]
\end{defn}

\begin{prop}\label{prop:archimedean-whittaker-factor}\label{lem:global-Haar-phase}
	For every totally positive definite
	\(T\in\Sym_{\nInLZorth}(\FGlobal)\) and
	\(\sz=\sx+i\sy\in\mathbb H_{\nInLZorth}\),
	\begin{equation}\label{eq:archimedean-whittaker-factor}
		q^{T/2}(\sz)
		=c_{\infty,\psi}^{\operatorname{tr}}\,
		\prod_{\RealPlace\mid\infty}
		\left((\det\sy_\RealPlace)^{-\dimofV/4}
		W_{T/2,\RealPlace}(g_{\sz,\RealPlace},0,\varphi_{\RealPlace,\psi})\right),
	\end{equation}
	with
	\begin{equation}\label{eq:archimedean-whittaker-global-constant}
		c_{\infty,\psi}^{\operatorname{tr}}
		=\left(
		\frac{\Gamma_{\nInLZorth}((\nInLZorth+1)/2)}
		{2^{\nInLZorth/2}
			\pi^{\nInLZorth(\nInLZorth+1)/2}}
		\right)^{[\FGlobal:\BQ]}
		\prod_{\RealPlace\mid\infty}(\gamma_{\Vincoherent_\RealPlace,\psi_\RealPlace}^{\nInLZorth}
		\kappa_\RealPlace^{-\nInLZorth(\nInLZorth+1)/4}).
	\end{equation}
\end{prop}
\begin{proof}
	Follows from
	\Cref{prop:archimedean-whittaker-factor-0} by taking products over
	the real places. 
\end{proof}
\begin{defn}\label{lem:global-Haar-phase-a}
	Set $a_{\infty,\psi}^{\operatorname{tr}}
	:=\left|c_{\infty,\psi}^{\operatorname{tr}}\right|
	=
	\left(
	\frac{\Gamma_{\nInLZorth}(\dimofV/2)}
	{2^{\nInLZorth/2}\pi^{\nInLZorth\dimofV/2}}
	\right)^{[\FGlobal:\BQ]}
	\prod_{\RealPlace\mid\infty}
	\kappa_\RealPlace^{-\nInLZorth\dimofV/4}$. Then $\prod_v\gamma_{\Vincoherent_v,\psi_v}=-1$ implies
	\begin{equation}\label{eq:Weil-phase-product-exact}
		a_{\infty,\psi}^{\operatorname{tr}}=(-1)^{\nInLZorth}
		c_{\infty,\psi}^{\operatorname{tr}}
		\prod_{v<\infty}
		\gamma_{\Vincoherent_v,\psi_v}^{\nInLZorth}.
	\end{equation}
\end{defn}

\begin{defn}\label{defn:matrix-local-density}
	For a finite place \(v\nmid2\), a sign \(\epsilon_v\in\{\pm1\}\), and a nonsingular
	\(T\in\Sym_{\nInLZorth}(\FGlobal_v)\), let \(L_{T,v}\) denote
	\(O_{\FGlobal_v}^{\nInLZorth}\) endowed with Gram matrix \(T\), and
	write
	\(\pDen_v^{\epsilon_v}(T):=\pDen^{\epsilon_v}(L_{T,v})\).
\end{defn}

\begin{thm}
	\label{thm:Siegel-Eisenstein}
	Let \(\varphi_{K_2}\in\sS(\Vincoherent_f^{\nInLZorth})^{K_2}\) be factorizable,
	put \(\varphi=\varphi_{K_2}\otimes\varphi_{\infty,\psi}\), and let
	\(T\in\Sym_{\nInLZorth}(\FGlobal)\) be a nonsingular Gram
	matrix with \(\Diff(T,\Vincoherent)=\{v\}\).
	
	Suppose first that \(v<\infty\). Let \(g=g_\infty1_f\), where
	\(g_\infty=(g_\RealPlace)_{\RealPlace\mid\infty}\in\Mp_{2\nInLZorth}(\FGlobal_\infty)\)
	is arbitrary. With
	\(\operatorname{Orb}_T(\varphi_{K_2}^v)\) denoting the finite-adelic orbital integral of
	\Cref{thm:global-int}, one has
	\begin{equation}\label{eq:orth-Luo-10-1-2-finite-0}
		\begin{aligned}
			(-1)^{\nInLZorth}\pEis_{T/2,v}(g,\varphi)
			={}&c_{\infty,\psi}^{\operatorname{tr}}
			\gamma_{\Vincoherent_v,\psi_v}^{\nInLZorth}\, W'_{T/2,v}(1,0,\varphi_v)\,
			\operatorname{Orb}_T(\varphi_{K_2}^v)\prod_{\RealPlace\mid\infty}W_{T/2,\RealPlace}(g_\RealPlace,0,\varphi_{\RealPlace,\psi}).
		\end{aligned}
	\end{equation}
	Here \(\gamma_{\Vincoherent_v,\psi_v}\) is the local Weil index of
	\(\Vincoherent_v\) with respect to \(\psi_v\).
	For every \(\sz\in\mathbb H_{\nInLZorth}\),
	\begin{equation}\label{eq:finite-fourier-normalized}
		(-1)^{\nInLZorth}q^{-T/2}(\sz)\pEis_{T/2,v}(\sz,\varphi)
		=\gamma_{\Vincoherent_v,\psi_v}^{\nInLZorth}\, W'_{T/2,v}(1,0,\varphi_v)\,
		\operatorname{Orb}_T(\varphi_{K_2}^v).
	\end{equation}
	Suppose moreover that \(v\nmid2\), that
	\(v\notin S_\psi\), and that
	\(\Lambda_v\subset\Vincoherent_v\) is self-dual with
	\(K_{2,v}=\SO(\Lambda_v)(O_{\FGlobal_v})\) and
	\(\varphi_v=\mathbf1_{\Lambda_v^{\nInLZorth}}\).  Put
	\(\epsilon_v:=\chi(\Lambda_v)=\chi(\Vincoherent_v)\).  One has
	\begin{equation}\label{eq:orth-Luo-10-1-2-finite}
		\begin{aligned}
			\frac{(-1)^{\nInLZorth}}{\log q_v}\,
			\pEis_{T/2,v}(g,\varphi)
			={}&c_{\infty,\psi}^{\operatorname{tr}}\vol^{\mathrm{GL}}_v(K_{2,v})\,\pDen_v^{\epsilon_v}(T)\,
			\operatorname{Orb}_T(\varphi_{K_2}^v)\prod_{\RealPlace\mid\infty}W_{T/2,\RealPlace}(g_\RealPlace,0,\varphi_{\RealPlace,\psi}).
		\end{aligned}
	\end{equation}
	
	If \(v=\RealPlace\in\Sigma_{F,\infty}\), then
	\begin{equation}\label{eq:orth-Luo-10-1-2-infty}
		(-1)^{\nInLZorth}q^{-T/2}(\sz)
		\pEis_{T/2,\RealPlace}(\sz,\varphi)
		=\frac{1}{\operatorname{Vol}^{\mathrm{GS}}(\SOShim\times_{\FGlobal,\RealPlace}\C,\Omega_{\mathcal E})}
		\int_{[\SOShim\times_{\FGlobal,\RealPlace}\C]}g_\RealPlace^{\mathrm K}(T,\sy,\varphi_{K_2}).
	\end{equation}
	The complex fiber is the base change of \(\SOShim\) along the real
	embedding \(\RealPlace\), as in \Cref{subsec:changing-iota}.  The current
	\(g_\RealPlace^{\mathrm K}(T,\sy,\varphi_{K_2})\) is defined in
	\eqref{eq:higher-Kudla-Green-current}, using the convention of
	\cite[Theorem~1.1.2]{garcia_sankaran_2019_green_forms}.
	
	For any nonsingular \(T\in\Sym_{\nInLZorth}(\FGlobal)\) with
	\(\#\Diff(T,\Vincoherent)\ge2\), one has
	\(\pEis_{T/2}(\sz,\varphi)=0\).
\end{thm}

\begin{proof}
	By \eqref{eq:eissum} and the central Whittaker vanishing at difference
	places \cite[Proposition~1.4]{Kudla1997CentralDerivatives},
	\(\pEis_{T/2}=\pEis_{T/2,v}\) for a singleton difference set, and
	\(\pEis_{T/2}=0\) when there are at least two difference places.
	
	Let \(v<\infty\).  
	For each finite \(w\ne v\), choose
	\(\mathbf x_w\in\Vincoherent_w^{\nInLZorth}\) with Gram matrix
	\(T\), and define
	\[
	\operatorname{Orb}_{T,w}(\varphi_w)
	=\int_{\SO(\Vincoherent_w)(\FGlobal_w)}
	\varphi_w(h^{-1}\mathbf x_w)\,dh_w.
	\]
	The exact local comparison \Cref{lem:local-Whittaker-GL-exact} gives
	\[
	W_{T/2,w}(1,0,\varphi_w)
	=\gamma_{\Vincoherent_w,\psi_w}^{-\nInLZorth}
	\operatorname{Orb}_{T,w}(\varphi_w).
	\]
	The definition of \(dg_f^v\) in
	\Cref{conv:haar-global-constant} gives
	$
	\prod_{\substack{w<\infty\\w\ne v}}
	\operatorname{Orb}_{T,w}(\varphi_w)
	=a_{\infty,\psi}^{\operatorname{tr}}
	\operatorname{Orb}_T(\varphi_{K_2}^v).
	$ 
	By \eqref{eq:Weil-phase-product-exact},
	\begin{equation*}
		\begin{aligned}
			\frac{a_{\infty,\psi}^{\operatorname{tr}}}
			{c_{\infty,\psi}^{\operatorname{tr}}}
			\prod_{\substack{w<\infty\\w\ne v}}
			\gamma_{\Vincoherent_w,\psi_w}^{-\nInLZorth}
			&=(-1)^{\nInLZorth}
			\gamma_{\Vincoherent_v,\psi_v}^{\nInLZorth}.
		\end{aligned}
	\end{equation*}
	Substituting these identities into \eqref{eq:eisenfactor} gives
	\eqref{eq:orth-Luo-10-1-2-finite-0}.
	Taking \(g=g_\sz\) and applying
	\eqref{eq:classical-eisenstein-normalization} and
	\eqref{eq:archimedean-whittaker-factor} gives \eqref{eq:finite-fourier-normalized}.
	Under the additional self-dual hyperspecial hypotheses, the
	unramified Weil index formula gives
	\(\gamma_{\Vincoherent_v,\psi_v}=1\).
	For nonintegral \(L_{T,v}\), both sides of
	\eqref{eq:orth-Luo-10-1-2-finite} vanish by
	\Cref{prop:pden-whittaker-ratio}. For integral \(L_{T,v}\),
	\Cref{lem:GL-hyperspecial-volume} gives
	\(\operatorname{Nor}^{\epsilon_v}(1,L_{T,v})=\vol^{\mathrm{GL}}_v(K_{2,v})\)
	and \Cref{prop:pden-whittaker-ratio} gives
	\begin{equation*}
		W'_{T/2,v}(1,0,\varphi_v)
		=\pDen_v^{\epsilon_v}(T)\,
		\operatorname{Nor}^{\epsilon_v}(1,L_{T,v})\,\log q_v,
	\end{equation*}
	so \eqref{eq:orth-Luo-10-1-2-finite} follows. 
	For an archimedean difference place,
	\eqref{eq:orth-Luo-10-1-2-infty} will be proved below in
	\Cref{thm:garcia-sankaran-archimedean}.
\end{proof}

\subsection{Associated Fourier--Jacobi series}
\label{subsec:diagonal-FJ}
We define the genus-one Fourier--Jacobi series associated to \(\pEis\) by the adelic operation below. The factor
\(2(-1)^{\nInLZorth}/\vol(K_2)\) matches the intersection normalization in
\Cref{thm:global-int}.
\begin{defn}\label{defn:Fourier-Jacobi-series}
	Fix lower-block data
	\[
	T_\flat\in\Sym_{\nInLZorth-1}(\FGlobal)_{>0},\qquad
	Y_\flat\in\prod_{\mu\mid\infty}\Sym_{\nInLZorth-1}(\BR)_{>0}.
	\]
	Let \(\phi_1\in\sS(\Vincoherent_f)^{K_2}\) and
	\(\phi_\flat\in\sS(\Vincoherent_f^{\nInLZorth-1})^{K_2}\).
	Assume \(\dimofV\ge3\).  Write
	\[
	\begin{gathered}
		\tau=u+iy\in\mathbb H_{\FGlobal},\quad
		\sz=\operatorname{diag}(\tau,iY_\flat),\quad
		\sy=\operatorname{diag}(y,Y_\flat),\\
		T_b(\xi)=\begin{pmatrix}\xi&b\\{}^tb&T_\flat\end{pmatrix},
		\quad \xi\in\FGlobal,\quad b\in\Mat_{1,\nInLZorth-1}(\FGlobal).
	\end{gathered}
	\]
	Let \(i_1:\Mp_2(\BA_{\FGlobal})\to
	\Mp_{2\nInLZorth}(\BA_{\FGlobal})\) be the lift of the
	symplectic embedding acting on the first symplectic plane and fixing
	its orthogonal complement. This map sends \((1,t)\) to \((1,t)\)
	for \(t\in\{\pm1\}\). Give
	\(m(\operatorname{diag}(1,Y_\flat^{1/2}))\) identity finite
	components, using the positive definite square root at each real place. For
	\(g\in\Mp_2(\BA_{\FGlobal})\), define
	\begin{equation}\label{eq:FJ-adelic-operation}
		\begin{aligned}
			&\FJ_{T_\flat}(g,Y_\flat;\phi_1\otimes\phi_\flat):=2\frac{(-1)^{\nInLZorth}}{\vol(K_2)}q^{-T_\flat/2}(iY_\flat)
			\prod_{\RealPlace\mid\infty}(\det Y_{\flat,\RealPlace})^{-\dimofV/4}\\
			&\qquad\cdot\int_{\Sym_{\nInLZorth-1}(\FGlobal)\backslash
				\Sym_{\nInLZorth-1}(\BA_{\FGlobal})}
			\pEis\!\left(n(\operatorname{diag}(0,B))i_1(g)
			m(\operatorname{diag}(1,Y_\flat^{1/2})),\phi_1\otimes\phi_\flat\right)\cdot
			\psi\!\left(-\tfrac12\operatorname{tr}(T_\flat B)\right)\,dB.
		\end{aligned}
	\end{equation}
	The Haar measure \(dB\) is normalized so that $\vol(\Sym_{\nInLZorth-1}(\FGlobal)\backslash
	\Sym_{\nInLZorth-1}(\BA_{\FGlobal}))=1$. 
	The function belongs to \(\CA_\infty(\SL_2(\FGlobal)\backslash\Mp_2(\BA_{\FGlobal}),\dimofV/2)\). 
	For \(\tau=u+iy\in\mathbb H_{\FGlobal}\), define its classical form by
	\[
	\begin{aligned}
		\FJ_{T_\flat}(\tau,Y_\flat;\phi_1\otimes\phi_\flat)
		=&\prod_{\RealPlace\mid\infty}y_\RealPlace^{-\dimofV/4} \cdot\FJ_{T_\flat}(n(u)m(y^{1/2})1_f,Y_\flat;
		\phi_1\otimes\phi_\flat).
	\end{aligned}
	\]
	Fourier expansion in \eqref{eq:FJ-adelic-operation} gives
	\begin{equation*}
		\begin{aligned}
			\FJ_{T_\flat}(\tau,Y_\flat;\phi_1\otimes\phi_\flat)
			=2\frac{(-1)^{\nInLZorth}}{\vol(K_2)}q^{-T_\flat/2}(iY_\flat)\cdot\sum_{\xi\in\FGlobal}
			\sum_{b\in\Mat_{1,\nInLZorth-1}(\FGlobal)}
			\pEis_{T_b(\xi)/2}(\sz,\phi_1\otimes\phi_\flat).
		\end{aligned}
	\end{equation*}
\end{defn}
\begin{defn}\label{defn:Fourier-Jacobi-coefficients}
	Assume  \(\phi_1\otimes\phi_\flat\) is
	factorizable and non-degenerate at two distinct finite places, as in
	\eqref{eq:Schwartz-nondegenerate}. 
	For \(v\in\Sigma_F\) and \(\xi\in\FGlobal\), define
	\begin{equation}\label{eq:FJ-local-coefficient}
		\FJ_{T_\flat,v,\xi}(y,Y_\flat;\phi_1\otimes\phi_\flat)
		:=2\frac{(-1)^{\nInLZorth}}{\vol(K_2)} 
		\sum_{\substack{b\in\Mat_{1,\nInLZorth-1}(\FGlobal)\\
				\det T_b(\xi)\ne0\\
				\Diff(T_b(\xi),\Vincoherent)=\{v\}}}
		\pEis_{T_b(\xi)/2,v}(i\sy,\phi_1\otimes\phi_\flat)
		q^{-T_b(\xi)/2}(i\sy).
	\end{equation}
	For \(\sz=\operatorname{diag}(\tau,iY_\flat)\), the unipotent
	transformation law removes the real part after division by
	\(q^{T_b(\xi)/2}(\sz)=q^{\xi/2}(\tau)q^{T_\flat/2}(iY_\flat)\).
	Thus by \eqref{eq:eissum} and \Cref{thm:two-place-singular} we have 
	\begin{equation}\label{eq:FJ-total-coefficient}
		\begin{aligned}
			[\FJ_{T_\flat}(\tau,Y_\flat;\phi_1\otimes\phi_\flat)]_{q^{\xi/2}}
			&=\sum_{v\in\Sigma_F}
			\FJ_{T_\flat,v,\xi}(y,Y_\flat;\phi_1\otimes\phi_\flat).
		\end{aligned}
	\end{equation}
\end{defn}
We check that for finite places \eqref{eq:FJ-local-coefficient} is a finite sum, and for archimedean places it converges absolutely and locally uniformly. Moreover, \eqref{eq:FJ-total-coefficient} is a finite sum.
\begin{prop}\label{prop:finite-contributing-indices}
	Fix \(\xi\in\FGlobal\),
	\(T_\flat\in\Sym_{\nInLZorth-1}(\FGlobal)_{>0}\), and a compact subset 
	\(C\subset\Sym_{\nInLZorth}(\BA_{\FGlobal,f})\).
	There are only finitely many
	\(b\in\Mat_{1,\nInLZorth-1}(\FGlobal)\) for which
	\(T_b(\xi)=\left(\begin{smallmatrix}
		\xi&b\\{}^tb&T_\flat\end{smallmatrix}\right)\) is totally positive
	definite and satisfies
	$
	T_b(\xi)\in C$.
	
\end{prop}

\begin{proof}
	The projection of \(C\) to its upper-right block is a compact
	subset of \(\BA_{\FGlobal,f}^{\nInLZorth-1}\) containing the finite
	components of every such \(b\).
	If \(\xi\) is nonpositive at a real embedding, there are no totally
	positive definite matrices \(T_b(\xi)\). Otherwise, for every real embedding
	\(\sigma:\FGlobal\hookrightarrow\R\), positive definiteness implies
	\[
	\sigma(b)\,\sigma(T_\flat)^{-1}\,{^t\sigma(b)}\leq\sigma(\xi).
	\]
	This inequality defines a closed ellipsoid. The product of
	these ellipsoids over the real embeddings bounds all archimedean 
	components of \(b\).  Consequently every such \(b\) lies in a fixed
	compact subset of \(\BA_{\FGlobal}^{\nInLZorth-1}\).  The diagonal
	image of \(\FGlobal^{\nInLZorth-1}\) is closed and discrete, so its
	intersection with this compact subset is finite.
\end{proof}
By \eqref{eq:local-standard-section-Iwasawa} and smoothness, the standard
section defining \(E(g,s,\phi_1\otimes\phi_\flat)\) is right invariant
under \(n(U)\) for every \(s\), for some compact open additive subgroup
\(U\subset\Sym_{\nInLZorth}(\BA_{\FGlobal,f})\). 
Since \(g_\sz\) has finite component \(1_f\), translating the Fourier
integral \eqref{eq:eisenstein-fourier-coefficient} by \(u\in U\) gives
\[
\pEis_{T/2}(g_\sz,\phi_1\otimes\phi_\flat)
=\psi_f\!\left(\tfrac12\operatorname{tr}(Tu)\right)
\pEis_{T/2}(g_\sz,\phi_1\otimes\phi_\flat).
\]
Thus $\pEis_{T/2}(g_\sz,\phi_1\otimes\phi_\flat)\neq 0$ only when $T_f$ lies in the lattice
\[
\left\{T_f\in\Sym_{\nInLZorth}(\BA_{\FGlobal,f}):
\psi_f\!\left(\tfrac12\operatorname{tr}(T_fu)\right)=1
\text{ for every }u\in U\right\}\subseteq \Sym_{\nInLZorth}(\BA_{\FGlobal,f}).
\]
For nonsingular \(T\) with \(\Diff(T,\Vincoherent)=\{v\}\) and \(v<\infty\), \(T\) is
totally positive definite and \eqref{eq:eissum} gives
\(\pEis_{T/2,v}=\pEis_{T/2}\).
Hence \Cref{prop:finite-contributing-indices} leaves only finitely many
\(b\) contributing at finite places. Each corresponding matrix has at
most one such difference place, so \eqref{eq:FJ-local-coefficient} is
finite for every finite \(v\), and only finitely many finite places
contribute to \eqref{eq:FJ-total-coefficient}. Since there are finitely
many real places, the latter sum is finite as well.

At a real place \(\RealPlace\), the singleton difference condition,
\eqref{eq:eissum}, and \Cref{rem:DiffTV} identify
\(\pEis_{T_b(\xi)/2,\RealPlace}\) with \(\pEis_{T_b(\xi)/2}\).
The finite-adelic support restriction above confines the contributing
\(b\) to a fixed fractional \(O_\FGlobal\)-lattice.
The smoothness of \(\pEis\) and integration by parts on the compact
additive quotient give rapid decay of these Fourier coefficients
in the archimedean entries of \(T\), locally uniformly together with
all derivatives in \(\sz\). For
\(\sz=\operatorname{diag}(\tau,iY_\flat)\), the factor
\(q^{T_b(\xi)/2}(\sz)\) is independent of \(b\).
Summing the rapid decay over this lattice proves absolute and locally
uniform convergence of the sum over \(b\) in
\eqref{eq:FJ-local-coefficient}, together with its derivatives in
\((\tau,Y_\flat)\).

\subsection{Adelic generating function}
\label{subsec:adelic-generating-function}
We use the following adelization for the special-cycle generating series of
\Cref{thm:Bruinier-modularity}, the arithmetic pairings of
\Cref{thm:Bruinier-modularity:cor} under their admissibility hypotheses,
and the archimedean series of \Cref{defn:modularity-I}.
As explained in \cite[\S2.5]{YuanZhangZhang2023Triple}, the central derivative of an incoherent Eisenstein series, and hence the corresponding Fourier--Jacobi series, generally fails to satisfy the Weil equivariance relation \eqref{eq:Weil-equiv}.

\begin{defn}\label{defn:adelic-generating-function}
	For adelization, put
	\(\mathbf i=(i)_{\RealPlace\mid\infty}\).  For
	\(g_\infty=(g_\RealPlace)_{\RealPlace\mid\infty}\in
	\Mp_2(\FGlobal_\infty)\), write
	\(\operatorname{pr}(g_\RealPlace)=
	\begin{psmallmatrix}a_\RealPlace&b_\RealPlace\\c_\RealPlace&d_\RealPlace\end{psmallmatrix}\).
	Let \(j_{1/2}(g_\RealPlace,\tau_\RealPlace)\) be the standard metaplectic automorphy
	factor, normalized by
	\(j_{1/2}(g_\RealPlace,\tau_\RealPlace)^2=c_\RealPlace\tau_\RealPlace+d_\RealPlace\), and put
	\[
	j_{\dimofV/2}(g_\infty,\mathbf i)
	:=\prod_{\RealPlace\mid\infty}
	j_{1/2}(g_\RealPlace,i)^{\dimofV},
	\qquad
	g_\infty\mathbf i:=(\operatorname{pr}(g_\RealPlace)i)_{\RealPlace\mid\infty}.
	\]
	Suppose for each 	\(\phi_1\in\sS(\Vincoherent_f)\) we have a generating series   \(F(\tau,\phi_1)=\sum_{\xi\in\FGlobal}
	\mathcal F(\xi,y,\phi_1)q^{\xi/2}\), where \(y=\operatorname{Im}\tau\) and the coefficients \(\mathcal F(\xi,y,\phi_1)\in\C\) depend linearly on \(\phi_1\).  We define its associated generating series on $\Mp_2(\BA_{\FGlobal})$ by
	\begin{equation}\label{eq:adelization-general}
		\mathcal F_{\phi_1}(g_\infty g_f)
		:=j_{\dimofV/2}(g_\infty,\mathbf i)^{-1}
		F\!\left(g_\infty\mathbf i,\omega_f(g_f)\phi_1\right),
		\qquad g_f\in\Mp_2(\BA_{\FGlobal,f}), g_\infty\in\Mp_2(\BA_{\FGlobal,\infty}).
	\end{equation} 
	We say the generating series is \emph{modular} of weight $m/2$ if
	\(\mathcal F_{\phi_1}\) converges and  belongs to
	\(\CA_\infty(\SL_2(\FGlobal)\backslash\Mp_2(\BA_\FGlobal),m/2)\).
	If the generating series is modular and $\mathcal F(\xi,y,\phi_1)$ is independent of $y$, then $\mathcal F_{\phi_1}$ is holomorphic by \Cref{subsec:nota-hilbert}.
	For
	\(h_f\in\Mp_2(\BA_{\FGlobal,f})\), set
	\(R(h_f)\mathcal F(g):=\mathcal F(gh_f)\).  By definition,
	\begin{equation}\label{eq:Weil-equiv}
		R(h_f)\mathcal F_{\phi_1}
		=\mathcal F_{\omega_f(h_f)\phi_1}.
	\end{equation}
\end{defn}

\section{Archimedean height pairings}\label{sec:height}

We recall the Kudla and Bruinier Green functions and the archimedean height formula of \cite{garcia_sankaran_2019_green_forms}. These results compare archimedean arithmetic intersections with Fourier coefficients of the central derivative of the incoherent Eisenstein series.

The results cited below are usually stated for the GSpin Shimura variety $\ShG(\Res_{\FGlobal/\Q}\GSpin(\Vcoherent), X)$. We apply them to our SO Shimura variety $\ShG_{K_2}$ using the following lemma and pushforward. 
\begin{lem}\label{lem:GSpin-SO-surj}
	For compatible levels $\GSpin(\Vcoherent)(\BA_{\FGlobal,f})\supset K_{\GSpin}\to K_2\subset \SO(\Vcoherent)(\BA_{\FGlobal,f})$, the finite \'etale map $\ShG_{K_{\GSpin}}\to \ShG_{K_2}$ is a cover. 
\end{lem}

\begin{proof}
	The kernel of
	\(\Res_{\FGlobal/\Q}\GSpin(\Vcoherent)
	\to\Res_{\FGlobal/\Q}\SO(\Vcoherent)\)
	is the induced torus \(\Res_{\FGlobal/\Q}\Gm\).  Hilbert's theorem 90 implies the claim. 
\end{proof}

\subsection{The Kudla Green functions}\label{subsec:Kudla}

For a real embedding \(\nu:\FGlobal\hookrightarrow\BR\), work on the
complex fiber \(\SOShim\times_{\FGlobal,\nu}\C\), which can be identified with the Shimura variety associated to the nearby space
\(V^{(\nu)}\) by \Cref{subsec:changing-iota}. Its real quadratic space at \(\nu\) has
signature \((\dimofV-2,2)\); write \(\mathcal D_\nu\) for its
orthogonal symmetric domain.  We transport the level $K_2$ and all finite
Schwartz functions along fixed finite-adelic isometries
\[
 V^{(\nu)}\otimes_\FGlobal\BA_{\FGlobal,f}
 \simeq\Vincoherent_f,
\]
as in \Cref{subsec:changing-iota}.  Finite adelic arguments below are
evaluated through these identifications.  We use
\(d^c=(\partial-\bar\partial)/(4\pi i)\), with the differentials acting
on \(\mathcal D_\nu\), and the Hodge curvature \(\Omega_{\mathcal E}\)
of \Cref{defn:GS-curvature-volume}.

\begin{defn}
\label{defn:Kudla-Green-functions}
For \(x\in V^{(\nu)}_\nu\setminus\{0\}\), let
\(\mathcal Z_\nu(x)=\{z\in\mathcal D_\nu:z\perp x\}\).
For \(z\in\mathcal D_\nu\), write
\(x=x_z+x_{z^\perp}\), where \(x_z\) is the projection to the
negative two-plane represented by \(z\), and set
\[
 R_\nu(x,z):=-(x_z,x_z)_\nu.
\]
This is positive away from \(\mathcal Z_\nu(x)\), which is empty
when \((x,x)_\nu\leq0\).
For \(r>0\), put
\[
 \beta_1(r):=\int_1^\infty e^{-rt}\frac{dt}{t}.
\]
For \(y_\nu>0\) and \(\psi_\nu(t)=e^{2\pi i\kappa_\nu t}\), the
Kudla Green function on \(\mathcal D_\nu\setminus\mathcal Z_\nu(x)\) is
\begin{equation}\label{eq:Kudla-Green-function}
 \CG^{\mathrm K}_\nu(x,y_\nu,z)
 :=\beta_1\!\left(2\pi\kappa_\nu y_\nu R_\nu(x,z)\right).
\end{equation}
The normalized Kudla--Millson coefficient form is the smooth extension
to \(V^{(\nu)}_\nu\times\mathcal D_\nu\) of the following expression
on \(R_\nu(x,z)>0\):
\begin{equation}\label{eq:KM-coefficient-form}
 \omega_{\mathrm{KM},\nu}(x,y_\nu,z)
 =e^{-2\pi\kappa_\nu y_\nu R_\nu(x,z)}
 \left(
 \frac{i\kappa_\nu y_\nu}{R_\nu(x,z)}
 \partial R_\nu(x,z)\wedge\bar\partial R_\nu(x,z)
 -\Omega_{\mathcal E}\right).
\end{equation}
In particular, \(\omega_{\mathrm{KM},\nu}(0,y_\nu)=-\Omega_{\mathcal E}\).

Let \(\xi\in\FGlobal^\times\) and
\(\phi_1\in\sS(\Vincoherent_f)^{K_2}\).  For
\(y=(y_v)_{v\mid\infty}\) with \(y_v>0\), on a component represented
by \(h_f\in\SO(V^{(\nu)})(\BA_{\FGlobal,f})\), define
\begin{equation*}
 \CG^{\mathrm K}_\nu(\xi,y,\phi_1)([z,h_f])
 =\sum_{\substack{x\in V^{(\nu)}(\FGlobal)\\(x,x)=\xi}}
 \phi_1(h_f^{-1}x_f)\,
 \CG^{\mathrm K}_\nu(x,y_\nu,z).
\end{equation*}
We define
\(\CG^{\mathrm K}_\nu(0,y,\phi_1)=-\phi_1(0)\log y_\nu\).
\end{defn}

\begin{prop}\label{prop:Kudla-Green-equation}
For \(x\ne0\), the function \(\CG^{\mathrm K}_\nu(x,y_\nu)\) is
locally integrable, has a logarithmic singularity along
\(\mathcal Z_\nu(x)\) when \((x,x)_\nu>0\), and is smooth when
\((x,x)_\nu\leq0\).  The form
\(\omega_{\mathrm{KM},\nu}(x,y_\nu)\) is closed of type \((1,1)\), and
\begin{equation}\label{eq:Kudla-Green-equation}
 dd^c[\CG^{\mathrm K}_\nu(x,y_\nu)]
 +\delta_{\mathcal Z_\nu(x)}
 =[\omega_{\mathrm{KM},\nu}(x,y_\nu)].
\end{equation}
Here \([\alpha]\) denotes the current associated to a locally
integrable form \(\alpha\), and \(\delta_Z\) denotes integration over
\(Z\), with \(\delta_\varnothing=0\).
\end{prop}
\begin{proof}
The smooth extension in \eqref{eq:KM-coefficient-form} and the current
identity are the rank-one case of
\cite[\S2.3 and (2.6.1)--(2.6.2)]{garcia_sankaran_2019_green_forms},
transported by \Cref{lem:additive-character-transport}.
\end{proof}

\begin{defn}\label{defn:higher-Kudla-Green}
Let \(1\leq r\leq\dimofV-1\), and let
\(\bx=(x_1,\ldots,x_r)\in(V^{(\nu)}_\nu)^r\) have nonsingular
Gram matrix \(((x_i,x_j)_\nu)_{i,j}\).  Put
\(\mathcal Z_\nu(\bx)=\bigcap_{i=1}^r\mathcal Z_\nu(x_i)\).
For a positive definite matrix \(Y_\nu\in\Sym_r(\BR)\), define the
natural higher Kudla Green form away from \(\mathcal Z_\nu(\bx)\) by
\begin{equation}\label{eq:higher-Kudla-Green-form}
 \CG^{\mathrm K}_\nu(\bx,Y_\nu,z)
 :=\int_1^\infty\sum_{i=1}^r
 e^{-2\pi\kappa_\nu tR_\nu((\bx Y_\nu^{1/2})_i,z)}
 \bigwedge_{\substack{1\leq j\leq r\\j\ne i}}
 \omega_{\mathrm{KM},\nu}((\bx Y_\nu^{1/2})_j,t,z)\,\frac{dt}{t}.
\end{equation}
Here \(Y_\nu^{1/2}\) is the positive definite symmetric square root,
\((\bx Y_\nu^{1/2})_i\) is the \(i\)-th vector of the resulting
tuple, and an empty wedge product equals \(1\).  This is a form of
type \((r-1,r-1)\); for \(r=1\) it gives
\eqref{eq:Kudla-Green-function}.

For nonsingular \(T\in\Sym_r(\FGlobal)\),
\(Y=(Y_v)_{v\mid\infty}\in\prod_{v\mid\infty}\Sym_r(\BR)_{>0}\), and
\(\phi_f\in\sS(\Vincoherent_f^r)^{K_2}\), its weighted current on the
component represented by \(h_f\) is
\begin{equation}\label{eq:higher-Kudla-Green-current}
 g_\nu^{\mathrm{K}}(T,Y,\phi_f)\big|_{h_f}
 :=\sum_{\substack{\bx\in V^{(\nu)}(\FGlobal)^r\\
                    ((x_i,x_j))_{i,j}=T}}
 \phi_f(h_f^{-1}\bx_f)
 [\CG^{\mathrm K}_\nu(\bx,Y_\nu)].
\end{equation}
\end{defn}

\begin{defn}
	For properly intersecting special cycles \(Z_1,Z_2\), let \(g_i\) be their
	Kudla Green currents defined above, and let \(\omega_i\) be the smooth forms
	satisfying \(dd^cg_i+\delta_{Z_i}=[\omega_i]\), for \(i=1,2\). Define their star-product by
	\begin{equation}\label{eq:star-product-height}
		g_1\star g_2:=g_1\wedge\delta_{Z_2}+\omega_1\wedge g_2.
	\end{equation}
	Then $dd^c(g_1\star g_2)+\delta_{Z_1\cap Z_2}=[\omega_1\wedge \omega_2]$, so $g_1\star g_2$ is a Green current for the cycle $Z_1\cap Z_2$. 
\end{defn}

\begin{prop}
\label{prop:higher-Kudla-Green}
The form in \eqref{eq:higher-Kudla-Green-form} is smooth away from
\(\mathcal Z_\nu(\bx)\), locally integrable on \(\mathcal D_\nu\), and
satisfies
\begin{equation*}
 dd^c[\CG^{\mathrm K}_\nu(\bx,Y_\nu)]
 +\delta_{\mathcal Z_\nu(\bx)}
 =\left[\bigwedge_{j=1}^r
 \omega_{\mathrm{KM},\nu}((\bx Y_\nu^{1/2})_j,1)\right].
\end{equation*}
The series \eqref{eq:higher-Kudla-Green-current} converges as a current
on \(\SOShim\times_{\FGlobal,\nu}\C\); summing this equation gives
its Green equation for the weighted special cycle. 
Let \(r=r_1+r_2\) with \(r_1,r_2\geq1\), and let
\(Y_{1,\nu}\in\Sym_{r_1}(\BR)_{>0}\) and
\(Y_{2,\nu}\in\Sym_{r_2}(\BR)_{>0}\).
Assume that the tuples \((x_1,\ldots,x_{r_1})\) and
\((x_{r_1+1},\ldots,x_r)\) have nonsingular Gram matrices and that
their special cycles intersect properly. For
\(Y_\nu=\operatorname{diag}(Y_{1,\nu},Y_{2,\nu})\), one has
\[
 [\CG^{\mathrm K}_\nu(\bx,Y_\nu)]
 \equiv[\CG^{\mathrm K}_\nu((x_1,\ldots,x_{r_1}),Y_{1,\nu})]
 \star[\CG^{\mathrm K}_\nu((x_{r_1+1},\ldots,x_r),Y_{2,\nu})]
 \pmod{\operatorname{im}\partial+\operatorname{im}\bar\partial},
\]
where the images are taken in the space of currents.  The weighted
sum of this congruence over tuples of Gram matrix \(T\) descends
to the Shimura fiber.
\end{prop}
\begin{proof}
The integral is the transgression construction of
\cite[(2.3.3), Proposition~2.4.5(a), and (2.6.1)]
{garcia_sankaran_2019_green_forms}.
The nonsingular Gram matrix makes \(\mathcal Z_\nu(\bx)\) empty or
smooth of codimension \(r\), so its equation~(2.6.2) gives local
integrability and the Green equation after
\Cref{lem:additive-character-transport}.
The estimates in its \S4.3.1 justify the weighted current sum and
termwise application of \(dd^c\).
Under the proper intersection hypothesis, Theorem~2.7.2 gives the
star-product congruence; the equivariance and convergence of the
correction currents in the proof of Theorem~4.5.1 allow its weighted
sum to descend to the Shimura fiber.
\end{proof}

\begin{defn}\label{defn:KM-theta}
For \(\tau\in\mathbb H_{\FGlobal}\), with imaginary part
\(y=(y_v)_{v\mid\infty}\), and \(\phi_1\in\sS(\Vincoherent_f)^{K_2}\),
the genus-one Kudla--Millson theta form is
\begin{align}
 \Theta_{\mathrm{KM},\nu}(\tau,\phi_1)
 &=\sum_{\xi\in\FGlobal}
       \Theta_{\mathrm{KM},\nu,\xi}(y,\phi_1)q^{\xi/2},
       \label{eq:KM-theta-series}\\
 \Theta_{\mathrm{KM},\nu,\xi}(y,\phi_1)(z,h_f)
 &:=\sum_{\substack{x\in V^{(\nu)}(\FGlobal)\\(x,x)=\xi}}
       \phi_1(h_f^{-1}x_f)\,
       \omega_{\mathrm{KM},\nu}(x,y_\nu,z).
       \notag
\end{align} 
\end{defn}

\begin{thm}\label{thm:KM-theta-modularity}
The series \eqref{eq:KM-theta-series}, together with all derivatives,
converges locally uniformly and defines a smooth metaplectic Hilbert
modular form of parallel weight \(\dimofV/2\), valued in closed
\((1,1)\)-forms on \(\SOShim\times_{\FGlobal,\nu}\C\).
It is holomorphic in \(\tau_v\) for \(v\ne\nu\) and generally
non-holomorphic in \(\tau_\nu\).
More precisely, for \(\gamma\in\SL_2(\FGlobal)\), choose compatible
components \(\gamma_f\) and \((\gamma_v)_{v\mid\infty}\) of its
canonical adelic metaplectic lift.  Then
\begin{equation}\label{eq:KM-theta-transformation}
 \Theta_{\mathrm{KM},\nu}(\gamma\tau,\phi_1)
 =\left(\prod_{v\mid\infty}
       j_{1/2}(\gamma_v,\tau_v)^{\dimofV}\right)
   \Theta_{\mathrm{KM},\nu}
       (\tau,\omega_f(\gamma_f^{-1})\phi_1),
\end{equation}
where \(j_{1/2}\) is the automorphy factor of
\Cref{subsec:adelic-generating-function}, and \(\omega_f\) is the
finite Weil representation for \(\psi_f\) on \(\sS(\Vincoherent_f)\).
Thus a finite level \(K_\bH\) fixing \(\phi_1\) under \(\omega_f\)
is a level of this modular form.
For \(\xi\ne0\), the Fourier coefficients satisfy
\[
 dd^c[\CG^{\mathrm K}_\nu(\xi,y,\phi_1)]
 +\delta_{\CZ_\xi(\phi_1)_\nu}
 =[\Theta_{\mathrm{KM},\nu,\xi}(y,\phi_1)],
\]
where \(\CZ_\xi(\phi_1)_\nu\) is the complex special divisor on the
\(\nu\)-fiber.  For \(\xi>0\),
\(\Theta_{\mathrm{KM},\nu,\xi}(y,\phi_1)\) therefore represents
the cohomology class of \(\CZ_\xi(\phi_1)_\nu\).
\end{thm}
\begin{proof}
At \(\nu\), the Schwartz form is
\(e^{-\pi\kappa_\nu(x,x)_\nu}\omega_{\mathrm{KM},\nu}(x,1)\);
at each \(v\ne\nu\), use the Gaussian
\(x\mapsto e^{-\pi\kappa_v(x,x)_v}\).
Each factor has Weil weight \(\dimofV/2\), and their tensor product
is closed on \(\mathcal D_\nu\).  Schwartz decay gives the asserted
convergence, and Poisson summation gives
\eqref{eq:KM-theta-transformation}; see
\cite[\S3.4]{Bruinier2012Regularized} and
\cite[\S3 and Theorem~5.2(i)]{KudlaMillson1990}.
Apply \Cref{lem:additive-character-transport} for the character \(\psi\).
The coefficient identity follows by summing
\eqref{eq:Kudla-Green-equation} with the weights
\(\phi_1(h_f^{-1}x_f)\).
\end{proof}

\begin{defn}
Assume that \(\phi_1\) is admissible. For \(\xi>0\), define the \emph{arithmetic special divisor}
\begin{equation*}
 \widehat{\CZ}^{\mathrm K}(\xi,\tau,\phi_1)
 :=\left(\CZ_\xi(\phi_1),
       (\CG^{\mathrm K}_\nu(\xi,y,\phi_1))_{\nu\mid\infty}\right).
\end{equation*}
Let \(\widehat\omega\) be the metrized Hodge line bundle, with underlying
line bundle \(\omega\) and curvature equal to the invariant K\"ahler form in the Kudla--Millson
normalization.  At \(\xi=0\) our convention is
\begin{equation*}
 \widehat{\CZ}^{\mathrm K}(0,\tau,\phi_1)
 =-\phi_1(0)\widehat\omega
   +\left(0,(-\phi_1(0)\log y_\nu)_{\nu\mid\infty}\right).
\end{equation*}
For nonzero \(\xi\) not totally positive, put
\[
 \widehat{\CZ}^{\mathrm K}(\xi,\tau,\phi_1)
 :=\left(0,(\CG^{\mathrm K}_\nu(\xi,y,\phi_1))_{\nu\mid\infty}\right).
\]
\end{defn}

\subsection{The Bruinier Green functions}

For \(\xi>0\), let \(\CG^{\mathrm B}_\nu(\xi,\phi_1)\) be the
Bruinier Green function from the regularized theta lift.  It is
independent of \(y\), and has the same logarithmic singularity as
\(\CG^{\mathrm K}_\nu(\xi,y,\phi_1)\). By \cite[Corollary~5.16 and Theorem~6.4]{Bruinier2012Regularized}  the curvature of the Bruinier
Green function is harmonic. We use the representatives and arithmetic Chow
convention in \cite[\S3]{Sankaran2025ArithmeticSpecialCycles}, with our
\(\widehat\omega\) equal to the negative of the arithmetic Hodge class denoted
by \(\widehat\omega\) there. For nonzero \(\xi\notin\FGlobal_{>0}\), set
\(\CG^{\mathrm B}_\nu(\xi,\phi_1)=0\).
Write
\[
 \widehat{\CZ}^{\mathrm B}(\xi,\phi_1)
 :=\bigl(\CZ_\xi(\phi_1),
       (\CG^{\mathrm B}_\nu(\xi,\phi_1))_{\nu\mid\infty}\bigr),
 \qquad
 \widehat{\CZ}^{\mathrm B}(0,\phi_1)=-\phi_1(0)\widehat\omega.
\]
The Bruinier series is
\begin{equation*}
 \widehat{\CZ}^{\mathrm B}(\tau,\phi_1)
 =\widehat{\CZ}^{\mathrm B}(0,\phi_1)
  +\sum_{\xi>0}\widehat{\CZ}^{\mathrm B}(\xi,\phi_1)q^{\xi/2}.
\end{equation*}

Put
\(\CG^{\mathrm{B-K}}_\nu(\xi,y,\phi_1)
 :=\CG^{\mathrm B}_\nu(\xi,\phi_1)-\CG^{\mathrm K}_\nu(\xi,y,\phi_1)\)
for \(\xi\ne0\). At \(\xi=0\), the common Hodge term of the
arithmetic divisors cancels, and we set
\[
 \CG^{\mathrm{B-K}}_\nu(0,y,\phi_1):=\phi_1(0)\log y_\nu.
\]
We use the following comparison theorem of Sankaran.
\begin{thm}[{\cite[Theorem~3.2]{Sankaran2025ArithmeticSpecialCycles}}]
 \label{thm:Kudla-Bruinier-modularity}
 There are smooth functions 
 \(\mathcal S^{\mathrm{B-K}}_\nu(\tau,\phi_1)\in C^\infty(\SOShim\times_{\FGlobal,\nu}\C)\) such that
$
 [\mathcal S^{\mathrm{B-K}}_\nu]_{q^{\xi/2}}
 -\CG^{\mathrm{B-K}}_\nu(\xi,y,\phi_1)
 $ is constant on each connected component of $\SOShim\times_{\FGlobal,\nu}\C$ for every $\xi\in\FGlobal$, and for each $x\in \SOShim\times_{\FGlobal,\nu}\C$ and each fixed $\phi_1$, the function
 $\mathcal S^{\mathrm{B-K}}_\nu(\tau,\phi_1)(x)$ is a smooth
 metaplectic Hilbert modular form of parallel weight $\dimofV/2$
 at a finite level.
 Consequently, for every zero-cycle \(C_\nu\) of degree zero on
 each component, the series
 \[
 \sum_{\xi\in\FGlobal}
 \CG^{\mathrm{B-K}}_\nu(\xi,y,\phi_1)(C_\nu)q^{\xi/2}
 \]
 is a smooth metaplectic Hilbert modular form of that weight.
 Its adelization intertwines finite right translation by
 \(h_f\in\Mp_2(\BA_{\FGlobal,f})\) with replacing \(\phi_1\) by
 \(\omega_f(h_f)\phi_1\).  This follows from the
 Schwartz-dual-valued transformation law in
 \cite[\S2.6, (2.40)--(2.41), and Theorem~3.2]{Sankaran2025ArithmeticSpecialCycles};
 the componentwise constants vanish on \(C_\nu\).
\end{thm}

\begin{thm}[{\cite[Theorem~1.3]{YuanZhangZhang2009}}]
	\label{thm:Bruinier-modularity}
	For an $\ov{\Q}$-valued finite Schwartz function $\phi_1\in \sS(\Vincoherent_f)^{K_2}$, the
	generating series of special cycles on the generic fiber is a holomorphic metaplectic Hilbert modular
	form of parallel weight $\dimofV/2$, i.e., it satisfies
	\[
	-\phi_1(0)[\omega]
	+\sum_{\xi>0}[Z_\xi(\phi_1)]q^{\xi/2}
	\in\CA_{\hol}(\SL_2(\FGlobal)\backslash\Mp_2(\BA_\FGlobal),\dimofV/2)_{\ov{\Q}}
	\otimes_{\ov{\Q}}\Pic(\SOShim)_{\ov{\Q}}.
	\]
\end{thm}
\begin{proof}
	The cited theorem gives the Picard-valued generating series after pullback to 
	the GSpin Shimura variety, which is a finite \'etale cover of $\ShG_{K_2}$ by \Cref{lem:GSpin-SO-surj}.  Pullback followed by pushforward acts on Picard classes as
	multiplication by the degree of the finite \'etale cover. The claim follows. 
\end{proof}
\subsection{The Garcia--Sankaran formula and its normalization}

\begin{defn}
\label{defn:GS-curvature-volume}
For a real place \(\nu\), the realization of the domain by isotropic
lines is
\[
 \mathcal D_\nu
 =\bigl\{[w]\in\mathbb P(V^{(\nu)}_\nu\otimes_{\mathbb R}\mathbb C):
       (w,w)_\nu=0,\ (w,\overline w)_\nu<0\bigr\},
\]
where \((\ ,\ )_\nu\) is extended complex bilinearly and conjugation
is taken with respect to \(V^{(\nu)}_\nu\).  This domain has two
connected components and complex dimension \(\dimofV-2\).
The tautological holomorphic line bundle \(\mathcal E\) on
\(\mathcal D_\nu\) has fiber \(\mathcal E_{[w]}=\mathbb Cw\), with
invariant Hermitian metric
\[
 \|v\|_{\mathcal E}^2=-(v,\overline v)_\nu,
 \qquad v\in\mathcal E_{[w]}.
\]
For a nowhere-vanishing local holomorphic section \(s\) of
\(\mathcal E\), define its first Chern form by
\begin{equation}\label{eq:GS-curvature-definition}
 \Omega_{\mathcal E}:=c_1(\mathcal E,\|\ \|_{\mathcal E})
 =-\frac{i}{2\pi}\partial\bar\partial\log\|s\|_{\mathcal E}^2
 =-dd^c\log\|s\|_{\mathcal E}^2.
\end{equation}
The expression is independent of the choice of \(s\) and gives the
invariant positive K\"ahler form on \(\mathcal D_\nu\); see
\cite[\S1.2, (2.2.6)--(2.2.8)]{garcia_sankaran_2019_green_forms}.
The metrized line bundle descends to the \(\SO\)-Shimura fiber as an
orbifold line bundle and identifies with the metrized Hodge line
\(\widehat\omega\) under
the period description.  We use \(\Omega_{\mathcal E}\) also for
the descended form, so \(c_1(\omega)=\Omega_{\mathcal E}\). 
Following \cite[(5.3.1)]{garcia_sankaran_2019_green_forms}, we define
\begin{equation}\label{eq:GS-volume-definition}
 \operatorname{Vol}^{\mathrm{GS}}
 (\SOShim\times_{\FGlobal,\nu}\C,\Omega_{\mathcal E})
 :=\int_{[\SOShim\times_{\FGlobal,\nu}\C]}
          \Omega_{\mathcal E}^{\dimofV-2}.
\end{equation}
\end{defn}

Garcia--Sankaran use the standard trace character $\psi_{\mathrm{std}}=\psi_{\BQ}\circ\operatorname{Tr}_{\FGlobal/\BQ}$.  The following
transport identity lets us apply their formula with the character
\(\psi\) fixed in this article.
\begin{lem}
	\label{lem:additive-character-transport}
	Let \(r\ge1\), and let \(U\) be a quadratic space over \(\FGlobal\), or a
	coherent or incoherent quadratic space over \(\BA_{\FGlobal}\).
	Let \(U[a_\psi]\) denote the
	same underlying vector space with bilinear form
	\((x,x')_{U[a_\psi]}:=a_\psi(x,x')_U\).  The identity on the underlying
	Schwartz space gives a canonical equality of genus-\(r\) Weil representations
	\begin{equation}\label{eq:additive-character-Weil-transport}
		\omega_{U,\psi}=\omega_{U[a_\psi],\psi_{\mathrm{std}}}.
	\end{equation}
	If \(\dim U=r+1\), \(T\in\Sym_r(\FGlobal)\), and \(\phi\) is adelic
	Schwartz data for \(U^r\), this equality identifies their standard
	Siegel--Weil sections and gives, for \(g\in\Mp_{2r}(\BA_{\FGlobal})\) and \(s\in\BC\),
	\begin{equation}\label{eq:additive-character-coefficient-transport}
		E^{\psi,U}_{T/2}(g,s,\phi)
		=E^{\psi_{\mathrm{std}},U[a_\psi]}_{a_\psi T/2}(g,s,\phi),
		\qquad
		q_\psi^{T/2}=q_{\psi_{\mathrm{std}}}^{a_\psi T/2}.
	\end{equation}
\end{lem}

\begin{proof}
	Equation \eqref{eq:additive-character-scalar} gives
	\(\psi((x,x')_U)=
	\psi_{\mathrm{std}}((x,x')_{U[a_\psi]})\).  Thus the unipotent
	characters, Fourier transforms, Weil indices, and Levi actions in the two
	Weil representations agree, proving
	\eqref{eq:additive-character-Weil-transport} and the assertion about
	standard sections.  On the \(N_r\)-coordinate, of dimension
	\(r(r+1)/2\), let \(db_{\psi,v}\) and
	\(db_{\psi_{\mathrm{std}},v}\) be the measures self-dual for the
	half-trace pairing with respect to \(\psi_v\) and
	\(\psi_{\mathrm{std},v}\), respectively.  They satisfy
	\[
	db_{\psi,v}
	=|a_\psi|_v^{r(r+1)/4}db_{\psi_{\mathrm{std}},v}.
	\]
	Their global products agree because
	\(\prod_v|a_\psi|_v=1\).  The Fourier characters also agree:
	\(\psi(-\tfrac12\operatorname{tr}(Tb))
	=\psi_{\mathrm{std}}(-\tfrac12
	\operatorname{tr}(a_\psi Tb))\).  This proves
	\eqref{eq:additive-character-coefficient-transport}.
\end{proof}

\begin{thm}
	\label{thm:garcia-sankaran-archimedean}
	Let \(\nu\) be an archimedean place of \(\FGlobal\).  Let \(T\) be a nonsingular Gram matrix in
	\(\Sym_{\dimofV-1}(\FGlobal)\) such that
	$
	\Diff(T,\Vincoherent)=\{\nu\}.
	$ 
	Fix a \(K_2\)-invariant function
	$
	\phi_f\in\sS(\Vadl_f^{\dimofV-1}).
	$ 
	For \(\sz=\sx+i\sy\in\mathbb H_{\dimofV-1}\), let
	\(g_\nu^{\mathrm K}(T,\sy,\phi_f)\) be the weighted current
	defined in \eqref{eq:higher-Kudla-Green-current}, using the finite-adelic
	identifications fixed in \Cref{subsec:Kudla}.  We then have
	\begin{equation}\label{eq:garcia-sankaran-archimedean}
		\frac{(-1)^{\dimofV-1}}
		{\operatorname{Vol}^{\mathrm{GS}}(\SOShim\times_{\FGlobal,\nu}\C,\Omega_{\mathcal E})}
		\int_{[\SOShim\times_{\FGlobal,\nu}\C]}g_\nu^{\mathrm{K}}(T,\sy,\phi_f)\;q^{T/2}
		=\pEis_{T/2,\nu}(\sz,\phi_f).
	\end{equation}
\end{thm}
\begin{proof}
	Apply \cite[Theorem~5.3.1, specialized in (1.1.19)]
	{garcia_sankaran_2019_green_forms} with genus \(\dimofV-1\),
	positive-signature parameter \(\dimofV-2\), and Siegel--Weil point
	\(0\).  Their normalization constant \(\kappa_0\) equals \(2\).
	The Chern power and the scalar in their formula become
	\[
	\Omega_{\mathcal E}^{(\dimofV-2)+1-(\dimofV-1)}=1,
	\qquad
	\frac{(-1)^{\dimofV-1}\kappa_0}{2}=(-1)^{\dimofV-1}.
	\]

	Multiply the bilinear forms on \(\Vincoherent\) and \(V^{(\nu)}\)
	by \(a_\psi\).  Their moment matrix is \(a_\psi T/2\), and
	\Cref{lem:additive-character-transport} identifies their standard
	sections and Fourier monomial with ours.  The rescaling gives
	\(g_\nu^{\mathrm K}(T,\sy,\phi_f)\) as in
	\eqref{eq:higher-Kudla-Green-form} and preserves
	\(\Omega_{\mathcal E}\) and its volume.  The local factors relating
	their trace-self-dual measures to our half-trace-self-dual measures
	cancel globally by the product formula.

	On a compatible finite \'etale \(\GSpin\) cover, the integral and
	volume both scale by the covering degree, so the normalized formula
	descends to \(\SOShim\).  Under
	\(\Diff(T,\Vincoherent)=\{\nu\}\), only the \(\nu\)-term
	contributes to the sums over complex fibers and derivative places,
	by \eqref{eq:higher-Kudla-Green-current}, \eqref{eq:eissum}, and
	\Cref{rem:DiffTV}.  This gives
	\eqref{eq:garcia-sankaran-archimedean}.
\end{proof}

\subsection{Comparison of measures}

The Garcia--Sankaran K\"ahler volume and the finite Haar volume
\(\vol(K_2)\) have different normalizations.  Their product in
\eqref{eq:GS-finite-mass} fixes the scalar that makes the archimedean
comparison compatible with \Cref{thm:global-int}.

\begin{thm}
	\label{thm:GS-finite-mass}
	Assume \(\dimofV\ge3\).  For
	\(\nu\in\Sigma_{F,\infty}\), assume that
	\(V^{(\nu)}\) is anisotropic over \(\FGlobal\).  Then the finite
	Haar measure fixed in
	\Cref{defn:Haar-measure} and the Garcia--Sankaran K\"ahler measure
	on \(\SOShim\times_{\FGlobal,\nu}\C\), with volume defined in
	\eqref{eq:GS-volume-definition}, satisfy
	\begin{equation}\label{eq:GS-finite-mass}
		\operatorname{Vol}^{\mathrm{GS}}
		(\SOShim\times_{\FGlobal,\nu}\C,\Omega_{\mathcal E})\,\vol(K_2)=4.
	\end{equation}
	The volume $\operatorname{Vol}^{\mathrm{GS}}
	(\SOShim\times_{\FGlobal,\nu}\C,\Omega_{\mathcal E})$ is a positive
	rational number, since it is the top self-intersection number of the
	Hodge line bundle with the orbifold convention of
	\Cref{defn:GS-curvature-volume}.  Hence $\vol(K_2)\in\Q_{>0}$.
\end{thm}
\begin{proof}
	Let \(d\tau\) be the Tamagawa measure on
	\(\SO(V^{(\nu)})(\BA_{\FGlobal})\).  We show that it has the following decomposition
	\[
	d\tau=dh_\infty\,dh_f,
	\qquad
	dh_\infty=\prod_{\mu\mid\infty}dh_\mu,
	\qquad
	dh_f=\prod_{v<\infty}dh_v,
	\]
	where the local Gelfand--Leray measures are defined in
	\Cref{defn:Gelfand--Leray}.    
	Choose a nonzero invariant top form
	\(\omega_{\SO(V^{(\nu)})}\) defined over \(\FGlobal\).
	The local measures \(|\omega_{\SO(V^{(\nu)})}|_v\) use
	coordinate additive measures with
	\(\operatorname{vol}(O_{\FGlobal_v})=1\) at finite places
	and Lebesgue measure at real places.
	For the standard trace character, the self-dual volume of
	\(O_{\FGlobal_v}\) is the inverse square root of the norm of the
	local different.  By \eqref{eq:GL-trace-density}, the product
	formula for \(2\), \(a_\psi\), the determinant of the trace pairing
	in an \(\FGlobal\)-basis of \(\mathfrak{so}(V^{(\nu)})\), and the
	value of \(\omega_{\SO(V^{(\nu)})}\) at \(1\) on that basis gives
	\[
	\prod_v dh_v
	=|D_{\FGlobal}|^{-\dim\SO(V^{(\nu)})/2}
	\prod_v|\omega_{\SO(V^{(\nu)})}|_v,
	\]
	where \(D_{\FGlobal}=D_{\FGlobal/\Q}\) is the absolute
	discriminant.  The right-hand side is the Tamagawa measure
	by definition.
	
	The stabilizer of an oriented negative two-plane
	in \(\mathcal D_\nu\)  is
	\[
	\bigl(\SO(\dimofV-2)\times\SO(2)\bigr)
	\times\prod_{\substack{\mu\mid\infty\\\mu\ne\nu}}
	\SO(V^{(\nu)}_\mu)(\mathbb R)
	\subset\prod_{\mu\mid\infty}\SO(V^{(\nu)}_\mu)(\mathbb R).
	\]
	Give the stabilizer factor at each real place \(\mu\) its probability 
	Haar measure \(dk_\mu\), and put
	\(dk_\infty=\prod_{\mu\mid\infty}dk_\mu\).
	The quotient of the real group by this stabilizer is \(\mathcal D_\nu\),
	including its two connected components.  The quotient Haar
	measure \(dh_\infty/dk_\infty\) and the positive curvature
	volume form $\Omega_{\mathcal E}^{\dimofV-2}$ on $\mathcal D_\nu$ are invariant under $\prod_{\mu\mid\infty}\SO(V^{(\nu)}_\mu)(\mathbb R)$, hence there is a
	unique constant \(C_\infty>0\) satisfying
	\begin{equation*}
		\Omega_{\mathcal E}^{\dimofV-2}
		=C_\infty\,\frac{dh_\infty}{dk_\infty}.
	\end{equation*}
	Recall the complex uniformization
	\[
	\SOShim\times_{\FGlobal,\nu}\C
	\simeq\SO(V^{(\nu)})(\FGlobal)\backslash
	\bigl(\mathcal D_\nu\times
	\SO(V^{(\nu)})(\BA_{\FGlobal,f})\bigr)/K_2;
	\] 
	The two invariant measures $\Omega_{\mathcal E}^{\dimofV-2}$, $\frac{dh_\infty}{dk_\infty}$
	descend to $\SOShim\times_{\FGlobal,\nu}\C$. Taking volumes gives
	\begin{equation*}
		\operatorname{Vol}^{\mathrm{GS}}
		(\SOShim\times_{\FGlobal,\nu}\C,\Omega_{\mathcal E})
		=\frac{C_\infty\,\tau(\SO(V^{(\nu)}))}
		{\vol^{\mathrm{GL}}(K_2)},
	\end{equation*}
	where \(\tau(\SO(V^{(\nu)}))\) is the \(d\tau\)-volume of
	\(\SO(V^{(\nu)})(\FGlobal)\backslash
	\SO(V^{(\nu)})(\BA_{\FGlobal})\), which equals \(2\) by
	\cite{Weil1965SiegelFormula}. 
	By \Cref{defn:Haar-measure}, we have $\vol^{\mathrm{GL}}(K_2)
	=a_{\infty,\psi}^{\operatorname{tr}}\vol(K_2)$. 
	It remains to prove
	\(C_\infty=2a_{\infty,\psi}^{\operatorname{tr}}\). We compute the ratio of the densities of $\Omega_{\mathcal E}^{\dimofV-2}$ and $dh_\infty/dk_\infty$ at the base point $e=1(\SO(\dimofV-2)\times\SO(2))\in \mathcal D_\nu\cong \SO(\dimofV-2,2)/(\SO(\dimofV-2)\times\SO(2))$. 
	
	The compact dual $\mathcal D_\nu^\vee$ of either connected component of
	\(\mathcal D_\nu\) is
	\[
	\SO(\dimofV)/(\SO(\dimofV-2)\times\SO(2))
	\simeq
	\{[z]\in\mathbb P^{\nInLZorth}_{\mathbb C}:
	\textstyle\sum_{j=1}^{\dimofV}z_j^2=0\},
	\]
	where \(z_1,\ldots,z_{\dimofV}\) are homogeneous coordinates and we assume $e\in \mathcal D_\nu$ maps to $[0,\ldots,0,i,1]\in \mathcal D_\nu^\vee$. 
	For each real place \(\mu\), endow \(\SO(\dimofV)\) with the
	Gelfand--Leray measure \(dh_{\mu,\mathrm c}\) of
	\Cref{defn:Gelfand--Leray}, defined using \(\psi_\mu\).
	Use matrix realizations with bilinear forms
	\(1_{\dimofV}\) and
	\(\operatorname{diag}(1_{\dimofV-2},-1,-1)\), respectively.
	Their common block-diagonal stabilizer
	\(\SO(\dimofV-2)\times\SO(2)\) carries the same probability
	measure \(dk_\nu\).  Identify the tangent complements of its
	Lie algebra by
	\[
	 \begin{pmatrix}0&B\\-B^t&0\end{pmatrix}
	 \longleftrightarrow
	 \begin{pmatrix}0&B\\B^t&0\end{pmatrix},
	 \qquad B\in M_{\dimofV-2,2}(\mathbb R).
	\] 
	Use the entries of \(B\) and the independent upper-triangular
	entries in the stabilizer Lie algebra as coordinates.
	The pairing \((A_1,A_2)\mapsto-\operatorname{tr}(A_1A_2)/2\)
	has absolute determinant \(1\) in both full coordinate bases.
	Thus by \eqref{eq:GL-trace-density} $dh_{\nu,\mathrm c}$ and $dh_\nu$ give equal density at $1$.  Dividing by the common stabilizer density
	\(dk_\nu\) shows that the quotient densities
	\(dh_{\nu,\mathrm c}/dk_\nu\) and \(dh_\nu/dk_\nu\)
	have equal values at $1$ under this tangent-space
	identification.

	On the compact quadric, equip the tautological bundle \(\mathcal O(-1)\)
	with the metric induced by the Euclidean norm on
	\(\mathbb C^{\dimofV}\), and equip \(\mathcal O(1)\) with the
	dual metric, whose Chern form is denoted by \(c_1(\mathcal O(1))\).
	On both quadrics, use the affine chart \(z_{\dimofV}=1\) and the
	local holomorphic coordinates \(z_1,\ldots,z_{\dimofV-2}\) near
	\([0,\ldots,0,i,1]\).  Direct differentiation of the Hermitian
	metrics gives
	\[
	\left.\Omega_{\mathcal E}\right|_{e}
	=\frac{i}{4\pi}\sum_{a=1}^{\dimofV-2}
	dz_a\wedge d\bar z_a=
	\left.c_1(\mathcal O(1))\right|_{e}.
	\]
	The Chern form \(c_1(\mathcal O(1))\) satisfies
	\[
	\int_{\SO(\dimofV)/(\SO(\dimofV-2)\times\SO(2))}
	c_1(\mathcal O(1))^{\dimofV-2}=2,
	\]
	because the quadric has degree \(2\).
	Both
	\(dh_{\nu,\mathrm c}/dk_\nu\) and
	\(c_1(\mathcal O(1))^{\dimofV-2}\) are positive invariant
	densities on the compact quotient, so they differ by a constant.   Integration therefore gives
	\begin{equation*}
		\frac{dh_{\nu,\mathrm c}}{dk_\nu}
		=\frac{\operatorname{vol}_{dh_{\nu,\mathrm c}}
			(\SO(\dimofV))}{2}
		c_1(\mathcal O(1))^{\dimofV-2}.
	\end{equation*}
	We conclude that
	\[\frac{dh_{\nu}}{dk_\nu}=\frac{\operatorname{vol}_{dh_{\nu,\mathrm c}}
		(\SO(\dimofV))}{2} \Omega_{\mathcal E}^{\dimofV-2}\]
	holds at $e\in \mathcal D_\nu$, hence on both connected components of
	$\mathcal D_\nu$ by invariance under $\SO(\dimofV-2,2)$.
	Taking the product over the real places gives
	$$\frac{dh_\infty}{dk_\infty}=\frac{1}{2}\prod_{\mu\mid \infty}\operatorname{vol}_{dh_{\mu,\mathrm c}}
	(\SO(\dimofV)) \, \Omega_{\mathcal E}^{\dimofV-2}, \qquad C_\infty=2\prod_{\mu\mid \infty}\operatorname{vol}_{dh_{\mu,\mathrm c}}
	(\SO(\dimofV))^{-1}.$$
	By \eqref{eq:GL-trace-density}, \(dh_{\mu,\mathrm c}\) equals
	\(2^{-\nInLZorth/2}\kappa_\mu^{\nInLZorth\dimofV/4}\) times the gauge measure
	in \cite[Corollary~3.5]{BruinierYang2021Arithmetic}. That corollary therefore gives
	\begin{equation*}
		\begin{aligned}
			\operatorname{vol}_{dh_{\mu,\mathrm c}}(\SO(\dimofV))=\frac{2^{\nInLZorth/2}\pi^{\nInLZorth\dimofV/2}}
			{\Gamma_{\nInLZorth}(\dimofV/2)}
			\kappa_\mu^{\nInLZorth\dimofV/4},
		\end{aligned}
	\end{equation*}
	Their product equals $\bigl(a_{\infty,\psi}^{\operatorname{tr}}\bigr)^{-1}$ by
	\Cref{lem:global-Haar-phase-a}. Hence $C_\infty=2a_{\infty,\psi}^{\operatorname{tr}}$.
\end{proof}
\section{Arithmetic modularity and the error series}\label{sec:modularity}

The standing assumption $\FGlobal\ne\Q$ gives the properness of
$\SOIntModel$ by \Cref{prop:properness}.
We write
\[
\BR_{\fkd}
:=\BR\Big/\sum_{\ell\mid\fkd}\BQ\log\ell, \qquad \BR_{\fkd,\ov{\BQ}}
:=\BR_{\fkd}\otimes_\Q \ov{\Q},\qquad \C_{\fkd}:=\C\Big/\sum_{\ell\mid\fkd}\ov{\Q}\log\ell
\]
for the target of the truncated arithmetic degree.  
The multiplication map $\BR \otimes_\Q \ov{\Q}\twoheadrightarrow \C$ induces a \(\ov{\BQ}\)-linear map $\BR_{\fkd,\ov{\BQ}}\twoheadrightarrow \C_{\fkd}$. 
In this section we construct a metaplectic Hilbert modular form
\[
\CE_{T_\flat}(\tau,\phi_1,\phi_\flat)\in \CA_{\hol}(\SL_2(\FGlobal)\backslash\Mp_2(\BA_\FGlobal),\dimofV/2)_{\ov{\BQ}}\otimes_{\ov{\Q}}\C_{\fkd}
\]
whose coefficients are essentially the difference between local intersection
numbers and derivatives of representation densities; see \Cref{prop:modularity}.  

\subsection{Admissible arithmetic pairing}
Let $E$ be a totally real number field, let $\fkd>1$, and let $\CX$ be a smooth proper
scheme of finite type over $O_E[\fkd^{-1}]$, with generic fiber
$\CX_\eta$. Fix a K\"ahler metric on every complex fiber, using
the invariant metric for the Shimura varieties in this article.
Let $\widehat{\Pic}^{\adm}(\CX)$ be the group of hermitian
line bundles whose first Chern form is harmonic at every archimedean
place, and let $\Pic^{\adm}(\CX_\eta)$ be its image in
$\Pic(\CX_\eta)$.  
Let $\wt{\CZ}_1(\CX)_\Q$ denote rational one-cycles modulo rational
equivalence within closed fibers, and let
$\wt{\CZ}_1(\CX)_{\Q,\deg=0}$ be the subgroup whose generic zero-cycle
has degree zero on every connected component.
\begin{defn}
	As explained in \cite[Lemma~4.1]{MihatschZhang2024}, the truncated arithmetic intersection pairing of
	\cite[Proposition~2.3.1]{BGS94-JAMES} $\widehat{\Pic}^{\adm}(\CX)_\Q\times
	\wt{\CZ}_1(\CX)_{\Q, \deg=0}\to \BR_\fkd$ induces a bilinear pairing
	\begin{equation}\label{eq:admissible-degree-zero-pairing}
		(-,-):\Pic^{\adm}(\CX_\eta)_\Q\times
		\wt{\CZ}_1(\CX)_{\Q, \deg=0}\longrightarrow \BR_\fkd.
	\end{equation}
	Write a rational one-cycle on $\CX$ as $C=\sum_i a_i C_i$, where
	$a_i\in\Q$ and the $C_i$ are integral curves. If
	$\widehat D=(\mathcal{D},(g_\RealPlace)_{\RealPlace\mid\infty})$ is an arithmetic
	Cartier divisor on $\CX$ with $\mathcal{D}$ effective, and the generic fiber of $C$ is disjoint
	from $\mathcal{D}_\eta$, then the pairing is given by 	\begin{equation*}
		(\widehat D,C)
		=\sum_{w\nmid\fkd}(\mathcal{D},C)_w\log q_w
		+\frac{1}{2}\sum_{\RealPlace\mid\infty}
		\langle g_\RealPlace,C_\RealPlace\rangle_\RealPlace\in \BR_{\fkd}
	\end{equation*}where $
	(\mathcal{D},C)_w:=\sum_i a_i\deg_w
	\bigl(\mathcal O_{\mathcal{D}}\otimes^{\mathbb L}\mathcal O_{C_i}\bigr)$. 
	The formula extends by linearity to arbitrary arithmetic Cartier divisors.
\end{defn}

\begin{cor}\label{thm:Bruinier-modularity:cor}
 Let $C\in\wt{\CZ}_1(\SOIntModel)_{\Q,\deg=0}$, and let $\phi_1\in \sS\bigl((\Vadl\otimes_\FGlobal \BA_{\FGlobal,f})\bigr)^{K_2}$ be an admissible
 $\ov{\Q}$-valued Schwartz function.  Then the arithmetic pairing gives  
 \[
 (\widehat{\CZ}^{\mathrm B}(\tau,\phi_1),C)
 \in\CA_{\hol}(\SL_2(\FGlobal)\backslash\Mp_2(\BA_\FGlobal),\dimofV/2)_{\ov{\Q}}
       \otimes_{\ov{\Q}}\BR_{\fkd,\ov{\BQ}}.
 \]
 Finite right translation by \(h_f\in\Mp_2(\BA_{\FGlobal,f})\)
 agrees with replacing \(\phi_1\) by \(\omega_f(h_f)\phi_1\) whenever
 the latter remains admissible on the fixed integral model.
\end{cor}
\begin{proof}
	The claim follows from  \eqref{eq:admissible-degree-zero-pairing} and  \Cref{thm:Bruinier-modularity}.
\end{proof}

\subsection{Auxiliary special cycles}
Throughout this subsection, assume $\dimofV\ge 3$. Recall that a special cycle of rank $\dimofV-2=\nInLZorth-1$ on the generic fiber $\SOShim$
has dimension zero.

Let $T_\flat\in\Sym_{\nInLZorth-1}(\FGlobal)_{>0}$ and consider an admissible Schwartz function $\phi_\flat\in \sS\bigl((\Vadl\otimes_\FGlobal \BA_{\FGlobal,f})^{\dimofV-2}\bigr)^{K_2}$. 
The generic fiber of $\CZ^{\Lder}_{T_\flat}(\phi_\flat)$ is $Z_{T_\flat}(\phi_\flat)$ by
\Cref{prop:generic-special-cycle-comparison}. We therefore have a map
\[\Fil^{\dimofV-2} K_0^{\CZ_{T_\flat}(\supp \phi_\flat)}(\SOIntModel)_{\ov{\Q}}\to \wt{\CZ}_1(\SOIntModel)_{\ov{\Q}}\]
and write $\CZ^{\Lder}_{T_\flat}(\phi_\flat)\in \wt{\CZ}_1(\SOIntModel)_{\ov{\Q}}$ for the corresponding image.   We now modify it to get a degree $0$-cycle. 
\begin{lem}\label{lem:SO-Sh-connected}
	The Shimura variety $\SOShim/\FGlobal$ is connected. 
\end{lem}

\begin{proof}
 By \Cref{lem:GSpin-SO-surj} the GSpin Shimura variety $\ShG_{K_{\GSpin}}$ maps surjectively to the SO variety, 
 so it suffices to prove connectedness for the former.  Consider the
 spinor similitude character \(\nu:\Res_{\FGlobal/\Q}\GSpin(\Vcoherent)\to T_{\FGlobal}\) of \Cref{subsec:SO-shimura-datum}.  Strong approximation for the Spin group identifies the geometric component set of $\ShG_{K_{\GSpin}}$ with the
 finite quotient
 \[
 \FGlobal_+^\times\backslash\BA_{\FGlobal,f}^\times/
                  \nu(K_{\GSpin}).
 \]
For the cocharacter of \(T_\FGlobal\)
 induced by \(\widetilde h_z\) in \Cref{defn:GSpin-datum}, the reflex norm
 \[
 r_\nu:\Res_{\FGlobal/\Q}\Gm\longrightarrow
 T_\FGlobal=\Res_{\FGlobal/\Q}\Gm
 \]
 is the identity. Thus \(\Gal_\FGlobal\) acts transitively on the
 geometric components.  The canonical model is therefore connected
 over \(\FGlobal\).
\end{proof}

We next choose an auxiliary derived special one-cycle of nonzero generic
degree. Subtracting a rational multiple of this cycle makes the generic
degree zero before applying the admissible arithmetic intersection pairing.

\begin{prop}
	\label{prop:aux-special-cycle}
	After enlarging $\fkd$ by finitely many primes away from $p$ and
	shrinking the product level $K_2$ at places dividing $\fkd$ if necessary, there exist
	\[
	T_\aux\in\Sym_{\nInLZorth-1}(\FGlobal)_{>0}
	\quad\text{and}\quad
	\phi_\aux=\phi_\aux^\fkd\otimes\phi_{\aux,\fkd}
	\in\sS\bigl((\Vadl\otimes_\FGlobal \BA_{\FGlobal,f})^{\nInLZorth-1}\bigr)^{K_2}
	\]
	with the following properties.
	\begin{enumerate}
		\item $\phi_\aux$ is admissible, rational-valued, and nonnegative;
		\item the weighted special cycle $Z_{T_\aux}(\phi_\aux)$ on
		$\SOShim$ is a nonzero effective zero-cycle, so
		\begin{equation*}
			\deg Z_{T_\aux}(\phi_\aux)>0;
		\end{equation*}
		\item For every $v\nmid\fkd$, one has
		$T_\aux\in\Sym_{\nInLZorth-1}(O_{\FGlobal_v})$, and its
		$(\nInLZorth-1,\nInLZorth-1)$-entry lies in $O_{\FGlobal_v}^{\times}$.
	\end{enumerate}
\end{prop}

\begin{proof}
	The real signatures allow positive definite tuples of length
	\(\nInLZorth-1\), and each self-dual lattice \(\Lambda_v\),
	\(v\mid p\), contains a vector of unit norm because \(p\neq2\).
	Weak approximation for these nonempty open conditions gives
	\(\mathbf x_\aux\in\Vcoherent(\FGlobal)^{\nInLZorth-1}\) with
	\begin{gather*}
	T_\aux:=T(\mathbf x_\aux)\in\Sym_{\nInLZorth-1}(\FGlobal)_{>0},\\
	\mathbf x_{\aux,v}\in\Lambda_v^{\nInLZorth-1},
	\quad Q(x_{\aux,\nInLZorth-1})\in O_{\FGlobal_v}^{\times}
	\quad(v\mid p).
	\end{gather*}
	We may simultaneously impose prescribed nonempty open conditions at
	finitely many additional finite places.
	Enlarge \(\fkd\), only at primes away from \(p\), so that for every
	\(v\nmid\fkd\) the tuple \(\mathbf x_{\aux,v}\) is integral,
	\(T_\aux\in\Sym_{\nInLZorth-1}(O_{\FGlobal_v})\), the level is the standard self-dual one,
	and \(Q(x_{\aux,\nInLZorth-1})\in O_{\FGlobal_v}^{\times}\).  
	
	For \(v\nmid\fkd\), set
	\(\phi_{\aux,v}=\mathbf1_{\Lambda_v^{\nInLZorth-1}}\).  At each
	\(v\mid\fkd\), shrink $K_{2,v}$ and choose a $K_{2,v}$-invariant compact open neighborhood
	\(U_v\) of \(\mathbf x_{\aux,v}\), and set
	\(\phi_{\aux,v}=\mathbf1_{U_v}\).  Thus
	\(\phi_\aux:=\bigotimes_v\phi_{\aux,v}\)
	is \(K_2\)-invariant, admissible, rational-valued and nonnegative, and
	\(\phi_\aux(\mathbf x_\aux)>0\).
	
	By the definition of the weighted special cycle, the summand
	corresponding to \((\mathbf x_\aux,1)\) occurs with positive
	coefficient. The orthogonal complement of \(\mathbf x_\aux\) has dimension \(2\) and is
	negative definite at the distinguished real place, so the associated
	Shimura variety is zero-dimensional and nonempty.  Therefore
	\(Z_{T_\aux}(\phi_\aux)\) is a nonzero effective zero-cycle, and hence
	\(\deg Z_{T_\aux}(\phi_\aux)>0.\)
	The integrality and unit condition in (3) were built into the choice of
	\(\fkd\).  
\end{proof}

\begin{cor}\label{cor:aux-degree-zero}
	Let $T_\flat\in\Sym_{\nInLZorth-1}(\FGlobal)_{>0}$ and let $\phi_\flat$ be a
	rational-valued admissible Schwartz function. There is a rational number
	$c$ such that
	\begin{equation*}
		\CZ^{\Lder}_{T_\flat}(\phi_\flat)-c\CZ^{\Lder}_{T_\aux}(\phi_\aux)
		\in \wt{\CZ}_1(\SOIntModel)_{\Q,\deg=0}.
	\end{equation*}
	More precisely, one may take
	\begin{equation}\label{eq:aux-c}
		c=
		\frac{\deg Z_{T_\flat}(\phi_\flat)}
		{\deg Z_{T_\aux}(\phi_\aux)}\in\Q.
	\end{equation}
	The same $c$ makes the degree zero on every geometric connected component.
\end{cor}

\begin{proof}
	The denominator in \eqref{eq:aux-c} is nonzero by
	\Cref{prop:aux-special-cycle}, so $c$ is well defined and rational. 
	By \Cref{lem:SO-Sh-connected}, the Galois group of \(\FGlobal\)
	acts transitively on the geometric connected components of the
	generic fiber.  Each of the two weighted cycles is defined over
	\(\FGlobal\), so its degree is the same on every geometric
	connected component.  Equation \eqref{eq:aux-c} makes their
	difference have total degree zero, and hence degree zero on every
	geometric connected component.
\end{proof}

\subsection{Modularity of the archimedean error series}
\label{subsec:arch-degree-zero}

Choose \(T_\aux,\phi_\aux\) and \(c\in\BQ\) as in
\Cref{prop:aux-special-cycle,cor:aux-degree-zero}, and define the degree zero $1$-cycle 
\begin{equation*}
 C_{T_\flat}
 :=\CZ^{\Lder}_{T_\flat}(\phi_\flat)-c\CZ^{\Lder}_{T_\aux}(\phi_\aux)\in \wt{\CZ}_1(\SOIntModel)_{\Q,\deg=0}.
\end{equation*} 
Fix a common lower
parameter \(Y_\flat>0\).
By \Cref{prop:higher-Kudla-Green}, the pair
\[
\left(\CZ^{\Lder}_{T_\flat}(\phi_\flat),
\bigl(g_\nu^{\mathrm K}(T_\flat,Y_\flat,\phi_\flat)\bigr)_{\nu\mid\infty}\right)
\]
represents an arithmetic Chow class of codimension $\dimofV-2$, and
the same construction applies to $\CZ^{\Lder}_{T_\aux}(\phi_\aux)$.
The admissible pairing uses the one-cycle components; the Green
currents enter \Cref{defn:modularity-I} below.
We use the Kudla--Millson theta form and its Fourier coefficients
from \Cref{defn:KM-theta}.
\begin{defn}\label{defn:modularity-I}
 For every \(\phi_1\in\sS(\Vincoherent_f)^{K_2}\), define
 \begin{align*}
 I_{C_{T_\flat},\nu}(\tau,Y_\flat;\phi_1)
 :={}&-\frac12\sum_{\xi\in\FGlobal}
 \CG^{\mathrm{B-K}}_\nu(\xi,y,\phi_1)((C_{T_\flat})_\nu)q^{\xi/2}
 \\
 &+\frac12\int_{[\SOShim\times_{\FGlobal,\nu}\C]}
 \Theta_{\mathrm{KM},\nu}(\tau,\phi_1)\wedge
 \bigl(g_\nu^{\mathrm{K}}(T_\flat,Y_\flat,\phi_\flat)
       -c\,g_\nu^{\mathrm{K}}(T_\aux,Y_\flat,\phi_\aux)\bigr).
 \end{align*}
 The difference \(\CG^{\mathrm{B-K}}_\nu\) is smooth, so its evaluation
 on \((C_{T_\flat})_\nu\) is defined without a disjointness assumption.
 The first term is smooth and modular by
 \Cref{thm:Kudla-Bruinier-modularity}, since \((C_{T_\flat})_\nu\) has
 degree zero on each component.  By the modularity of the Kudla--Millson theta form (cf. \Cref{thm:KM-theta-modularity}) and the compactness of
 $\SOShim\times_{\FGlobal,\nu}\C$, the second term is smooth and modular. 
Thus \(I_{C_{T_\flat},\nu}(\tau,Y_\flat;\phi_1)\) belongs to
 \(\CA_{\infty}(\SL_2(\FGlobal)\backslash\Mp_2(\BA_\FGlobal),\dimofV/2)\).
\end{defn}

For the remainder of this section, assume that
\(\phi_1\otimes\phi_\flat\) and \(\phi_1\otimes\phi_\aux\)
are admissible factorizable Schwartz functions in $\sS((\Vadl\otimes_\FGlobal\BA_{\FGlobal,f})^{\dimofV-1})$ and are non-degenerate at the same two distinct finite
places dividing \(\fkd\).  Every tuple whose first vector is zero has
singular Gram matrix, so \(\phi_1(0)\phi_\aux=0\).
Since \(\phi_\aux\ne0\) by \Cref{prop:aux-special-cycle}, we have
\(\phi_1(0)=0\).
Since \(\FGlobal\ne\Q\), the space \(V^{(\nu)}\) is positive definite
at a real place different from \(\nu\), hence is anisotropic over
\(\FGlobal\).  Therefore the zero coefficients of
\(\Theta_{\mathrm{KM},\nu}\), \(\CG^{\mathrm{B-K}}_\nu\), and
\(\widehat{\CZ}^{\mathrm B}\) vanish.
For an arithmetic divisor \((D,G)\) and a derived one-cycle \(C\) with
\(\supp(C_\nu)\cap\supp(D_\nu)=\varnothing\), the local pairing is
\begin{equation*}
 ((D,G),C)_\nu:=\frac12
 \int_{[\SOShim\times_{\FGlobal,\nu}\C]}G\wedge\delta_{C_\nu}.
\end{equation*}
Non-degeneracy ensures disjointness for the special-cycle pairings below:
an intersection would give \(\dimofV-1\) vectors perpendicular to a
negative two-plane, with singular Gram matrix.
Pairings with generating series are taken coefficientwise.

\begin{prop}
 \label{prop:arch-degree-zero-comparison}
 Assume that \(\phi_1\otimes\phi_\flat\) and
 \(\phi_1\otimes\phi_\aux\) are factorizable and non-degenerate at
 the same two distinct finite places dividing \(\fkd\), as above.
 For every \(\nu\mid\infty\) and every \(\xi\in\FGlobal\), the Fourier coefficient is given by
 \begin{align*}
 [I_{C_{T_\flat},\nu}(\tau,Y_\flat;\phi_1)]_{q^{\xi/2}}
 ={}&\FJ_{T_\flat,\nu,\xi}(y,Y_\flat;\phi_1\otimes\phi_\flat)
 \notag\\
 &-c\,\FJ_{T_\aux,\nu,\xi}(y,Y_\flat;\phi_1\otimes\phi_\aux)
 -(\widehat{\CZ}^{\mathrm B}(\xi,\phi_1),C_{T_\flat})_\nu. 
 \end{align*}
 All terms in this identity vanish at \(\xi=0\).
\end{prop}
\begin{proof} 
 First let \(\xi=0\).  The zero coefficients just described give
 \([I_{C_{T_\flat},\nu}(\tau,Y_\flat;\phi_1)]_{q^0}=0\).
 A nonsingular \(T_b(0)\) cannot be represented at a positive definite
 real place \(\mu\ne\nu\).
 Equation \eqref{eq:FJ-local-coefficient} therefore makes both
 Fourier--Jacobi terms vanish, proving the assertion at \(\xi=0\).

 Assume now \(\xi\ne0\).
 We work with \((T_\flat,\phi_\flat)\); the same calculation applies to \((T_\aux,\phi_\aux)\).
We use the matrices \(T_b(\xi)\) defined in 
\Cref{subsec:diagonal-FJ}, with
\(b\in\Mat_{1,\nInLZorth-1}(\FGlobal)\) and
\(\sy=\operatorname{diag}(y,Y_\flat)\).  
Combining \eqref{eq:FJ-local-coefficient} with
\Cref{thm:garcia-sankaran-archimedean} and simplifying the volume factor
using \Cref{thm:GS-finite-mass}\eqref{eq:GS-finite-mass}, we obtain
\begin{align*}
	&\FJ_{T_\flat,\nu,\xi}(y,Y_\flat;\phi_1\otimes\phi_\flat)
	\quad=\frac12
	\sum_{\substack{b\in\Mat_{1,\nInLZorth-1}(\FGlobal)\\
			\det T_b(\xi)\ne0\\
			\Diff(T_b(\xi),\Vincoherent)=\{\nu\}}}
	\int_{[\SOShim\times_{\FGlobal,\nu}\C]}
	g_\nu^{\mathrm{K}}(T_b(\xi),\sy,\phi_1\otimes\phi_\flat).	
\end{align*}

By \cite[Theorem~4.5.1]{garcia_sankaran_2019_green_forms}, we have
\begin{multline*}
	\sum_{\substack{b\in\Mat_{1,\nInLZorth-1}(\FGlobal)\\
	\det T_b(\xi)\ne0}}
	g_\nu^{\mathrm{K}}(T_b(\xi),\sy,\phi_1\otimes\phi_\flat)\\
	\equiv\CG^{\mathrm K}_\nu(\xi,y,\phi_1)
	\wedge\delta_{Z_{T_\flat}(\phi_\flat)_\nu}
	+\Theta_{\mathrm{KM},\nu,\xi}(y,\phi_1)
	\wedge g_\nu^{\mathrm{K}}(T_\flat,Y_\flat,\phi_\flat)
	\quad\bmod(\operatorname{im}\partial+\operatorname{im}\bar\partial).
\end{multline*}
Moreover, by \eqref{eq:higher-Kudla-Green-current} only terms with \(\Diff(T_b(\xi),\Vincoherent)=\{\nu\}\) contribute. 
Integrating and applying Stokes' formula gives
\begin{align*}
	&\sum_{\substack{b\in\Mat_{1,\nInLZorth-1}(\FGlobal)\\
			\det T_b(\xi)\ne0\\
			\Diff(T_b(\xi),\Vincoherent)=\{\nu\}}}
	\int_{[\SOShim\times_{\FGlobal,\nu}\C]}
	g_\nu^{\mathrm{K}}(T_b(\xi),\sy,\phi_1\otimes\phi_\flat)\\
	&\qquad=\CG^{\mathrm{K}}_\nu(\xi,y,\phi_1)(Z_{T_\flat}(\phi_\flat)_\nu)+\int_{[\SOShim\times_{\FGlobal,\nu}\C]}
	\Theta_{\mathrm{KM},\nu,\xi}(y,\phi_1)\wedge
	g_\nu^{\mathrm{K}}(T_\flat,Y_\flat,\phi_\flat).
\end{align*}The claim follows.  
\end{proof}

\subsection{A holomorphic correction of the archimedean series}
\label{subsec:arch-error-holomorphy}
\begin{defn}\label{def:error-series-J}
Using the Fourier--Jacobi series of \Cref{subsec:diagonal-FJ}, define
 \begin{align*}
 J_{T_\flat}(\tau,Y_\flat;\phi_1,\phi_\flat)
 :={}&\sum_{\nu\mid\infty}I_{C_{T_\flat},\nu}(\tau,Y_\flat;\phi_1)
 \notag\\
 &-\FJ_{T_\flat}(\tau,Y_\flat;\phi_1\otimes\phi_\flat)
 +c\,\FJ_{T_\aux}(\tau,Y_\flat;\phi_1\otimes\phi_\aux).
 \end{align*}
By \Cref{defn:modularity-I} and
\Cref{defn:Fourier-Jacobi-series}, $J_{T_\flat}(\tau,Y_\flat;\phi_1,\phi_\flat)\in \CA_{\infty}(\SL_2(\FGlobal)\backslash\Mp_2(\BA_\FGlobal),\dimofV/2)$.
By \Cref{prop:arch-degree-zero-comparison}, its $\xi$-Fourier coefficients for any $\xi\in \FGlobal$ are
\begin{align}
	[J_{T_\flat}(\tau,Y_\flat;\phi_1,\phi_\flat)]_{q^{\xi/2}}
	=&-\sum_{\nu\mid\infty}
	(\widehat{\CZ}^{\mathrm B}(\xi,\phi_1),C_{T_\flat})_\nu
	\notag\\
	&-\sum_{v<\infty}\Bigl(
	\FJ_{T_\flat,v,\xi}(y,Y_\flat;\phi_1\otimes\phi_\flat)
	-c\,\FJ_{T_\aux,v,\xi}(y,Y_\flat;\phi_1\otimes\phi_\aux)
	\Bigr).
	\label{eq:J-holomorphic-decomposition}
\end{align}

\end{defn}
The following is the analogue of
\cite[Proposition~10.4.2]{Luo2025KudlaRapoport}.
\begin{thm}\label{thm:arch-comparison}
 Assume that \(\phi_1\otimes\phi_\flat\) and
 \(\phi_1\otimes\phi_\aux\) are factorizable and non-degenerate at
 the same two distinct finite places \(w_1,w_2\) dividing \(\fkd\).
 Then
 \[
 J_{T_\flat}(\tau,Y_\flat;\phi_1,\phi_\flat)
 \in\CA_{\hol}(\SL_2(\FGlobal)\backslash\Mp_2(\BA_\FGlobal),\dimofV/2)\cong \CA_{\hol}(\SL_2(\FGlobal)\backslash\Mp_2(\BA_\FGlobal),\dimofV/2)_{\ov{\Q}}
       \otimes_{\ov{\Q}}\C.
 \]
 Its coefficients vanish for
 \(\xi\notin\FGlobal_{>0}\), and its \(q^{\xi/2}\)-coefficient
 for \(\xi>0\) is independent of \(y=\Im\tau\).
\end{thm}
\begin{proof}
	For a finite place \(v\), the condition
	\(\Diff(T_b(\xi),\Vincoherent)=\{v\}\) forces
	\(T_b(\xi)\) to be positive definite at every real place, because
	\(\Vincoherent\) is positive definite there.  Its first diagonal
	entry \(\xi\) is therefore totally positive.  By
	\eqref{eq:FJ-local-coefficient} and
	\Cref{thm:Siegel-Eisenstein}\eqref{eq:finite-fourier-normalized}, the coefficients
	\(\FJ_{T_\flat,v,\xi}(y,Y_\flat;\phi_1\otimes\phi_\flat)\)
	and \(\FJ_{T_\aux,v,\xi}(y,Y_\flat;\phi_1\otimes\phi_\aux)\)
	are independent of \(y\) and vanish unless \(\xi>0\).  The Bruinier divisor is independent of \(y\) and
	vanishes for \(\xi\notin\FGlobal_{>0}\).
	Thus
	$[J_{T_\flat}(\tau,Y_\flat;\phi_1,\phi_\flat)]_{q^{\xi/2}}$ is \(y\)-independent
	for every \(\xi\in\FGlobal\) and \(y>0\), so \(J_{T_\flat}\) is holomorphic
	on \(\mathbb H_{\FGlobal}\). As explained in
	\Cref{subsec:nota-hilbert}, it follows that
	 $J_{T_\flat}(\tau,Y_\flat;\phi_1,\phi_\flat)
	\in\CA_{\hol}(\SL_2(\FGlobal)\backslash\Mp_2(\BA_\FGlobal),\dimofV/2)$. 
\end{proof}

\subsection{Modularity of the error series}
Retain \(C_{T_\flat}\), the parameter \(Y_\flat\), and the two-place
support assumptions of \Cref{subsec:arch-degree-zero}.

\begin{defn}\label{def:error-series}
 Fix two distinct finite places \(w_1,w_2\) dividing \(\fkd\) and
 \(\phi_{1,w_i}\in\sS(\Vadl\otimes_{\FGlobal}\FGlobal_{w_i})^{K_{2,w_i}}\)
 such that \(\phi_{1,w_i}\otimes\phi_{\flat,w_i}\) and
 \(\phi_{1,w_i}\otimes\phi_{\aux,w_i}\) are non-degenerate for
 \(i=1,2\).  Define
 \begin{align}
 \CE_{T_\flat}(\tau,\phi_1,\phi_\flat)
 :=\vol(K_2)\Bigl(&
     (\widehat{\CZ}^{\mathrm B}(\tau,\phi_1),C_{T_\flat})
     +J_{T_\flat}(\tau,Y_\flat;\phi_1,\phi_\flat)\Bigr) \in\CA_{\hol}(\SL_2(\FGlobal)\backslash\Mp_2(\BA_\FGlobal),\dimofV/2)_{\ov{\BQ}}
     \otimes_{\ov{\BQ}}\C_\fkd.
 \label{eq:error-series-decomposition}
 \end{align}
 Both terms are holomorphic and modular by
 \Cref{thm:Bruinier-modularity:cor,thm:arch-comparison}. 
 We use the \(\ov{\BQ}\)-structure of the displayed holomorphic space;
 the maps \(\BR_{\fkd,\ov{\BQ}}\to\C_\fkd\) and
 \(\C\to\C_\fkd\) act on the coefficient factors.
 Finite right translations act on the holomorphic factor. Point values and Whittaker
 coefficients below are understood after applying a
 \(\ov{\BQ}\)-linear functional on \(\C_\fkd\) to the coefficient factor.
 At finite component \(1\), the classical series is recovered by
 \eqref{eq:adelization-general} coefficientwise.
\end{defn}

\begin{lem}\label{lem:bad-finite-whittaker-algebraicity}
For \(\ov{\Q}\)-valued \(K_2\)-invariant finite Schwartz data \(\varphi_f\)
that are factorizable and non-degenerate at two distinct finite places,
for every \(\xi\in\FGlobal\) and $w\in \Sigma_{\FGlobal,f}$ one has
\begin{equation*}
 \frac{\vol(K_2)}{\log q_w}
 \FJ_{T_\flat,w,\xi}(y,Y_\flat;\varphi_f)
 \in\ov{\Q}.
\end{equation*}
\end{lem}
\begin{proof}
	Write
	\(\varphi_f=\varphi_w\otimes\varphi^w\).
	By \Cref{thm:GS-finite-mass} we have \(\vol(K_2)\in\BQ_{>0}\).
	The Gelfand--Leray measure \(dh_w\) of
	\Cref{defn:Gelfand--Leray} takes algebraic values on compact open
	subsets, as follows from its density in \eqref{eq:GL-trace-density}.
	In particular, \(\vol^{\mathrm{GL}}_w(K_{2,w})\in\ov{\Q}\cap\BR_{>0}\), and
	\Cref{defn:Haar-measure} gives
	\[
	\vol(K_2^w)
	=\frac{\vol(K_2)}{\vol^{\mathrm{GL}}_w(K_{2,w})}\in\ov{\Q}\cap\BR_{>0}.
	\]
	For a nonsingular \(T\) with \(\Diff(T,\Vincoherent)=\{w\}\),
	\(\operatorname{Orb}_T(\varphi^w)/\vol(K_2^w)\) is a finite
	\(\ov{\Q}\)-valued weighted sum over \(K_2^w\)-cosets.
	Consequently \(\operatorname{Orb}_T(\varphi^w)\in\ov{\Q}\).
	By \Cref{lem:bad-finite-whittaker-rational},
	\(W'_{T/2,w}(1,0,\varphi_w)/\log q_w\in\ov{\Q}\),
	and the Weil index \(\gamma_{\Vincoherent_w,\psi_w}\) is a root of unity.
	It follows that
	\[
	 \gamma_{\Vincoherent_w,\psi_w}^{\nInLZorth}
	 \frac{W'_{T/2,w}(1,0,\varphi_w)}{\log q_w}
	 \operatorname{Orb}_T(\varphi^w)\in\ov{\Q}.
	\]
	By \eqref{eq:FJ-local-coefficient} and
	\Cref{thm:Siegel-Eisenstein}\eqref{eq:finite-fourier-normalized},
	\(\frac{\vol(K_2)}{\log q_w}
	\FJ_{T_\flat,w,\xi}(y,Y_\flat;\varphi_f)\) is twice a finite sum
	of the displayed algebraic numbers, so the claim follows.
\end{proof}

\begin{prop}\label{prop:modularity}
 For every
 \(\xi\in\FGlobal^\times\) the \(q^{\xi/2}\)-coefficient of
 \(\CE_{T_\flat}\) is
 \begin{align*}
 \CE_{T_\flat,\xi}(\phi_1,\phi_\flat)
 ={}&\sum_{w\in\Sigma_{F,f}^{\fkd}}
       2\vol^{\mathrm{GL}}_w(K_{2,w})\log q_w\\
 &\quad\Big(
 \sum_{\substack{
       T=\left(\begin{smallmatrix}\xi&*\\ *&T_\flat\end{smallmatrix}\right)
          \in\Sym_{\dimofV-1}(\FGlobal)_{>0}\\
       \Diff(T,\Vincoherent)=\{w\}}}
       \bigl(\Int_w^{\epsilon_w}(T)-\pDen_w^{\epsilon_w}(T)\bigr)
       \operatorname{Orb}_T((\phi_1\otimes\phi_\flat)^w)\\
\quad-c &
 \sum_{\substack{
       T=\left(\begin{smallmatrix}\xi&*\\ *&T_\aux\end{smallmatrix}\right)
          \in\Sym_{\dimofV-1}(\FGlobal)_{>0}\\
       \Diff(T,\Vincoherent)=\{w\}}}
       \bigl(\Int_w^{\epsilon_w}(T)-\pDen_w^{\epsilon_w}(T)\bigr)
       \operatorname{Orb}_T((\phi_1\otimes\phi_\aux)^w)
       \Big)\in \C_\fkd.
 \end{align*}
Moreover, the constant coefficient of $\CE_{T_\flat}$ vanishes.
The coefficient formula shows that $\CE_{T_\flat}$ is independent of
$Y_\flat$; the auxiliary data remain fixed.
\end{prop}

\begin{proof}
 Combining \eqref{eq:error-series-decomposition} and
 \eqref{eq:J-holomorphic-decomposition} gives
 \begin{align*}
[ \CE_{T_\flat}(\tau,\phi_1,\phi_\flat)]_{q^{\xi/2}}
 =\vol(K_2)\Bigl(&
 (\widehat{\CZ}^{\mathrm B}(\xi,\phi_1),C_{T_\flat})
 -\sum_{\nu\mid\infty}
 (\widehat{\CZ}^{\mathrm B}(\xi,\phi_1),C_{T_\flat})_\nu\\
 &-\sum_{v<\infty}\bigl(
 \FJ_{T_\flat,v,\xi}(y,Y_\flat;\phi_1\otimes\phi_\flat)
 -c\,\FJ_{T_\aux,v,\xi}(y,Y_\flat;\phi_1\otimes\phi_\aux)\bigr)\Bigr)
 \end{align*}
 for any $\xi\in \FGlobal$. 
 For \(\xi>0\), associativity of derived tensor products and
 \Cref{prop:prod-special} give
 \[
 \CZ^{\Lder}_{\xi}(\phi_1)\otimes^{\Lder}
 \CZ^{\Lder}_{T_\flat}(\phi_\flat)
 =\sum_{\substack{T\in\Sym_{\dimofV-1}(\FGlobal)\\
 T=\left(\begin{smallmatrix}\xi&*\\ *&T_\flat\end{smallmatrix}\right)}}
 \CZ^{\Lder}_T(\phi_1\otimes\phi_\flat).
 \]
 Applying
 \Cref{thm:global-int} to each summand computes
 the finite pairing with the derived lower cycle, so we obtain
 \[\sum_{w\in\Sigma_{F,f}^{\fkd}}
 \bigl(\CZ_\xi(\phi_1),\CZ^{\Lder}_{T_\flat}(\phi_\flat)\bigr)_w\log q_w
 =\sum_{w\in\Sigma_{F,f}^{\fkd}}
 \frac{2\log q_w}{\vol(K_2^w)}
 \sum_{\substack{T=\left(\begin{smallmatrix}\xi&*\\ *&T_\flat\end{smallmatrix}\right)>0\\
 		\Diff(T,\Vincoherent)=\{w\}}}
 \Int_w^{\epsilon_w}(T)\operatorname{Orb}_T((\phi_1\otimes\phi_\flat)^w).
 \]
  For $v\in \Sigma_{F,f}^{\fkd}$, applying \Cref{thm:Siegel-Eisenstein}\eqref{eq:orth-Luo-10-1-2-finite}
  at \(g=g_\sz\), with the normalization \eqref{eq:classical-eisenstein-normalization}, gives
 \[\FJ_{T_\flat,v,\xi}(y,Y_\flat;\phi_1\otimes\phi_\flat)=\sum_{\substack{T=\left(\begin{smallmatrix}\xi&*\\ *&T_\flat\end{smallmatrix}\right)>0\\
 		\Diff(T,\Vincoherent)=\{v\}}} \frac{2\log q_v}{\vol(K_2^v)}
 \pDen_v^{\epsilon_v}(T)\operatorname{Orb}_T((\phi_1\otimes\phi_\flat)^v).\]
 For a finite place $v\mid\fkd$, by \Cref{lem:bad-finite-whittaker-algebraicity} $\FJ_{T_\flat,v,\xi}(y,Y_\flat;\phi_1\otimes\phi_\flat)= 0$  
 as an equality in $\C_\fkd$. 
 The same calculation applies to the auxiliary lower data
 \(T_\aux,\phi_\aux\).

 For nonzero \(\xi\) not totally positive, the Bruinier arithmetic
 divisor and every finite-place Fourier--Jacobi coefficient are zero, as in the
 proof of \Cref{thm:arch-comparison}; hence the same formula, with
 empty sums, holds.  The constant coefficient vanishes by
 \Cref{prop:arch-degree-zero-comparison}.  
\end{proof}
\begin{prop}\label{prop:dual-relation}
	Let \(v_0\nmid\fkd\) be a finite place with
	\(v_0\notin S_\psi\).  Assume that
	\(\VlattiveLocalo \subset V_{v_0}\) is self-dual and
	\(\phi_{1,v_0}=\mathbf1_{\VlattiveLocalo }\).  Then the adelic tensor satisfies
	\begin{align}
		R(n_{v_0}(b))\,
		\CE_{T_\flat}(\,\cdot\,,\phi_1,\phi_\flat)
		&=
		\CE_{T_\flat}(\,\cdot\,,\phi_1,\phi_\flat),
		\qquad b\in O_{\FGlobal_{v_0}},
		\label{eq:error-NO-invariance}\\
		R(\widetilde w_{v_0})\,
		\CE_{T_\flat}(\,\cdot\,,\phi_1,\phi_\flat)
		&=
		\gamma_{V_{v_0},\psi_{v_0}}\,
		\CE_{T_\flat}(\,\cdot\,,\phi_1,\phi_\flat).
		\label{eq:error-w-eigen}
	\end{align}
	In fact, the adelizations of all the summands $I_{C_{T_\flat},\nu}(\tau,Y_\flat;\phi_1)$, $\FJ_{T_\flat}(\tau,Y_\flat;\phi_1\otimes\phi_\flat)$, $\FJ_{T_\aux}(\tau,Y_\flat;\phi_1\otimes\phi_\aux)$, $(\widehat{\CZ}^{\mathrm B}(\tau,\phi_1),C_{T_\flat})$ satisfy both displayed identities.
\end{prop}

\begin{proof}
	Because \(v_0\notin S_\psi\), the local character \(\psi_{v_0}\) has
	conductor \(O_{\FGlobal_{v_0}}\).  Since
	\(\VlattiveLocalo \) is self-dual, the
	\(\psi_{v_0}\)-self-dual measure gives
	\[
	\widehat{\mathbf1_{\VlattiveLocalo }}^{\,\psi_{v_0}}
	=
	\mathbf1_{\VlattiveLocalo }.
	\]
	Moreover \(Q(\VlattiveLocalo )\subset O_{\FGlobal_{v_0}}\), so
	\eqref{eq:weil-n} and \eqref{eq:weil-w}, applied with \(r=1\), give
	\[
	\omega_{v_0}(n(b))\mathbf1_{\VlattiveLocalo }
	=
	\mathbf1_{\VlattiveLocalo }
	\quad (b\in O_{\FGlobal_{v_0}}),
	\qquad
	\omega_{v_0}(\widetilde w_{v_0})\mathbf1_{\VlattiveLocalo }
	=
	\gamma_{V_{v_0},\psi_{v_0}}
	\mathbf1_{\VlattiveLocalo }.
	\]
	These local Weil operators multiply the chosen local Schwartz function
	by a scalar and leave the other local factors unchanged.  Thus the
	hypotheses of \Cref{thm:Bruinier-modularity:cor} remain valid.
	The Weil equivariance of that pairing and of
	\Cref{defn:modularity-I} proves \eqref{eq:error-NO-invariance} and
	\eqref{eq:error-w-eigen} for the arithmetic pairing and the
	\(I_{C_{T_\flat},\nu}\)-summands.
	
	Let $k_{v_0}=n_{v_0}(b)$ or $\widetilde w_{v_0}$, so $k_{v_0}$ belongs to
	$ K^{\mathrm{std}}_{1,v_0}$.  The image of $k_{v_0}$ under the
	first-block embedding lies in $K^{\mathrm{std}}_{\nInLZorth,v_0}$, so
	\eqref{eq:local-standard-section-Iwasawa} gives, for every $s$,
	\[
	R(i_1(k_{v_0}))\Phi_{\phi_1\otimes\phi_j}(-,s)
	=\Phi_{(\omega_{v_0}(k_{v_0})\phi_1)\otimes\phi_j}(-,s),
	\qquad j\in\{\flat,\aux\}.
	\]
	Here the local Weil operator acts on the first Schwartz factor
	at $v_0$ and leaves all other factors fixed.  Differentiating
	the Eisenstein series at $s=0$ and applying the compact integral
	\eqref{eq:FJ-adelic-operation} gives the same transformation
	law for each complete Fourier--Jacobi series.  The claim follows. 
\end{proof}

\section{Proof of the main theorem}\label{sec:proof}
\begin{thm}\label{thm:main2}
	Assume that the local quadratic space $\Vlocal$ admits a self-dual
	$O_{\FLocal}$-lattice, as in \Cref{thm:main}, and put
	$\epsilon=\chi(\Vlocal)\in\{\pm1\}$.
	Let $\bu=(u_1,\ldots,u_{\dimofV-1})\in\Vspecial^{\dimofV-1}$
	have a nonsingular Gram matrix. Then
	$\Int^\epsilon(\bu)=\pDen^\epsilon(\bu).$
\end{thm}
We prove Theorem~\ref{thm:main2} by induction on \(\dimofV=\dim \Vlocal\).

\begin{prop}\label{prop:base-case}
	Theorem~\ref{thm:main2} holds for \(\dimofV=2\).
\end{prop}

\begin{proof}
	If \(\chi(\Vlocal)=1\), then \(\Vspecial=0\) by
	\Cref{prop:V_b-SO-2}. Thus no one-tuple in \(\Vspecial\) has a nonsingular Gram matrix,
	and the statement is vacuous. Hence we may assume \(\chi(\Vlocal)=-1\). 
	Let $x\in\Vspecial$ satisfy $(x,x)\ne0$. If $(x,x)\notin O_\FLocal$, then both sides vanish. We may assume $(x,x)\in O_\FLocal$. 
	Then $\valFLocal((x,x))\in 2\BZ_{\geq0}+1$. 
	By \Cref{prop:n-two-special},
	\begin{equation}\label{eq:base-case-geometric}
		\Int^\epsilon(x)
		=\length_{O_{\breve \FLocal}}O_{\CZ_{\mathrm{SO}}(x)}
		=\frac{\valFLocal((x,x))+1}{2}.
	\end{equation}
	
	We now compute the local density. The integral lattices \(L'\) of rank one in \(\FLocal x\) with
	\(O_\FLocal x\subseteq L'\subseteq L'^\vee\) are exactly
	\[
	L_j:=O_\FLocal\varpi^{-j}x,
	\qquad 0\leq j\leq\frac{\valFLocal((x,x))-1}{2}.
	\]
	Each \(L_j\) has type one and positive odd valuation, so its reduction
	modulo \(\varpi\) is totally radical. Thus \(U_0=0\) in
	\eqref{eq:LZ-sign-convention}, and
	\(\operatorname{sgn}_2(L_j)=1\) for every \(j\).
	Substituting this value and \(\epsilon=-1\) into the density formula of
	\cite[Theorem~3.4.5]{LZ-orth} gives
	\[
	\operatorname{Den}^{\epsilon}(X,O_\FLocal x)
	=(1-X)\sum_{j=0}^{(\valFLocal((x,x))-1)/2}X^{2j}.
	\]
	It follows that
	\begin{equation}\label{eq:base-case-analytic}
		\pDen^\epsilon(x)
		=-\left.\frac{d}{dX}\operatorname{Den}^{\epsilon}(X,O_\FLocal x)
		\right|_{X=1}
		=\frac{\valFLocal((x,x))+1}{2}.
	\end{equation}
	Comparing \eqref{eq:base-case-geometric} and
	\eqref{eq:base-case-analytic} proves the proposition.
\end{proof}

\subsection{Cancellation at boundary cases}
Assume \(\dimofV\ge3\) and that
\Cref{thm:main2} holds in dimension \(\dimofV-1\) for every local
datum satisfying its hypotheses. 

\begin{lem}
	\label{lem:local-discrepancy-support}
	Let \(\bu\in\Vspecial^{\dimofV-1}\) have a nonsingular Gram matrix.
	If \(T(\bu)\notin\Sym_{\dimofV-1}(O_\FLocal)\), then
	$\Int^\epsilon(\bu)=\pDen^\epsilon(\bu)=0.$
\end{lem}
\begin{proof}
	On a GSpin cover, the vectors in \(\bu\) act on the integral
	Dieudonn\'e module at every geometric point of the special cycle.
	Their squares are scalar endomorphisms, so their norms lie in
	\(O_\FLocal\).  Since the residue characteristic is odd,
	polarization gives \(T(\bu)\in\Sym_{\dimofV-1}(O_\FLocal)\).
	Thus a nonintegral Gram matrix gives an empty special cycle and
	zero intersection number.  The density derivative is zero by the
	nonintegral-lattice convention in \Cref{defn:local-density}.
\end{proof}

\begin{lem}\label{lem:local-unit-cancellation}
	Assume Theorem~\ref{thm:main2} in dimension \(\dimofV-1\).
	Let
	\(\mathbf y=(y_1,\ldots,y_{\dimofV-1})\in
	\Vspecial^{\dimofV-1}\) have a nonsingular Gram matrix, and suppose that
	\((y_i,y_i)\in O_\FLocal^\times\) for some \(i\).  Then
	$\Int^\epsilon(\mathbf y)=\pDen^\epsilon(\mathbf y).$
\end{lem}
\begin{proof}
	If the Gram matrix has a nonintegral entry, both sides are zero by
	\Cref{lem:local-discrepancy-support}.  Otherwise permute the tuple
	so that \(i=1\), and define
	\[
	 y_j^\flat=y_j-\frac{(y_j,y_1)}{(y_1,y_1)}y_1,
	 \qquad 2\le j\le\dimofV-1,
	\]
	with \(\mathbf y^\flat=(y_2^\flat,\ldots,y_{\dimofV-1}^\flat)\).
	Its Gram matrix is the Schur complement of the unit entry
	\((y_1,y_1)\), hence is integral and nonsingular.
	By \Cref{lem:unit-special-vector-gauge}, the lower coherent space
	\(\Vlocal^\flat\) has the self-dual lattice
	\(\VlattiveLocaltot^\flat\), and its space of special
	quasi-endomorphisms is identified with \(y_1^\perp\).
	Thus the lower tuple satisfies
	the hypotheses of the induction statement.
	Geometric cancellation
	\Cref{cor:cancellation} gives
	\(\Int^\epsilon(\mathbf y)=\Int^{\chi(\Vlocal^\flat)}(\mathbf y^\flat)\).
	On the analytic side, the cancellation formula
	\cite[(3.4.8.1)]{LZ-orth} gives
	\(\pDen^\epsilon(\mathbf y)=\pDen^{\chi(\Vlocal^\flat)}(\mathbf y^\flat)\).
	The induction hypothesis identifies the two lower-dimensional terms.
\end{proof}

\subsection{Globalization and isolation}
Recall that we have chosen a local field $\FLocal$ of odd residue characteristic, and a quadratic space $\Vlocal$ admitting a self-dual
$O_{\FLocal}$-lattice $\VlattiveLocaltot$.
Choose a tuple \(\bu=(u_1,\ldots,u_{\dimofV-1})\in\Vspecial^{\dimofV-1}\) with a nonsingular Gram matrix. The goal of this subsection is to globalize our local data and isolate \(\Int^\epsilon(\bu)-\pDen^\epsilon(\bu)\) in a Fourier coefficient of the error series. 
\begin{prop}\label{prop:globalization}
	The following global data exist. 
	\begin{enumerate}[(F1)]
		\item a totally real field \(\FGloball\neq \Q\) and a finite place \(v_0\) with
		\(\FGloball_{v_0}\simeq \FLocal\);
		
		\item a global additive character \(\psi:\FGlobal\backslash\BA_{\FGlobal}\to\BC^\times\) that is unramified at every \(v\mid p\);
		
		\item an even positive integer \(\fkd\), prime to \(p\), divisible by
		all rational primes lying below \(S_\psi\);
		
		\item a quadratic space \(\Vcoherent/\FGloball\) of dimension
		\(\dimofV\), of signature \((\dimofV-2,2)\) at one real place $\RealPlace_0$ and
		positive definite at every other real place, whose localization at $v_0$ is $\Vlocal$. It contains a global lattice $\VlattiveGlobal$ which is self-dual away from the
		finite set of places dividing \(\fkd\), with \(\VlattiveLocalo\cong \VlattiveLocaltot\).
		
		\item a sufficiently small neat product level
		\(K_2=\prod_{v<\infty}K_{2,v}\), with
		\(K_{2,v}=\SO(\VlattiveGlobal\otimes_{O_\FGlobal}O_{\FGlobal_v})(O_{\FGlobal_v})\)
		for \(v\nmid\fkd\);
		
		\item \label{item:globalization-isolation-I1}
		a tuple
		\(\bu'\in\VspecialGlobal^{\dimofV-1}\), together with an isometry
		\((\VspecialGlobal)_{v_0}\simeq\Vspecial\), such that
		\[
		\Int^\epsilon(\bu)=\Int^\epsilon(\bu'_{v_0}),\qquad \pDen^\epsilon(\bu)=\pDen^\epsilon(\bu'_{v_0}),
		\]
		where the local tuple \(\bu'_{v_0}\) is viewed in \(\Vspecial^{\dimofV-1}\) through this isometry.
		The Gram matrix has the block decomposition
		\[
		\RT=T(\bu')=
		\begin{pmatrix}t&*\\ *&T_\flat\end{pmatrix},
		\]
		where \(t=(u'_1,u'_1)\) and
		\(T_\flat=T(u'_2,\ldots,u'_{\dimofV-1})\).
		Its last diagonal entry is a unit at every place
		\(v\nmid\fkd\), \(v\neq v_0\).  
		
		\item \label{item:globalization-isolation-I2}
		an auxiliary pair \((T_\aux,\phi_\aux)\) satisfying
		\Cref{prop:aux-special-cycle}, with \(\phi_\aux\) factorizable,
		and admissible factorizable rational-valued Schwartz functions
		$\phi_1\in \sS\bigl((\Vadl\otimes_\FGlobal \BA_{\FGlobal,f})\bigr)^{K_2}$ and $\phi_\flat\in \sS\bigl((\Vadl\otimes_\FGlobal \BA_{\FGlobal,f})^{\nInLZorth-1}\bigr)^{K_2}$, such that both $\phi=\phi_1\otimes \phi_\flat$ and $\phi_1\otimes \phi_\aux$ are non-degenerate at the same two 	places dividing $\fkd$, and 
		$\operatorname{Orb}_{\RT}(\phi^{v_0})\ne0,$
		and such that, for every positive definite
		\(\widetilde T\ne\RT\) of the form
		\[
		\widetilde T=
		\begin{pmatrix}t&*\\ *&T_\flat\end{pmatrix},
		\qquad
		\Diff(\widetilde T,\Vincoherent)=\{v_0\},
		\]
		one has
		\[	\bigl(
		\Int^\epsilon_{v_0}(\widetilde T)
		-\pDen^\epsilon_{v_0}(\widetilde T)
		\bigr)
		\operatorname{Orb}_{\widetilde T}(\phi^{v_0})
		=0;\]		
	\end{enumerate}
\end{prop}
\begin{proof}
First choose a totally real field \(\FGloball\neq\Q\) with a finite
place \(v_0\) and an identification \(\FGloball_{v_0}\simeq\FLocal\).
When \([\FLocal:\Q_p]>1\), weak approximation and Krasner's lemma give
a monic polynomial over \(\Q\) defining \(\FLocal\) over \(\Q_p\)
and with all roots real and distinct. The field generated by a root has
the required properties. When \(\FLocal=\Q_p\), use weak approximation
to choose a monic quadratic over \(\Q\) sufficiently close to \(X(X-1)\) over \(\Q_p\),
to \(X^2-1\) over \(\R\), and to an Eisenstein quadratic at another
rational prime. Its quadratic field is totally real and \(p\) splits in it.
We fix a real place \(\RealPlace_0\) of \(\FGloball\).

Identify the given lattice \(\VlattiveLocaltot\) with
\(\VlattiveLocalo\) under \(\FGloball_{v_0}\simeq\FLocal\).
At every other place \(v\mid p\),
choose a quadratic space of dimension \(\dimofV\) with a self-dual
lattice \(\VlattiveLocalv\).  Apply weak approximation to Gram matrices,
imposing signature \((\dimofV-2,2)\) at \(\RealPlace_0\), positive
definiteness at the other real places, and the prescribed lattice
isometry classes at the places above \(p\).  This gives a global
quadratic lattice \(\VlattiveGlobal\subset\Vcoherent\) whose
completion at \(v_0\) is isomorphic to \(\VlattiveLocalo\), and
whose completion at each \(v\mid p\), \(v\ne v_0\), is isomorphic to
\(\VlattiveLocalv\).  This lattice is
self-dual at every place above \(p\) and at all but finitely many
finite places. 

 Choose the global additive character \(\psi\) to be
unramified at every \(v\mid p\); weak approximation for the totally
positive scalar in \eqref{eq:additive-character-scalar} allows this
choice.  We now choose an even positive integer \(\fkd\), prime to
\(p\), containing all rational primes below \(S_\psi\) and below the
places where \(\VlattiveGlobal\) fails to be self-dual.  The associated incoherent space
\(\Vincoherent\) and its nearby space \(\VspecialGlobal\) at
\(v_0\) are then as in \Cref{sec:notat-quadr-spac}; in particular,
\((\VspecialGlobal)_{v_0}\simeq\Vspecial\), and
\((\VspecialGlobal)_v\simeq\Vcoherent_v\) at every finite
\(v\neq v_0\).

With these choices fixed, apply weak approximation in
\(\VspecialGlobal^{\dimofV-1}\) to choose \(\bu'\) sufficiently
close to \(\bu\) at \(v_0\).  At each other \(v\mid p\), require
\(\bu'\) to lie in \(\VlattiveLocalv^{\dimofV-1}\), with
nonsingular Gram matrix and last vector of unit norm.  Such local
tuples exist because \(\VlattiveLocalv\) is self-dual and \(p\)
is odd, and these conditions are open.  The local constancy results
\Cref{thm:local-constancy,prop:analytic-locally-constant} ensure
\[
\Int^\epsilon(\bu'_{v_0})=\Int^\epsilon(\bu),\qquad
\pDen^\epsilon(\bu'_{v_0})=\pDen^\epsilon(\bu).
\]
Since \(\VspecialGlobal\) is totally positive definite,
\(\RT=T(\bu')\) is totally positive definite and satisfies
\(\Diff(\RT,\Vincoherent)=\{v_0\}\).  Enlarge \(\fkd\) to contain
the rational primes below the finitely many places \(v\neq v_0\)
where \(\bu'\) fails to lie in the completion of
\(\VlattiveGlobal^{\dimofV-1}\), or where the last diagonal entry
of \(\RT\) fails to be a unit.  This proves \eqref{item:globalization-isolation-I1}.

Choose a sufficiently small neat product level \(K_2\) with
\(K_{2,v}=\SO(\VlattiveGlobal\otimes_{O_\FGlobal}O_{\FGlobal_v})(O_{\FGlobal_v})\)
for \(v\nmid\fkd\).
Choose preliminary factorizable rational-valued Schwartz
data using the characteristic functions of the chosen self-dual
lattices at \(v\nmid\fkd\), and characteristic
functions of compact open product neighborhoods of \(\bu'\) at
the places dividing \(\fkd\).  Their orbital integral at \(\RT\)
away from \(v_0\) is nonzero.  
Choose two finite places away from \(p\fkd\) where the global data
are unramified, and enlarge \(\fkd\) to contain the rational primes
below them.  Shrink the
product neighborhoods around the localizations of \(\bu'\) at these two
places so that all their Gram matrices are nonsingular.
The preliminary Schwartz data are then non-degenerate at these two places,
and their orbital integral at \(\RT\) away from \(v_0\) remains nonzero.

In the construction of \Cref{prop:aux-special-cycle}, impose the
lower neighborhoods just chosen at the same two auxiliary places.
Weak approximation gives the auxiliary tuple in these neighborhoods,
with totally positive Gram matrix and with its last vector of unit norm at every place above \(p\).
After shrinking the product level at the places dividing \(\fkd\), choose the auxiliary lower Schwartz
support inside these neighborhoods. Then \(\phi_1\otimes\phi_\aux\)
is non-degenerate at the same two places, as required in
\Cref{subsec:arch-degree-zero}. Applying
\Cref{prop:isolate-contributing-indices} gives \eqref{item:globalization-isolation-I2} by restricting the
support at a further auxiliary place away from \(p\) and these two
places, shrinking its local level, and enlarging \(\fkd\) accordingly. 
\end{proof}

\begin{prop}
	\label{prop:isolate-contributing-indices}
	Assume that
	\[
	\RT=
	\begin{pmatrix}t&*\\ *&T_\flat\end{pmatrix}
	\in\Sym_{\dimofV-1}(\FGloball)_{>0}
	\]
	satisfies \(\Diff(\RT,\Vincoherent)=\{v_0\}\).  Suppose that a preliminary
	admissible factorizable rational-valued Schwartz function
	\(\phi^\circ=\phi_1^\circ\otimes\phi_\flat^\circ\) has
	$\operatorname{Orb}_{\RT}((\phi^\circ)^{v_0})\ne0,$
	and has the standard self-dual component at \(v_0\).  After enlarging
	\(\fkd\) away from \(p\) and replacing \(K_2\) by a smaller product level
	at one further auxiliary finite place \(w_3\mid\fkd\), one may alter the
	Schwartz factors at \(w_3\) and obtain an admissible
	factorizable rational-valued Schwartz function
	\(\phi=\phi_1\otimes\phi_\flat\) such that
	$\operatorname{Orb}_{\RT}(\phi^{v_0})\ne0$
	and, for every
	\[
	\widetilde T=
	\begin{pmatrix}t&*\\ *&T_\flat\end{pmatrix}
	\in\Sym_{\dimofV-1}(\FGloball)_{>0},
	\qquad
	\widetilde T\ne\RT,
	\]
	with \(\Diff(\widetilde T,\Vincoherent)=\{v_0\}\), one has
	\[
	\bigl(\Int^\epsilon_{v_0}(\widetilde T)-\pDen^\epsilon_{v_0}(\widetilde T)\bigr)\,
	\operatorname{Orb}_{\widetilde T}(\phi^{v_0})=0.
	\]
\end{prop}

\begin{proof}
	Consider the totally positive definite matrices
	\(\widetilde T=\left(\begin{smallmatrix}
	t&z\\{}^tz&T_\flat\end{smallmatrix}\right)\), with
	\(z\in\Mat_{1,\dimofV-2}(\FGlobal)\), satisfying
	\begin{equation}\label{eq:isolation-support-condition}
	(\widetilde T_v)_{v<\infty}
	\in\Sym_{\dimofV-1}(O_{\FGlobal_{v_0}})
	\times T\!\left(\supp((\phi^\circ)^{v_0})\right).
	\end{equation}
	The Gram map in the second factor is evaluated on
	\((\Vcoherent\otimes_\FGlobal\BA_{\FGlobal,f}^{v_0})^{\dimofV-1}\).
	The right-hand side is compact, so
	\Cref{prop:finite-contributing-indices} gives only finitely many
	totally positive definite completions of \(t\) and \(T_\flat\)
	satisfying \eqref{eq:isolation-support-condition}.
	
	Choose a finite place \(w_3\nmid p\fkd\) where the global data are unramified.
	Enlarge \(\fkd\) to contain the rational prime below \(w_3\).
	The nonzero preliminary orbital integral supplies a local tuple
	\(\mathbf x^\circ_{w_3} = (x^\circ_{1,w_3},\mathbf x^\circ_{\flat,w_3})\)
	with Gram matrix \(\RT_{w_3}\) at which both local factors of
	\(\phi^\circ_{w_3}\) are nonzero.
	
	Since \(\FGlobal\hookrightarrow\FGlobal_{w_3}\) is injective, choose
	a compact open neighborhood
	\(B\subset\Sym_{\dimofV-1}(\FGlobal_{w_3})\)
	of \(\RT_{w_3}\) containing none of the \(w_3\)-components of
	the finitely many \(\widetilde T\ne\RT\) satisfying
	\eqref{eq:isolation-support-condition}.  By continuity of the Gram map
	and local constancy of the factors of \(\phi^\circ_{w_3}\), choose
	compact open neighborhoods
	$U_1\ni x^\circ_{1,w_3}$, $U_\flat\ni\mathbf x^\circ_{\flat,w_3}$ 
	contained in their respective nonzero loci and so small that
	\(T(U_1\times U_\flat)\subset B.\)
	After shrinking the level at \(w_3\), we may assume that \(U_1\) and \(U_\flat\) are
	stable under the new local level.  Replace the two local factors at
	\(w_3\) by
	\[
	\phi_{1,w_3}=\mathbf1_{U_1},
	\qquad
	\phi_{\flat,w_3}=\mathbf1_{U_\flat},
	\]
	and leave all other local factors unchanged.
	
	The target local orbit meets \(U_1\times U_\flat\), so its local orbital
	integral at \(w_3\) is nonzero; hence
	\(\operatorname{Orb}_{\RT}(\phi^{v_0})\ne0\).
	Now suppose a totally positive definite completion \(\widetilde T\)
	with \(\Diff(\widetilde T,\Vincoherent)=\{v_0\}\) satisfies
	\[
	\bigl(\Int^\epsilon_{v_0}(\widetilde T)-\pDen^\epsilon_{v_0}(\widetilde T)\bigr)
	\operatorname{Orb}_{\widetilde T}(\phi^{v_0})\ne0.
	\]
	By \Cref{lem:local-discrepancy-support},
	\(\widetilde T_{v_0}\) is integral.  The new support is contained
	in the preliminary support, so this nonzero orbital integral gives
	\eqref{eq:isolation-support-condition}. The matrix \(\widetilde T_{w_3}\)
	lies in \(T(U_1\times U_\flat)\subset B\), and the choice of \(B\)
	therefore forces \(\widetilde T=\RT\), proving the required vanishing.
\end{proof}

\subsection{Vanishing of the error series}
For the remainder of the proof, fix global data satisfying
\Cref{prop:globalization}. We use the following
one-place boundary-vanishing criterion.

\begin{prop}\label{lem:boundary-vanishing}
	Let \(\psi:\FGloball\backslash\BA_{\FGloball}\to\BC^\times\) be the fixed global
	additive character from \Cref{subsec:Siegel-Eisenstein}.  Assume that
	its local component \(\psi_{v_0}\) is unramified:
	$\operatorname{cond}(\psi_{v_0})=O_{\FGlobal_{v_0}}.$
	Let
	\[
	w=\begin{pmatrix}0&1\\-1&0\end{pmatrix},
	\qquad
	\widetilde w_{v_0}\in\Mp_2(\FGlobal_{v_0})
	\]
	be the standard Weyl element and the lift used in the local Weil
	representation.  Let \(f\) be a continuous scalar-valued automorphic
	function on \(\SL_2(\FGloball)\backslash\Mp_2(\BA_{\FGloball})\),
	smooth at the finite places.
	For \(\xi\in \FGloball\), put
	\[
	W_{\xi,f}(h)
	=
	\int_{\FGloball\backslash\BA_{\FGloball}}
	f(n(t)h)\psi(-\xi t)\,dt,
	\]
	where \(dt\) is the \(\psi\)-self-dual Haar measure.  Assume:
	\begin{enumerate}[(i)]
		\item \label{lem:boundary-vanishing-1}\(f(hn_{v_0}(x))=f(h)\) for every \(x\in O_{\FGlobal_{v_0}}\);
		\item \label{lem:boundary-vanishing-2}\(f(h\widetilde w_{v_0})=a\,f(h)\) for some \(a\in\BC^\times\) and
		every \(h\);
		\item \label{lem:boundary-vanishing-3}\(W_{\xi,f}(h_\infty)=0\) for every \(\xi\in \FGloball^\times\) with
		\(v_0(\xi)=0\) and every \(h_\infty\in\Mp_2(\FGloball\otimes_\BQ\BR)\).
	\end{enumerate}
	Then
	\[
	W_{\xi,f}(h_\infty)=0
	\qquad(\xi\in \FGloball^\times,\ h_\infty\in\Mp_2(\FGloball\otimes_\BQ\BR)).
	\]
	If in addition \(W_{0,f}(h_\infty)=0\) for every \(h_\infty\), then
	\(f=0\).
\end{prop}

\begin{proof}
	We first record the Fourier-uniqueness principle used below.  If an
	automorphic function \(g\) of the same regularity satisfies
	\[
	W_{\xi,g}(h_\infty)=0
	\qquad(\xi\in \FGloball,\ h_\infty\in\Mp_2(\FGloball\otimes_\BQ\BR)),
	\]
	then \(g=0\).  Indeed, for fixed \(h_\infty\), Fourier uniqueness on the
	compact group \(\FGloball\backslash\BA_{\FGloball}\) gives \(g(h_\infty)=0\).  Finite
	smoothness gives a compact open subgroup \(K_f\) fixing \(g\) on the
	right, and strong approximation for \(\SL_2\) propagates the vanishing
	to arbitrary finite components.
	
	Because \(\psi_{v_0}\) has conductor \(O_{\FGlobal_{v_0}}\), the annihilator of \(O_{\FGlobal_{v_0}}\) under
	$(\xi,x)\longmapsto\psi_{v_0}(\xi x)$
	is again \(O_{\FGlobal_{v_0}}\).  Hence, if \(v_0(\xi)<0\), there is \(x\in O_{\FGlobal_{v_0}}\) with
	\(\psi_{v_0}(\xi x)\ne1\). By (i), for every \(h_\infty\), viewed as an adelic element
	with identity finite component,
	$W_{\xi,f}(h_\infty)
	=
	\psi_{v_0}(\xi x)\,W_{\xi,f}(h_\infty),$
	and therefore \(W_{\xi,f}(h_\infty)=0\) for \(v_0(\xi)<0\).
	
	Now take \(b\in\varpi^{-1}O_{\FGlobal_{v_0}}\), where \(\varpi\) is a
	uniformizer.  The automorphic function
	\(h\longmapsto f(hn_{v_0}(b))-f(h)\) has Fourier coefficients
	\[
	\bigl(\psi_{v_0}(\xi b)-1\bigr)W_{\xi,f}(h_\infty).
	\]
	They vanish for \(\xi=0\); for \(v_0(\xi)<0\) by the preceding paragraph;
	for \(v_0(\xi)=0\) by (iii); and for \(v_0(\xi)>0\) because
	\(\xi b\in O_{\FGlobal_{v_0}}\), so \(\psi_{v_0}(\xi b)=1\).  Fourier uniqueness therefore gives
	\begin{equation}\label{eq:boundary-N-minus-one}
		f(hn_{v_0}(b))=f(h)
		\qquad(b\in\varpi^{-1}O_{\FGlobal_{v_0}}).
	\end{equation}
	
	For \(x\in\FGlobal_{v_0}\), write
	\(n^-(x)=\begin{psmallmatrix}1&0\\x&1\end{psmallmatrix}\)
	and use the same notation for its canonical metaplectic lift.
	For an additive subgroup \(A\subseteq\FGlobal_{v_0}\), let
	\(N^-_{v_0}(A)=\{n^-(x):x\in A\}\).
	Conditions (i) and (ii) give invariance under the opposite
	unipotent subgroup \(N_{v_0}^-(O_{\FGlobal_{v_0}})\), as
	$w^{-1}n(x)w=n^-(-x).$
	In \(\SL_2(\FGlobal_{v_0})\), one has
	\begin{align*}
	&\bigl(n(\varpi^{-1})n^-(-\varpi)n(\varpi^{-1})\bigr)
	 \bigl(n(1)n^-(-1)n(1)\bigr)^{-1}\\
	&\qquad=
	\begin{pmatrix}0&\varpi^{-1}\\-\varpi&0\end{pmatrix}
	\begin{pmatrix}0&1\\-1&0\end{pmatrix}^{-1}
	=m(\varpi^{-1}).
	\end{align*}
	Here all upper-unipotent parameters belong to
	\(\varpi^{-1}O_{\FGlobal_{v_0}}\), and all opposite-unipotent
	parameters belong to \(O_{\FGlobal_{v_0}}\).
	Take the canonical unipotent lifts in the product on the left.
	Every factor fixes \(f\), so the product fixes \(f\).
	Conjugation by this product sends \(n(x)\) to
	\(n(\varpi^{-2}x)\).  Using this identity
	and \eqref{eq:boundary-N-minus-one}, we obtain invariance under
	\(N_{v_0}(\varpi^{-(2j+1)}O_{\FGlobal_{v_0}})\) for every \(j\geq0\).  Their union is
	\(N_{v_0}(\FGlobal_{v_0})\).  Hence \(f\) is right invariant under \(N_{v_0}(\FGlobal_{v_0})\).
	
	Finally, for any \(\xi\ne0\), choose \(x\in \FGlobal_{v_0}\) with
	\(\psi_{v_0}(\xi x)\ne1\). Translation by \(n_{v_0}(x)\) gives, for every \(h_\infty\),
	$W_{\xi,f}(h_\infty)
	=
	\psi_{v_0}(\xi x)W_{\xi,f}(h_\infty),$
	so \(W_{\xi,f}(h_\infty)=0\). If \(W_{0,f}(h_\infty)=0\) for every \(h_\infty\),
	the Fourier-uniqueness principle gives \(f=0\).
\end{proof}

\begin{lem}\label{lem:Fourier-E}
	The modular form $\CE_{T_\flat}(\phi_1,\phi_\flat)\in\CA_{\hol}(\SL_2(\FGlobal)\backslash\Mp_2(\BA_\FGlobal),\dimofV/2)_{\ov{\BQ}}
	\otimes_{\ov{\BQ}}\C_\fkd$ defined in \Cref{def:error-series} satisfies for each $\xi\in\FGlobal$, 
	\[\CE_{T_\flat,\xi}(\phi_1,\phi_\flat)=2\vol^{\mathrm{GL}}_{v_0}(K_{2,v_0})\log q_{v_0}
	\sum_{\substack{
			T=\left(\begin{smallmatrix}\xi&*\\ *&T_\flat\end{smallmatrix}\right)
			\in\Sym_{\dimofV-1}(\FGlobal)_{>0}\\
			\Diff(T,\Vincoherent)=\{{v_0}\}}}
	\bigl(\Int_{v_0}^{\epsilon_{v_0}}(T)-\pDen_{v_0}^{\epsilon_{v_0}}(T)\bigr)
	\operatorname{Orb}_T((\phi_1\otimes\phi_\flat)^{v_0})\]
\end{lem}

\begin{proof}
	By \Cref{prop:aux-special-cycle,lem:local-unit-cancellation}, the auxiliary terms
	in \Cref{prop:modularity} vanish. For
	\(w\in\Sigma_{F,f}^{\fkd}\setminus\{v_0\}\), condition~\ref{item:globalization-isolation-I1}
	gives a unit diagonal entry in \(T_\flat\), so \Cref{lem:local-unit-cancellation}
	also eliminates the \(w\)-terms. Only the \(v_0\)-terms remain.
\end{proof}

\begin{thm}\label{thm:vanishing-generating-series}
	The modular form
	$\CE_{T_\flat}(\tau,\phi_1,\phi_\flat)$ vanishes.  
\end{thm}

\begin{proof}
	For every \(\xi\in\FGlobal_{>0}\) with \(v_0(\xi)=0\) and every matrix
	\[
	T=\begin{pmatrix}\xi&*\\ *&T_\flat\end{pmatrix}
	\in\Sym_{\dimofV-1}(\FGlobal)_{>0},
	\qquad \Diff(T,\Vincoherent)=\{v_0\}, 
	\]
	the entry \(T_{11}=\xi\) is a unit at \(v_0\). 
	So \Cref{lem:local-unit-cancellation} gives
	\(\Int_{v_0}^{\epsilon_{v_0}}(T)=\pDen_{v_0}^{\epsilon_{v_0}}(T)\). Thus \Cref{lem:Fourier-E}  yields
	\begin{equation}\label{eq:boundary-coefficients}
	\CE_{T_\flat,\xi}(\phi_1,\phi_\flat)=0
	\qquad
	\bigl(\xi\in\FGlobal_{>0},\ v_0(\xi)=0\bigr).
	\end{equation}

	Let \(\lambda:\C_{\fkd}\to\BC\) be any
	\(\ov{\BQ}\)-linear functional, acting on the coefficient factor
	as in \Cref{def:error-series}.
	By \Cref{prop:globalization} and the admissibility of
	\(\phi_1\), \(\psi_{v_0}\) is unramified and
	\(\phi_{1,v_0}=\mathbf1_{\VlattiveLocalo}\) for the
	self-dual lattice \(\VlattiveLocalo\).
	By \Cref{def:error-series},
	\(\lambda\circ\CE_{T_\flat}
	(\,\cdot\,,\phi_1,\phi_\flat)\) is a continuous scalar-valued automorphic function,
	smooth at the finite places and holomorphic of parallel weight
	\(\dimofV/2\) at infinity.
	We verify \eqref{lem:boundary-vanishing-1}--\eqref{lem:boundary-vanishing-3} of \Cref{lem:boundary-vanishing} in order.

	For \eqref{lem:boundary-vanishing-1} and \eqref{lem:boundary-vanishing-2}, apply \Cref{prop:dual-relation};
	\(\ov{\BQ}\)-linearity preserves the Weyl eigenvalue
	\(\gamma_{V_{v_0},\psi_{v_0}}\in\ov{\BQ}^{\times}\).
	For \eqref{lem:boundary-vanishing-3},
	\eqref{eq:adelization-general}, Fourier orthogonality,
	\Cref{prop:modularity}, and \eqref{eq:boundary-coefficients} give
	\begin{equation}\label{eq:error-series-Whittaker-coefficients}
	W_{\xi/2,\lambda\circ\CE_{T_\flat}
	(\,\cdot\,,\phi_1,\phi_\flat)}(h_\infty)=0
	\end{equation}
	for every \(\xi\in\FGlobal\) with \(\xi=0\) or \(v_0(\xi)=0\)
	and every
	\(h_\infty\in\Mp_2(\FGlobal\otimes_\BQ\BR)\).
	Since \(v_0\) has odd residue characteristic,
	\(v_0(\xi/2)=v_0(\xi)\) for \(\xi\ne0\).
	This proves \eqref{lem:boundary-vanishing-3} and the vanishing
	of the constant Whittaker coefficient.
	Thus \Cref{lem:boundary-vanishing} yields
	\(\lambda\circ\CE_{T_\flat}
	(\,\cdot\,,\phi_1,\phi_\flat)=0\).
	Since \(\lambda\) is arbitrary and these functionals separate points of
	\(\C_{\fkd}\), every Fourier coefficient of
	\(\CE_{T_\flat}(\tau,\phi_1,\phi_\flat)\) vanishes.
\end{proof}

\begin{proof}[Proof of Theorem~\ref{thm:main2}]
	By \Cref{thm:vanishing-generating-series,lem:Fourier-E}
	and condition~\ref{item:globalization-isolation-I2} of
	\Cref{prop:globalization}, we have
	\[
	2\vol^{\mathrm{GL}}_{v_0}(K_{2,v_0})\log q_{v_0}
	\bigl(\Int^\epsilon_{v_0}(\RT)
	-\pDen^\epsilon_{v_0}(\RT)\bigr)
	\operatorname{Orb}_{\RT}(\phi^{v_0})=0\quad\text{in }\C_\fkd.
	\]
	By the algebraicity established in the proof of
	\Cref{lem:bad-finite-whittaker-algebraicity} and
	condition~\ref{item:globalization-isolation-I2} of \Cref{prop:globalization},
	\[
	 \vol^{\mathrm{GL}}_{v_0}(K_{2,v_0})\operatorname{Orb}_{\RT}(\phi^{v_0})
	 \in\ov{\BQ}^{\times}.
	\]
	The discrepancy
	\(\Int^\epsilon_{v_0}(\RT)-\pDen^\epsilon_{v_0}(\RT)\) is rational.
	By Baker's theorem \cite{Baker1966LinearForms}, the numbers $\log\ell$, with $\ell$ prime, are linearly independent over $\ov{\BQ}$.
	Since \(v_0\nmid\fkd\), the class of \(\log q_{v_0}\) in
	$\C_\fkd$ is nonzero.
	It follows that
	\(\Int^\epsilon_{v_0}(\RT) =\pDen^\epsilon_{v_0}(\RT).\)
	By condition~\ref{item:globalization-isolation-I1} of \Cref{prop:globalization}, this is equivalent to
	\(\Int^\epsilon(\bu)=\pDen^\epsilon(\bu)\), which proves the theorem.
\end{proof}

\subsection{Arithmetic Siegel--Weil formula: the semi-global identity}

\begin{thm}[Semi-global arithmetic Siegel--Weil formula]\label{thm:semi-global-ASW}
	Assume \(\dimofV\ge3\).  Let $v\nmid\fkd$
	be a finite place and let
	$T\in\Sym_{\dimofV-1}(\FGloball)$ be nonsingular with
	$\Diff(T,\Vincoherent)=\{v\}$.  Choose
	\(\mathbf x\in(\VspecialGlobal^{(v)})^{\dimofV-1}\) with
	\(T(\mathbf x)=T\).  For an admissible factorizable Schwartz
	function $\phi$, one has
	\begin{equation*}
		\deg_v\CZ_T^{\BL}(\phi)
		=\frac{2}{\vol(K_2^v)}\,
		\pDen_v^{\epsilon_v}(\mathbf x)\,
		\operatorname{Orb}_{\mathbf x}(\phi^v).
	\end{equation*}
\end{thm}

\begin{proof}
	By \Cref{thm:global-int,thm:main2} we have 
	\[
	 \deg_v\CZ_T^{\BL}(\phi)=\frac{2}{\vol(K_2^v)}\,
	 \Int_v^{\epsilon_v}(\mathbf x)\operatorname{Orb}_{\mathbf x}(\phi^v)=\frac{2}{\vol(K_2^v)}\,
	 \pDen_v^{\epsilon_v}(\mathbf x)\,
	 \operatorname{Orb}_{\mathbf x}(\phi^v).\qedhere
	\]
\end{proof}

	\appendix
	\section{Prismatic crystals}\label{sec:prism}
The purpose of this appendix is to construct the exact crystalline
realization functor $\omega_\cris$ of \Cref{prop:cris-tensor}. Applied to
$V_{\Z_p}\hookrightarrow\End(W_{\Z_p})$, it produces the subcrystal
$\CV_\cris$ used to define special endomorphisms. 

We first recall some basic notions from prismatic theory and the results of
\cites{guo2022prismatic,du2022completed,AnschutzLeBras2023}. 
In \Cref{thm:para-equiv}, we generalize
\cite[Theorem~2.29 and Corollary~2.31]{IKY26} from reductive group
schemes over $\Z_p$ to connected parahoric group schemes. Using quasi-affineness of $\GL(\Lambda)/\CG$, we reduce the extension of a $\CG_\Prism$-torsor to extension over punctured
spectra of rings of the form $W(\bar\kappa)[[u]]$, where Ansch\"utz's
extension theorem applies; see \cite[Proposition~10.3]{Ans22}. 

\Cref{prop:iwahori-obstruction} gives a counterexample to the parahoric
analogue of the regular-surface extension corollary
\cite[p.6, corollary]{IKY26} over a ramified regular surface.

\subsection{Prismatic crystals and \'etale realization}

Throughout this article, a morphism of locally Noetherian formal schemes
is called smooth if it is formally smooth and locally of pseudo-finite type.
Let $K$ be a complete discretely valued field of mixed characteristic $(0,p)$
with ring of integers $O_K$ and perfect residue field, and let $X$ be a smooth
$p$-adic (i.e., the structure morphism $X\to\Spf O_K$ is adic) formal scheme over $O_K$.  
\begin{defn}
	\begin{enumerate}
		\item A $\delta$-ring is a commutative ring $A$ with a map
		$\delta:A\to A$ satisfying $\delta(0)=\delta(1)=0$ and
		\[
		\begin{aligned}
			\delta(a+b)&=\delta(a)+\delta(b)
			-\sum_{i=1}^{p-1}\frac{1}{p}\binom pi a^i b^{p-i},\\
			\delta(ab)&=a^p\delta(b)+b^p\delta(a)+p\delta(a)\delta(b).
		\end{aligned}
		\]
		The map $\varphi_A(a)=a^p+p\delta(a)$ is then a ring endomorphism
		lifting Frobenius modulo $p$. A \emph{prism} is a
		pair $(A,I)$ consisting of a $\delta$-ring $A$ and an ideal $I\subset A$ such
		that $I$ defines a Cartier divisor, $A$ is derived $(p,I)$-complete, and
		\(p\in I+\varphi_A(I)A.\)
		A prism is called \emph{bounded} if the $p^\infty$-torsion of $A/I$ is
		annihilated by a power of $p$.
		
		\item  The absolute prismatic site $X_\Prism$ is the opposite of the category of
		bounded prisms $(A,I)$ equipped with a morphism
		\(\Spf(A/I)\longrightarrow X.\)
		Its coverings are the families of morphisms of prisms whose underlying maps
		of rings are $(p,I)$-completely faithfully flat. The site is equipped with the
		prismatic structure sheaf and the prismatic ideal sheaf
		\(\CO_\Prism(A,I)=A\), \(\CI_\Prism(A,I)=I.\)
		The Frobenius endomorphisms $\varphi_A$ assemble to an endomorphism
		\(\varphi:\CO_\Prism\longrightarrow\CO_\Prism.\)
		
		\item  A prismatic crystal in vector bundles assigns a finite projective
		$A$-module $\CE_A$ to each $(A,I)$ and a base-change isomorphism
		$B\otimes_A\CE_A\simeq\CE_B$ to each morphism
		$(A,I)\to(B,J)$, subject to the identity and cocycle conditions.
		Morphisms are compatible families of module homomorphisms. Equivalently,
		$$\Vect(X_\Prism,\CO_\Prism) = \varprojlim_{(A,I)\in X_\Prism}\Vect(\Spec A).$$
		Here and below, limits of categories are taken as $2$-limits along pullback.
		A sequence in $\Vect(X_\Prism,\CO_\Prism)$ is called \emph{exact} if for each $(A,I)\in X_\Prism$ the corresponding sequence of $A$-vector bundles is exact.

		For a bounded prism $(A,I)$, put
		\(U_{A,I}:=\Spec A\setminus V(p,I).\)
		The category of analytic prismatic crystals is
		\[
		\Vect^{\an}(X_\Prism,\CO_\Prism)
		:=
		\varprojlim_{(A,I)\in X_\Prism}\Vect(U_{A,I}).
		\]
		
		\item For a prismatic crystal $\CE$, write
		\(\varphi^*\CE := \CO_\Prism\otimes_{\varphi,\CO_\Prism}\CE.\)
		A \emph{prismatic $F$-crystal} is a pair $(\CE,\varphi_\CE)$ with
		$\CE\in\Vect(X_\Prism,\CO_\Prism)$ and an isomorphism
		\[
		\varphi_\CE:
		(\varphi^*\CE)[1/\CI_\Prism]
		\xrightarrow{\sim}
		\CE[1/\CI_\Prism].
		\]
		The resulting category is denoted by
		\(\Vect^\varphi(X_\Prism,\CO_\Prism).\)

		\item An \emph{analytic prismatic $F$-crystal} is a pair
		$(\CE,\varphi_\CE)$ with
		$\CE\in\Vect^{\an}(X_\Prism,\CO_\Prism)$ and an isomorphism
		\[
		\varphi_\CE:
		(\varphi^*\CE)[1/\CI_\Prism]
		\xrightarrow{\sim}
		\CE[1/\CI_\Prism].
		\]
		Its category is denoted by
		\(\Vect^{\an,\varphi}(X_\Prism).\)
		By \cite[Proposition~3.7]{guo2022prismatic}, restriction to the analytic loci gives a fully faithful functor
		\begin{equation}\label{eq:analytic-fully-faithful}
			r: \Vect^\varphi(X_\Prism,\CO_\Prism)
			\longrightarrow
			\Vect^{\an,\varphi}(X_\Prism).
		\end{equation}
	\end{enumerate}
\end{defn}

The category $\Loc_{\Z_p}(X_\eta)$ consists of finite-rank $\Z_p$-local
systems on the pro-\'etale site of the rigid generic fiber $X_\eta$.
Its full subcategory $\Loc^{\cris}_{\Z_p}(X_\eta)$ consists of those
whose rationalizations are crystalline in the sense of the local
crystalline period comparison; see \cite[\S2.3.3]{IKY26}.
Exactness of local systems means exactness of the underlying sheaves.

The following equivalence is proved in \cite{guo2022prismatic}*{Theorem~A}
and \cite{du2022completed}*{Theorem~3.46 and Remark~1.5}.
\begin{thm}\label{thm:prism-etale}
	We have a natural symmetric monoidal
	equivalence, compatible with pullback,
	\[
	T:
	\Vect^{\an,\varphi}(X_\Prism)
	\xrightarrow{\sim}
	\Loc_{\Z_p}^{\cris}(X_\eta)
	\]
	between analytic prismatic $F$-crystals on $X$ and crystalline
	$\Z_p$-local systems on the rigid generic fiber $X_\eta$. 
\end{thm}
This generalizes the
case $X=\Spf O_K$ in \cite[Theorem 1.2]{BS23}.

Put $T_{\et}:=T\circ r$. An analytic prismatic $F$-crystal need not
extend to a crystal in vector bundles on all of $X_\Prism$; see
\cite[Example~3.36]{du2022completed}. Prismatic good reduction below
requires such an extension.

\begin{defn}\label{def:prismatic-good}
	A local system $\BL\in\Loc_{\Z_p}(X_\eta)$ has \emph{prismatic good
		reduction relative to $X$} if there exist
	$\CE\in\Vect^\varphi(X_\Prism,\CO_\Prism)$ and an isomorphism $T_{\et}(\CE)\simeq\BL$.
	Let $\Loc^{\Prism\text{-}\mathrm{gr}}_{\Z_p}(X_\eta)$ be the full
	subcategory of such objects, equipped with the exact sequences inherited
	from $\Loc_{\Z_p}(X_\eta)$. This is the convention of
	\cite[Definition~2.28]{IKY26}.
\end{defn}

\begin{defn}\label{def:admissible-prismatic-dieudonne}
An $F$-crystal $(\CE,\varphi_\CE)$ is \emph{effective of Frobenius height
at most one} if its Frobenius linearization extends to a morphism
$\varphi^*\CE\to\CE$ such that
\[
\CI_\Prism\cdot
\operatorname{coker}\bigl(\varphi_\CE:\varphi^*\CE\to\CE\bigr)=0.
\]
We write $\Vect_{[0,1]}^\varphi(X_\Prism,\CO_\Prism)$ for the category
of such crystals that are admissible in the sense of
\cite[Definition~4.5]{AnschutzLeBras2023}. This condition is local for the
quasisyntomic topology by \cite[Proposition~4.9]{AnschutzLeBras2023}.
\end{defn}

\begin{thm}[{\cite[Theorem~4.74]{AnschutzLeBras2023}}]\label{thm:prism-dieudonne}
	Let $\pdivgroup(X)$ be the category of $p$-divisible groups over $X$.
	The contravariant prismatic Dieudonn\'e functor is an antiequivalence
	\[
	\mathbb D_\Prism^{\mathrm{contra}}:
	\pdivgroup(X)^{\mathrm{op}}\xrightarrow{\sim}
	\Vect_{[0,1]}^\varphi(X_\Prism,\CO_\Prism).
	\]
	For $\mathscr H\in\pdivgroup(X)$, the \'etale realization of $\mathbb D_\Prism^{\mathrm{contra}}(\mathscr H)$ is the dual of the homological Tate module by \cite[Proposition~5.25]{AnschutzLeBras2023}, namely
	$T_p(\mathscr H_{X_\eta})^\vee$. We use the covariant functor
	\[
	\mathbb D_\Prism(\mathscr H)
	:=\bigl(\mathbb D_\Prism^{\mathrm{contra}}(\mathscr H)\bigr)^\vee,
	\qquad
	T_{\et}(\mathbb D_\Prism(\mathscr H))\simeq T_p(\mathscr H_{X_\eta}).
	\]
This covariant functor takes values in the category of prismatic
$F$-crystals. Its crystalline realization is the covariant crystalline
Dieudonn\'e crystal, denoted
	$\BD(\mathscr H)$ in the main text. If $\varphi_{\BD(\mathscr H)}$
	denotes its rational Frobenius linearization on the characteristic-$p$
	fiber, then both $p\varphi_{\BD(\mathscr H)}$ and
	$\varphi_{\BD(\mathscr H)}^{-1}$ are integral, since the contravariant
	crystal has Frobenius height at most one.
\end{thm}
\subsection{Parahoric tensor functors and the main theorem}

Let $k$ be the perfect residue field of $O_K$, and fix a
uniformizer $\pi\in O_K$. Let $E(u)\in W(k)[u]$ be the monic Eisenstein
minimal polynomial of $\pi$ over $W(k)[1/p]$, so that $O_K=W(k)[u]/(E(u))$.
Let $\CG$ be a parahoric group scheme over $\Z_p$, i.e.,  $\CG$
is the neutral component of the Bruhat--Tits group scheme attached to a
facet in the enlarged building of a connected reductive $\Q_p$-group $G$.

\begin{defn}\label{def:parahoric-objects}
	For a flat affine group scheme $H$ of finite type over $\Z_p$, let
	$\Rep_{\Z_p}(H)$ be the category of finite free $\Z_p$-modules
	equipped with algebraic $H$-actions and $H$-equivariant linear maps. Its
	short exact sequences are those that are exact on underlying $\Z_p$-modules.
	A representation $\rho:H\to\GL(\Lambda)$ is \emph{faithful} if $\rho$ is
	a closed immersion of group schemes. If $\mathcal C$ is an exact
	$\Z_p$-linear symmetric monoidal category, write
	\[
	H\text{-}\mathcal C
	:=\operatorname{Fun}^{\otimes,\mathrm{ex}}_{\Z_p}
	\bigl(\Rep_{\Z_p}(H),\mathcal C\bigr)
	\]
	for the groupoid of exact $\Z_p$-linear symmetric monoidal functors and
	monoidal natural isomorphisms.
	
	For a $\Z_p$-scheme $S$, an $H$-torsor on $S$ means a right fppf
	$H_S$-torsor. For such a torsor $P$ and $V\in\Rep_{\Z_p}(H)$,
	put
	\[
	P(V):=P\times^{H}(V\otimes_{\Z_p}\CO_S),
	\]
	where the contracted product is defined by $(ph,v)\sim(p,hv)$.
	It is a vector bundle on $S$. We identify vector bundles with finite
	locally free sheaves, and on affine schemes with finite projective modules.
\end{defn}

By \cite[Chapter~VIII, Aside~9.4]{MilneAGS}, there is a faithful
representation $\rho:\CG\hookrightarrow\GL(\Lambda)$. We prove the following
generalization of \cite[Theorem~2.29 and Corollary~2.31]{IKY26}.
\begin{thm}\label{thm:para-equiv}
	Let $\omega\in\CG\text{-}\Loc_{\Z_p}(X_\eta)$, and fix a faithful
	representation $\rho:\CG\hookrightarrow\GL(\Lambda)$. If
	$\omega(\Lambda)$ has prismatic good reduction, there exist an exact
	$\Z_p$-linear symmetric monoidal functor
	\[
	\nu:\Rep_{\Z_p}(\CG)\longrightarrow\Vect^\varphi(X_\Prism,\CO_\Prism)
	\]
	and a monoidal natural isomorphism $\alpha:T_{\et}\circ\nu\simeq\omega$.
	By the full faithfulness of $T_{\et}$, such a pair $(\nu,\alpha)$ is unique up to a unique isomorphism. 
\end{thm}

\subsection{\texorpdfstring{Proof of \Cref{thm:para-equiv}}{Proof of the main theorem}}
Work with $\omega$ and $\Lambda$ as in \Cref{thm:para-equiv}.
Since $\omega(\Lambda)$ has prismatic good reduction, it is crystalline.
By \cite[Proposition~2.21]{IKY26}, $\omega$ therefore takes values in
crystalline local systems. Thus the analytic tensor equivalence of
\cite[Proposition~2.23]{IKY26} applies to $\omega$.

By \cite[Remark A.15]{IKY26}, the closed subgroup
$\CG\subset\GL(\Lambda)$ is the common scheme-theoretic stabilizer of
finitely many tensors $s_\alpha\in\Lambda^{\otimes(i_\alpha,j_\alpha)}$.

The prismatic good reduction hypothesis gives a crystal
$\CE_\Lambda\in\Vect^\varphi(X_\Prism,\CO_\Prism)$ together with an
identification $T_{\et}(\CE_\Lambda)\simeq\omega(\Lambda)$.
By \Cref{thm:prism-etale} and the full faithfulness of
\eqref{eq:analytic-fully-faithful}, the tensors lift uniquely to morphisms
$s_{\alpha,\Prism}:\CO_\Prism\to
\CE_\Lambda^{\otimes(i_\alpha,j_\alpha)}$.
Define
\[
\CQ_\omega:=\Isom_{\CO_\Prism}
\bigl((\Lambda,s_\alpha)\otimes_{\Z_p}\CO_\Prism,
      (\CE_\Lambda,s_{\alpha,\Prism})\bigr).
\]
This is a pseudo-torsor for $\CG_\Prism$ on $X_\Prism$.
By \cite[Proposition~2.23]{IKY26}, its restriction to the analytic locus
$U_{A,I}=\Spec A\setminus V(p,I)$ of each prism $(A,I)$ is a torsor.
By \cite[Proposition~1.29]{IKY26}, it remains to prove that
$\CQ_\omega$ is a torsor on the entire site.

Choose an open covering $\{\Spf R_i\}$ of $X$ by small, and hence base,
$O_K$-algebras in the sense of \cite[\S1.1.5]{IKY26}, i.e., each $R_i$ has the form 
\[
R_i\simeq R_{0,i}\widehat\otimes_{W(k)}O_K,
\]
where $R_{0,i}$ is a $p$-adically complete $W(k)$-algebra admitting a
$p$-adically \'etale toric framing. The completed tensor product is taken
with respect to the $p$-adic topology. After refining the covering, we
may assume that $\Spec R_{0,i}$ is connected.
The Witt Frobenius and the $p$th-power maps on the framing coordinates
extend uniquely to a Frobenius lift on $R_{0,i}$ by formal \'etaleness.
Extend this lift to $R_{0,i}[[u]]$ by $\varphi(u)=u^p$.
The associated Breuil--Kisin prism is
\[
(\fkS_{R_i},(E(u)))=(R_{0,i}[[u]],(E(u))),
\qquad \fkS_{R_i}/(E(u))\simeq R_i.
\]
Here $E(u)\in W(k)[u]$ is the Eisenstein polynomial fixed above.
Since $E(u)$ reduces to a power of $u$ modulo $p$, one has
$V(p,E(u))=V(p,u)$. By \cite[Proposition~1.16]{IKY26}, it suffices to
show that the evaluation of $\CQ_\omega$ on each of these
Breuil--Kisin prisms is a $\CG_{\fkS_{R_i}}$-torsor.

\begin{lem}
	\label{lem:parahoric-quotient}
	Let $H\subset\GL(\Lambda)$ be smooth over $\Z_p$ with reductive generic
	fiber. If $S$ is a Noetherian affine $\Z_p$-scheme and $M$ is a vector
	bundle of rank $\operatorname{rank}_{\Z_p}\Lambda$ on $S$, put
	\[
	F=\Isom_{\CO_S}(\Lambda\otimes_{\Z_p}\CO_S,M).
	\]
	The fppf quotient $Y=F/H$ is a quasi-affine $S$-scheme, and $F\to Y$
	is an $H$-torsor. In particular, $\GL(\Lambda)/H$ is quasi-affine.
\end{lem}
\begin{proof}
	Choose a finite family of defining tensors
	$s_\alpha\in\Lambda^{\otimes(i_\alpha,j_\alpha)}$ for $H$, with
	$i_\alpha,j_\alpha\ge0$, as in
	\cite[Proposition~1.3.2]{Kisin2010}. Put
	\[
	N=\bigoplus_\alpha M^{\otimes(i_\alpha,j_\alpha)},\qquad
	\mathbb V(N)=\underline{\Spec}_S(\Sym(N^\vee)).
	\]
	Thus $\mathbb V(N)$ is the total space of $N$.
	The action of $H$ on $F$ is free, so
	\cite[\href{https://stacks.math.columbia.edu/tag/06PH}{Tag~06PH}]{stacks-project} makes $Y$ an algebraic space and
	$F\to Y$ an $H$-torsor. The morphism
	\[
	F\longrightarrow\mathbb V(N),\qquad
	f\longmapsto\bigl(f(s_\alpha)\bigr)_\alpha
	\]
	descends to $\iota:Y\to\mathbb V(N)$.
	The morphism $\iota$ is a monomorphism of fppf sheaves. Indeed, fppf locally two sections
	of $Y$ lift to frames; if their ordered tensor tuples agree, the change
	of frame lies in their common scheme-theoretic stabilizer $H$.
	
	The morphism $\iota$ is of finite type, since $Y$ is of finite type
	over $S$. As a monomorphism, it has geometric fibers
	with at most one point; hence it is quasi-finite and separated. By
	\cite[\href{https://stacks.math.columbia.edu/tag/0418}{Tag~0418}]{stacks-project}, $Y$ is a scheme. By \cite[\href{https://stacks.math.columbia.edu/tag/05K0}{Tag~05K0}]{stacks-project}, $\iota$ factors as an open immersion followed by a
	finite morphism. Since $\mathbb V(N)$ is affine, the intermediate scheme
	is affine. Hence $Y$ is quasi-affine over $S$.
\end{proof}
Fix $i$, and apply \Cref{lem:parahoric-quotient} with
$S=\Spec\fkS_{R_i}$, $H=\CG$, and
$M=\CE_\Lambda(\fkS_{R_i},(E(u)))$. Write $F$ for the frame torsor
in that lemma and set $U=S\setminus V(p,u)$.
The torsor $\CQ_\omega|_U$ defines a reduction of $F|_U$ to $\CG$,
and hence a section
\begin{equation}\label{eq:parahoric-reduction-section}
s_U:U\longrightarrow Y:=F/\CG.
\end{equation}
The pseudo-torsor $\CQ_\omega$ is a torsor over $S$ if and only if
$s_U$ extends to a section of $Y$ over $S$. Indeed, an extension gives
a $\CG$-reduction of $F$ whose defining tensors agree with
$s_{\alpha,\Prism}$ on $U$ and hence on $S$, since $U$ is dense in the
integral scheme $S$ by \Cref{lem:parahoric-slices}.
We check the extension using the following two lemmas.

\begin{lem}\label{lem:section-test}
	Let $S$ be a Noetherian normal integral scheme, let $U\subset S$ be
	open with complement of codimension at least two, and let $Y\to S$
	be quasi-affine. Fix a section $s_U$ of $Y$ over $U$.
	Suppose that for each $z\in S\setminus U$ there are a morphism
	$a_z:T_z\to S$ and a point $t_z\in T_z$ such that
	\begin{enumerate}[(i)]
		\item $a_z(t_z)=z$;
		\item\label{lem:section-test:ii} $U_z:=a_z^{-1}(U)$ is schematically dense in $T_z$, meaning that
		the restriction map $\CO_{T_z}(V)\to\CO_{U_z}(V\cap U_z)$ is injective
		for every open subset $V\subseteq T_z$;
		\item\label{lem:section-test:iii}  the pullback of $s_U$ to $U_z$ extends to a section $s_z$ of
		$Y\times_S T_z$ over $T_z$.
	\end{enumerate}
	Then $s_U$ extends uniquely to a section of $Y$ over $S$.
\end{lem}
\begin{proof}
	Uniqueness follows from separatedness of $Y\to S$. 
	Choose an open immersion $Y\hookrightarrow Z$ with $Z$ affine over $S$.
	By the Hartogs property the composite section over $U$
	extends uniquely to $\bar s:S\to Z$. 
	The composite of $s_z$ with
	$Y\times_S T_z\hookrightarrow Z\times_S T_z$ and the pullback
	$a_z^*\bar s$ are sections of the same affine $T_z$-scheme and agree
	over $U_z$. Schematic density of $U_z$ in $T_z$ implies that they agree
	on $T_z$. In particular,
	$\bar s(z)=\bar s(a_z(t_z))\in Y$. Thus $\bar s$ factors through $Y$.
\end{proof}
\begin{lem}\label{lem:parahoric-slices}
	Let $R_0$ be a nonzero $p$-adically complete $W(k)$-algebra with
	connected spectrum. Assume that it admits a $p$-adically \'etale framing
	\[
	W(k)\langle t_1^{\pm1},\ldots,t_d^{\pm1}\rangle
	\longrightarrow R_0.
	\]
	Put $S=\Spec R_0[[u]]$, $Z=V(p,u)$, and
	$U=S\setminus Z$. Then $R_0[[u]]$ is a Noetherian regular normal
	domain, and $Z$ has pure codimension two. For every point $z\in Z$,
	there are an algebraically closed field $\bar\kappa$ of characteristic
	$p$ and a $W(k)$-morphism
	\[
	a_z:T_z:=\Spec A_z\longrightarrow S,
	\qquad A_z:=W(\bar\kappa)[[u]],
	\]
	whose closed point maps to $z$ and for which
	\[
	a_z^{-1}(U)=U_z:=T_z\setminus\{(p,u)\}.
	\]
\end{lem}
\begin{proof}
	The framing implies that $R_0$ is Noetherian and regular, and the
	same holds for $R_0[[u]]$. Since $\Spec R_0$ is connected
	and regular, $R_0$ is an integral domain, and so is $R_0[[u]]$.
	The sequence $p,u$ is regular in $R_0[[u]]$, and
	$R_0[[u]]/(p,u)=R_0/p$ is smooth over $k$. Thus $Z$ has pure
	codimension two.

	Choose an algebraic closure $\bar\kappa$ of $\kappa(z)$. The point
	$z$ gives a $k$-algebra homomorphism
	$\bar a:R_0/p\to\bar\kappa$. Let $\bar t_i\in\bar\kappa^\times$
	be the image of $t_i$. Using the map $W(k)\to W(\bar\kappa)$
	induced by $k\to\bar\kappa$, Teichm\"uller lifts define a map
	\[
	W(k)\langle t_1^{\pm1},\ldots,t_d^{\pm1}\rangle
	\longrightarrow W(\bar\kappa),\qquad
	t_i\longmapsto[\bar t_i].
	\]
	Formal \'etaleness lifts $\bar a$ successively and uniquely to
	compatible maps $R_0\to W(\bar\kappa)/p^n$. Taking the inverse
	limit gives a $W(k)$-algebra map
	$\widetilde a:R_0\to W(\bar\kappa)$. Applying it coefficientwise
	gives
	\begin{equation*}
		R_0[[u]]\longrightarrow A_z,\qquad
		\sum_{n\geq0}b_nu^n\longmapsto
		\sum_{n\geq0}\widetilde a(b_n)u^n.
	\end{equation*}
	The induced residue map is $\bar a$, so the closed point of $T_z$
	maps to $z$. The ideal $(p,u)$ extends to $(p,u)A_z$, which proves
	the assertion about $a_z^{-1}(U)$.
\end{proof}

\begin{proof}[Proof of \Cref{thm:para-equiv}]
	Return to $S=\Spec\fkS_{R_i}$ and the reduction section $s_U$ of
	\eqref{eq:parahoric-reduction-section}, where $U=S\setminus V(p,u)$.
	Each slice $T_z=\Spec W(\bar\kappa)[[u]]$ is integral, so its
	nonempty open subset $U_z$ is schematically dense. By
	\Cref{lem:section-test,lem:parahoric-slices}, it suffices to extend
	$a_z^*s_U$ to $T_z$ for each $z\in S\setminus U$.

	The base change of $\CG$ to $W(\bar\kappa)$ is a connected
	parahoric group scheme. Apply \cite[Proposition~10.3]{Ans22} to
	the power-series ring $A_z=W(\bar\kappa)[[u]]$, whose coefficient
	ring is a complete discrete valuation ring with algebraically
	closed residue field. This extends the pulled-back torsor
	$a_z^*(\CQ_\omega|_U)$ on $U_z$ to a $\CG_{A_z}$-torsor
	$\CQ'$ on $T_z$.

	Let
	\[
	\CE'_\Lambda:=\CQ'\times^{\CG_{A_z}}
	(\Lambda\otimes_{\Z_p}A_z)
	\]
	be the associated vector bundle, equipped with its defining
	tensors $s'_{\alpha,\Prism}$. The chosen extension of the torsor
	gives an isomorphism on $U_z$,
	\[
	(a_z^*\CE_\Lambda,a_z^*s_{\alpha,\Prism})|_{U_z}
	\simeq
	(\CE'_\Lambda,s'_{\alpha,\Prism})|_{U_z}.
	\]
	Since $T_z$ is normal and $T_z\setminus U_z$ has codimension two,
	the Hartogs property extends this isomorphism and its inverse
	uniquely to $T_z$. The extended isomorphism preserves the tensors,
	because it does so on $U_z$. It therefore identifies $\CQ'$ with
	a reduction of $a_z^*F$ and gives the required section of
	$Y\times_S T_z$ over $T_z$.

	The section criterion now extends $s_U$ over $S$. Thus
	$\CQ_\omega$ is a torsor on every chosen Breuil--Kisin prism,
	and hence on $X_\Prism$. The Frobenius isomorphism of
	$\CE_\Lambda$ preserves the defining tensors after inverting
	$\CI_\Prism$, and therefore induces an isomorphism
	\[
	\varphi^*\CQ_\omega[1/\CI_\Prism]
	\xrightarrow{\sim}\CQ_\omega[1/\CI_\Prism].
	\]
	The torsor--tensor-functor correspondence of
	\cite[Proposition~1.29]{IKY26} consequently gives a tensor functor
	$\nu$ with values in prismatic $F$-crystals. It is exact because,
	fppf locally where the torsor is trivial, it is the forgetful
	functor followed by extension of scalars. By construction and
	\cite[Proposition~2.23]{IKY26}, the identifications
	$T_{\et}(\nu(V))\simeq\omega(V)$ for $V\in\Rep_{\Z_p}(\CG)$ are
	natural in $V$ and compatible with tensor products. They give the required monoidal
	identification $\alpha:T_{\et}\circ\nu\simeq\omega$.
\end{proof}

\subsection{Crystalline realization}

\begin{prop}\label{prop:compare-pris-crys}
	For a quasisyntomic $\BF_p$-scheme $Y$, there is a natural exact
	symmetric monoidal equivalence
	\[
	\Vect(Y_\Prism,\CO_\Prism)\xrightarrow{\sim}
	\Vect((Y/\Z_p)_\cris,\CO_\cris),
	\]
	compatible with pullback and Frobenius. This is
	\cite[Example~4.7]{BS23}; see also
	\cite[Theorem~6.4]{guo2022prismatic}.
\end{prop}

For $X_0:=X\times_{\Spf\Z_p}\Spec\BF_p$, restrict a prismatic crystal
on $X$ to $(X_0)_\Prism$ and apply this comparison. The scheme $X_0$
is quasisyntomic because $X$ is smooth over $O_K$. This defines
\[
\realization_\cris:\Vect^\varphi(X_\Prism,\CO_\Prism)
\longrightarrow\Vect^\varphi((X_0/\Z_p)_\cris,\CO_\cris).
\]
The target consists of finite locally free crystals $\CE$ equipped with an
isomorphism $F^*\CE[1/p]\simeq\CE[1/p]$, where $F$ is the absolute
Frobenius of $X_0$. 
The functor is exact, symmetric monoidal, and compatible with duals and
pullback. 
\begin{prop}\label{prop:crys-mod-p}
	Let $Y$ be a formal $\Z_p$-scheme and put $Y_0:=Y\times_{\Spf\Z_p}\Spf\BF_p$.
	For the crystalline sites relative to $\Z_p$ with its canonical divided
	powers on $(p)$, restriction induces an exact tensor equivalence
	\[
	\Vect((Y/\Z_p)_\cris,\CO_\cris)
	\xrightarrow{\sim}\Vect((Y_0/\Z_p)_\cris,\CO_\cris).
	\]
\end{prop}
\begin{proof}
	For the closed immersion $Y_0\hookrightarrow Y$, crystalline restriction
	and crystalline pushforward are quasi-inverse exact tensor equivalences
	on finite locally free crystals; see \cite[\S2.3]{Kim2018HodgeTypeRZ}.
	For formal schemes, apply this result modulo $p^n$ and pass to the
	compatible inverse systems of crystals.
\end{proof}
On $(Y/\Z_p)_\cris$, the functor $F^*$ means restriction to $Y_0$,
pullback by the absolute Frobenius of $Y_0$, and crystalline pushforward
along $Y_0\hookrightarrow Y$, as in \cite[\S2.3]{Kim2018HodgeTypeRZ}.

\begin{proof}[Proof of \Cref{prop:cris-tensor}]
	By \Cref{prop:crys-mod-p}, restriction gives an equivalence between
	the categories of vector-bundle crystals on
	$(\widehat{\Sh}_K/\Z_p)_\cris$ and on
	$(\widehat{\Sh}_K\times_{\Spf\Z_p}\Spec\BF_p/\Z_p)_\cris$. Let
	\[
	\omega_p:\Rep_{\Z_p}(\CG)\longrightarrow
	\Loc_{\Z_p}((\widehat{\Sh}_K)_\eta)
	\]
	be the exact tensor functor defined by the $p$-adic level tower. It sends
	$W_{\Z_p}$ to $\CW_p=T_p(\CA_{(\widehat{\Sh}_K)_\eta})$ and each defining tensor
	$s_\alpha$ to $s_{\alpha,p}$. Define
	\[
	\CW_\Prism:=\mathbb D_\Prism(\CA[p^\infty])
	\in\Vect^\varphi((\widehat{\Sh}_K)_\Prism,\CO_\Prism).
	\]
	By \Cref{thm:prism-dieudonne},
	$T_{\et}(\CW_\Prism)\simeq\CW_p$. The representation $W_{\Z_p}$ is
	faithful, so \Cref{thm:para-equiv} gives an exact tensor lift
	\[
	\begin{gathered}
	\omega_\Prism:\Rep_{\Z_p}(\CG)\longrightarrow
	\Vect^\varphi((\widehat{\Sh}_K)_\Prism,\CO_\Prism),\\
	T_{\et}\circ\omega_\Prism\simeq\omega_p.
	\end{gathered}
	\]
	Full faithfulness identifies $\omega_\Prism(W_{\Z_p})$ with
	$\CW_\Prism$, with the specified \'etale identification.
	
	Compose $\omega_\Prism$ with $\realization_\cris$ and a quasi-inverse of the equivalence in
	\Cref{prop:crys-mod-p}. The result is the exact tensor functor
	$\omega_\cris$ in the statement. Compatibility of prismatic and
	crystalline Dieudonn\'e theory identifies its value on $W_{\Z_p}$ with
	$\CW_\cris$. For each $\alpha$, apply the functor to the morphism
	$s_\alpha:\Z_p\to W_{\Z_p}^\otimes$. This gives the integral
	crystalline tensor $s_{\alpha,\cris}$, whose rationalization is
	Frobenius invariant.
	
	Let $x:\Spec\Omega\to\widehat{\Sh}_K$ be a geometric special-fiber
	point. Formal smoothness gives a lift
	$\widetilde x:\Spf O_L\to\widehat{\Sh}_K$, where $O_L$ is a complete
	discrete valuation ring over $\widehat O_{\reflexHodge}$ with residue field
	$\Omega$. By \cite[Remark~7.11]{BS23}, evaluation of
	$\widetilde x^*\omega_\Prism$ on a Breuil--Kisin prism agrees with
	Kisin's tensor functor, compatibly with \'etale realization.
	Under the crystalline comparison of
	\cite[Proposition~5.3]{BS23}, both
	$\realization_\cris(\widetilde x^*\omega_\Prism(s_\alpha))(W(\Omega))$
	and $s_{\alpha,x,0}$ correspond to $s_{\alpha,\widetilde x,\et}$.
	Consequently,
	\[
	(x^*s_{\alpha,\cris})(W(\Omega))
	=\realization_\cris(\widetilde x^*\omega_\Prism(s_\alpha))(W(\Omega))
	=s_{\alpha,x,0}.\qedhere
	\]
\end{proof}

\subsection{The limitation of extension over regular surfaces}

The following example concerns a ramified regular surface to which the
slice construction in \Cref{lem:parahoric-slices} does not apply. It shows
that the extension corollary for reductive groups stated in the
introduction to \cite{IKY26} cannot be generalized to parahoric groups
on arbitrary regular surfaces. 

\begin{prop}
	\label{prop:iwahori-obstruction}
	Let $k$ be a perfect field of characteristic $p$, and put
	\[
	R=W(k)[[x,y]]/(xy-p),\qquad
	S=\Spec R,\qquad U=S\setminus\{(x,y)\}.
	\]
	There are a connected Iwahori group scheme $\mathcal H/\Z_p$ for $\GL_{2, \Q_p}$, a faithful
	representation $\mathcal H\hookrightarrow\GL_{4,\Z_p}$, and an
	$\mathcal H$-torsor $P$ on $U$ with trivial associated rank-four vector
	bundle, such that $P$ admits no extension to $S$.
	The ring $R$ is a complete regular local domain of dimension two.
\end{prop}
\begin{proof}
	Let $V_0=V_1=\Z_p^2$, with fixed standard bases, and consider
	\[
	f_0=\begin{pmatrix}1&0\\0&p\end{pmatrix}:V_0\longrightarrow V_1,
	\qquad
	g_0=\begin{pmatrix}p&0\\0&1\end{pmatrix}:V_1\longrightarrow V_0.
	\]
	The automorphism group scheme $\mathcal H$ of this diagram is an Iwahori model for $\GL_{2,\Q_p}$.
	Explicitly, for every $\Z_p$-algebra $B$ the elements of $\CH(B)$ are pairs
	\[
	(h_0,h_1)=
	\left(\begin{pmatrix}a&pb\\c&d\end{pmatrix},
	\begin{pmatrix}a&b\\pc&d\end{pmatrix}\right),
	\qquad a,b,c,d\in B,\quad ad-pbc\in B^\times.
	\]
	These are exactly the pairs of invertible matrices satisfying
	$f_0h_0=h_1f_0$ and $g_0h_1=h_0g_0$. 
	The group scheme $\mathcal H$ is smooth, its generic fiber is $\GL_{2,\Q_p}$, and its special
	fiber is connected. The action on $V_0\oplus V_1$ gives a closed
	immersion into $\GL_{4,\Z_p}$: its image is cut out by preservation
	of the two summands and of the maps $f_0,g_0$.
	
	The maximal ideal of $R$ is $\mathfrak m=(x,y)$. 
	The ring has dimension two and its maximal ideal is generated by two
	elements; hence it is regular, and in particular a normal domain.
	Moreover $U=D(x)\cup D(y)$. On $U$ let $E_0=E_1=\CO_U^2$ and set
	\[
	f=\begin{pmatrix}x&0\\0&y\end{pmatrix}:E_0\longrightarrow E_1,
	\qquad
	g=\begin{pmatrix}y&0\\0&x\end{pmatrix}:E_1\longrightarrow E_0.
	\]
	Let $P$ be the sheaf of isomorphisms from $(V_0,V_1,f_0,g_0)_U$
	to $(E_0,E_1,f,g)$. The following pairs of matrices give such
	isomorphisms on the indicated opens:
	\[
	\begin{array}{c|cc}
		& A:V_0\otimes\CO_U\longrightarrow E_0
		& B:V_1\otimes\CO_U\longrightarrow E_1\\[2pt]\hline
		D(x)& I_2 & \operatorname{diag}(x,x^{-1})\\[2pt]
		D(y)& \begin{pmatrix}0&1\\1&0\end{pmatrix}
		& \begin{pmatrix}0&y^{-1}\\y&0\end{pmatrix}
	\end{array}
	\]
	Indeed, all these matrices are invertible on their specified opens and
	satisfy $fA=Bf_0$ and $gB=Ag_0$, using $xy=p$. 
	This defines a (Zariski) $\CH$-torsor on $U$. 
	
	Suppose that $P$ extended to an $\mathcal H$-torsor $Q$ on $S$. 
	The associated rank $2$ bundles $Q(V_i)$ on $S$ satisfy $Q(V_i)|_U\cong E_i$. By
	the Hartogs property, these identifications extend uniquely
	to $Q(V_i)\simeq\CO_S^2$. Under these identifications, the morphism $Q(f_0):Q(V_0)\to Q(V_1)$ is
	the matrix $\operatorname{diag}(x,y)$, because both morphisms agree on
	the dense open subset $U$. The fiber of the cokernel of $Q(f_0)$ at
	$\mathfrak m$ therefore has dimension $2$ over $k$. Since $Q$ is an $\mathcal H$-torsor,
	this fiber is fppf-locally isomorphic to $\mathrm{coker}(f_0)\otimes_{\Z_p}k$,
	which has dimension $1$ over $k$, a contradiction.
\end{proof}

	\bibliography{ref}
\end{document}